\documentclass[12pt]{article}
\usepackage{amsfonts,amsmath,amsxtra}
\usepackage{amssymb,amscd}
\usepackage[matrix,arrow]{xy}
\usepackage[mathscr]{eucal}
\def\hybrid{\topmargin 0pt      \oddsidemargin 0pt
        \headheight 0pt \headsep 0pt
        \textwidth 16.5cm
        \textheight 23cm
        \marginparwidth 0.0in
        \parskip 5pt plus 1pt   \jot = 1.5ex}
\catcode`\@=11
\def\marginnote#1{}
\newcount\hour
\newcount\minute
\newtoks\amorpm
\hour=\time\divide\hour by60 \minute=\time{\multiply\hour by60
\global\advance\minute by-\hour}
\edef\standardtime{{\ifnum\hour<12 \global\amorpm={am}%
        \else\global\amorpm={pm}\advance\hour by-12 \fi
        \ifnum\hour=0 \hour=12 \fi
      \number\hour:\ifnum\minute<10 0\fi\number\minute\the\amorpm}}
\edef\militarytime{\number\hour:\ifnum\minute<10 0\fi\number\minute}

\def\draftlabel#1{{\@bsphack\if@filesw {\let\thepage\relax
   \xdef\@gtempa{\write\@auxout{\string
      \newlabel{#1}{{\@currentlabel}{\thepage}}}}}\@gtempa
   \if@nobreak \ifvmode\nobreak\fi\fi\fi\@esphack}
        \gdef\@eqnlabel{#1}}
\def\@eqnlabel{}
\def\@vacuum{}
\def\draftmarginnote#1{\marginpar{\raggedright\scriptsize\tt#1}}

\def\draft{\oddsidemargin -0.1truein
        \def\@oddfoot{\sl preliminary draft \hfil
        \rm\thepage\hfil\sl\today\quad\militarytime}
        \let\@evenfoot\@oddfoot \overfullrule 3pt
        \let\label=\draftlabel
        \let\marginnote=\draftmarginnote
\def\@eqnnum{{\rm (\theequation)}
\rlap{\kern\marginparsep\tt\@eqnlabel}%
\global\let\@eqnlabel\@vacuum}  }

\newcommand{\RR}{{\mathbb{R}}}
\newcommand{\CC}{{\mathbb{C}}}

\newfont{\Bbbb}{msbm7 scaled 1\@ptsize00}
\newcommand{\zs}{\raise-1pt\hbox{$\mbox{\Bbbb Z}$}}

scaled 1\@ptsize00

\font\sevenmsa=msam6 
\newfam\msafam
\textfont\msafam=\sevenmsa
\def\hexnumber@#1{\ifnum#1<10 \number#1\else
\ifnum#1=10 A\else\ifnum#1=11 B\else\ifnum#1=12 C\else \ifnum#1=13
D\else\ifnum#1=14 E\else\ifnum#1=15 F\fi\fi\fi\fi\fi\fi\fi}
\def\msa@{\hexnumber@\msafam}
\def\llcorner{\delimiter"4\msa@78\msa@78 }
\def\lrcorner{\delimiter"5\msa@79\msa@79 }
\mathchardef\blacktriangleright="3\msa@49
\mathchardef\blacktriangleleft="3\msa@4A \font\tenmsb=msbm10 scaled
1\@ptsize00
\newfam\msbfam
\textfont\msbfam=\tenmsb \scriptfont\msbfam=\tenmsb

\newdimen\linethick  \linethick=0.4pt
\newdimen\hboxitspace    \hboxitspace=5pt
\newdimen\vboxitspace    \vboxitspace=5pt

\def\fr#1{%
\be\new \vcenter{ \hrule height\linethick
           \hbox{\vrule width\linethick
                 \kern\hboxitspace
                 \vbox{\kern\vboxitspace
                       \hbox{$\begin{array}{c}\displaystyle#1
          \end{array}$}%
                       \kern\vboxitspace}%
                 \kern\hboxitspace
                 \vrule width\linethick}%
           \hrule height\linethick}%
\ee}

\newdimen\Squaresize \Squaresize=14pt
\newdimen\Thickness \Thickness=0.5pt

\def\Square#1{\hbox{\vrule width \Thickness
   \vbox to \Squaresize{\hrule height \Thickness\vss
      \hbox to \Squaresize{\hss#1\hss}
   \vss\hrule height\Thickness}
\unskip\vrule width \Thickness} \kern-\Thickness}

\def\Vsquare#1{\vbox{\Square{$#1$}}\kern-\Thickness}

\def\numberbysection{\@addtoreset{equation}{section}
        \def\theequation{\thesection.\arabic{equation}}}
\numberbysection

\renewcommand{\theequation}{\thesection.\arabic{equation}}
\def\titlepage{\@restonecolfalse\if@twocolumn\@restonecoltrue\onecolumn
     \else \newpage \fi \thispagestyle{empty}\c@page\z@
        \def\thefootnote{\fnsymbol{footnote}} }

\def\endtitlepage{\if@restonecol\twocolumn \else  \fi
        \def\thefootnote{\arabic{footnote}}
        \setcounter{footnote}{0}}  
\relax

\hybrid
\makeatletter
\newdimen\normalarrayskip            
\newdimen\minarrayskip               
\normalarrayskip\baselineskip \minarrayskip\jot
\newif\ifold             \oldtrue            \def\new{\oldfalse}
\def\arraymode{\ifold\relax\else\displaystyle\fi}
\def\eqnumphantom{\phantom{(\theequation)}} 
\def\@arrayskip{\ifold\baselineskip\z@\lineskip\z@
     \else
     \baselineskip\minarrayskip\lineskip1\baselineskip\fi}

\def\@arrayclassz{\ifcase \@lastchclass \@acolampacol \or
\@ampacol \or \or \or \@addamp \or
   \@acolampacol \or \@firstampfalse \@acol \fi
\edef\@preamble{\@preamble
  \ifcase \@chnum
     \hfil$\relax\arraymode\@sharp$\hfil
     \or $\relax\arraymode\@sharp$\hfil
     \or \hfil$\relax\arraymode\@sharp$\fi}}

\def\@array[#1]#2{\setbox\@arstrutbox=\hbox{\vrule
     height\arraystretch \ht\strutbox
     depth\arraystretch \dp\strutbox
width\z@}\@mkpream{#2}\edef\@preamble{\halign \noexpand\@halignto
\bgroup \tabskip\z@ \@arstrut \@preamble \tabskip\z@ \cr}%
\let\@startpbox\@@startpbox \let\@endpbox\@@endpbox
  \if #1t\vtop \else \if#1b\vbox \else \vcenter \fi\fi
  \bgroup \let\par\relax
  \let\@sharp##\let\protect\relax
  \@arrayskip\@preamble}
\def\eqnarray{\stepcounter{equation}%
              \let\@currentlabel=\theequation
              \global\@eqnswtrue
              \global\@eqcnt\z@
              \tabskip\@centering              
              \let\\=\@eqncr
              $$%
            \halign to \displaywidth  \bgroup
             \eqnumphantom \@eqnsel
      \hskip\@centering                               
    $\displaystyle  \tabskip\z@ {##}$%
    &\global\@eqcnt\@ne \hskip 2\arraycolsep
         $ \displaystyle  \arraymode{##}$\hfil
    &\global\@eqcnt\tw@ \hskip 2\arraycolsep
         $\displaystyle\tabskip\z@{##}$\hfil
         \tabskip\@centering
    &{##}\tabskip\z@\cr}
\makeatother
\newtheorem{te}{Theorem}[section]

\newtheorem{prop}{Proposition}[section]           
\newtheorem{cor}{Corollary}[section]
\newtheorem{lem}{Lemma}[section]
\newtheorem{ex}{Example}[section]
\newtheorem{rem}{Remark}[section]
\newcommand{\beq}[1]{\begin{equation}\label{#1}}
\newcommand\eeq{\end{equation}}
\newcommand\bqa{\begin{eqnarray}}
\newcommand\eqa{\end{eqnarray}}
\def\be{\begin{eqnarray}\new\begin{array}{cc}}
\def\ee{\end{array}\end{eqnarray}}

\def\beq{\begin{equation}}
\def\eeq{\end{equation}}
\def\bse{\begin{subequations}}                
\def\ese{\end{subequations}}
\def\bp{\begin{pmatrix}}
\def\ep{\end{pmatrix}}

\def\i{\imath}

\newcommand{\R}{\mathbb{R}}

\def\stack#1#2{\raise0.7pt\hbox{$\mathrel{\mathop{#2}\limits^{#1}}$}}
\def\tr{\triangleright}
\def\tl{\triangleleft}
\def\sem{\mathsurround=0pt \raise1pt
\hbox{$\scriptscriptstyle>\!\!$}\:\!\!\tl}
\def\mes{\mathsurround=0pt \tr\!\:\!\raise0.8pt
\hbox{$\scriptscriptstyle\!\!<$}\,}
\def\]{\mathsurround=0pt ]\raise-2pt\hbox{$_\ast$}}

\def\la{\lambda}

\def\<{\langle}
\def\>{\rangle}

\def\frak{\mathfrak}

\def\CH{\mathcal{H}}

\def\we{\raise-1pt\hbox{$\,\stackrel{\wedge}{,}\,$}}

\def\pr {\partial}

0
\begin{document}

\footnotesize
\normalsize

\newpage

\thispagestyle{empty}

\begin{center}

\phantom.
\bigskip
{\hfill{\normalsize hep-th/yymmnnn}\\
\hfill{\normalsize ITEP-TH-XX/07}\\
\hfill{\normalsize HMI-07-06}\\
\hfill{\normalsize TCD-MATH-07-12}\\
[10mm]\Large\bf  New Integral Representations  of Whittaker
Functions for Classical Lie Groups} \vspace{0.5cm}

\bigskip\bigskip
{\large A. Gerasimov}
\\ \bigskip
{\it Institute for Theoretical and
Experimental Physics, 117259, Moscow,  Russia,} \\
{\it  School of Mathematics, Trinity
College, Dublin 2, Ireland, } \\
{\it  Hamilton
Mathematics Institute, TCD, Dublin 2, Ireland},\\
\bigskip
{\large D. Lebedev\footnote{E-mail: lebedev@itep.ru}}
\\ \bigskip
{\it Institute for Theoretical and Experimental Physics, 117259,
Moscow, Russia},\\{\it Max-Planck-Institut f\"ur Mathematik,
Vivatsgasse 7, D-53111
Bonn, Germany},\\
\bigskip
{\large S. Oblezin} \footnote{E-mail: Sergey.Oblezin@itep.ru}\\
\bigskip {\it Institute for Theoretical and Experimental Physics,
117259, Moscow, Russia},\\{\it Max-Planck-Institut f\"ur Mathematik,
Vivatsgasse 7, D-53111
Bonn, Germany}.\\
\end{center}

\vspace{0.5cm}

\begin{abstract}
\noindent

We propose   integral representations of the Whittaker functions for
the classical Lie algebras  $\mathfrak{sp}_{2\ell}$, $\mathfrak{so}_{2\ell}$
and $\mathfrak{so}_{2\ell+1}$. These integral
representations generalize the integral representation of
$\mathfrak{gl}_{\ell+1}$-Whittaker functions first introduced by
Givental.  One of the salient features of the Givental
representation is its recursive structure with respect to the rank
$\ell$ of the Lie algebra $\mathfrak{gl}_{\ell+1}$. The proposed
generalization of the Givental representation to the classical Lie algebras
retains this property.
 It was shown elsewhere that the integral recursion operator
for $\mathfrak{gl}_{\ell+1}$-Whittaker function in the Givental
representation coincides with  a degeneration of the Baxter
$\mathcal{Q}$-operator for $\widehat{\mathfrak{gl}}_{\ell+1}$-Toda
chains. We construct $\mathcal{Q}$-operator for  affine Lie algebras
$\widehat{\mathfrak{so}}_{2\ell}$,
$\widehat{\mathfrak{so}}_{2\ell+1}$ and a twisted form of
$\widehat{\mathfrak{gl}}_{2\ell}$.
We  demonstrate that the relation between recursion
integral operators of the generalized Givental representation and
degenerate $\mathcal{Q}$-operators remains valid for all classical
Lie algebras.

\end{abstract}

\vspace{1cm}

\clearpage \newpage



\section{Introduction}

A remarkable  integral representation  for  common eigenfunctions of
$\mathfrak{gl}_{\ell+1}$-Toda chain  Hamiltonian operators
 was proposed by A.~Givental  \cite{Gi} (see also \cite{JK}). The integral
representation appears naturally in a  construction of a mirror dual
of the theory of Type A topological closed strings on
$\mathfrak{gl}_{\ell+1}$-flag manifolds.  The Givental integral
representation has many interesting properties. It has an explicit
recursive structure over the rank $\ell$ of the corresponding Lie
algebra $\mathfrak{gl}_{\ell+1}$. The integrand in the integral
representation allows for purely combinatorial description in terms of a
simple graph. This graph captures a flat degenerations of  flag
manifolds to Gorenstein toric Fano varieties (torification)
\cite{Ba}, \cite{BCFKS}.

In \cite{GKLO}, the Givental integral representation was reconsidered in
the framework of the representation theory approach to quantum
integrable systems.  According to B.~Kostant \cite{Ko1},\cite{Ko2}
the common eigenfunctions of $\mathfrak{g}$-Toda chain  Hamiltonian
operators are given by generalizations of classical Whittaker
functions and can be expressed in terms of  the matrix elements of
infinite-dimensional representations of the universal enveloping algebra
$\mathcal{U}(\mathfrak{g})$. It was demonstrated in \cite{GKLO} that
 the  Givental representation of $\mathfrak{gl}_{\ell+1}$-Toda
eigenfunctions coincides with an integral representation of the
relevant matrix elements obtained by using a particular parametrization
of an open part of the $\mathfrak{gl}_{\ell+1}$-flag manifold. A
conceptual explanation for  the particular choice of coordinates on
flag manifolds was proposed using a relation with the Baxter
$\mathcal{Q}$-operator formalism. It was noticed that the Givental
integral representation has  a recursive structure connecting 
the $\mathfrak{gl}_{\ell}$ -  and $\mathfrak{gl}_{\ell+1}$-Whittaker
functions by  simple integral transformations. The corresponding
integral operator coincides with a particular degeneration of the
Baxter $\mathcal{Q}$-operator  for
$\widehat{\mathfrak{gl}}_{\ell+1}$-Toda chain \cite{PG}. It is
well-known that $\mathcal{Q}$-operators realise the quantum B\"acklund
transformations in quantum integrable systems. On the other hand,  in
the classical limit, the formalism of $\mathcal{Q}$-operators allows
us to define special coordinate system on the phase space. Thus
degenerate $\mathcal{Q}$-operators define particular coordinates on
an open part of flag manifolds and therefore  lead to the Givental
integral representation of $\mathfrak{gl}_{\ell+1}$-Whittaker
function.

Until now a generalization of the Givental integral representation of
$\mathfrak{gl}_{\ell+1}$-Whittaker functions to Lie algebras
other then  $\mathfrak{gl}_{\ell+1}$ was not known. The only
known generalization \cite{BCFKS}, \cite{Ri} of the Givental
construction is an integral representation for common eigenfunctions
of certain degenerations $\mathfrak{gl}_{\ell+1}$-Toda chains
\cite{STS}.  It is  based on  flat degenerations
of partial flag  manifolds $G/P$ for $G=GL(\ell+1,\mathbb{C})$, $P$
being a parabolic subgroup of $G$ \cite{BCFKS}. In this paper we
propose a generalization of the Givental construction 
for classical Lie algebras 
$\mathfrak{sp}_{2\ell}$, $\mathfrak{so}_{2\ell}$ and $\mathfrak{so}_{2\ell+1}$. Our
construction possesses all  characteristic properties of the original
Givental integral representation. The integral representations for
 the classical Lie algebras have recursive structure. The integrands of the
integral representations have  combinatorial descriptions in terms
of graphs. The proposed generalization to the classical Lie algebras
is based on a modification of  a well-known factorized
representation of  generic elements of maximal
 unipotent subgroups of the corresponding Lie groups
 \cite{Lu} (see  also  \cite{FZ}, \cite{BZ}). The construction of
 the modified factorized representation essentially
uses the realization of  maximal unipotent subgroups of classical
Lie groups as explicitly described subgroups of upper-triangular
matrices (see e.g. \cite{DS}). We define Baxter
$\mathcal{Q}$-operators associated with  the classical affine Lie
algebras  $\widehat{\mathfrak{so}}_{2\ell}$,  $\widehat{\mathfrak{so}}_{2\ell+1}$
and a twisted form of $\widehat{\mathfrak{gl}}_{2\ell}$.
We  demonstrate that the relation between recursion
integral operators of the generalized Givental representation and
degenerate $\mathcal{Q}$-operators remains valid for all classical
Lie algebras.

 The novel feature of the constructed integral representation is that,
in contrast with the $\mathfrak{gl}_{\ell+1}$ case (where the kernel of
the recursion operator is a simple function), the integral kernels
of the recursion  operators for all  other classical Lie groups are given by
non-trivial integrals. This suggests that the recursion
operators can be obtained as a composition of elementary
operators. Indeed, for zero eigenvalues, recursion operators
relating the Toda chain eigenfunctions of the Lie algebras 
 with adjacent ranks can be represented as
compositions of elementary recursion operators relating the Toda chain
eigenfunctions of {\it different} classical series. This might not be
so surprising due to the fact that we  essentially use a
realisation of all classical Lie groups as subgroups  of general Lie
groups of large enough rank.

Let us stress that the construction of  integral representations of
$\mathfrak{g}$-Whittaker functions presented in this paper also has 
a natural interpretation  in terms of  torification of flag manifolds
associated with classical Lie groups. The graph encoding the
integrand of the Givental representation for a  classical Lie group
allows us to describe toric degeneration of the corresponding flag
manifold explicitly (thus  generalizing the results of \cite{BCFKS}
to all classical Lie groups).

One of the interesting applications  of the Givental integral
representation of $\mathfrak{g}$-Whittaker functions for classical
Lie algebras might be a construction of mirror duals  for closed
strings on  flag spaces associated with classical Lie groups $G$,
$\mathfrak{g}={\rm Lie}(G)$.
 According to Givental \cite{Gi} the mirror dual to
Type A topological string theories  on flag manifolds 
associated with Lie groups $G$ should
be Landau-Ginzburg models associated with Langlands dual Lie groups
$G^{\vee}$ such that the generating function of the genus zero
correlators is a $\mathfrak{g}^{\vee}$-Whittaker function,
$\mathfrak{g}^{\vee}={\rm Lie}(G^{\vee})$. In the case of
$\mathfrak{g}=\mathfrak{gl}_{\ell+1}$ Givental provides a
description of the  dual
 Landau-Ginzburg model in terms of the integrand of the integral
 representation of the corresponding Whittaker function. Moreover
the interpretation  of the integral representation of
a $\mathfrak{g}$-Whittaker function in terms of a torification of flag
manifolds \cite{Ba}, \cite{BCFKS} allows us to construct the  mirror
map  explicitly. Thus using the same reasoning, the generalization
of the Givental integral representation  for classical Lie groups
allows us  to infer the superpotential of the 
corresponding Landau-Ginzburg model from the integrand. 
Moreover, similar  to the case of 
$\mathfrak{g}=\mathfrak{gl}_{\ell+1}$,  the interpretation  of the
integrand in terms of a toric degeneration of the flag manifold
provides an explicit construction of the mirror map. We will 
discuss this construction elsewhere.

Finally, note that some of the results presented in this paper was
announced previously in \cite{GLO}, \cite{GLO1}.

The plan of this paper  is as follows. In Part I we formulate the
results for the classical Lie algebras $\mathfrak{sp}_{2\ell}$,
$\mathfrak{so}_{2\ell}$ and $\mathfrak{so}_{2\ell+1}$. The main results
are formulated in the Theorems \ref{teintmfgl}, \ref{teintmfso2l1},
\ref{teintmfsp2l}, \ref{teintmfso2l} respectively.
In Part II we collect the proofs of the results presented in Part I.

{\em Acknowledgments}:  The authors are grateful to S.~Kharchev for
discussions at the initial stage of this project.  The 
research of the first author  was  partly supported by the
Enterprise Ireland Basic Research Grant. The second author is
grateful to School of Mathematics, Trinity College Dublin and
Max-Planck-Institut f\"ur Mathematik for hospitality. The third
author thanks to  Max-Planck-Institut f\"ur Mathematik for
hospitality. The research of the third author been partly supported
by the RF President Grant MK-134.2007.1.

\section{ Part I:  Results}

\subsection{Toda chain eigenfunctions as  matrix elements}

Eigenfunctions of $\mathfrak{g}$-Toda chain are given by particular
matrix elements of infinite-dimensional representations of Lie
algebra $\mathfrak{g}$ \cite{Ko1}, \cite{Ko2} (for detailed
exposition see e.g. \cite{Et}). In this section we provide integral
representations of these matrix elements with integrands  being
expressed in terms of matrix elements of finite-dimensional
representations of $\mathfrak{g}$. In the following sections we
derive explicit expressions for the  relevant matrix elements of
finite-dimensional representations  and thus obtain integral
representations of $\mathfrak{g}$-Toda chain eigenfunctions
generalizing the results of Givental for
$\mathfrak{g}=\mathfrak{gl}_{\ell+1}$. The construction will be
given for all classical Lie algebras.  We start with standard
definitions in the theory of Lie algebras mostly following \cite{K}
(for a discussion of root data of reductive groups see e.g.
\cite{S}).

\subsubsection{Root data for reductive groups}

Root datum is a quadruple $(X,\Phi,X^{\vee},\Phi^{\vee})$ where $X$
is a lattice of a finite rank, $X^{\vee}$ is a dual lattice, $\Phi$
and $\Phi^{\vee}$ are subsets of $X$ and $X^{\vee}$ supplied with  a
bijection $\alpha\mapsto\alpha^{\vee}$ of $\Phi$ onto $\Phi^{\vee}$
and the following conditions hold.  One has
$\langle\alpha,\alpha^{\vee}\rangle=2$ for any $\alpha\in \Phi$.
Subsets $\Phi$ and $\Phi^{\vee}$ should be stable with respect to
any automorphisms $s_{\alpha}$, $s_{\alpha^{\vee}}$:   $$
s_{\alpha}(x)=x-\langle x,\alpha^{\vee}\rangle\alpha,\qquad
s_{\alpha^{\vee}}(y)=y-\langle y,\alpha\rangle\alpha^{\vee},\qquad
 x\in X,\,\,\, y\in X^{\vee},\,\,\,\alpha\in \Phi.$$ Let $Q\subset X$ be  a sublattice
generated by elements of $\Phi$, and $P$ be a lattice defined as
$$P=\{x\in X\otimes\mathbb{Q}|\,\langle
x,\,\alpha^{\vee}\rangle\in\mathbb{Z},\,\alpha\in\Phi\}.$$ One has
$Q\subset P$ and $P/Q$ is a finite group. Let $X_0\subset X$ be
sublattice defined as
$$
X_0=\{x\in X|\,\langle x, y\rangle=0, \, y\in \Phi^{\vee}\}.
$$

With any reductive Lie group one can associate root datum.
 Let $G$ be a connected reductive complex Lie group
 and $H\subset G$ be a maximal torus. We associate to a pair $(G,H)$ a root
datum $(X,\Phi,X^{\vee},\Phi^{\vee})$ as follows. Here $X$ is a free
abelian finite rank group of $\mathbb{Q}$-characters of $H$,
$X^\vee={\rm Hom}(\mathbb{C}^*,H)$ is a dual group of one-parameter
multiplicative subgroups of $H$. The pairing
 $\langle\, ,\,\rangle: X\times X^{\vee}\to \mathbb{Z}$
is defined as
$$\lambda(u(t))=t^{\langle\lambda,\,u\rangle},\,\hspace{2cm}
 \lambda\in X,\,\,\,\, u\in X^{\vee},\,\,\,\,\, t\in\mathbb{C}^*.$$
Then $\Phi$ and $\Phi^\vee$ are finite subsets of $X$ and
$X^\vee$ respectively, and there is a bijection
$\alpha\mapsto\alpha^\vee$ of $\Phi$ onto $\Phi^\vee$.

Adjoint action of $H$ on a Lie algebra $\mathfrak{g}={\rm Lie}(G)$
defines a decomposition
$$\mathfrak{g}=\mathfrak{h}\oplus
\sum_{\alpha\in\Phi}\mathbb{C}e_\alpha,\qquad \mathfrak{h}={\rm
Lie}(H),$$ and thus defines a subset $\Phi\subset X$. Let $B$ be Borel subgroup
containing $H$. There is a unique ordering $>$ of $\Phi$ such that
$\mathfrak{b}={\rm Lie}(B)$ is generated by $\mathfrak{h}={\rm
Lie}(H)$ and $e_\alpha$ with $\alpha>0$. One fixes a basis
$\Pi=\{\alpha_i\}$ of $\Phi$ compatible with  the ordering of $\Phi$
associated to $B$.

There is a decomposition $G=Z_0\cdot G'$ where $Z_0$ is the identity
component of the center $Z$ of $G$ and $G'$ is a semisimple group
(derived group of $G$). We have $H=Z_0\cdot H'$ where $H'$ is a
maximal torus of $G'$. The root datum associated with $(G',H')$ is
$(X/X_0,\Phi,Q^\vee,\Phi^\vee)$, with $Q\subset X/X_0$. Given a
basis $\{\alpha_i^\vee\}$, $i\in I$  in $Q^\vee$ and a basis
$\{\omega_j\}$ $j\in J$  in $X$, one can choose a basis of
representatives of the form $\{\omega'_i=\omega_i+X_0\}$, $i\in
I\subset J$  in $X/X_0$ such that $\{\omega'_i\}$, $i\in I$ form a
basis dual to $\{\alpha_i^{\vee}\}$, $i\in I$.

From now on if not explicitly mentioned $\mathfrak{g}$ be a
semisimple  Lie algebra. Let $\mathfrak{h}\subset \mathfrak{g}$ be a
Cartan subalgebra and $\mathfrak{b}_{\pm}$ be a pair of opposite
Borel subalgebras of $\mathfrak{g}$ containing $\mathfrak{h}$. We
have a decomposition
$\mathfrak{g}=\mathfrak{n}_{-}\oplus\mathfrak{h}\oplus\mathfrak{n}_+$
where $\mathfrak{n}_{\pm}\subset\mathfrak{b}_{\pm}$ is a nilpotent
radical. Denote by $\Gamma$ the set of vertexes of Dynkin graph
 associated with the root system of $\mathfrak{g}$. Let
$\Pi=\{\alpha_i\in\mathfrak{h}^{*},\,\,i\in\Gamma\}$ be the set of
simple roots, $\{\omega_i\in\mathfrak{h}^*,\,\,i\in\Gamma\}$ be the
set of fundamental weights and
$\Pi^{\vee}\{\alpha_i^{\vee}\in\mathfrak{h},\,\,i\in\Gamma\}$ with
be the set of the co-roots defined by
$\<\alpha_i^{\vee},\omega_j\>=\delta_{ij}$. Let
$A=\|a_{ij}\|,\,\,i,j=1,\ldots,\ell$ be the Cartan matrix of
$\mathfrak{g}$ defined by $a_{ij}=\<\alpha_i^{\vee},\alpha_j\>$.
Denote $R_+$  the set of positive roots of $\mathfrak{g}$ and let
$\rho$ be a half of the sum of the positive roots
$\rho=\frac{1}{2}\sum_{\alpha \in R_+}\alpha$.  There exist co-prime
positive rational numbers $d_1,\ldots,d_{\ell}$ such that the matrix
$\|b_{ij}\|=\|d_i a_{ij}\|$ is symmetric.  Define a symmetric
bilinear form on $\mathfrak{h}^{*}$ by $(\alpha_i,\alpha_j)=b_{ij}$.
This form defines a non-degenerate pairing
$\nu:\,\mathfrak{h}\widetilde{\rightarrow}\mathfrak{h}^*$ given by
$\nu(\alpha_i^\vee)=d_i^{-1}\alpha_i$.

Let $W$ be a  Weyl group of root system associated with
 Lie algebra $\mathfrak{g}$. It  is  generated by
simple reflections $s_1,\ldots,s_{\ell}$ acting by linear
transformations in $\mathfrak{h}^*$: \be
s_i(\lambda)=\lambda-\langle\lambda,\alpha_i^{\vee}\rangle\alpha_i,\,\,\,
\lambda\in\mathfrak{h}^{*}.\ee  Defining relations can be
represented as: \be
s_i^2=1,\,\,\,\,\,(s_is_j)^{m_{ij}}=1,\,\,\,i,j=1,\ldots,\ell, \ee
where $m_{ij}$ are equal to
$$m_{ij}=2,\\,3,\,\,4,\,\,6,\,\,\infty,$$
for
$$ a_{ij}a_{ji}=0,\,\,1,\,\,2,\,\,3,\,\,\geq 4, $$
 respectively.
For any $w\in W$ a reduced word is a sequence of indexes $I_w
=(i_1,\ldots,i_{l(w)})$, $i_k\in \Gamma$,  of shortest possible
length such that $w=s_{i_1}s_{i_2}\cdots s_{i_{l(w)}}$. The integer
$l(w)$ is called the length of $w$. Denote by $w_0$ the unique
element of maximal length in Weyl group and let $m=l(w_0)$. In the
following we fix a lift $\dot{w}\in G$, $\mathfrak{g}={\rm Lie}(G)$
of an element $w\in W$ such that $w(u)={\rm Ad}_{\dot{w}}u$,
$u\in\mathfrak{g}$. For simple reflections $s_i$ we define
$$\dot{s_i}=e^{e_i}e^{-f_i}e^{e_i},$$
and for $w=s_{i_1}s_{i_2}\cdots s_{i_{l(w)}}$ we take
$\dot{w}=\dot{s}_{i_1}\dot{s}_{i_2}\cdots \dot{s}_{i_{l(w)}}$. Thus
defined $\dot{w}$ does not depend on the decomposition into the
product of simple reflections (see e.g. \cite{K}, Lemma 3.8).

Denote by $e_i,f_i,h_i,\,\,i=1,\ldots,\ell$ the set of standard
generators of a semisimple Lie algebra $\mathfrak{g}$ satisfying the
following relations: \be [h_i,h_j]=0,\ee \be [h_i,e_j ] =a_{ij} e_j
,\,\,[h_i,f_j ] =-a_{ij}f_j ,\,\,
[e_i,f_j ] =\delta_{ij} h_i,\\
\ee \be ({\rm ad}\,\,\, e_i)^{1-a_{ij}}e_j=0, \,\,\,({\rm ad}\,\,\,
f_i)^{1-a_{ij}}f_j=0,\,\,
 for\,\,\, i\neq j. \ee
The invariant symmetric bilinear form on $\mathfrak{g}$ is  given by
$$(h_i,h_j)=b_{ij}d^{-1}_id^{-1}_j,\,\,\,
(e_i,f_j)=\delta_{ij}d^{-1}_i,\,\,\, (e_i,h_j)=(f_i,h_j)=0.$$

The only example of a non-semisimple reductive Lie algebra that will
be considered in this paper is the reductive Lie algebra
$\mathfrak{gl}_{\ell+1}$. In this case we explicitly define Lie
algebra by generators and relations as follows. Introduce the set of
generators $$\{e_{i,i\pm1},\,i=1,\ldots,\ell;\quad
e_{k,k},\,k=1,\ldots,\ell+1\},$$
 of $\mathfrak{gl}_{\ell+1}$. They satisfy the following
relations
 \bqa[e_{i,i},e_{j,j}]=0,\hspace{2cm}
[e_{i,i+1},e_{i+1,i}]=e_{i,i}-e_{i+1,i+1},\nonumber\\
\nonumber [e_{i,i},e_{i,i+1}]=e_{i,i+1},\hspace{2cm}
[e_{i+1,i+1},e_{i,i+1}]=-e_{i,i+1},\\ \,
[e_{i,i},e_{i+1,i}]=-e_{i+1,i},\hspace{2cm}
[e_{i+1,i+1},e_{i+1,i}]=e_{i+1,i}, \eqa \bqa ({\rm
ad}_{e_{i,i+1}})^{2}e_{j,j+1}=0, \qquad ({\rm
ad}_{e_{i+1,i}})^{2}e_{j+1,j}=0,\,\, \hspace{2cm}|i-j|=1.
\nonumber\eqa

\subsubsection{Whittaker model of principal series representations}

 Let $\mathcal{U}(\frak{g})$ be
a universal enveloping of $\mathfrak{g}$ and  $V$, $V'$  be
$\mathcal{U}(\frak{g})$-modules. Modules $V$ and $V'$ are called
dual if there exists a non-degenerate pairing $\<.\,,.\>: V'\times
V\to\CC$ such that $\<v',Xv\>=-\<X v',v\>$ for all $v\in V$, $v'\in
V'$ and $X\in\frak{g}$. We will  assume that the action of the
Cartan subalgebra on $V$, $V'$ is integrated to the action of the
Cartan torus.

Let $B_{-}=N_{-}H$ and $B_+=HN_+$ be a pair of opposed Borel
subgroups where $H$ is a maximal torus, and $N_{\pm}$ are opposite
maximal unipotent subgroups of $G$.  Characters of
$\mathfrak{n_{\pm}}={\rm Lie}(N_{\pm})$
 are defined by their values on simple root generators. Let
$\chi_{\pm}\!:\mathfrak{n}_{\pm}\!\rightarrow\CC$ be the characters
of $\mathfrak{n}_{\pm}$ defined by $\chi_+(e_{i}):=-1$ and
$\chi_-(f_{i}):=-1$ for all $i=1,\ldots,\ell$.  A vector $\psi_R\in
V$ is called a Whittaker vector with respect to
 $\chi_+$ if \be\label{ewt} \hspace{2cm} e_i\psi_R=-\psi_{R}\,,
 \hspace{1cm}(i=1,\ldots,\ell), \ee and a vector
$\psi_L\in V'$ is called a Whittaker vector with respect to $\chi_-$
if \be\label{fwt} \hspace{2cm} f_{i} \psi_L=-\psi_L\,, \hspace{1cm}
(i=1,\ldots,\ell). \ee

A Whittaker vector $\psi$ is called cyclic in  $V$ if
$\mathcal{U}(\mathfrak{g})\psi=V$,
 and  $\mathcal{U}(\mathfrak{g})$-module $V$ is a Whittaker
module if it contains a cyclic Whittaker vector. The Whittaker
$\mathcal{U}(\mathfrak{g})$-module $V$ admits an infinitesimal
character $\xi$ i.e. there exists a homomorphism of the center
$\mathcal{Z}(\mathfrak{g})\subset \mathcal{U}(\mathfrak{g})$
$\xi:\mathcal{Z}(\mathfrak{g})\rightarrow\CC$ such that $zv=\xi(z)v$
for all $z\in \mathcal{Z}(\mathfrak{g})$ and $v\in V$.

Consider the principal series representation ${\rm
Ind}_{B_-}^G\,\chi_{\mu}$ of $G$, induced from the character
$\chi_{\mu}$ of $B_-=HN_-$ trivial on $N_-$.
 It is realized in the space of functions $f\in
L^2(G)$ satisfying 
\bqa\label{twistf}
f(bg)=\chi_{\mu}(b)f(g).
\eqa
The action of $G$
 is given by the right action $(g_1\cdot f)(g_2)=f(g_2g_1)$.
 We will be interested in the infinitesimal  form
${\rm Ind}_{U(\mathfrak{b}_-)}^{U(\mathfrak{g})}\,\chi_{\mu}$ of
this representation.  The action of the Lie algebra
$\mathfrak{g}={\rm Lie}(G)$ is given by  the infinitesimal form  of
the right action \bqa
\label{infact}(Xf)(g)=\frac{d}{d\epsilon}f(ge^{\epsilon
X})|_{\epsilon\rightarrow 0}\,\,. \eqa Denote $V_{\mu}$ the
corresponding $U(\mathfrak{g})$-module.

Let $G(\mathbb{R})$ be a totally split real form of a reductive Lie
group $G$, $\mathfrak{g}_{\mathbb{R}}={\rm Lie}(G(\mathbb{R}))$ be a
corresponding Lie algebra and $N(\mathbb{R})_+\subset N_+ \cap
G(\mathbb{R})$ be a nilpotent subalgebra of $G(\mathbb{R})$. Let
$d\mu_{G(\mathbb{R})}$ be a bi-invariant (Haar) measure on
$G(\mathbb{R})$.  We have the Bruhat decomposition
$G(\mathbb{R})=\coprod_{w\in
  W}\,B_-(\mathbb{R})wB_+(\mathbb{R})$.
 Let $G_0(\mathbb{R})=B_-(\mathbb{R})N_+(\mathbb{R})$
be a $w=1$ component  in this decomposition. Restriction of the
measure $d\mu_{G(\mathbb{R})}$ on $G_0(\mathbb{R})$
 up to normalization has the following form
\cite{He}: \bqa\label{factmeas}
 d\mu_{G(\mathbb{R})}(g)=\delta_{B_+(\mathbb{R})}(b)\,
d\mu_{B_+(\mathbb{R})}(b)\wedge d\mu_{N_+(\mathbb{R})}(x). \eqa
 Here $\delta_{B_+(\mathbb{R})}$ is the
 modular function on $B_+(\mathbb{R})$. For any $b=ng_0\in
N_+(\mathbb{R})\,H$ it is equal  to
$\delta_{B_+(\mathbb{R})}(b)=\exp\,2 \<\rho,\ln g_0\>$.

 Let $\mu=i\lambda-\rho$. Consider the following non-degenerate pairing
$V_{\mu}\times V_{\mu} \rightarrow\mathbb{C}$:
$$
\langle f_1,f_2\rangle= \int_{N_+(\mathbb{R})}d\mu_{N_+(\mathbb{R})}
(x)\, f_1(x)\,\overline{f_2(x)},
$$
where $d\mu_{N_+(\mathbb{R})}$ is a restriction of (\ref{factmeas})
on $N_+(\mathbb{R})$. It defines  on $V_{\mu}$ a  structure of a
unitary representation $\pi_{\lambda}$ of
$U(\mathfrak{g}_{\mathbb{R}})$
 and we have $\langle f_1,X f_2\rangle=-\langle Xf_1,f_2\rangle$ for any
$X\in \mathfrak{g}_{\mathbb{R}}$.

We shall consider a slightly more general pairing defined as
follows. Note that $N_+(\mathbb{R})\subset N_+$ is a real
non-compact  middle dimension subspace. One has a natural
holomorphic structure on a Lie algebra
$\mathfrak{g}=\mathfrak{g}_{\mathbb{R}}\otimes \mathbb{C}$ which
induces the holomorphic structure on $N_+$.
 Consider the space of  holomorphic functions on $N_+$.
It is a module with respect to the holomorphic action of a
corresponding holomorphic subalgebra of $\mathfrak{g}$.  The
right-invariant measure $d\mu_{N_+(\mathbb{R})}$ can be extended to
a holomorphic  top-dimensional form $d\mu_{N_+}^{hol}$ on $N_+$. Let
$C\subset N_+$ be an arbitrary non-compact middle-dimensional
submanifold. Consider the following pairing
$$
\langle f_1,f_2\rangle_C= \int_{C}d\mu^{hol}_{N_+} (x)\,
f_1(x)\,\overline{f_2(\overline{x})},
$$
on the space $\mathcal{S}_C^{hol}(N_+)$ of holomorphic functions on
$N_+$ exponentially decreasing with all its derivatives when
restricted to  $C$. This pairing satisfies
 $\langle f_1,X f_2\rangle_C=-\langle Xf_1,f_2\rangle_C$ for any
holomorphic  $X \in \mathfrak{g}$.

\subsubsection{Whittaker function as Toda wave function}

According to B.~Kostant \cite{Ko1},\cite{Ko2} eigenfunctions of
$\mathfrak{g}$-Toda chain can  be written in terms of the invariant
pairing on Whittaker modules as follows
 \bqa\label{pairing} \Psi^{\mathfrak{g}}_{\lambda} (x)=e^{-
\<\rho,x\>}\<\psi_L\,, \pi_{\lambda}(e^{-h_x}) \,\psi_R\>,\, \qquad
\qquad x\in \mathfrak{h}, \eqa  where
$h_x:=\sum\limits_{i=1}^{\ell}\<\omega_i,x\>\,h_i$.  In the special
case of $\mathfrak{g}=\mathfrak{sl}_2$ the function
$\Psi^{\mathfrak{g}}_{\lambda}(x)$, $x\in \mathbb{R}$ coincides with
the classical Whittaker function. In the following we will use the
term $\mathfrak{g}$-Whittaker function  for (\ref{pairing}) ( see
e.g. \cite{Et}).  A slightly different  notion of the 
Whittaker functions was used in  \cite{Ja}, \cite{Ha}.

One can  introduce a set of commuting differential operators
$\CH_k\in{\rm Diff} (\mathfrak{h})$, $k=1,\cdots, \ell$
corresponding to a set $\{c_k\}$ of  generators of the center
$\mathcal{Z}\subset \mathcal{U}(\mathfrak{g})$  as follows: \bqa
\label{hamdef} \CH_k
\Psi^{\mathfrak{g}}_{\lambda}(x)=e^{-\<\rho,x\>}
\<\psi_L\,,\pi_{\lambda}(e^{-h_x})\,c_k\,\psi_R\>.\qquad \eqa

Operators $\CH_k$  provide a complete set of commuting Hamiltonians
of $\mathfrak{g}$-Toda chain \cite{Ko1}. Connection with  Toda
chains can be seen as follows. Quadratic generator of the center
$\mathcal{Z}(\mathfrak{g})$ (Casimir element) is given by
\bqa\label{seccas}
c_2=\frac{1}{2}\sum\limits_{i,j=1}^{\ell}c_{ij}h_ih_j+\frac{1}{2}
\sum\limits_{\alpha\in
R_+}(e_{\alpha}f_{\alpha}+f_{\alpha}e_{\alpha}),\eqa where the
matrix $\|c_{ij}\|=\|d_id_j(b^{-1})_{ij}\|$ is  inverse to the
matrix $\|(\alpha_i^{\vee},\alpha_j^{\vee})\|$. Let $\{\epsilon_i\}$
be an orthogonal bases $(\epsilon_i,\epsilon_j)=\delta_{ij}$ in
$\mathfrak{h}$ and $x=\sum_{i=1}^{\ell} x_i\epsilon_i$ be  a
decomposition of $x\in \mathfrak{h}$ in this bases. Then the
projection  (\ref{hamdef}) of (\ref{seccas}) gives the well-known
Hamiltonian operator of $\mathfrak{g}$-Toda chain \cite{STS}
\be\CH^{\mathfrak{g}}_2=-\frac{1}{2}\sum_{i=1}^{\ell}
\frac{\partial^2}{{\partial x_i}^2} +\sum_{i=1}^\ell d_i
e^{\<\alpha_i,x\>}.\ee

The eigenfunctions (\ref{pairing}) of $\mathfrak{g}$-Toda chain are
written in an abstract form. To get explicit integral
representations we start with representations of  matrix elements
(\ref{pairing}) of  infinite-dimensional representations in terms of
matrix elements of finite-dimensional representations of
$\mathcal{U}(\mathfrak{g})$. Let $\pi_i$ be a set of fundamental
representations corresponding to all  fundamental weights $\omega_i$
of $\mathfrak{g}$ and $\xi^{+/-}_{\omega_i}$ be  highest/lowest
vectors in these representations such that
$\langle\xi_{\omega_i}^-|\xi_{\omega_i}^+\rangle =1$. For highest
weight vector $\xi^{+}_{\omega_i}$ in a fundamental representation
$V_{\omega_i}$ we have
$\dot{s}_i^{-1}\xi^{+}_{\omega_i}=f_i\xi^{+}_{\omega_i}$. Consider
following  matrix elements in fundamental finite-dimensional
representations \be
\Delta_{\omega_i,\dot{w}}(g)=
\langle\xi_{\omega_i}^{-}|\,\pi_i(g)\pi_i(\dot{w})
|\,\xi_{\omega_i}^{+}\rangle,\qquad w\in W, \,\,\,g\in G.  \ee

\begin{lem}\label{Whittakervectlem}
The left/right Whittaker vectors defined by (\ref{ewt}) and (\ref{fwt})
are given by: \bqa\label{psir}
\psi_R(v)=\exp\Big\{-\sum\limits_{i=1}^{\ell}\frac{\Delta_{\omega_i,
\dot{s}_{i}^{-1}}(v)}{ \Delta_{\omega_i,1}(v)}\Big\},\eqa
\bqa\label{psil} \psi_L(v)=\prod_{i=1}^{\ell}
(\Delta_{\omega_i,\dot{w}_0^{-1}}(v))^{\imath\<\lambda,\alpha_i^{\vee}\>-1}
\times \exp\Big\{\sum\limits_{i=1}^{\ell}
\frac{\Delta_{\omega_i,\dot{w}_0^{-1}\dot{s}_i^{-1}}(v)}
{\Delta_{\omega_i,\dot{w}_0^{-1}}(v)}\Big\},
 \eqa
\end{lem}
The proof in given in Part II,  Section \ref{Witvec}.

\begin{prop}\label{propwhitgen}
Common eigenfunctions (\ref{pairing}) of $\mathfrak{g}$-Toda chain
can be represented in the following integral form: \bqa\label{int}
\Psi^{\mathfrak{g}}_{\lambda}(x)=e^{\imath \<\lambda,x\>}\,
\int_{C}d\mu^{hol}_{N_+}(v)\prod\limits_{i=1}^{\ell}
(\Delta_{\omega_i,\dot{w}_0^{-1}}(v))^{\imath\<\lambda,\alpha_i^{\vee}\>-1}\times\\
\times \exp\Big\{\sum\limits_{i=1}^{\ell}
\left(\frac{\Delta_{\omega_i,\dot{w}_0^{-1}\dot{s}_i^{-1}}(v)}
{\Delta_{\omega_i,\dot{w}_0^{-1}}(v)}-e^{\<\alpha_i,x\>}
\frac{\Delta_{\omega_i,\dot{s}_{i}^{-1}}(v)}{
\Delta_{\omega_i,1}(v)}\right)\Big\}.\nonumber \eqa Here  $C \subset
N_+$ is a  middle-dimensional non-compact cycle such that the
integrand decreases exponentially  at the boundaries and infinities.
The measure of the integration  is the restriction on $C$ of the
holomorphic continuation $d\mu^{hol}_{N_+}$ of the right-invariant
measure $d\mu_{N_+(\mathbb{R})}$ on $N_+(\mathbb{R})$.
\end{prop}

 The first example of this type of integral representation for
$\mathfrak{gl}_n$-Whittaker function was considered in
\cite{GKMMMO}. Its generalization given above is straightforward.
The proof of the Proposition is given in Part II, Section \ref{Witvec}.

The expression (\ref{int}) for  a Whittaker function is much more
detailed then (\ref{pairing}) but does not yet provide explicit
integral representation. To obtain explicit integral representations
of Whittaker functions one should choose a parameterization of $N_+$
(or an open part of it) and express the measure $d\mu^{hol}_{N_+}$
and  various matrix elements entering (\ref{int}) in terms of the
coordinates on $N_+$. Natural choice would be a  factorized
representation of the elements of an  open part of a maximal
unipotent subgroup of an arbitrary Lie group \cite{Lu} (see also
\cite{BZ}, \cite{FZ}). For each $i=1,\ldots, \ell$  let
$X_i(t)=\exp\{te_i\}$
 be a one-parameter subgroup in $N_+$. Pick a decomposition
of the longest element $w_0$ in the Weyl group $W$ corresponding to
a reduced word $I_{w_0}=(i_1,\ldots,i_m)$, $l(w_0)=m={\rm
dim}\,N_+$.
 Then the following map  \bqa\label{factparam}
 \mathbb{C}^m\,\longrightarrow\, N_+^{(0)},
\hspace{16mm} (t_1,\ldots,t_m)\,\longmapsto\,v(t_1,\cdots,
t_m)=X_{i_1}(t_1)\cdots X_{i_m}(t_m),\eqa is a birational
isomorphism. This provides a parametrization of an open part
$N_+^{(0)}$ of $N_+$.  Parametrizations corresponding to different
choices of the reduced word $I_{w_0}$ are related by birational
transformations described explicitly by G.~Lustzig \cite{Lu}. The
right-invariant measure has the following description in the
factorized representation.
\begin{lem}\label{lemmeasure}
The right-invariant measure $d\mu^{hol}_{N_+}$ in the factorized
parametrization is given by:
 \be\label{meas} d\mu^{hol}_{N_+}(v)=\prod\limits_{i=1}^{\ell}
\Delta_{\omega_i,\dot{w}_0^{-1}}(v)\bigwedge\limits_{k=1}^{m}
\frac{d t_k}{t_k}. \ee
\end{lem}
The proof is given in Part II, Section \ref{measureN}.

Thus the problem of finding explicit integral representations of
Whittaker functions in the factorized parametrization
(\ref{factparam}) is reduced to a calculation of the matrix elements
of finite-dimensional representations  of $\mathfrak{g}$ in this
parametrization. In the following we provide explicit expressions
for finite-dimensional matrix elements for classical Lie groups and
give corresponding integral representations  of Whittaker functions.
Let us stress however that thus obtained integral representation for
$\mathfrak{g}=\mathfrak{gl}_{\ell+1}$ does not coincide with
Givental representation \cite{Gi}. Note that for
 classical series  of Lie algebras the factorized
parametrization (\ref{factparam}) has a recursive structure over the
rank $\ell$  reflecting the recursive structure of the reduced
decomposition of $w_0\in W$. This recursive structure is not
translated, however, into a simple recursive structure of the
infinite-dimensional matrix element in the factorized parametrization
and does not reproduce the recursive structure of the Givental
integral representation.

In \cite{GKLO} a modification of the factorized parametrization
(\ref{factparam})  for $\mathfrak{g}=\mathfrak{gl}_{\ell+1}$ was
proposed and it was shown that the integral representation
(\ref{int}) in this parametrization exactly reproduces the Givental
integral representation of $\mathfrak{gl}_{\ell+1}$-Whittaker
functions. In particular for this parametrization the recursive
structure of the reduced decomposition of $w_0\in W$ directly
translates into the recursive structure of the integral
representation of the corresponding Whittaker function.

Below we generalize the results of \cite{GKLO} to
 all classical series of Lie algebras.  We propose a modification of factorized
parametrization (\ref{factparam}) based on a particular realization
of  maximal unipotent subgroups $N_+\subset G$ of  classical Lie
groups as explicitly defined  subgroups of the maximal unipotent
subgroups of general linear groups.  For any classical simple Lie
group, the maximal unipotent subgroup  can be realized as a subgroup
of a group  of upper-triangular matrices of appropriate size with
units on diagonal (see e.g. \cite{DS}). The corresponding subset of
upper-triangular matrices for  classical Lie group can be describe
explicitly. We define a parametrization of maximal unipotent
subgroups of classical Lie groups by constructing a particular form
of  the parametrization  of  the corresponding subset of
upper-triangular matrices. Using this  parametrization we derive
explicit integral representations of Whittaker functions associated
with all classical groups and demonstrate that these integral
representations have all characteristic properties of the Givental
integral representation for $\mathfrak{gl}_{\ell+1}$-Whittaker
functions. In particular the recursive structure of Whittaker
functions  is  explicit in  this new parametrization.

\subsection{ Integral representations of
  $\mathfrak{gl}_{\ell+1}$- and $ \mathfrak{sl}_{\ell+1}$-Toda chain
  eigenfunctions}

In this section we recall the construction of integral
representations of $\mathfrak{gl}_{\ell+1}$- and
$\mathfrak{sl}_{\ell+1}$-Toda eigenfuctions  using  factorized
parametrization (\ref{factparam}) of a maximal unipotent subgroup
$N_+\subset GL(\ell+1)$ and its modification introduced in
\cite{GKLO}. The second parametrization leads to an integral
representation obtained earlier by Givental \cite{Gi} using
different approach. In the following these constructions will be
generalized  to $\mathfrak{g}$-Toda theory for arbitrary classical
Lie algebras $\mathfrak{g}$.

We start with the case of the reductive Lie algebra
$\mathfrak{gl}_{\ell+1}$. Let $(\epsilon_1,\ldots
,\epsilon_{\ell+1})$ be an orthogonal basis in $\RR^{\ell+1}$,
$(\epsilon_i,\epsilon_j)=\delta_{ij}$. Roots and fundamental wights
of $\mathfrak{gl}_{\ell+1}$ considered as vectors in $\RR^{\ell+1}$
are given by:
 \bqa\label{RootDataGLn}
\alpha_i=\epsilon_{i+1}-\epsilon_i,\,\,\,i=1,\ldots,\ell,\qquad
\omega_j=\epsilon_j,\,\,\,\,j=1,\ldots,(\ell+1). \eqa Coroots
$\alpha_i^{\vee}$ can be identified with the corresponding roots
$\alpha_i$ with respect to the pairing in $\RR^{\ell+1}$.
 To this root/weight system one associates $\mathfrak{gl}_{\ell+1}$-Toda
quantum integrable system having a set of $(\ell+1)$ mutually
commuting functionally independent quantum Hamiltonians
$H^{\mathfrak{gl}_{\ell+1}}_k$, $k=1,\cdots, (\ell+1)$. We are
interested in the explicit integral representations for common
eigenfunctions of the full set of quantum Hamiltonian operators for
$\mathfrak{gl}_{\ell+1}$. For instance linear and quadratic quantum
Hamiltonians of $\mathfrak{gl}_{\ell+1}$-Toda chain are  given by
\bqa \CH_1^{\mathfrak{gl}_{\ell+1}}&=&-\i\sum\limits_{i=1}^{\ell+1}
\frac{\partial}{\partial x_i},\\
\CH_2^{\mathfrak{gl}_{\ell+1}}&=&-\frac{1}{2}\sum\limits_{i=1}^{\ell+1}
\frac{\partial^2}{{\partial x_i}^2}+ \sum\limits_{i=1}^{\ell}
e^{x_{i+1}-x_{i}},  \eqa
 and the  eigenfunction should satisfy the
following equation \bqa
\CH_1^{\mathfrak{gl}_{\ell+1}}(x)\,\,\,\Psi^{\mathfrak{gl}_{\ell+1}}_{\la_1,\cdots
,\la_{\ell+1}} (x_1,\ldots,x_{{\ell}+1})&=&
\sum\limits_{i=1}^{\ell+1}\lambda_{i}\,\,\,
\Psi^{\mathfrak{gl}_{\ell+1}}_{\la_1,\cdots ,\la_{\ell+1}}
(x_1,\ldots,x_{\ell+1}),\\
\CH_2^{\mathfrak{gl}_{\ell+1}}(x)\,\,\,\Psi^{\mathfrak{gl}_{\ell+1}}_{\la_1,\cdots
,\la_{\ell+1}} (x_1,\ldots,x_{{\ell}+1})&=&
\frac{1}{2}\sum\limits_{i=1}^{\ell+1}\lambda_{i}^2\,\,\,
\Psi^{\mathfrak{gl}_{\ell+1}}_{\la_1,\cdots ,\la_{\ell+1}}
(x_1,\ldots,x_{\ell+1}). \eqa

Common eigenfunction of the quantum Hamiltonians has the
 following representation as a matrix element \bqa\label{pairgl}
\Psi_{\lambda}^{\mathfrak{gl}_{\ell+1}}(x)=e^{- \sum
x_i\rho_i}\<\psi_L\,, \pi_{\lambda}(e^{-\sum x_i E_{i,i}})
\,\psi_R\>,\eqa where $\rho_i=\frac{1}{2}(\ell-2i+2)$ are the
components of $\rho$ in the standard basis $\{\epsilon_i\}$ in
$\RR^{\ell+1}.$

The construction for the  semisimple Lie algebra
$\mathfrak{sl}_{\ell+1}$ is  quite similar to that for reductive Lie
algebra $\mathfrak{gl}_{\ell+1}$. The roots and fundamental wights
for semisimple Lie algebra $\mathfrak{sl}_{\ell+1}$ can be written
in the following form  (see \cite{Bou}): \bqa\label{weightsSL}
\alpha_i=\epsilon_{i+1}-\epsilon_i, \hspace{2cm}
\omega_i=-(\epsilon_1+\ldots+\epsilon_i)+
\frac{i}{\ell+1}(\epsilon_1+\ldots+ \epsilon_{\ell+1}), \eqa for
$i=1,\ldots,\ell$. This representation of the $A_{\ell}$ root/weight
system can be obtained from the root/weight system of the reductive
Lie algebra $\mathfrak{gl}_{\ell+1}$ as follows. Let us pick an
orthogonal basis of fundamental weights of $\mathfrak{gl}_{\ell+1}$:
$$\omega'_i=-\epsilon_1-\ldots-\epsilon_i, $$ such that
$\langle\omega'_i,\,\alpha_j^\vee\rangle=\delta_{ij}$ for
$i,j=1,\ldots \ell$, and
$\langle\omega'_{\ell+1},\,\alpha_j^\vee\rangle=0$ for $j=1,\ldots
\ell$. Then $\omega'_{\ell+1}$ can be identified as a
 generator of $X_0$. Introducing
$$\omega_i=\omega'_i-\frac{i}{\ell+1}\omega'_{\ell+1},$$ one readily
obtains the set (\ref{weightsSL}) of fundamental weights for
$\mathfrak{sl}_{\ell+1}$.

To this root/weight system one  associates
$\mathfrak{sl}_{\ell+1}$-Toda quantum integrable system possessing a
set of $\ell$ mutually commuting functionally independent
Hamiltonians $\CH^{\mathfrak{sl}_{\ell+1}}_k$, $k=1,\cdots, \ell$.
It is convenient however to consider $\mathfrak{sl}_{\ell+1}$-Toda
chain Hamiltonians as a subset $\CH^{\mathfrak{gl}_{\ell+1}}_k$,
$k=2,\cdots, (\ell+1)$ of $\mathfrak{gl}_{\ell+1}$-Toda chain
Hamiltonians acting on the kernel of the linear Hamiltonian
$\CH_1^{\mathfrak{gl}_{\ell+1}}$. For instance the eigenfunction of
a quadratic quantum Hamiltonian of $\mathfrak{sl}_{\ell+1}$-Toda
chain should satisfy the equation \bqa
\CH_2^{\mathfrak{sl}_{\ell+1}}\Psi^{\mathfrak{sl}_{\ell+1}}_{\la_1,\cdots
,\la_{\ell+1}} (x_1,\ldots,x_{{\ell}+1})=\nonumber \eqa \bqa
=(-\frac{1}{2}\sum\limits_{i=1}^{\ell+1} \frac{\partial^2}{{\partial
x_i}^2}+ \sum\limits_{i=1}^{\ell}
e^{x_{i+1}-x_{i}}))\Psi^{\mathfrak{sl}_{\ell+1}}_{\la_1,\cdots
,\la_{\ell+1}} (x_1,\ldots,x_{{\ell}+1})=\\ \nonumber
=\frac{1}{2}\sum\limits_{i=1}^{\ell+1}\lambda_{i}^2\,\,
\Psi^{\mathfrak{sl}_{\ell+1}}_{\la_1,\cdots ,\la_{\ell+1}}
(x_1,\ldots,x_{{\ell}+1}), \eqa with an additional constraint
$\la_1+\ldots+\la_{\ell+1}=0$. The eigenfunctions for
$\mathfrak{sl}_{\ell+1}$-Toda chain  can be also written using a
reduced set of variable \bqa \label{slgl}
\Psi_{\nu_1,\ldots,\nu_\ell}^{\mathfrak{sl}_{\ell+1}}
(y_1,\ldots,y_\ell)\equiv
\Psi_{\lambda_1,\ldots,\lambda_{\ell+1}}^{\mathfrak{sl}_{\ell+1}}
(x_1,\ldots,x_{\ell+1})\qquad \nu_i=\lambda_{i+1}-\lambda_i,\,\,\,\,
y_i=x_{i+1}-x_i \eqa Note that without imposing the constraint
$\la_1+\ldots+\la_{\ell+1}=0$, the eigenfunctions of
$\mathfrak{sl}_{\ell+1}$-Toda chain  can be expressed through
eigenfunctions  of $\mathfrak{gl}_{\ell+1}$-Toda theory in the
following simple way \bqa\label{psisl}
\Psi_{\nu_1,\ldots,\nu_\ell}^{\mathfrak{sl}_{\ell+1}}
(y_1,\ldots,y_\ell)=\exp\Big\{-\frac{\imath}{\ell+1}\sum_{i=1}^{\ell+1}\lambda_i\cdot
\sum_{i=1}^{\ell+1}x_i\Big\}\cdot
\Psi_{\lambda_1,\ldots,\lambda_{\ell+1}}^{\mathfrak{gl}_{\ell+1}}
(x_1,\ldots,x_{\ell+1}), \eqa where we use notations (\ref{slgl}).
In the following we will  consider  mostly
$\mathfrak{gl}_{\ell+1}$-Toda chain eigenfunctions making comments
on the corresponding modifications for $\mathfrak{sl}_{\ell+1}$ case
(we  will mostly use the non-reduced form
$\Psi_{\lambda}^{\mathfrak{sl}_{\ell+1}}(x)$).

\subsubsection{ $\mathfrak{gl}_{\ell+1}$-Whittaker function:
 factorized   parametrization}

To make the integral representation (\ref{int}) for
$\mathfrak{gl}_{\ell+1}$-Whitaker functions explicit one should pick
a particular parametrization of $N_+\subset GL(\ell+1)$. Let $w_0$
be an element of maximal length  of the Weyl group $W=S_{\ell+1}$ of
$\mathfrak{gl}_{\ell+1}$. Consider the  reduced decomposition of
$w_0$ corresponding to the following  reduced word $I_{\ell}$
$$I_{\ell}=(i_1,i_2,\ldots,i_m):=(1,21,321,\ldots,(\ell\ldots 321)).$$
The reduced word $I_\ell$ has an obvious  recursive structure:
$I_{\ell+1}=I_{\ell}\sqcup(\ell+1\ldots 321)$.
 Thus the corresponding parametrization of  unipotent elements $v^{(\ell)}$
in an open part $N_+^{(0)}$ of $N_+$ can be written in a recursive
form: \bqa\label{rec1A}
v^{(\ell)}=v^{(\ell-1)}\cdot\mathfrak{X}^{\ell}_{\ell-1},\eqa where
\bqa\label{rec2A}
\mathfrak{X}^{A_\ell}_{A_{\ell-1}}=X_{\ell}(y_{\ell,1})\cdots
X_2(y_{2,\ell-1})X_1(y_{1,\ell}),\eqa and   $X_i(y)=\exp(ye_i)$.
Parameters $y_{ik}$ of one-parametric subgroups will be called
factorization parameters. The action of Lie algebra
$\mathfrak{gl}_{\ell+1}$ on $G/B_-$ considered at the beginning of
the previous section defines an action of the Lie algebra on the
space of functions $V_{\mu}$ restricted to $N_+^{(0)}$.

\begin{prop}\label{propgengl}
The following differential operators define a realization of the
representation $\pi_{\lambda}$ of
$\mathcal{U}(\mathfrak{gl}_{\ell+1})$ in $V_{\mu}$ in terms of
factorized  parametrization (\ref{rec1A}), (\ref{rec2A}): \bqa
E_{i,i}&=&\mu_i-\sum_{l=1}^{\ell+1-i}
y_{i,l}\frac{\partial}{\partial y_{i,l}}+ \sum_{l=1}^{\ell+2-i}
y_{i-1,l}\frac{\partial}{\partial y_{i-1,l}},\nonumber \\
E_{i,i+1}&=&\sum_{k=0}^{i-1}\prod_{s=0}^k
\frac{y_{i-s,\ell+2-i}}{y_{i+1-s,\ell+1-i}} \frac{\partial}{\partial
y_{i-k,\ell+1-i}}-
\prod_{s=0}^k\frac{y_{i-(s+1),\ell+2-i}}{y_{i-s,\ell+1-i}}
\frac{\partial}{\partial y_{i-(k+1),\ell+2-i}},\eqa \bqa
E_{i+1,i}&=&\sum_{k=1}^{\ell}\left[(\mu_{i+1}-\mu_i)
y_{i,k+1-i}-y_{i,k+1-i}\Big(y_{i,k+1-i}\frac{\partial}{\partial
y_{i,k+1-i}}-
y_{i+1,k-i}\frac{\partial}{\partial y_{i+1,k-i}}\Big)\right.+\nonumber \\
 &+& \left. y_{i,k+1-i}\sum_{s=1}^{k-1}\Big(
y_{i-1,s+2-i}\frac{\partial}{\partial y_{i-1,s+2-i}}-
2y_{i,s+1-i}\frac{\partial}{\partial y_{i,s+1-i}}+
y_{i+1,s-i}\frac{\partial}{\partial
y_{i+1,s-i}}\Big)\right],\nonumber \eqa where
$E_{i,j}=\pi_{\lambda}(e_{i,j})$,  $\mu_k=\imath\lambda_k-\rho_k$
and $\rho_k=\frac{1}{2}(\ell-2k+2).$
\end{prop}
{\it Proof.} The proof is given in Part II, Section \ref{gengl}.

 The calculation of  matrix elements
entering the integral (\ref{int}) in the factorized parametrization
(\ref{rec1A}), (\ref{rec2A}) can be done following  \cite{BZ} and
\cite{FZ} (see Section 3.3 for details). Another,  more
straightforward approach is to find left and right Whittaker vectors
solving the equations (\ref{ewt})-(\ref{fwt}) directly. In the
following we will use the  convention: $\sum_{i=k}^{j}=0,$ when
$k>j$ and $\prod_{i=k}^{j}=1,$ when $k>j$.
\begin{lem}\label{Whitfacrep}
 The following expressions for the left/right
Whittaker vectors in terms of factorization parameters hold:
\bqa\psi_R(y)=\exp\Big\{-\sum_{i=1}^\ell
\sum\limits_{n=1}^{\ell+1-i}y_{i,n}\Big\},\eqa
\bqa\psi_L(y)=\prod_{i=1}^\ell\Big(
\prod\limits_{k=1}^{i}\prod_{n=i+1-k}^{\ell}
y_{k,n}\Big)^{(\mu_{i+1}-\mu_i)}\times\\ \nonumber \times
\exp\Big\{-\sum_{k=1}^\ell\frac{1}{y_{\ell+1-k,k}}\Big(1+
\sum_{n=1}^{\ell-k}\prod_{i=1}^n
\frac{y_{\ell+1-k-i,k+1}}{y_{\ell+1-k-i,k}} \Big)\Big\}.\eqa
\end{lem}

Using (\ref{meas}) we have the following expression for
$\mathfrak{gl}_{\ell+1}$-Whittaker function in the factorized
parametrization.

\begin{te} \label{teintfgl}
  Eigenfunctions of the $\mathfrak{gl}_{\ell+1}$-Toda
chain (\ref{pairgl}) admit the integral representation:
\bqa\label{intlugl}
\Psi_{\lambda_1,\ldots,\lambda_{\ell+1}}^{\mathfrak{gl}_{\ell+1}}
(x_1,\ldots,x_{\ell+1})=e^{\imath\sum_{k=1}^{\ell+1}\lambda_kx_k}
\int\limits_C\,\,\bigwedge_{i=1}^\ell\bigwedge_{n=1}^{\ell+1-i}
\frac{dy_{i,n}}{y_{i,n}}\prod_{i=1}^\ell\Big(
\prod\limits_{k=1}^{i}\prod_{n=i+1-k}^{\ell}
y_{k,n}\Big)^{\imath(\lambda_{i+1}-\lambda_i)} \nonumber\\
\exp\Big\{ -\Big(\sum_{k=1}^\ell\frac{1}{y_{\ell+1-k,k}}\Big(1+
\sum_{n=1}^{\ell-k}\prod_{i=1}^n\frac{y_{\ell+1-k-i,k+1}}{y_{\ell+1-k-i,k}}
\Big)\,\,+\,\,\sum_{i=1}^\ell\, e^{x_{i+1}-x_{i}}
\sum\limits_{n=1}^{\ell+1-i}y_{i,n}\Big)\,\,\Big\}.\eqa Here $C
\subset N_+$ is a  middle-dimensional non-compact submanifold such
that the integrand decreases exponentially  at the boundaries and
infinities. In particular one can take $C=\RR_{+}^{\ell(\ell+1)/2}$.
\end{te}
The proof is given in  Part II, Section \ref{melgl}.

\subsubsection{ $\mathfrak{gl}_{\ell+1}$-Whittaker function:
modified factorized   parametrization}

Now we consider a modification of the factorized  parametrization
(\ref{rec1A}), (\ref{rec2A})  leading  to the  Givental integral
representation of  $\mathfrak{gl}_{\ell+1}$-Whittaker function.
This modified factorized parametrization was first introduced in
\cite{GKLO}. There is an important difference between factorized and
modified factorized parametrizations. Note that the parametrization
(\ref{rec1A}), (\ref{rec2A}) is defined in terms of group elements
of $N_+$. To define a modified factorized parametrization of $N_+$
we shall consider the image of a group element in a faithful
finite-dimensional representation of $G$. In the case of
$\mathfrak{gl}_{\ell+1}$ and $\mathfrak{sl}_{\ell+1}$ we use a
tautological  representation $\pi_{\ell+1}:
\mathfrak{gl}_{\ell+1}\to End(\mathbb{C}^{\ell+1})$. Let
$\epsilon_{i,j}$ be a set of  elementary
$(\ell+1)\times(\ell+1)$-matrices with units at $(i,j)$-places and
zeros, otherwise. Consider the following set of diagonal matrices
$$U_k=\sum\limits_{i=1}^{k}e^{-x_{k,i}}\epsilon_{i,i}+\sum\limits_{i=k+1}^{N}
\epsilon_{i,i}.$$
 Define the following upper-triangular
 deformation of $U_k$
\be\label{deformed} \tilde{U}_k=\sum\limits_{i=1}^{k}
e^{-x_{k,i}}\epsilon_{i,i} + \sum\limits_{i=k+1}^{N} \epsilon_{i,i}
+ \sum\limits_{i=1}^{k-1} e^{-x_{k-1,i}}\epsilon_{i,i+1}.\ee

The modified factorized parametrization of $N_+$ is then defined as
follows.

\begin{te}\label{Mpargl}

i)  The image of any generic unipotent element
  $v\in N_+$ in the tautological representation
$\pi_{\ell+1}:\mathfrak{gl}_{\ell+1}\to End(\mathbb{C}_{\ell+1})$
 can be represented in the form
 \be\label{givpar}
\pi_{\ell+1}(v)=\tilde{U}_{2}U_{2}^{-1}\tilde{U}_{3}U_{3}^{-1}\cdots
\tilde{U}_{N-1}U_{N-1}^{-1}\tilde{U}_{N}, \ee where we assume that
$x_{\ell+1,i}=0,\,\,\,\,i=1,\ldots,\ell+1$.

ii) This defines a parametrization of an open part $N_+^{(0)}$ of
$N_+$.
\end{te}

{\it Proof.} Let  $v(y)$ be elements of $N_+$  parametrized
according to (\ref{rec1A}), (\ref{rec2A}). Let us now change the
variables in the following way: \be\label{anzacgl}
y_{i,n}=e^{x_{n+i,i+1}-x_{n+i-1,i}},\ee where $x_{\ell+1,n}=0,\,\,\,
n=1,\ldots,\ell+1$ are assumed. By elementary operations it is easy
to check that after the change of variables, the image
$\pi_{\ell+1}(v)$ of $v$ defined by (\ref{rec1A}), (\ref{rec2A})
transforms into  (\ref{givpar}).
 Taking into account that the change of variables (\ref{anzacgl}) is invertible
 we get a parametrization of $N_+^{(0)}\subset N_+$ $\Box$

Considering the image of the factorized group element (\ref{rec1A}),
(\ref{rec2A}) in the tautological representation $\pi_{\ell+1}$ we
obtain the following relations between factorization and modified
factorization parameters: \be\label{givparGL}
y_{i,n}=e^{x_{n+i,i+1}-x_{n+i-1,i}},\ee where $x_{\ell+1,n}=0,\,\,\,
n=1,\ldots,\ell+1$ are assumed.
Applying the change of variables (\ref{givparGL}) to the expressions in 
Proposition \ref{propgengl} one obtains the realization in the modified
factorized parametrization. 

\begin{prop}
The following differential operators define a realization of
representation $\pi_{\lambda}$ of $\mathfrak{gl}_{\ell+1}$ in
$V_{\mu}$ in terms of modified factorized parametrization
(\ref{givpar}), (\ref{givparGL}) of $N_+$: \bqa E_{i,i}&=& \mu_i -
\sum_{k=1}^{i-1}\frac{\partial}{\partial x_{\ell+1+k-i,k}} +
\sum_{k=i}^{\ell}\frac{\partial}{\partial x_{k,i}},\nonumber \\
E_{i,i+1}&=& -\sum_{k=1}^i\left(\,
\sum_{s=k}^ie^{x_{\ell+1+s-i,s}-x_{\ell+s-i,s}}\right)
\left(\frac{\partial}{\partial x_{\ell+k-i,k}}-
\frac{\partial}{\partial x_{\ell+k-i,k-1}}\right), \label{GGrep}\\
E_{i+1,i}&=& -\sum_{k=1}^{\ell}e^{(x_{k,i}-x_{k+1,i+1})}
\left(\mu_i-\mu_{i+1}+\sum_{s=1}^k \left(\frac{\partial}{\partial
x_{s,i+1}}-\frac{\partial}{\partial x_{s,i}}\right)\right),
\nonumber \eqa where $E_{i,j}=\pi_{\lambda}(e_{i,j})$,
$\mu_k=\imath\lambda_k-\rho_k,\,$ and
$\rho_k=\frac{1}{2}(\ell+2k-2).$ We let $x_{\ell+1,k}=0,$
$(k=1,\ldots,\ell+1).$
\end{prop}

 This realization of the principal series representation
of $\mathfrak{gl}_{\ell+1}$  by differential operators is based on a
particular parametrization of the maximal unipotent subgroup $N_+$
entering the Gauss decomposition of the group $G$
 and was inspired by the Givental integral formula. In \cite{GKLO} we coined the
term Gauss-Givental representation for this realization
 of the principal series representation.
Applying the change of variables (\ref{givparGL}) to the expressions in 
Lemma \ref{Whitfacrep} one obtains Whitaker vectors 
in the modified factorized parametrization. 

\begin{lem}\label{Witmfgl}
The following expressions for  the left/right Whittaker vectors hold:
\bqa\psi_R(x)&=&
\exp\Big\{-\sum\limits_{i=1}^{\ell}\sum\limits_{n=1}^{\ell+1-i}
e^{x_{n+i,i+1}-x_{n+i-1,i}}\Big\},
\\
\psi_L(x)&=&\exp\Big\{\sum\limits_{k=1}^{\ell}\sum\limits_{i=1}^k(\mu_{k+1}-
\mu_k)x_{k,i}\Big\}\,
\exp\Big\{-\sum\limits_{i=1}^{\ell}\sum\limits_{k=1}^{\ell+1-i}
e^{x_{k+i-1,k}-x_{k+i,k}} \Big\},\nonumber  \eqa where we set
$x_{\ell+1,i}=0,\,\,\,\, i=1,\ldots,\ell+1$.
\end{lem}

Now we are ready to write down  the integral representation of the
pairing (\ref{pairing}) using the modified factorized
representation. Going from (\ref{pairing}) to (\ref{int}),
(\ref{intlugl})  we chose to act by an  element of the Cartan torus
to the right in (\ref{pairing}). Different choice (for example the
action to the left) leads to the integrand that differs by total
derivative. The choice made in (\ref{int}), (\ref{intlugl}) is not
the most symmetric one. One of the special features of
Gauss-Givental representation is that up to a simple exponential
term in $\psi_L(x)$ the left and right Whittaker vectors are very
similar  (compare in this respect with the 
case of factorized parameterization (\ref{psir}),
 (\ref{psil})).  We would like to maintain this symmetry
in the integrand of the integral representation. Let us represent
the Cartan group element in the following way:
$$e^{H}=e^{H_L}e^{H_R},$$ where
\bqa
 e^{H}&=&e^{-\sum x_i
E_{i,i}}=\exp\Big\{\sum_{i=1}^{\ell} x_{\ell+1,i}\Big( \mu_i -
\sum_{k=1}^{i-1}\frac{\partial}{\partial x_{\ell+1+k-i,k}} +
\sum_{k=i}^{\ell}\frac{\partial}{\partial x_{k,i}}\Big)\Big\},\\
e^{H_L}&=&\exp\Big\{\sum\limits_{i=1}^{\ell+1}x_{\ell+1,i}\sum_{k=i}^{\ell}
\frac{\partial}{\partial x_{k,i}}\},\\
e^{H_R}&=& \exp\{\sum\limits_{i=1}^{\ell+1}x_{\ell+1,i}\mu_i
-\sum\limits_{i=1}^{\ell+1}x_{\ell+1,i}\sum_{k=1}^{i-1}\frac{\partial}{\partial
x_{\ell+1+k-i,k}}\}.\label{HRgl}\eqa In the calculation of the matrix element we
will chose the differential operator $H_L$ acting on the left vector
and $H_R$ acting on the right vector in (\ref{pairing}). This way
we obtain the following integral formula for eigenfunctions of the
$\mathfrak{gl}_{\ell+1}$-Toda chain.

\begin{te}\label{teintmfgl}
Eigenfunctions of the $\mathfrak{gl}_{\ell+1}$-Toda chain
(\ref{pairgl}) admit the integral representation: \bqa\label{giv}
\Psi_{\lambda_1,\ldots,\lambda_{\ell+1}}^{\mathfrak{gl}_{\ell+1}}
(x_1,\ldots,x_{\ell+1})= \int_{C}
\bigwedge_{k=1}^{\ell}\bigwedge_{i=1}^kdx_{k,i}\,\,
e^{\mathcal{F}^{\mathfrak{gl}_{\ell+1}}(x) }, \eqa where the
function $\mathcal{F}^{{\mathfrak{gl}}_{\ell+1}}(x)$ is given by
\bqa\label{intrep}
\mathcal{F}^{{\mathfrak{gl}}_{\ell+1}}(x)=\imath\sum\limits_{k=1}^{\ell+1}
\lambda_k\Big(\sum\limits_{i=1}^{k}
x_{k,i}-\sum\limits_{i=1}^{k-1}x_{k-1,i}\Big)-\sum\limits_{k=1}^{\ell}
\sum\limits_{i=1}^{k-1}
\Big(e^{x_{k-1,i}-x_{k,i}}+e^{x_{k,i+1}-x_{k-1,i}}\Big).\eqa Here
$x_i=-x_{\ell+1,i},\,\,\,i=1,\ldots,\ell+1$ and  $C \subset N_+$ is
a middle-dimensional non-compact submanifold such that the integrand
decreases  exponentially at the boundaries and at infinities. In
particular  $C$ can be chosen to be
$C=$$\RR^{\frac{(\ell+1)\ell}{2}}$.
\end{te}

\begin{cor} The $\mathfrak{sl}_{\ell+1}$-Whittaker function has the
following integral representation:

\bqa\Psi_{\lambda_2-\lambda_1,\ldots,\lambda_{\ell+1}-\lambda_{\ell}}
^{\mathfrak{sl}_{\ell+1}} (x_2-x_1,\ldots,x_{\ell+1}-x_{\ell})=\\
\nonumber=\exp\Big\{\frac{\imath}{\ell+1}\sum_{i=1}^{\ell+1}\,
x_i\,\Big(\,\ell\lambda_i-\sum_{j\neq i}
\lambda_j\,\Big)\,\,\Big\}\,\,\, \int_{C}
\bigwedge_{k=1}^{\ell}\bigwedge_{i=1}^kdx_{k,i}\,\,
e^{\mathcal{F}^{\mathfrak{sl}_{\ell+1}}(x)},\eqa where
\bqa\mathcal{F}^{\mathfrak{sl}_{\ell+1}}(x)=\imath\sum_{k=1}^\ell
(\lambda_{k+1}-\lambda_k)\,
\sum_{i=1}^k(x_k-x_{k,i})\,-\\ -\nonumber\sum_{k=1}^{\ell}
\sum_{i=1}^{k-1}
\Big(e^{x_{k-1,i}-x_{k,i}}+e^{x_{k,i+1}-x_{k-1,i}}\Big).\eqa
\end{cor}


This integral representation of the $\mathfrak{gl}_{\ell+1}$-Toda
chain eigenfunctions  was first obtained by A.~Givental in his study
of quantum cohomology of the $\mathfrak{gl}_{\ell+1}$ flag manifold
\cite{Gi} (see also \cite{JK}). The description of the  Givental
integral formula in terms of the matrix element (\ref{pairgl}) was
first obtained in \cite{GKLO}.

The function $\mathcal{F}^{{\mathfrak{gl}}_{\ell+1}}(x)$ allows
simple description in terms of the following  diagram introduced by
A.~Givental

\bqa
 \label{AnDiag}
 \xymatrix{
 x_{\ell,1}\ar[d] &&&&\\
 x_{\ell-1,1}\ar[d]\ar[r] & x_{\ell,2}\ar[d] &&&\\
 \vdots\ar[d] & \ddots & \ddots &&\\
 x_{21}\ar[d]\ar[r] & \ddots\ar[d] & \ddots\ar[r] &
 x_{\ell,\ell-1}\ar[d]\\
 x_{11}\ar[r] & x_{22}\ar[r] &
 \ldots\ar[r] & x_{\ell-1,\ell-1}\ar[r] & x_{\ell,\ell}
}\eqa

We assign variables $x_{k,i}$ to  vertexes $(k,i)$ and functions
$e^{y-x}$ to  arrows $(x\longrightarrow y)$ of the diagram
(\ref{AnDiag}). Then the potential function
$\mathcal{F}^{{\mathfrak{gl}}_{\ell+1}}(x)$ (\ref{intrep}) at zero
spectrum $\la_i=0$ is given by the  sum of the functions assigned
to all arrows.

As it was demonstrated in the Theorem \ref{Mpargl} variables
$\{x_{k,i}\}$ provide a parametrization of  an open part $N_+^{(0)}$
of the flag manifold $X=SL(\ell+1,\mathbb{C})/B$. The non-compact
manifold $N_+^{(0)}$ has a natural action of the torus $T^{l(w_0)}$
and can be compactified to a (singular) toric variety. The set of
the monomial relations defining this compactification can be
described as follows. Introduce new variables
$$a_{k,i}=e^{x_{k,i}-x_{k+1,i}}
,\,\,\,\,\,b_{k,i}=e^{x_{k+1,i+1}-x_{k,i}},\,\,\,\,\,\,\quad 1\leq
k\leq\ell,\,\,1\leq i\leq k,\,$$ assigned to  arrows of the diagram
(\ref{AnDiag}). Then the following defining relations hold
\be\label{defrelAn}
 a_{k,i}\cdot b_{k,i}\,=\,
b_{k+1,i}\cdot a_{k+1,i+1},\qquad 1\leq k< \ell,\,\,1\leq i\leq k,\\
a_{\ell,i}\cdot b_{\ell,i}=e^{x_{\ell,i+1}-x_{\ell,i}}.\ee They can
be interpreted as relations between various compositions of
elementary paths having the same initial and final vertexes. The set
of relations between  more general paths (following from
(\ref{defrelAn}))  provides a toric embedding of the degeneration of
flag manifold (see \cite{BCFKS} for details).

\subsubsection{ Relation with $\widehat{\mathfrak{gl}}_{\ell+1}$-Toda
chain  Baxter  $\mathcal{Q}$-operator}

Integral representation (\ref{giv}), (\ref{intrep})  of
$\mathfrak{gl}_{\ell+1}$-Whittaker function  has a recursive
structure over the rank $\ell$ of the Lie algebra. Indeed the
integral representation can be rewritten in the following form
 \bqa\label{giviter}
\Psi_{\lambda_1,\ldots,\lambda_{\ell+1}}^{\mathfrak{gl}_{\ell+1}}
(\underline{x}_{\ell+1})= \int_{C}
\bigwedge_{k=1}^{\ell}\bigwedge_{i=1}^kdx_{k,i}\,\,
\prod_{k=1}^{\ell+1}
Q^{\mathfrak{gl}_{k}}_{\mathfrak{gl}_{k-1}}(\underline{x}_{k};\underline{x}_{k-1};
\lambda_k),\eqa where \bqa\label{QRecA}
Q^{\mathfrak{gl}_{k+1}}_{\mathfrak{gl}_{k}}
(\underline{x}_{k+1};\underline{x}_{k};\lambda_{k+1})=\eqa \bqa
\hspace{-0.3cm}=\exp\left\{\i\lambda_{k+1}\Big(\sum_{i=1}^{k+1}
x_{k+1,i}- \sum_{i=1}^{k}x_{k,i}\Big)-
\sum_{i=1}^{k}\Big(e^{x_{k,i}-x_{k+1,i}}+e^{x_{k+1,i+1}-x_{k,i}}\Big)
\right\}.\nonumber\eqa Here  we denote
$\underline{x}_{k}=(x_{k,1},\ldots,x_{k,k})$ and assume that
$Q^{\mathfrak{gl}_{1}}_{\mathfrak{gl}_{0}}=e^{\imath\lambda_1x_{1,1}}.$

Let us chose linear coordinates
$\underline{x}_{k}=(x_{k,1},\ldots,x_{k,k})$ in
$\mathbb{C}^{k}$. Let  $C_k$ be a non-compact
middle-dimensional submanifold  in $\mathbb{C}^{k}$ such that (\ref{QRecA})
as a function of $\underline{x}_{k}$ decreases exponentially at
possible  boundaries and infinities of $C_k$. Consider the following
integral operator
 $$
(Q^{\mathfrak{gl}_{k+1}}_{\mathfrak{gl}_{k}}f)(\underline{x}_{k})=\int_{C_k}
Q^{\mathfrak{gl}_{k+1}}_{\mathfrak{gl}_{k}}
(\underline{x}_{k+1};\underline{x}_{k};\lambda_{k+1})
f(\underline{x}_k)d\underline{x}_k.$$ acting on  functions not
growing too fast at possible boundaries and infinities of $C_k$.
Integral operators $Q^{\mathfrak{gl}_{k+1}}_{\mathfrak{gl}_{k}}$
provide a recursive  construction of
 $\mathfrak{gl}_{\ell+1}$-Whittaker functions:
 \bqa\label{recursgl}
\Psi_{\lambda_1,\ldots,\lambda_{\ell+1}}^{\mathfrak{gl}_{\ell+1}}
(\underline{x}_{\ell+1})= \int_C
\bigwedge_{i=1}^{\ell}dx_{\ell,i}\,\,
Q^{\mathfrak{gl}_{\ell+1}}_{\mathfrak{gl}_{\ell}}
(\underline{x}_{\ell+1};\underline{x}_{\ell};\lambda_{\ell+1})
\Psi_{\lambda_1,\ldots,\lambda_{\ell}}^{\mathfrak{gl}_{\ell}}
(\underline{x}_{\ell}).\eqa There is a  natural oriented path in the
diagram (\ref{AnDiag}), which can be associated with   the recursive
operator $Q^{\mathfrak{gl}_{\ell+1}}_{\mathfrak{gl}_{\ell}}:$

\bqa\label{orpath}
\xymatrix{
 x_{\ell,1}\ar[ddr] && x_{\ell,2}\ar[ddr] && \ldots &
 x_{\ell,\ell}\\
 &&&&&\\
 & x_{\ell-1,1}\ar[uur] && \ldots & x_{\ell-1,\ell-1}\ar[uur] &
}
\eqa

Diagram (\ref{AnDiag}) can be considered as a
 collection of the oriented pathes (\ref{orpath}) and
thus the recursive construction of the integral representation is
encoded in the  diagram (\ref{AnDiag}) in an obvious way.

As a consequence of (\ref{recursgl}),  integral operators
$Q^{\mathfrak{gl}_{k+1}}_{\mathfrak{gl}_{k}}$
 with the kernels
 $Q^{\mathfrak{gl}_{k+1}}_{\mathfrak{gl}_{k}}(\underline{x}_{k+1},\,
\underline{x}_{k};\lambda_{k+1})$ satisfy braiding relations with
the Quantum Toda chain Hamiltonians. For example the following
 relation between quadratic Hamiltonians
$\CH_2^{\mathfrak{gl}_{k+1}}(\underline{x}_{k+1})$ and
$\CH_2^{\mathfrak{gl}_k}(\underline{x}_{k}),$  holds
\bqa\label{intertw} \CH_2^{\mathfrak{gl}_{k+1}}(\underline{x}_{k+1})
Q^{\mathfrak{gl}_{k+1}}_{\mathfrak{gl}_{k}}(\underline{x}_{k+1},\,
\underline{x}_{k};\lambda_{k+1})=
Q^{\mathfrak{gl}_{k+1}}_{\mathfrak{gl}_{k}}
(\underline{x}_{k+1},\,\underline{x}_{k};\lambda_{k+1})
\CH_2^{\mathfrak{gl}_{k}}(\underline{x}_k)
+\frac{1}{2}\lambda_{k+1}^2. \eqa   We shall assume that in the
relation above and similar ones,  Hamiltonian operators on l.h.s.
act to the right and  Hamiltonians on r.h.s. act to the left.
Similar braiding relations hold for  higher quantum Hamiltonian
operators (see \cite{GKLO}  for  details).

The recursion operators
$Q^{\mathfrak{gl}_{k+1}}_{\mathfrak{gl}_{k}}$
 appear to be related with an important object in the theory of Quantum
 Integrable
Systems,  $\mathcal{Q}$-operator. $\mathcal{Q}$-operator was
introduced by R.~Baxter \cite{B} for certain statistical models as a
tool  to solve  quantum integrable models explicitly. In the case of
$\widehat{\mathfrak{gl}}_{\ell+1}$-Toda chain, with the quadratic
Hamiltonian \be
\CH_2^{\widehat{\mathfrak{gl}}_{\ell+1}}=-\frac{1}{2}\sum\limits_{i=1}^{\ell+1}
\frac{\partial^2}{\partial x_i^2}+ \sum\limits_{i=1}^{\ell}
e^{x_{i+1}-x_{i}}+ge^{x_1-x_{\ell+1}}\,, \ee where $g$ is an
arbitrary coupling constant,  the $\mathcal{Q}$-operator has the
following  integral kernel \bqa\label{baff}
\mathcal{Q}^{\widehat{\mathfrak{gl}}_{\ell+1}}(\underline{x}^{(\ell+1)},
\underline{y}^{(\ell+1)};\lambda)=\exp\Big\{\,\imath
\lambda\sum\limits_{i=1}^{\ell+1}(x_i-y_i)-\\ \nonumber
-\Big(\sum_{i=1}^{\ell} \left(e^{x_{i}-y_{i}}+
e^{y_{i+1}-x_{i}}\right)\,+e^{x_{\ell+1}-y_{\ell+1}}+
ge^{y_1-x_{\ell+1}}\Big)\Big\}. \eqa Here we use notations
$\underline{x}^{(\ell+1)}=(x_1,\dots,x_{\ell+1})$ and
$\underline{y}^{(\ell+1)}=(y_1,\ldots,y_{\ell+1})$. This $\mathcal{Q}$-operator
was first constructed in \cite {PG}. It commutes with all
Hamiltonians  of $\widehat{\mathfrak{gl}}_{\ell+1}$-Toda chain and
generates quantum B\"{a}cklund transformations \cite{PG}. For
instance, for the quadratic Hamiltonians we have:
\bqa\label{intertwaff}
\CH_2^{\widehat{\mathfrak{gl}}_{\ell+1}}(\underline{x}^{(\ell+1)})
\mathcal{Q}^{\widehat{\mathfrak{gl}}_{\ell+1}}(\underline{x}^{(\ell+1)},\,
 \underline{y}^{(\ell+1)},\lambda)=
 \mathcal{Q}^{\widehat{\mathfrak{gl}}_{\ell+1}} (\underline{x}^{(\ell+1)},\,
\underline{y}^{(\ell+1)},\lambda)
\CH_2^{\widehat{\mathfrak{gl}}_{\ell+1}}(\underline{y}^{(\ell+1)}).
\eqa

To establish a relation between Baxter $\mathcal{Q}$-operator for
$\widehat{\mathfrak{gl}}_{k+1}$-Toda theory and a recursion operator
for $\mathfrak{gl}_{k+1}$-Toda theory it is useful to introduce a
 slightly modified recursion operator
 $Q^{\mathfrak{gl}_{k+1}}_{\mathfrak{gl}_{k}\oplus \mathfrak{gl}_1}$
 with the kernel:
\bqa\label{newitertw}
Q^{\mathfrak{gl}_{k+1}}_{\mathfrak{gl}_{k}\oplus \mathfrak{gl}_1}
(\underline{x}^{(k+1)},\underline{y}^{(k+1)},\lambda)=
\exp\left\{\i \la \,y_{k+1}\right\}\,
Q^{\mathfrak{gl}_{k+1}}_{\mathfrak{gl}_{k}}
(\underline{x}^{(k+1)},\underline{y}^{(k)},\lambda)=\eqa \bqa
\hspace{-0.3cm}\exp\left\{\i\lambda\Big(\sum_{i=1}^{k+1}
x_{i}- \sum_{i=1}^{k}y_{i}\Big)-
\sum_{i=1}^{k}\Big(e^{y_{i}-x_{i}}+e^{x_{i+1}-y_{i}}\Big)
\right\},\nonumber \eqa
where $\underline{x}^{(k+1)}=(x_1,\dots,x_{k+1})$,
$\underline{y}^{(k)}=(y_1,\ldots,y_{k})$ and 
$\underline{y}^{(k+1)}=(y_1,\ldots,y_k,y_{k+1})$.

 This modified operator  intertwines Hamiltonian operators of
$\mathfrak{gl}_{k+1}$- and $\mathfrak{gl}_{k}\oplus
\mathfrak{gl}_{1}$-Toda chains (the new variable $y_{k+1}$ enters
only $\mathfrak{gl}_{1}$-Toda chain). Thus for quadratic Hamiltonian
operators we have \bqa\nonumber
\CH_2^{\mathfrak{gl}_{k+1}}(\underline{x}^{(k+1)})
Q^{\mathfrak{gl}_{k+1}}_{\mathfrak{gl}_{k}\oplus \mathfrak{gl}_1}
(\underline{x}^{(k+1)},\,\underline{y}^{(k+1)},\lambda)
=Q^{\mathfrak{gl}_{k+1}}_{\mathfrak{gl}_{k}\oplus \mathfrak{gl}_1}
(\underline{x}^{(k+1)},\,\underline{y}^{(k+1)},\lambda)
(\CH_2^{\mathfrak{gl}_{k}}(\underline{y}^{(k)})+\CH_2^{\mathfrak{gl}_1}(y_{k+1})),
 \eqa
where
$\CH_2^{\mathfrak{gl}_1}(y_{k+1})=-\frac{1}{2}\frac{\pr^2}{\pr
y_{k+1}^2}$. Obviously the projection of the above relation  on the
subspace of functions $F(\underline{y},y_{k+1})=\exp(\i
\la\,y_{k+1})\,f(\underline{x})$ leads to (\ref{intertw}).

Now consider a one-parameter family of integral operators
\bqa\label{onepar}
\mathcal{Q}^{\widehat{\mathfrak{gl}}_{\ell+1}}(\underline{x}^{(k+1)},
\underline{y}^{(k+1)};\lambda;\varepsilon)=\varepsilon^{\i \la}
\exp\Big\{\,\imath \lambda\sum\limits_{i=1}^{\ell+1}(x_i-y_i)-\\
\nonumber -\Big(\sum_{i=1}^{\ell} \left(e^{x_{i}-y_{i}}+
e^{y_{i+1}-x_{i}}\right)\,+\varepsilon\,e^{x_{\ell+1}-y_{\ell+1}}+
\varepsilon^{-1}ge^{y_1-x_{\ell+1}}\Big)\Big\}. \eqa
 obtained from (\ref{baff}) by a shift of the variable
$x_{\ell+1}\to x_{\ell+1}+\ln\varepsilon$.  The limiting
 behavior of  (\ref{onepar}) when $\varepsilon\to 0$, $g
 \epsilon^{-1}\to 0$ can be described as follows
\bqa Q_{\mathfrak{gl}_{k+1}}^{\mathfrak{gl}_{k}\oplus
\mathfrak{gl}_1}
(\underline{x}^{(k+1)},\,\underline{y}^{(k+1)},\lambda)=\lim_{
\varepsilon\to 0,\, g  \varepsilon^{-1}\to 0}\, \varepsilon^{-\i
  \la}\,\mathcal{Q}^{\widehat{\mathfrak{gl}}_{k+1}}(\underline{x}^{(k+1)},
\underline{y}^{(k+1)},\lambda,\varepsilon). \eqa  This provides a relation
between the Baxter $\mathcal{Q}$-operator and the (modified) recursion
operator.
\newpage


\newpage
\subsection{ Integral representations of $\mathfrak{so}_{2\ell+1}$-Toda chain
eigenfunctions }

In this subsection we provide a generalization of the Givental
integral representation of $\mathfrak{gl}_{\ell+1}$-Whittaker
functions to the case of $\mathfrak{so}_{2\ell+1}$. We start with a
derivation of the integral representation of
$\mathfrak{so}_{2\ell+1}$-Whittaker functions using the factorized
representation. Then we consider a
modification of  the factorized representation that directly leads
to a Givental type integral representation.

Consider $B_{\ell}$ type root system corresponding to
 Lie algebra $\mathfrak{so}_{2\ell+1}$.
 Let $(\epsilon_1,\ldots,\epsilon_{\ell})$ be an
orthogonal basis in $\R^{\ell}.$ We realize $B_{\ell}$
 roots, coroots and fundamental weights as vectors in $\RR^{\ell}$
in the following way:

\bqa\label{RootBn}
\begin{array}{l}
\alpha_1=\epsilon_1,\\
\alpha_2=\epsilon_2-\epsilon_1,\\ \ldots\\
\alpha_{\ell}=\epsilon_{\ell}-\epsilon_{\ell-1},
\end{array}
\hspace{2cm}
\begin{array}{l}
\alpha_1^{\vee}=2\epsilon_1,\\
\alpha_2^{\vee}=\epsilon_2-\epsilon_1,\\ \ldots\\
\alpha_{\ell}^{\vee}=\epsilon_{\ell}-\epsilon_{\ell-1},
\end{array}
\hspace{2cm}
\begin{array}{l}
\omega_1=\frac{1}{2}(\epsilon_1+\ldots+\epsilon_{\ell}),\\
\omega_2=\epsilon_2+\ldots+\epsilon_{\ell},\\ \ldots\\
\omega_{\ell}=\epsilon_{\ell}.
\end{array} \eqa
The Cartan matrix is then given by
$a_{ij}=\<\alpha_i^{\vee},\alpha_j\>$ and positive rational numbers
$d_1=\frac{1}{2},d_2=1,\ldots,d_{\ell}=1$ are  such that the matrix
$\|b_{ij}\|=\|d_i a_{ij}\|$ is symmetric.
 One associates with these
data a Quantum Toda chain with a quadratic Hamiltonian
\bqa\label{BtwoH}
\CH_2^{B_{\ell}}&=&-\frac{1}{2}\sum\limits_{i=1}^{\ell}
\frac{\partial^2}{{\partial x_i}^2}+
\frac{1}{2}e^{x_1}+\sum\limits_{i=1}^{\ell-1} e^{x_{i+1}-x_{i}}
.\eqa One can complete (\ref{BtwoH}) to a full set of $\ell$
mutually commuting  functionally independent Hamiltonians
$H^{B_{\ell}}_k$ of the $\mathfrak{so}_{2\ell+1}$-Toda chain. We are
looking for common eigenfunction integral representations of
 the commuting set of the Hamiltonians. Corresponding eigenfunction problem for the
quadratic Hamiltonian can be written in the following form  \be
\CH_2^{B_{\ell}}(x)\,\,\,\Psi^{{B}_{\ell}}_{\la_1,\cdots
,\la_{\ell}} (x_1,\ldots,x_{\ell})=\frac{1}{2}
\sum\limits_{i=1}^{\ell}\lambda_{i}^2\,\,\,
\Psi^{{B}_{\ell}}_{\la_1,\cdots ,\la_{\ell}} (x_1,\ldots,x_{\ell}).
\ee

\subsubsection{ $\mathfrak{so}_{2\ell+1}$-Whittaker function:
factorized   parametrization}

The reduced word for the element $w_0$ of maximal length in the Weyl
group of $B_{\ell}$-type can be represented in the recursive form:
$$I=(i_1,i_2,\ldots,i_m):=(1,212,32123,\ldots,(\ell\ldots212\ldots\ell)),$$
where indexes  $i_k$ correspond to elementary reflections with
respect to  roots $\alpha_k$. Let $N_+\subset G$ be a maximal
unipotent subgroup of $G=SO(2\ell+1)$. One associates with the
reduced word  $I$ the following recursive parametrization of a
generic unipotent element $v^{B_{\ell}}\in N_+$: \bqa\label{rec1B}
v^{B_{\ell}}=v^{B_{\ell-1}}\,\mathfrak{X}^{B_{\ell}}_{B_{\ell-1}}
,\eqa where
\bqa\label{rec2B}\mathfrak{X}^{B_\ell}_{B_{\ell-1}}=X_{\ell}(y_{\ell,1})\cdots
X_k(y_{k,2(\ell+1-k)-1})\cdots X_2(y_{2,2\ell-3})\times\\
\nonumber \times X_1(y_{1,\ell})X_2(y_{2,2\ell-2})\cdots
X_k(y_{k,2(\ell+1-k)})\cdot X_{\ell}(y_{\ell,2}). \eqa Here
$X_i(y)=e^{ye_i}$ and $e_i\equiv e_{\alpha_i}$ are simple root
generators. The subset $N_+^{(0)}$ of elements allowing the
representation is an open part of $N_+$. The action of the Lie
algebra $\mathfrak{so}_{2\ell+1}$ on $N_+$ (\ref{infact}) considered
at the beginning of the previous section defines an action of the
Lie algebra on $N_+^{(0)}$. The following proposition explicitly
describes this action on the space $V_{\mu}$
 considered as a space of functions on $N_+^{(0)}$.

\begin{prop}\label{Bfactrep}
The following differential operators define a realization of a
principal series representation $\pi_{\lambda}$ of
$\mathcal{U}(\mathfrak{so}_{2\ell+1})$
 in terms of factorized parametrization  of $N_+^{(0)}$:
\bqa E_1&=& \frac{\partial}{\partial y_{1,\ell}}\,+\,
\sum_{n=1}^{\ell-1}\Big\{\left(\frac{\partial}{\partial y_{1,n}}-
\frac{\partial}{\partial y_{1,n+1}}\right) \prod_{j=n}^{\ell-1}
\frac{y_{2,2i}}{y_{2,2i-1}}\,+ \nonumber\\+
&&2\left(\frac{\partial}{\partial y_{2,2n-1}}-
\frac{\partial}{\partial y_{2,2n}}\right)
\frac{y_{2,2(n-1)}}{y_{1,n}}
\prod_{i=n+1}^{\ell-1}\frac{y_{2,2i}}{y_{2,2i-1}}\Big\},\nonumber\\
E_k&=& \frac{\partial}{\partial y_{k,2(\ell+1-k)}}\,+\,
\sum_{n=1}^{n-k}\Big\{\left(\frac{\partial}{\partial y_{k,2n}}-
\frac{\partial}{\partial y_{k,2n+1}}\right)
\prod_{i=n}^{\ell-k}\frac{y_{k+1,2i}}{y_{k+1,2i-1}}
\frac{y_{k,2(i+1)-1}}{y_{k,2(i+1)}}\,+\\+ \nonumber
&&\left(\frac{\partial}{\partial y_{k+1,2n-1}}-
\frac{\partial}{\partial y_{k+1,2n}}\right)
\frac{y_{k+1,2n}}{y_{k,2(n+1)}}
\prod_{i=n+1}^{\ell-k}\frac{y_{k+1,2i}}{y_{k+1,2i-1}}
\frac{y_{k,2(i+1)-1}}{y_{k,2(i+1)}}\Big\},\qquad 1<k<\ell,\nonumber \\
E_\ell&=&\frac{\partial}{\partial y_{\ell,2}},\nonumber \\
H_k&=&\<\mu\,,\alpha_k^\vee\>\,+\,\sum_{i=1}^\ell
a_{k,i}\sum_{j=1}^{n_i}y_{i,j}\frac{\partial}{\partial y_{i,j}},
\hspace{3cm}1\leq k\leq\ell,\eqa

\bqa F_1&=& \sum_{n=1}^\ell\,y_{1,n}\Big(\<\mu,\alpha_1^\vee\>+
\sum_{j=1}^{2(n-1)-1}2y_{2,j}\frac{\partial}{\partial y_{2,j}}
-2\sum_{j=1}^{n-1}y_{1,j}\frac{\partial}{\partial y_{1,j}}-
y_{1,n}\frac{\partial}{\partial y_{1,n}}\Big),\nonumber \\
F_2&=& \sum_{n=1}^{2(\ell-1)}\,y_{2,n}\Big(\<\mu,\alpha_2^\vee\>+
2\sum_{j=1}^{[n/2]+1}y_{1,j}\frac{\partial}{\partial y_{1,j}}
-2\sum_{j=1}^{n-1}y_{2,j}\frac{\partial}{\partial y_{2,j}}+\nonumber \\
&&+\sum_{j=1}^{2[(n+1)/2]-3} y_{3,j}\frac{\partial}{\partial
y_{3,j}}- y_{2,n}\frac{\partial}{\partial y_{2,n}}\Big),\\
F_k&=& \sum_{n=1}^{2(\ell+1-k)}\,y_{k,n}\Big(\<\mu,\alpha_k^\vee\>+
\sum_{j=1}^{2[n/2]+1}y_{k-1,j}\frac{\partial}{\partial y_{k-1,j}}
-2\sum_{j=1}^{n-1}y_{k,j}\frac{\partial}{\partial y_{k,j}}+\nonumber\\
&&+ \sum_{j=1}^{2[(n+1)/2]-3} y_{k+1,j}\frac{\partial}{\partial
y_{k+1,j}}- y_{k,n}\frac{\partial}{\partial y_{k,n}}\Big), \qquad
2<k<\ell,
\nonumber \\
F_\ell&=&(y_{\ell,1}+y_{\ell,2})\left( \<\mu,\alpha_\ell^\vee\>+
y_{\ell-1,1}\frac{\partial}{\partial y_{\ell-1,1}}+
y_{\ell-1,2}\frac{\partial}{\partial
y_{\ell-1,2}}\right)+\nonumber\\
&& +y_{\ell,2}\left( y_{\ell-1,3}\frac{\partial}{\partial
y_{\ell-1,3}}+
y_{\ell-1,4}\frac{\partial}{\partial y_{\ell-1,4}}\right)-
\left(y_{\ell,1}^2\frac{\partial}{\partial y_{\ell,1}}+
2y_{\ell,1}y_{\ell,2}\frac{\partial}{\partial y_{\ell,1}}+
y_{\ell,2}^2\frac{\partial}{\partial y_{\ell,2}}\right),\nonumber
\eqa where
$\pi_{\lambda}(e_i)=E_i,\,\,\,\pi_{\lambda}(f_i)=F_i,\,\,\,\pi_{\lambda}(h_i)=H_i\,\,\,
i=1,\ldots,\ell$,  $n_1=\ell$, $n_k=2(\ell+1-k)$ for $1<k\leq\ell$,
$a_{ij}$ is a Cartan matrix and we assume that the terms  containing
$y_{i,j}$ with the   indexes not in  the set $\{1\leq i,j\leq
\ell\}$ should be  omitted.

\end{prop}
For the proof see Part II, Section \ref{genso2l1}.

Left/right Whittaker vectors in the factorized parametrization
have the following expressions.

\begin{lem}\label{lemwitso2l1}
The following expressions for the left/right
Whittaker vectors hold:
\bqa\psi_R(y)=\exp\Big\{-\Big(\sum\limits_{n=1}^{\ell}y_{1,n}
+\sum\limits_{k=2}^{\ell}\sum\limits_{n=1}^{n_k}y_{k,n}\Big)\Big\},\nonumber\\
\psi_L(y)= \Big(\prod_{n=1}^{\ell}y_{1,n}
\prod_{i=2}^{\ell}\prod_{n=1}^{n_i/2}y_{i,2n-1}\Big)^{\<\mu,\alpha_1^{\vee}\>}
\times\nonumber\\
\times  \prod_{k=2}^\ell\Big(\prod_{n=2}^{\ell}y_{1,n}^2
\prod_{i=k+1}^{\ell}\prod_{n=1}^{n_i/2}y_{i,2n-1}^2
\prod_{i=2}^k\prod_{n=1}^{n_i/2}y_{i,2n-1}y_{i,2n}
\Big)^{\<\mu,\alpha_k^{\vee}\>}\times \\
\nonumber \times \exp\Big\{ -\Big(\sum_{n=1}^\ell\frac{1}{y_{1,n}}
\Big(1+\frac{y_{2,2(n-1)}}{y_{2,2(n-1)-1}}\Big)\prod_{i=n+1}^\ell
\frac{y_{2,2(i-1)}}{y_{2,2(i-1)-1}}+\\ \nonumber+\sum_{k=2}^\ell
\sum_{n=1}^{n_k/2} \frac{1}{y_{k,2n}}
\Big(1+\frac{y_{k+1,2(n-1)}}{y_{k+1,2(n-1)-1}}\Big)
\prod_{i=n+1}^{n_k/2} \frac{y_{k+1,2(i-1)}}{y_{k+1,2(-1)-1}}
\frac{y_{k,2i-1}}{y_{k,2i}}\Big)\Big\}, \eqa where $n_1=\ell$ and
$n_k=2(\ell+1-k),$ $k=2,\ldots,\ell.$
\end{lem}
For the proof see Part II, Section \ref{melso2l1}.

Using (\ref{pairing}) and (\ref{meas}) we obtain an integral
representation of $\mathfrak{so}_{2\ell+1}$-Whittaker functions in
the factorized parametrization.

\begin{te}\label{tewfso2l1}

The eigenfunctions of the  $\mathfrak{so}_{2\ell+1}$-Toda
chain (\ref{pairing}) admit the following integral representation:
\bqa\label{intlu}\Psi_{\lambda_1,\ldots,\lambda_{\ell}}^{B_{\ell}}
(x_1,\ldots,x_{\ell})=e^{\imath\lambda_1x_1+\ldots+\imath\lambda_{\ell}x_{\ell}}\int_C
\bigwedge_{i=1}^\ell\bigwedge_{k=1}^{n_i}\frac{dy_{i,k}}{y_{i,k}}
\Big(\prod_{n=1}^{\ell}y_{1,n}
\prod_{i=2}^{\ell}\prod_{n=1}^{n_i/2}y_{i,2n-1}\Big)^{2\imath\lambda_1}\times
\nonumber \\
\times \prod_{k=2}^\ell\Big(\prod_{n=2}^{\ell}y_{1,n}^2
\prod_{i=k+1}^{\ell}\prod_{n=1}^{n_i/2}y_{i,2n-1}^2
\prod_{i=2}^k\prod_{n=1}^{n_i/2}y_{i,2n-1}y_{i,2n}
\Big)^{\imath(\lambda_k-\lambda_{k-1})}\times \nonumber \\ \times
\exp\Big\{ -\Big(\sum_{n=1}^\ell\frac{1}{y_{1,n}}
\Big(1+\frac{y_{2,2(n-1)}}{y_{2,2(n-1)-1}}\Big)\prod_{i=n+1}^\ell
\frac{y_{2,2(i-1)}}{y_{2,2(i-1)-1}}+\\ \nonumber +\sum_{k=2}^\ell
\sum_{n=1}^{n_k/2} \frac{1}{y_{k,2n}}
\Big(1+\frac{y_{k+1,2(n-1)}}{y_{k+1,2(n-1)-1}}\Big)
\prod_{i=n+1}^{n_k/2} \frac{y_{k+1,2(i-1)}}{y_{k+1,2(-1)-1}}
\frac{y_{k,2i-1}}{y_{k,2i}} \,\,+\\ +\nonumber
e^{x_{1}}\sum_{n=1}^\ell y_{1,n}\,+\,\sum_{k=2}^{\ell}
e^{x_{k}-x_{k-1}} \sum\limits_{n=1}^{n_k}y_{k,n}\Big)\Big\},\eqa
where  $n_1=\ell$, $n_k=2(\ell+1-k),\,\,\,k=2,\ldots,\ell$ and $C
\subset N_+$ is a middle-dimensional non-compact submanifold  such
the integrand decays exponentially  at the boundaries and
infinities. In particular one can chose $C=\RR_{+}^{\ell^2}.$
\end{te}
The proof is given in Part II, Section \ref{melso2l1}.

\begin{ex} Let $\ell=2$. In this case, the general formula
(\ref{intlu}) acquirs the form \bqa\label{B2Lu}
\Psi^{B_2}_{\lambda_1,\lambda_2}(x_{1},x_{2})= e^{\imath\lambda_1
x_{1}+\imath\lambda_2 x_{2}}\int\limits_{C}
\bigwedge_{i=1}^{2}\bigwedge_{k=1}^{2}\frac{d y_{i,k}}{y_{i,k}}
(y_{11}y_{21}y_{12})^{2\imath\lambda_1}
(y_{21}y_{12}^2y_{22})^{\imath\lambda_2-\imath\lambda_1}\cdot\\
\nonumber \exp\Big\{-\Big(\Big\{\frac{1}{y_{12}}+
\frac{y_{22}}{y_{21}}\Big(\frac{1}{y_{11}}+
\frac{1}{y_{12}}\Big)\Big\}+\frac{1}{y_{22}}+e^{x_{1}}(y_{11}+y_{12})+
e^{x_{2}-x_{1}}(y_{21}+y_{22}) \Big)\Big\},\eqa  where one can chose
 $C=\RR_+^4.$
\end{ex}

\subsubsection{ $\mathfrak{so}_{2\ell+1}$-Whittaker function:
modified factorized   parametrization}

In this part we introduce a modified  factorized parametrization of
$N_+$.   We use this parametrization to construct the integral
representations for $\mathfrak{so}_{2\ell+1}$-Whittaker functions.
In contrast with the integral representations described above these
integral representations have  a simple recursive structure over the
rank $\ell$ and can be described in purely combinatorial terms using
suitable graphs. Thus these representations can be considered as a
generalization of Givental integral representations to
$\mathfrak{so}_{2\ell+1}$.

There exists a realization of a tautological representation
$\pi_{2\ell+1}:\mathfrak{so}_{2\ell+1}\to End(\mathbb{C}_{2\ell+1})$
such that Weyl generators corresponding to  Borel (Cartan)
subalgebra of $\mathfrak{so}_{2\ell+1}$ are realized by
 upper triangular (diagonal) matrices.
This defines an embedding $\mathfrak{so}_{2\ell+1}\subset
\mathfrak{gl}_{2\ell+1}$ such that Borel (Cartan) subalgebra maps
into Borel (Cartan) subalgebra (see e.g. \cite{DS}). To define the
corresponding embedding of the groups  consider the following
involution on $GL(2\ell+1)$: \bqa g\longmapsto g^*:=\dot{w}_0\cdot
(g^{-1})^{t}\cdot\dot{w}_0^{-1},\eqa where $a^{t}$ is induced by the
standard transposition of the  matrix $a$ and $\dot{w}_0$ is a lift
of the maximal length element $w_0$ of the Weyl group of
$\mathfrak{gl}_{2\ell+1}$. In a matrix form it can be written as 
$$ \dot{w}_0=S\cdot J,$$ where
 $S={\rm diag}(1,-1,\ldots,-1,1)$ and $J=\|J_{i,j}\|=\|\delta_{i+j,2\ell+2}\|$.
 The orthogonal group $G=SO(2\ell+1)$ then can be defined as a
 following  subgroup of   $GL(2\ell+1)$
$$ SO(2\ell+1)=\{g\in GL(2\ell+1):g^{*}=g\}.$$

Let $\epsilon_{i,j}$ be elementary $(2\ell+1)\times (2\ell+1)$
matrices with  units at the $(i,j)$ place and zeros, otherwise. For
any $n=2,\ldots,\ell$ introduce matrices $U_n,\,\widetilde{U}_n$ and
$V_n,\,\widetilde{V}_n:$ \bqa
U_n=\sum_{i=1}^{\ell-n}(\epsilon_{i,i}+\epsilon_{2\ell+2-i,2\ell+2-i})+
\sum_{i=1}^n\epsilon_{\ell-n+i,\ell-n+i}\,+
e^{-z_{n,1}}\epsilon_{\ell+1,\ell+1}+\\
\nonumber +e^{z_{n,1}}\epsilon_{\ell+2,\ell+2}\,+
\sum_{i=1}^{n-1}e^{-z_{\ell,\ell+1-i}}\epsilon_{\ell+n+2-i,\ell+n+2-i},\eqa

\bqa \widetilde{U}_n=
\sum_{i=1}^{\ell-n}(\epsilon_{i,i}+\epsilon_{2\ell+2-i,2\ell+2-i})+
\sum_{i=1}^n\epsilon_{\ell-n+i,\ell-n+i}\,+e^{-z_{n,1}}\epsilon_{\ell+1,\ell+1}+
e^{z_{n,1}}\epsilon_{\ell+2,\ell+2}\,+\\ \nonumber
+\sum_{i=1}^{n-1}e^{-z_{\ell,\ell+1-i}}\epsilon_{\ell+n+2-i,\ell+n+2-i}\,+\,
e^{x_{n-1,1}}\epsilon_{\ell+1,\ell+2}+
\sum_{i=2}^{n-2}e^{-x_{n-1,i}}\epsilon_{\ell+i,\ell+i+1}, \eqa

and \bqa
V_n=\sum_{i=1}^{\ell-n}(\epsilon_{i,i}+\epsilon_{2\ell+2-i,2\ell+2-i})+\\
\nonumber +\sum_{i=1}^ne^{x_{n,n+1-i}}\epsilon_{\ell-n+i,\ell-n+i}+
e^{-x_{n,1}}\epsilon_{\ell+1,\ell+1}+
\sum_{i=1}^n\epsilon_{\ell+i+1,\ell+i+1},\eqa

\bqa \widetilde{V}_n=\sum_{i=1}^{\ell-n}(\epsilon_{i,i}+
\epsilon_{2\ell+2-i,2\ell+2-i})+
\sum_{i=1}^ne^{x_{n,n+1-i}}\epsilon_{\ell-n+i,\ell-n+i}+
e^{-x_{n,1}}\epsilon_{\ell+1,\ell+1}+\\
\nonumber +\sum_{i=1}^n\epsilon_{\ell+i+1,\ell+i+1}
+\sum_{i=1}^{n-1}e^{z_{n,n+1-i}}\epsilon_{\ell-n+i,\ell-n+i+1}+
e^{-z_{n,1}}\epsilon_{\ell,\ell+1},\eqa

\bqa
U_1&=&\sum_{i=1}^{\ell-1}\epsilon_{i,i}+e^{-z_{11}}\epsilon_{\ell,\ell}+
e^{+z_{11}}\epsilon_{\ell+1,\ell+1}+\sum_{i=\ell+2}^{2\ell+1}\epsilon_{i,i},
\eqa
\bqa \widetilde{U}_1&=&\sum_{i=1}^{\ell-1}\epsilon_{i,i}+
e^{-z_{11}}\epsilon_{\ell,\ell}+
e^{z_{11}}\epsilon_{\ell+1,\ell+1}+\sum_{i=\ell+2}^{2\ell+1}\epsilon_{i,i}+
e^{x_{11}}\epsilon_{\ell,\ell+1},\eqa

\bqa
V_1&=&\sum_{i=1}^{\ell}\epsilon_{i,i}+e^{-z_{11}}\epsilon_{\ell+1,\ell+1}+
e^{z_{11}}\epsilon_{\ell+2,\ell+2}+\sum_{i=\ell+3}^{2\ell+1}\epsilon_{i,i},
\eqa \bqa
\widetilde{V}_1&=&U_1^*=\sum_{i=1}^{\ell}\epsilon_{i,i}+
e^{-z_{11}}\epsilon_{\ell+1,\ell+1}+
e^{z_{11}}\epsilon_{\ell+2,\ell+2}+\sum_{i=\ell+3}^{2\ell+1}\epsilon_{i,i}+
e^{x_{11}}\epsilon_{\ell+1,\ell+2},\eqa
where $x_{\ell,k}=0,\,\,\, k=1,\ldots,\ell$ are assumed. Note that
$\widetilde{V}_i$, $\widetilde{U}_i$ can be considered as
off-diagonal deformations of $V_i$, $U_i$. Now we can define a
modified factorized representation for $N_+\subset SO(2\ell+1)$.

\begin{te}\label{Mparso2l1}

i)  The image of any generic unipotent element
  $v^{B_{\ell}}\in N_+$ in the tautological representation
$\pi_{2\ell+1}:\mathfrak{so}_{2\ell+1}\to End(\mathbb{C}_{2\ell+1})$
 can be represented in the form \bqa\label{par1}
\pi_{2\ell+1}(v^{B_{\ell}})=\mathfrak{X}_1\mathfrak{X}_{2}\cdots\mathfrak{X}_{\ell},\eqa
where \bqa\label{par2}
\mathfrak{X}_1=\widetilde{U}_1U_1^{-1}\widetilde{V}_1V_1^{-1},\eqa
\bqa\label{par3}
\mathfrak{X}_n=\widetilde{U}_nU_n^{-1}[\widetilde{U}_nU_n^{-1}]^*
\widetilde{V}_nV_n^{-1}[\widetilde{V}_nV_n^{-1}]^*,\qquad
n=2,\ldots,\ell \nonumber \eqa and $x_{\ell,k}=0\,\,$for
$k=1,\ldots\ell$ are assumed.

ii) This defines a parametrization of an open part $N_+^{(0)}$ of
$N_+$.
\end{te}

{\it Proof.} Let  $v^{B_{\ell}}(y)$ be a parametrization of an open
part of  $N_+$ according to (\ref{rec1B})-(\ref{rec2B}).
 Let $\tilde{X}_i(y)=e^{y e_{i,i+1}}$
be a  one-parametric unipotent subgroup in $GL(2\ell+1),$  then
$\tilde{X}_i(y)^*=\tilde{X}_{2\ell+1-i}(y).$ Embed an elementary
unipotent element $X_i(y)$ of $SO(2\ell+1)$ into $GL(2\ell+1)$ as
follows:
$$X_i(y)=\tilde{X}_i(y)^*\cdot\tilde{X}_i(y).$$ This maps
 an arbitrary regular unipotent element $v^{B_\ell}$
into unipotent subgroup of $GL(2\ell+1).$ Let us now change the
variables in the following way: \bqa\label{1yBn}
y_{11}=e^{x_{11}-z_{11}},\qquad
y_{1,k}=\Big(e^{x_{k-1,1}-z_{k,1}}+ e^{x_{k,1}-z_{k,1}}\Big),\\
\nonumber \label{2yBn}
y_{k,2r-1}=e^{z_{k+r-1,k}-x_{k+r-2,k-1}},\hspace{2cm}
k=2,\ldots,\ell,\\ \nonumber
 y_{k,2r}=e^{z_{k+r-1,k}-x_{k+r-1,k-1}},\hspace{2cm}
 r=1,\ldots,\ell+1-k,\eqa where the conditions
$x_{\ell,k}=0,\,\,k=1,\ldots\ell$ are assumed.
 By elementary operations it is easy to check that after the change of variables,
the image $\pi_{2\ell+1}(v^{B_{\ell}})$ of $v^{B_{\ell}}$ defined by
(\ref{rec1B})-(\ref{rec2B})  transforms into the (\ref{par1}).
 Taking into account that the change of variables
(\ref{1yBn}) is invertible
we get a parametrization of $N_+^{(0)}\subset N_+$ $\Box$


The modified factorized parameterization  of a unipotent group
 $N_+\subset SO(2\ell+1)$ defines a particular realization of a principal
series representation of $U(\mathfrak{so}_{2\ell+1})$ by
differential operators. It can be obtained using the change of
variables (\ref{1yBn}) applied to the representation given in
Proposition \ref{Bfactrep}. We shall use the term Gauss-Givental
representation  for this realization.

\begin{prop}
The following differential operators define a representation
$\pi_{\lambda}$ of $\mathfrak{so}_{2\ell+1}$ in $V_{\mu}$ in terms
of the modified factorized parametrization:
 \bqa E_1=-2\sum_{n=1}^\ell e^{z_{n,1}}\left(\frac{1}{2}
\frac{\partial}{\partial z_{n,1}}+
\frac{e^{x_{n,1}}}{e^{x_{n-1,1}}+e^{x_{n,1}}}
\frac{\partial}{\partial z_{n,2}}+\right.\\ \nonumber +\left.
\sum_{n=1}^{\ell-1} \Big(\frac{\partial}{\partial z_{\ell,1}} +
\frac{\partial}{\partial z_{\ell,2}}+ \frac{\partial}{\partial
x_{\ell,1}}\Big)\right),\eqa

\bqa E_2=\left(\frac{\partial}{\partial z_{11}}+
\frac{\partial}{\partial x_{11}}\right)\Big(
e^{x_{22}-z_{22}}+\sum_{k=3}^\ell e^{x_{k-1,2}-z_{k,2}}+
e^{x_{k,2}-z_{k,2}}\Big)+\\ \nonumber +\sum_{n=2}^\ell
\frac{\partial}{\partial\,z_{n,1}}\Big(
e^{x_{n-1,1}-z_{n,1}}\frac{e^{x_{n-1,2}-z_{n,2}}+
e^{x_{n,2}-z_{n,2}}}{e^{x_{n-1,1}-z_{n,1}}+ e^{x_{n,1}-z_{n,1}}}+
\sum_{k=n+1}^\ell e^{x_{k-1,2}-z_{k,2}}+ e^{x_{k,2}-z_{k,2}}\Big)+\\
\nonumber+\sum_{n=2}^\ell\left(\frac{\partial}{\partial z_{n,2}}-
\frac{\partial}{\partial z_{n,3}}\right)\Big(e^{x_{n,2}-z_{n,2}}+
\sum_{k=n+1}^\ell e^{x_{k-1,2}-z_{k,2}}+ e^{x_{k,2}-z_{k,2}}\Big)+\\
\nonumber +\sum_{n=2}^\ell\left(\frac{\partial}{\partial x_{n,1}}-
\frac{\partial}{\partial x_{n,2}}\right) \sum_{k=n+1}^\ell\Big(
e^{x_{k-1,2}-z_{k,2}}+ e^{x_{k,2}-z_{k,2}}\Big),\eqa

\bqa  E_k=\sum_{n=k-1}^{\ell-1}\left( \frac{\partial}{\partial
x_{n,k-1}}- \frac{\partial}{\partial x_{n,k}} \right)
\sum_{i=n+1}^\ell\Big(e^{x_{i-1,k}-z_{i,k}}+e^{x_{i,k}-z_{i,k}}\Big)
+\\ \nonumber +\sum_{n=k}^{\ell}\left( \frac{\partial}{\partial
z_{n,k}}- \frac{\partial}{\partial
z_{n,k+1}}\right)\Big(e^{x_{n,k}-z_{n,k}}+ \sum_{i=n+1}^\ell
e^{x_{i-1,k}-z_{i,k}}+e^{x_{i,k}-z_{i,k}}\Big),\quad 3\leq k\leq \ell,\eqa

\bqa H_k=\<\mu,\alpha_k^\vee\>+ \sum_{n=1}^\ell
a_{k,n}\sum_{i=n}^\ell\frac{\partial}{\partial z_{i,n}},\quad
1\leq k \leq \ell,  \eqa

\bqa F_1=\<\mu,\alpha_1^\vee\>\Big( e^{x_{11}-z_{11}}+
\sum_{k=2}^\ell e^{x_{k-1,1}-z_{k,1}}+ e^{x_{k,1}-z_{k,1}}\Big)+\\
\nonumber +\sum_{n=1}^\ell\Big(e^{x_{n,1}-z_{n,1}}-
e^{x_{n-1,1}-z_{n,1}}\Big)\frac{\partial}{\partial z_{n,1}}-\\
\nonumber +2\sum_{n=1}^\ell \frac{\partial}{\partial\,x_{n,1}}
\sum_{k=n+1}^\ell\Big(e^{x_{k-1,1}-z_{k,1}}+
e^{x_{k,1}-z_{k,1}}\Big),\eqa

\bqa F_2=\left(\<\mu,\alpha_2^\vee\>+ \frac{\partial}{\partial
x_{11}}\right) \sum_{k=2}^\ell\Big(e^{z_{k,2}-x_{k-1,1}}+
e^{z_{k,2}-x_{k,1}} \Big)-\\ \nonumber- \sum_{n=2}^\ell
\frac{\partial}{\partial\,z_{n,2}}\Big(e^{z_{n,2}-x_{n,1}}+
\sum_{k=n+1}^\ell e^{z_{k,2}-x_{k-1,1}}+ e^{z_{k,2}-x_{k,1}}\Big)+\\
\nonumber+ \sum_{n=2}^\ell\left(
\frac{\partial}{\partial\,x_{n,1}}-\frac{\partial}{\partial\,x_{n,2}}
\right)\sum_{k=n+1}^\ell\Big( e^{z_{k,2}-x_{k-1,1}}+
e^{z_{k,2}-x_{k,1}} \Big),\eqa

\bqa F_k=\Big( \<\mu,\alpha_k^\vee\>+\frac{\partial}{\partial
x_{k-1,k-1}}+ \frac{\partial}{\partial z_{k-1,k-1}}\Big)
\sum_{n=k}^\ell\Big( e^{z_{n,k}-x_{n-1,k-1}}+
e^{z_{n,k}-x_{n,k-1}}\Big)-\nonumber
\\ -\nonumber
\sum_{n=k}^\ell\left( \frac{\partial}{\partial z_{n,k}}-
\frac{\partial}{\partial
z_{n,k-1}}\right)\Big(e^{z_{n,k}-x_{n,k-1}}+ \nonumber
\sum_{i=n+1}^\ell e^{z_{i,k}-x_{i-1,k-1}}+
e^{z_{i,k}-x_{i,k-1}}\Big)+\\ +\nonumber \sum_{n=k}^\ell\left(
\frac{\partial}{\partial x_{n,k-1}}- \frac{\partial}{\partial
x_{n,k}}\right) \sum_{i=n+1}^\ell\Big( e^{z_{i,k}-x_{i-1,k-1}}+
e^{z_{i,k}-x_{i,k-1}}\Big),\qquad 3\leq k\leq\ell,\eqa where
$E_i=\pi_{\la}(e_i)$, $F_i=\pi_{\la}(f_i)$, $H_i=\pi_{\la}(h_i)$,
$x_{\ell,k}=0, k=1,\ldots,\ell$ are assumed and the derivatives over
$x_{i,k},\,\,\, z_{i,k},\,\,\,i < k$,
$x_{\ell,n},\,\,\,n=1,\ldots,\ell$ are omitted.
\end{prop}

We are going to write down the matrix element (\ref{pairing})
explicitly  in Gauss-Givental representation.
 Whittaker vectors $\psi_R$ and $\psi_L$ in this representation
 satisfy the system of differential
\bqa\label{bwitt} E_i\psi_R(x)= -\psi_R(x),\hspace{2cm}
F_i\psi_L(x)=-\psi_L(x),\qquad  1\leq i\leq\ell.\eqa Its solution
has the following form.
\begin{lem}\label{BWvector}
The functions \bqa\psi_L(x,z) = e^{2\mu_1x_{1,1}}\prod_{n=2}^\ell
\Big(e^{x_{n,1}}+e^{x_{n-1,1}}\Big)^{2\mu_n}
\times  \\
\nonumber\times \prod_{n=1}^\ell\exp\Big\{\,-\mu_n\Big(\sum_{i=1}^nx_{n,i}+
2z_{n,1}-2\sum_{i=2}^{n}z_{n,i}+\sum_{i=1}^{n-1}x_{n-1,i}\Big)\Big\}
\,\,\times \\\times  \nonumber\exp\Big\{-\Big(\sum_{k=1}^\ell e^{z_{k,1}}+
\sum_{k=2}^\ell e^{x_{k,k}-z_{k,k}}+
\sum_{k=2}^{\ell}\sum_{n=k+1}^\ell\Big(e^{x_{n-1,k}-z_{n,k}}+
e^{x_{n,k}-z_{n,k}}\Big)\Big) \Big\},\eqa

\bqa\psi_R(x,z)=\exp\Big\{-\Big( e^{x_{11}-z_{11}}+\sum_{n=2}^\ell
\Big(e^{x_{n-1,1}-z_{n,1}}+e^{x_{n,1}-z_{n,1}}\Big)\Big)+\\
\nonumber +\sum_{k=2}^\ell\,\,
\sum_{n=k}^\ell\Big(e^{z_{n,k}-x_{n-1,k-1}}+e^{z_{n,k}-x_{n,k-1}}
\Big)\Big\},\eqa  are solutions of the linear differential equations
(\ref{bwitt}).  We let $x_{\ell,k} =0\,\,\,$ for $k=1,\ldots,\ell,$
and $\mu_k=\imath\lambda_k-\rho_k,$ where $\rho_k=\frac{2k-1}{2}$.
\end{lem}

Now we are ready to find the integral representation of the pairing
(\ref{pairing}) in terms of modified factorization parameters. To
get explicit expression for the integrand, one uses the same type of
a Cartan element decomposition as in the case of
$\mathfrak{gl}_{\ell+1}$:
$$e^{-H_x}=\pi_\lambda(\exp(-\sum_{i=1}^{\ell}
\langle\omega_i,x\rangle h_i))=e^{H_L}e^{H_R},$$ where \bqa\label{HRso2l1}
-H_x=H_L+H_R=-\sum_{i=1}^{\ell}\langle\omega_i,x\rangle \langle
\mu,\alpha^{\vee}_i\rangle+\\ \nonumber +x_{\ell,1}\sum_{n=1}^\ell
\frac{\partial}{\partial z_{n,1}}+
\sum_{k=1}^{\ell-1}(x_{\ell,i}-x_{\ell,i+1})\sum_{n=k}^\ell
\frac{\partial}{\partial z_{n,k}},\eqa with  \bqa
H_L=\sum_{k=1}^{\ell-1}\sum_{n=1}^{k}x_{\ell,n}
\frac{\partial}{\partial x_{k,n}}+ \sum_{k=2}^\ell \sum_{n=2}^k
x_{\ell,n}\frac{\partial}{\partial z_{k,n}},\eqa \bqa H_R=-H_x-H_L
.\eqa
 We imply  that the
differential operator $H_L$ acts on the left  vector, and $H_R$ acts
on the right vector in (\ref{pairing}). Taking into account the
results of the Proposition \ref{BWvector} one obtains the following
theorem.

\begin{te}\label{teintmfso2l1}
 The eigenfunctions of $\mathfrak{so}_{2\ell+1}$-Toda chain (\ref{pairing})
admit the integral representation:  \bqa\label{iterwa}
\Psi^{B_\ell}_{\lambda_1,\ldots,\lambda_\ell}(x_{\ell,1},\ldots,x_{\ell,\ell})\,=\,
\int\limits_{C}
\bigwedge_{k=1}^{\ell-1}\bigwedge_{i=1}^kdx_{k,i}\bigwedge_{k=1}^{\ell}
\bigwedge_{i=1}^kdz_{k,i}\,\, \nonumber e^{{\mathcal
F}^{B_{\ell}}},\eqa where \bqa
\mathcal{F}^{B_{\ell}}=-\imath\lambda_1(-x_{1,1}+2z_{1,1})-\nonumber\\-
\imath\sum\limits_{n=2}^{\ell}\lambda_n \Big(\sum_{i=1}^nx_{n,i}+
2z_{n,1}-2\sum_{i=2}^{n}z_{n,i}+\sum_{i=1}^{n-1}x_{n-1,i}
-2\ln(e^{x_{n,1}}+e^{x_{n-1,1}})\Big)-
\\-\Big\{\sum_{n=1}^{\ell}e^{z_{n,1}}+
\sum_{k=2}^{\ell}\sum_{n=k+1}^\ell\Big(e^{x_{n-1,k}-z_{n,k}}+
e^{x_{n,k}-z_{n,k}}\Big)+\nonumber\\+\sum_{n=k}^\ell\Big
(e^{z_{n,k}-x_{n-1,k-1}}+e^{z_{n,k}-x_{n,k-1}} \Big)
+\sum_{n=1}^{\ell}e^{x_{n,n}-z_{n,n}}\Big\},\nonumber \eqa where we
set $x_i:=x_{\ell,i},\,\,\, 1\le i\leq \ell $. Here ${C} \subset
N_+$ is a middle-dimensional non-compact submanifold such that the
integrand decreases  exponentially at the boundaries and at
infinities. In particular the domain of integration can be chosen to
be ${C}=$$\RR^{m},$ where $m=l(w_0).$
\end{te}

\begin{ex} For $\ell=2$  the general formula
(\ref{iterwa}) has the following form: \bqa\label{psigg}
\Psi^{B_2}_{\lambda_1,\lambda_2}(x_{2,1},x_{2,2})= \int_C
dz_{1,1}\wedge dx_{1,1}\wedge dz_{2,1}\wedge dz_{2,2}\times   \\\times
\nonumber e^{-\imath\lambda_1(-x_{1,1}+2z_{1,1})-\imath\lambda_2
(2z_{2,1}-2z_{2,2}+x_{1,1}+x_{2,1}+x_{22})}
\Big(e^{x_{2,1}}+e^{x_{1,1}}\Big)^{2\imath\lambda_2}\times  \\\times
\nonumber \exp\Big\{-\Big(e^{z_{1,1}}+e^{z_{2,1}}
+e^{x_{2,2}-z_{2,2}}+ e^{x_{1,1}-z_{1,1}}+e^{x_{1,1}-z_{2,1}}+\\+
\nonumber
e^{x_{2,1}-z_{2,1}}+e^{z_{2,2}-x_{1,1}}+e^{z_{2,2}-x_{2,1}}\Big)\Big\},
\eqa where we set  $x_1=x_{21},x_2=x_{22}$ and the contour of
integration $C=\mathbb{R}_+^4$.
\end{ex}

There is a simple combinatorial description of the potential
$\mathcal{F}^{B_{\ell}}$ for zero spectrum $\{\la_i=0\}$. Namely, it
can be presented as the sum over all the arrows in the following
diagram. The diagram for $B_{\ell}$ reads

$$
\xymatrix{
 &&&& \circ\ar[d] &&&&\\
 &&& \circ\ar[r]& z_{\ell,1}\ar[d]\ar[r] &
 x_{\ell,1}\ar[d] &&&\\
 && \circ\ar[d] & \ldots & x_{\ell-1,1}\ar[r] &
 z_{\ell,2}\ar[d]\ar[r] & \ddots &&\\
 & \circ\ar[d]\ar[r] & z_{21}\ar[d]\ar[r] & \ddots\ar[d] &
 \ddots & \ddots & \ddots\ar[d]\ar[r] & x_{\ell,\ell-1}\ar[d] &\\
 \circ\ar[r] & z_{11}\ar[r] & x_{11}\ar[r] &
 z_{22}\ar[r] & \ldots & \ldots\ar[r] & x_{\ell-1,\ell-1}\ar[r] &
 z_{\ell,\ell}\ar[r] & x_{\ell,\ell}
}
$$
\hspace{1cm}

Here we use the same rules for assigning variables to the arrows of
the diagram as in $A_{\ell}$ case. In addition we assign functions
$e^x$ to the arrows $\circ\longrightarrow x$.

Note that the diagram for $B_{\ell}$  can be obtained  by a
factorization of the diagram (\ref{AnDiag}) for $A_{2\ell}$.
Consider the  following involution \be \label{inv} \iota\,:\quad
X\longmapsto  \dot{w}_0^{-1}X^t \dot{w}_0, \ee where  $\dot{w}_0$ is
the longest element of $A_{2\ell}$ Weyl group and $X^t$ denotes the
standard transposition. Corresponding action on the modified
factorization parameters is given by  \be \label{w0invBn}
w_0\,:\qquad x_{k,i}\longleftrightarrow -x_{k,k+1-i}. \ee This
defines a factorization of $A_{2\ell}$-diagram  that gives  the
diagram for $B_{\ell}$.

 An analog of  $\mathfrak{gl}_{\ell+1}$-monomial relations
(\ref{defrelAn}) can be described as  follows. Associate  variables
$a_{k,i}$, $b_{k,i}$, $c_{k,i}$, $d_{k,i}$ to the arrows of the
Givental diagram as \be
\,a_{k,i}=e^{z_{k,i}-x_{k-1,i-1}},\,\,\,b_{k,i}=e^{z_{k,i}-x_{k,i-1}},\,\,
\,c_{k,i}=e^{z_{k,i}-x_{k,i}},\,\,\,d_{l,j}=e^{x_{l,j}-z_{l+1,j}},\,\,\,\\
\,1\leq k\leq \ell,\,\,1\leq i\leq k, \quad 1\leq l\leq \ell-1,\,\,
1\leq j\leq l.\,\ee Then the following relations hold:
\bqa\label{Brel} a_{k,1}&=&b_{k,1},\qquad \qquad
\,\,\,\,\,\,\,\,\,\,\,\,\,\,\,\,\,\,\,\,1\leq k\leq \ell,\nonumber
\\
d_{k,i}\cdot a_{k+1,i+1}\,&=&\,c_{k+1,i}\cdot
b_{k+1,i+1},\,\,\,\qquad 1\leq k<\ell-1,\,\,1\leq i\leq k, \\
b_{k,i}\cdot c_{k,i}\,&=&\,
a_{k+1,i}\cdot d_{k,i},\qquad \qquad \,1\leq k<\ell-1,\,\,1\leq i\leq k,\nonumber\\
b_{\ell,i}\cdot
c_{\ell,i}&=&e^{x_{\ell,i}-x_{\ell,i-1}}.\nonumber\eqa The above
relations can be considered as relations between elementary paths on
the Givental diagram. Using a set of relations for more general
paths that follows from (\ref{Brel}) one can define a toric
degeneration of the $\mathfrak{so}_{2\ell+1}$ flag manifolds thus
generalizing the results in \cite{BCFKS}.

\subsubsection{Recursion for $\mathfrak{so}_{2\ell+1}$-Whittaker
  functions and  $\mathcal{Q}$-operator for \\  $B^{(1)}_{\ell}$-Toda
  chain }
The integral representation (\ref{iterwa}) of
$\mathfrak{so}_{2\ell+1}$-Whittaker functions
 possesses a remarkable recursive structure over the rank $\ell$.
Let us introduce integral operators $Q^{B_n}_{B_{n-1}}$,
$n=2,\ldots,\ell$  with the  kernels
$Q^{B_n}_{B_{n-1}}(\underline{x}_{n};\,\underline{x}_{n-1};\lambda_n)$
defined as follows \bqa\label{BQB}
Q^{B_n}_{B_{n-1}}(\underline{x}_{n};\,\underline{x}_{n-1};\lambda_n)=
\int\bigwedge_{i=1}^n
dz_{n,i}\,\,\,\Big(e^{x_{n,1}}+e^{x_{n-1,1}}\Big)^{2\imath\lambda_n}
\times  \\
\nonumber\times
\exp\Big\{\,-\imath\lambda_n\Big(\sum_{i=1}^nx_{n,i}+
2z_{n,1}-2\sum_{i=2}^{n}z_{n,i}+\sum_{i=1}^{n-1}x_{n-1,i}\Big)\Big\}
\,\,\times \\ \nonumber \times
Q^{B_n}_{BC_{n}}(\underline{x}_{n};\,\underline{z}_{n})\,\,\,
Q^{BC_{n}}_{B_{n-1}}(\underline{z}_{n};\,\underline{x}_{n-1}),\quad\eqa
where \be Q^{BC_n}_{\,\,\,\,B_{n-1}}(\underline{z}_n
;\,\underline{x}_{n-1})= \exp\Big\{-\Big(
\frac{1}{2}e^{z_{n,1}}+\sum_{i=1}^{n-1}\Big(e^{x_{n-1,i}-z_{n,i}}+
e^{z_{n,i+1}-x_{n-1,i}}\Big)\,\Big)\Big\},\label{kerBCB}\ee
\be Q^{B_n}_{\,\,\,\,BC_{n}}(\underline{x}_n ;\,\underline{z}_n)=\\
=\exp\Big\{\,-\Big(
\frac{1}{2}e^{z_{n,1}}+\sum_{i=1}^{n-1}\Big(e^{x_{n,i}-z_{n,i}}+
e^{z_{n,i+1}-x_{n,i}}\Big)+e^{x_{n,n}-z_{n,n}}\,\Big)\Big\}.\label{kerBBC}\ee
We set for $n=1$
$$Q^{B_1}_{B_{0}}(x_{1,1};\lambda_1)=\int dz_{1,1}e^{\imath\lambda_1
x_{1,1}-2\imath\lambda_1
z_{1,1}}\exp\Big\{-\Big(e^{z_{1,1}}+e^{x_{1,1}-z_{1,1}}\Big)\Big\}.$$
Using integral operators $Q^{B_n}_{B_{n-1}}$
 the integral representation (\ref{iterwa}) can be written   in a
recursive form.
\begin{te} The eigenfunction of $B_{\ell}$-Toda chain can be
written as \bqa\label{iterwa1}
\Psi^{B_\ell}_{\lambda_1,\ldots,\lambda_\ell}(x_1,\ldots,x_\ell)\,=\,
\int\limits_{C} \bigwedge_{k=1}^{\ell-1}\bigwedge_{i=1}^kdx_{k,i}\,
\nonumber
\prod_{k=1}^{\ell}Q^{B_{k}}_{\,\,B_{k-1}}(\underline{x}_{k};\,
\underline{x}_{k-1};\lambda_{k}),\eqa or equivalently
\bqa\label{iterwa2}
\Psi^{B_\ell}_{\lambda_1,\ldots,\lambda_\ell}(x_{1},\ldots,x_{\ell})\,=\,
\int\limits_{C} \bigwedge_{i=1}^{\ell-1}dx_{\ell-1,i} \,
 Q^{B_{\ell}}_{\,\,B_{\ell-1}}(\underline{x}_{\ell};\,\underline{x}_{\ell-1};
\lambda_{\ell})
\Psi^{B_{\ell-1}}_{\lambda_1,\ldots,\lambda_{\ell-1}}
(\underline{x}_{\ell-1}),\eqa where we assume
$x_n:=x_{\ell,n},\,\,\,  1\le n\leq \ell $ and ${C} \subset N_+$ is
a middle-dimensional non-compact submanifold such that the integrand
decreases exponentially at possible  boundaries and at infinities.
In particular as the domain of integration one can chose
${C}=$$\RR^{m},$ where $m=l(w_0).$
\end{te}
Let us note that in contrast with the case of
$\mathfrak{gl}_{\ell+1}$ integral
 representations, kernels of $Q^{B_n}_{B_{n-1}}$, $n=1,\ldots,\ell$
have more complicated  form. Curious new structure appeares if we
consider the Whittaker functions for zero spectrum\footnote{\ Note
that the zero spectrum Whittaker functions are directly related to
the  quntum cohomology of flag manifolds in Givental description.}
 $\{\la_i=0\}$. As it is clear from (\ref{BQB})
the kernel of $Q^{B_n}_{B_{n-1}}$ is given by a  convolution of two
kernels
 $Q^{B_n}_{BC_{n}}(\underline{x}_{n};\,\underline{z}_{n})$ and
$Q^{BC_{n}}_{B_{n-1}}(\underline{z}_{n};\,\underline{x}_{n-1})$.
Corresponding integral operators  $Q^{B_n}_{BC_{n}}$,
$Q^{BC_{n}}_{B_{n-1}}$ can be regarded as  elementary intertwiners
relating Toda chains for $B_n$, $BC_n$ and $BC_n$, $B_{n-1}$ root
systems. $BC_{\ell}$-Toda chain\footnote{ $BC_{\ell}$-Toda chain can
be also considered as a most general form of $C_{\ell}$-Toda chain
(see e.g. \cite{RSTS}, Remark p.61). In
 the following we will use the term $BC_{\ell}$-Toda chain to
distinguish it from a more standard $C_{\ell}$-Toda chain that will
be  consider below.}  is defined in terms of the non-reduced root
system $BC_{\ell}$ in a standrad fashion. Let us recall the
construction of the non-reduced root system $BC_{\ell}$. Root system
of $BC_{\ell}$ type  can be realized in terms of an orthogonal bases
$\{\epsilon_i\}$ in $\R^{\ell}$ as \bqa \alpha_0=2\epsilon_1,\qquad
\alpha_1=\epsilon_1,\qquad
\alpha_{i+1}=\epsilon_{i+1}-\epsilon_i,\qquad 1\leq i\leq \ell-1,
\eqa and the corresponding Dynkin diagram is
$$
\xymatrix{
 {\alpha_0\atop\alpha_1}\ar@{<=>}[r] & \alpha_2\ar@{-}[r] &
 \ldots\ar@{-}[r] & \alpha_{\ell-1}\ar@{-}[r] & \alpha_\ell
}$$
where the first vertex from the left is a doubled vertex
corresponding to a
 reduced $\alpha_1=\epsilon_1$ and non-reduced
 $\alpha_0=2\epsilon_1$ roots. Then for example the quadratic
Hamiltonian operator of $BC_{\ell}$-Toda chain is given by
 \bqa \CH_2^{BC_{\ell}}(\underline{x}^{(\ell)})=-
\frac{1}{2}\sum_{i=1}^{\ell}\frac{\partial^2}{\partial x_i^2}+
\frac{1}{4}\Big(e^{x_1}+\frac{1}{2}e^{2x_1}\Big)+
\sum_{i=1}^{\ell-1}e^{x_{i+1}-x_i}. \eqa Integral operators
$Q^{B_n}_{BC_{n}}$ and $Q^{BC_{n}}_{B_{n-1}}$ intertwine Hamitonian
operators of different Toda chains. Thus for quadratic Hamiltonians
one can directly check the following relations.

\begin{prop}1. The operators
$Q^{B_n}_{BC_{n}}$ and $Q^{BC_{n}}_{B_{n-1}}$ defined by the kernels
(\ref{kerBCB}), (\ref{kerBBC}) intertwine quadratic Hamiltonians of
$B$ and $BC$ Toda chains:
\bqa\CH_2^{BC_{n}}(\underline{z}_n)Q^{BC_{n}}_{B_{n-1}}(\underline{z}_n,\,
\underline{x}_{n-1})= Q^{BC_n}_{B_{n-1}}
(\underline{z}_n,\,\underline{x}_{n-1})\CH_2^{B_{n-1}}(\underline{x}_{n-1}),
\eqa
\bqa\CH_2^{B_{n}}(\underline{x}_n)Q^{B_{n}}_{BC_{n}}(\underline{x}_n,\,
\underline{z}_{n})= Q^{B_n}_{BC_{n}}
(\underline{x}_n,\,\underline{z}_{n})\CH_2^{BC_{n}}(\underline{z}_{n}).
\eqa 2. Integral operator $Q^{B_{n}}_{B_{n-1}}$ at $\la_n=0$
intertwines  Hamiltonians ${\CH}_2^{B_n}$ and ${\CH}_2^{B_{n-1}}:$
\bqa\CH_2^{B_{n}}(\underline{x}_n)Q^{B_{n}}_{B_{n-1}}(\underline{x}_n,\,
\underline{x}_{n-1};\la_n=0)= Q^{B_n}_{B_{n-1}}
(\underline{x}_n,\,\underline{x}_{n-1};\la_n=0)\CH_2^{B_{n-1}}(\underline{x}_{n-1}).
\eqa
\end{prop}
The kernel $Q^{B_n}_{B_{n-1}}(\underline{x}_n,\,\underline{x}_{n-1};\la_n=0)$:
can be succinctly encoded into the following sub-diagramm of
$\mathfrak{so}_{2\ell+1}$ Givental diagram

\bqa
\label{Bsubdiag}
\xymatrix{
 & \circ\ar[d] &&&&\\
 \circ\ar[r] & z_{n,1}\ar[d]\ar[r] & x_{n,1}\ar[d] &&&\\
 & x_{n-1,1}\ar[r] & z_{n,2}\ar[d]\ar[r] & \ddots &&\\
 && \ddots & \ddots\ar[d]\ar[r] & x_{n,n-1}\ar[d] &\\
 &&& x_{n-1,n-1}\ar[r] & z_{n,n}\ar[r] & x_{n,n}
}\eqa

Here the upper and lower descending paths of the oriented diagram
correspond to the kernels of elementary intertwiners
$Q^{B_n}_{BC_n}$ and $Q_{B_{n-1}}^{BC_n}$ respectively. The
convolution of the kernels $Q^{B_n}_{BC_n}$ and
$Q_{B_{n-1}}^{BC_n}$ in (\ref{BQB}) at $\la_n=0$ corresponds to the
integration over the variables $z_{n,i}$ associated with the inner
vertexes of the subdiagram (\ref{Bsubdiag}).

Similarly to the case of $\mathfrak{gl}_{\ell+1}$, recursion
operators $Q^{B_n}_{B_{n-1}}$ can be considered as particular
degenerations of  Baxter $\mathcal{Q}$-operators for affine
$B^{(1)}_{\ell}$-Toda chains. Below we provide the integral
representations for these $\mathcal{Q}$-operator. Let us stress that
up to now $\mathcal{Q}$-operators were known only for
$\widehat{\mathfrak{gl}}_{\ell+1}$-case. We will not present here
the complete set of  properties characterizing the introduced
$\mathcal{Q}$-operators and only consider the commutation relations
with quadratic affine Toda chain Hamiltonians. The detailed account
will be given elsewhere.

We start with a description of $B_{\ell}^{(1)}$-Toda chain. The set
of simple roots of the affine root system $B^{(1)}_{\ell}$ can be
represented in the following form:
 \be \alpha_1=\epsilon_1,\qquad \alpha_{i+1}=\epsilon_{i+1}-\epsilon_i,\qquad
1\leq i\leq \ell-1\qquad
\alpha_{\ell+1}=-\epsilon_{\ell}-\epsilon_{\ell-1}. \ee The
corresponding Dynkin diagram is
$$
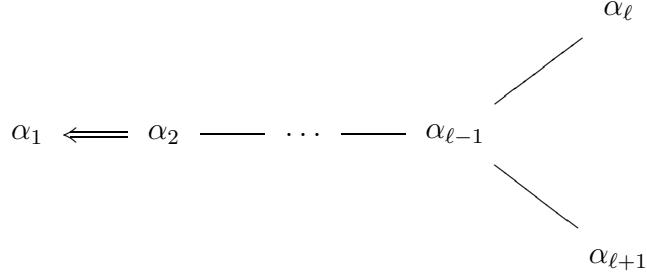
\xymatrix{
 &&&& \alpha_{\ell}\\
 \alpha_1 & \alpha_2\ar@{=>}[l]\ar@{-}[r] & \ldots\ar@{-}[r] &
 \alpha_{\ell-1}\ar@{-}[ur]\ar@{-}[dr] &\\
 &&&& \alpha_{\ell+1}
}
$$

These root data allows to define affine $B^{(1)}_{\ell}$-Toda chain
with a quadratic Hamiltonian given by \be
\CH_2^{B_{\ell}^{(1)}}=-\frac{1}{2}
\sum\limits_{i=1}^{\ell+1}\frac{\partial^2}{\partial x_i^2}+
\frac{1}{2}e^{x_1}+\sum\limits_{i=1}^{\ell-1}
e^{x_{i+1}-x_{i}}+ge^{-x_{\ell}-x_{\ell-1}}\,. \ee Here $g$ is an
arbitrary coupling constant.

Define the Baxter $\mathcal{Q}$-operator of
 $B_{\ell}^{(1)}$-Toda chain as  an integral operator with the following kernel
 \bqa\label{QBone}
\mathcal{Q}^{B_{\ell}^{(1)}}
(\underline{x}^{(\ell)},\underline{y}^{(\ell)}\,,\,\lambda)=\int\bigwedge_{i=1}^{\ell}
dz_{i}\,\,\,\Big(e^{x_{1}}+e^{y_{1}}\Big)^{2\imath\lambda}
\Big(e^{-x_{\ell}}+e^{-y_{\ell}}\Big)^{-2\imath\lambda}\times  \\
\nonumber\times \exp\Big\{\,-\imath\lambda\Big(\sum_{i=1}^\ell
x_{i}+
2z_{1}-2\sum_{i=2}^{\ell}z_{i}+\sum_{i=1}^{\ell}y_{i}\Big)\Big\}
\,\,  Q^{B^{(1)}_\ell}_{BC^{(1)}_\ell}(x_{i};\,z_{i})\,\,\,
Q^{BC^{(1)}_\ell}_{B^{(1)}_\ell}(z_{i};\,y_{i}),\quad\eqa where
 \be\label{inttwBn}
Q_{B^{(1)}_{\ell}}^{\,\,\,\,\,BC^{(1)}_{\ell}}(z_i,\,y_i)=\\
=\exp\Big\{-\Big(\,\frac{1}{2} e^{z_1}+ \sum_{i=1}^{\ell-1}\Big(
e^{y_i-z_i}+e^{z_{i+1}-y_i}\Big)+
e^{y_{\ell}-z_{\ell}}+ge^{-y_{\ell}-z_{\ell}}\Big)\,\Big\},\ee and
\be
Q^{\,\,\,\,B^{(1)}_{n}}_{BC^{(1)}_{n}}(x_i,\,z_i)=Q_{B^{(1)}_{n}}
^{\,\,\,\,BC^{(1)}_{n}}(z_i,\,x_i).\ee
Here we denote $\underline{x}^{(\ell)}=(x_1,\ldots ,x_{\ell})$, and 
$\underline{y}^{(\ell)}=(y_1,\ldots ,y_{\ell})$.

 The following Proposition can be proved by a direct check.
\begin{prop}
 The $\mathcal{Q}$-operator (\ref{QBone}) commutes with quadratic
Hamiltonian  of the  $B_\ell^{(1)}$ Toda chain, that is the kernel
intertwines the Hamiltonians $\CH_2^{B_{\ell}^{(1)}}$
\bqa\CH_2^{B_{\ell}^{(1)}}(\underline{x}^{(\ell)})
\mathcal{Q}^{^{B_{\ell}^{(1)}}}(\underline{x}^{(\ell)}\,;\,
\underline{y}^{(\ell)}\,,\,\lambda)= \mathcal{Q}^{^{B_{\ell}^{(1)}}}
(\underline{x}^{(\ell)}\,,\,\underline{y}^{(\ell)}\,;\,\lambda)\CH_2^{{B_{\ell}^{(1)}}}
(\underline{y}^{(\ell)}). \eqa
\end{prop}
Now we will demonstrate that recursion operator
$Q^{B_{\ell}}_{B_{\ell-1}}$ can be considered as  a degeneration of
Baxter $\mathcal{Q}$-operators for $B_{\ell}^{(1)}$. Let us
introduce a slightly modified recursion operator
$Q^{B_{\ell}}_{B_{\ell-1}\oplus B_{1}}$ with the kernel: \bqa\label{modrecop}
Q^{B_{\ell}}_{B_{\ell-1}\oplus B_{1}}
(\underline{x}^{(\ell)}\,,\,\underline{y}^{(\ell)}\,,
\,\lambda):=e^{\imath\lambda y_{\ell}}\,\,
Q^{B_{\ell}}_{B_{\ell-1}}(\underline{x}^{(\ell)}\,,\,\underline{y}^{(\ell-1)}\,,
\,\lambda),
\eqa
where $\underline{y}^{(\ell-1)}=(y_1,\ldots ,y_{\ell-1})$.
Operator (\ref{modrecop})  intertwines 
Hamiltonians of $\mathfrak{so}_{2\ell+1}$- and
 $\mathfrak{so}_{2\ell-1}\oplus \mathfrak{so}_2$-Toda chains. Thus
 for quadratic Hamiltonians we have
 \bqa\nonumber \CH_2^{B_{\ell}}(\underline{x}^{(\ell)})
Q^{B_{\ell}}_{B_{\ell-1}\oplus B_{1}}
(\underline{x}^{(\ell)}\,,\,\underline{y}^{(\ell)},\la)=
 Q^{B_{\ell}}_{B_{\ell-1}\oplus B_{1}}
(\underline{x}^{(\ell)}\,,\,\underline{y}^{(\ell)},\la)
\Big(\CH_2^{B_{\ell-1}}(\underline{y}^{(\ell-1)})+
\CH_2^{B_1}(y_{\ell})\Big), \eqa where
$\CH_2^{B_1}(y_{\ell})=-\frac{1}{2}\Big(\partial^2/\partial
y_{\ell}^2\Big).$ Obviously the projection of the above relation  on the
subspace of functions
$F(\underline{y}^{(\ell)},y_{\ell})=\exp(\imath\lambda y_{\ell})
f(\underline{y}^{(\ell-1)})$ recovers the  genuine recursion operator
satisfying:
\bqa\CH_2^{B_{\ell}}(\underline{x}^{(\ell)})
Q^{B_{\ell}}_{B_{\ell-1}}(\underline{x}^{(\ell)}\,,\,\underline{y}^{(\ell)}\,;\,\lambda)
= Q^{B_{\ell}}_{B_{\ell-1}}(\underline{x}^{(\ell)}\,,\,
\underline{y}^{(\ell-1)}\,,\,\lambda)
\Big(\CH_2^{B_{\ell-1}}
(\underline{y}^{(\ell-1)})+\frac{1}{2}\lambda^{2}\Big). \eqa Let us
introduce a one-parameter family of the kernels \bqa\nonumber
\mathcal{Q}^{^{B_{\ell}^{(1)}}}(\underline{x}^{(\ell)},\,
\underline{y}^{(\ell)};\lambda;\varepsilon)=
\varepsilon^{\imath\lambda}e^{\imath\lambda y_{\ell}}
\int\bigwedge_{i=1}^{\ell}
dz_{i}\,\,\,\Big(e^{x_{1}}+e^{y_{1}}\Big)^{2\imath\lambda}
\Big(\varepsilon e^{y_{\ell}-x_{\ell}}+1\Big)^{-2\imath\lambda}\times  \\
\times \exp\Big\{\,-\imath\lambda\Big(\sum_{i=1}^\ell x_{i}+
2z_{1}-2\sum_{i=2}^{\ell}z_{i}+\sum_{i=1}^{\ell-1}y_{i}\Big)\Big\}
\,\,  Q^{B^{(1)}_\ell}_{BC^{(1)}_\ell}(x_{i};\,z_{i})\,\,\,
Q^{BC^{(1)}_\ell}_{B^{(1)}_\ell}(z_{i};\,y_{i};\varepsilon),\quad\eqa
where \be\label{inttwBnde}
Q_{B^{(1)}_{\ell}}^{\,\,\,\,\,BC^{(1)}_{\ell}}(\underline{z}_{\ell},\,
\underline{y}_{\ell};\varepsilon) =\exp\Big\{-\Big(\,\frac{1}{2}
e^{z_1}+ \sum_{i=1}^{\ell-1}\Big(
e^{y_i-z_i}+e^{z_{i+1}-y_i}\Big)+\\+ \nonumber \varepsilon
e^{y_{\ell}-z_{\ell}}+\varepsilon^{-1}ge^{-y_{\ell}-z_{\ell}}\Big)\,\Big\},
\nonumber\ee obtained from the kernel of the operator 
$\mathcal{Q}^{B^{(1)}_{\ell}}$  by the change of the
variable $y_{\ell}=y_{\ell}+\ln\varepsilon.$ Consider limiting
behavior of (\ref{inttwBnde}) when $\varepsilon\rightarrow 0,\,\,\,
g\varepsilon^{-1}\rightarrow 0.$ Then the following relation between
$\mathcal{Q}$-operator for $B_{\ell}^{(1)}$-Toda chain and
(modified) recursion operator for
$\mathfrak{so}_{2\ell+1}$-Whittaker function holds \bqa
Q^{B_{\ell}}_{B_{\ell-1}\oplus
B_{1}}(\underline{x}^{(\ell)},\underline{y}^{(\ell)};\lambda)=
\lim_{\varepsilon\rightarrow 0,\,\,g\epsilon^{-1}\rightarrow 0}
\varepsilon^{-\imath\lambda}
\mathcal{Q}^{^{B_{\ell}^{(1)}}}(\underline{x}^{(\ell)},\,
\underline{y}^{(\ell)};\lambda;\varepsilon). \eqa


\newpage

\subsection{ Integral representations of $\mathfrak{sp}_{2\ell}$-Toda
  chain eigenfunctions }

In this subsection we provide an analog of  the Givental integral
representation of Whittaker functions for $\mathfrak{sp}_{2\ell}$
Lie algebras. As in the case of $\mathfrak{so}_{2\ell+1}$, we start
with a derivation of the integral representation of
$\mathfrak{sp}_{2\ell}$-Whittaker functions using the factorized
parametrization. Then we consider a
modification of  the factorized parametrization  leading  to a
Givental type integral representation of
$\mathfrak{sp}_{2\ell}$-Whittaker functions.

Consider $C_{\ell}$ type root system corresponding to a  Lie algebra
$\mathfrak{sp}_{2\ell}$. Let
 $(\epsilon_1,\ldots,\epsilon_{\ell})$
be an orthogonal basis in $\RR^{\ell}.$ We use the following
realization of simple roots, coroots and fundamental weights as
vectors in $\RR^{\ell}$:
 \bqa\label{rootc}
\begin{array}{l}
\alpha_1=2\epsilon_1,\\
\alpha_2=\epsilon_2-\epsilon_1,\\ \ldots\\
\alpha_{\ell}=\epsilon_{\ell}-\epsilon_{\ell-1},
\end{array}
\hspace{2cm}
\begin{array}{l}
\alpha^{\vee}_1=\epsilon_1,\\
\alpha^{\vee}_2=\epsilon_2-\epsilon_1,\\ \ldots\\
\alpha^{\vee}_{\ell}=\epsilon_{\ell}-\epsilon_{\ell-1},
\end{array}
\hspace{2cm}
\begin{array}{l}
\omega_1=\epsilon_1+\ldots+\epsilon_{\ell}\\
\omega_2=\epsilon_2+\ldots+\epsilon_{\ell},\\ \ldots\\
\omega_{\ell}=\epsilon_{\ell}.
\end{array} \eqa
Cartan matrix $\|a_{ij}\|=\|\<\alpha^{\vee}_i,\alpha_j\>\|$ can be
made symmetric $\|b_{ij}\|=\|d_ia_{ij}\|$ with $d_1=2$, $d_{i}=1$,
$i=2,\ldots, \ell$.  One associates with these root data a
$\mathfrak{sp}_{2\ell}$-Toda chain with a quadratic Hamiltonian
given by \bqa\label{CquadHam}
\CH_2^{C_{\ell}}&=&-\frac{1}{2}\sum\limits_{i=1}^{\ell}
\frac{\partial^2}{{\partial z_i}^2}+
2e^{2z_1}+\sum\limits_{i=1}^{\ell-1} e^{z_{i+1}-z_{i}} .\eqa One can
complete  (\ref{CquadHam}) to a full set of $\ell$
 mutually commuting operators $H^{C_{\ell}}_k$ of
 $C_{\ell}$-Toda chain. We are looking for integral representations of
 common eigenfunctions of the full commuting set of Hamiltonians.
The corresponding eigenfunction problem for quadratic Hamiltonian
can be written in the following form  \be
\CH_2^{C_{\ell}}\,\,\,\Psi^{{C}_{\ell}}_{\la_1,\cdots
,\la_{\ell}} (z_1,\ldots,z_{\ell})=
\frac{1}{2}\sum\limits_{i=1}^{\ell}\lambda_{i}^2\,\,\,
\Psi^{{C}_{\ell}}_{\la_1,\cdots ,\la_{\ell}} (z_1,\ldots,z_{\ell}).
\ee

\subsubsection{ $\mathfrak{sp}_{2\ell}$-Whittaker function:
  factorized   parametrization}

The reduced word for the  maximal length element $w_0$ in the Weyl
group of $\mathfrak{sp}_{2\ell}$ can be represented in the recursive
form:
$$I=(i_1,i_2,\ldots,i_m):=
(1,212,32123,\ldots,(\ell\ldots212\ldots\ell)),$$
where indexes $i_k$ correspond  to elementary reflections with
respect to the  roots $\alpha_k$. Let $N_+\subset G$ be a maximal
unipotent subgroup of $G=Sp(2\ell)$. One associates with the reduced
word $I$ the following recursive parametrization of a generic
element  $v^{C_{\ell}}\in N_+$: \bqa\label{rec1C}
v^{C_{\ell}}=v^{C_{\ell-1}}\cdot\mathfrak{X}^{C_{\ell}}_{C_{\ell-1}}
,\eqa where
\bqa\label{rec2C}\mathfrak{X}^{C_\ell}_{C_{\ell-1}}=X_{\ell}(y_{\ell,1})\,
X_k(y_{k,2(\ell+1-k)-1})\, X_2(y_{2,2\ell-3})\times\\ \times
\nonumber X_1(y_{1,\ell})X_2(y_{2,2\ell-2})\,
X_k(y_{k,2(\ell+1-k)})\cdot X_{\ell}(y_{\ell,2}) .\eqa Here
$X_i(y)=e^{ye_i}$ and  $e_i\equiv e_{\alpha_i}$ are simple root
generators. The subset  $N_+^{(0)}$ allowing representation
(\ref{rec1C}), (\ref{rec2C})  is an open part of $N_+$. The action
of the Lie algebra $\mathfrak{sp}_{2\ell}$ on $N_+$ given by
(\ref{infact}) defines an action on the space of functions on
$N_+^{(0)}$. The following proposition explicitly describes the
action on the space $V_{\mu}$ of (twisted) functions on $N_+^{(0)}$.

\begin{prop}\label{realCn}
The following differential operators define a realization of the
representation $\pi_{\lambda}$ of
$\mathcal{U}(\mathfrak{sp}_{2\ell})$ in $V_{\mu}$ in terms of
factorized parametrization of $N^{(0)}_+$: \bqa
E_1&=&\sum_{n=1}^\ell\left(\frac{\partial}{\partial y_{1,n}}-
\frac{\partial}{\partial y_{1,n+1}}\right)
\prod_{j={n}}^{\ell-1} \Big(\frac{y_{2,2j}}{y_{2,2j-1}}\Big)^2\,+\nonumber\\
&+&\sum_{n=1}^{\ell-1}\left(\frac{\partial}{\partial y_{2,2n-1}}-
\frac{\partial}{\partial y_{2,2n}}\right) \frac{y_{2,2n}}{y_{1,n}}
\Big(1+\frac{y_{2,2n}}{y_{2,2n-1}}\Big)
\prod_{j=n+1}^{\ell-1}\Big(\frac{y_{2,2j}}{y_{2,2j-1}}\Big)^2,\\
E_k&=&\sum_{n=1}^{\ell+1-k}\left( \frac{\partial}{\partial
y_{k,2n}}- \frac{\partial}{\partial y_{k,2n+1}}\right)
\prod_{i=n}^{\ell-k}\frac{y_{k+1,2j}}{y_{k+1,2j-1}}
\frac{y_{k,2(j+1)-1}}{y_{k,2(j+1)}}\,+\nonumber\\ \nonumber
&+&\sum_{n=1}^{\ell-k}\left( \frac{\partial}{\partial y_{k+1,2n-1}}-
\frac{\partial}{\partial y_{k+1,2n}}\right)
\frac{y_{k+1,2n}}{y_{k,2(n-1)}}
\prod_{i=n}^{\ell-k}\frac{y_{k+1,2j}}{y_{k+1,2j-1}}
\frac{y_{k,2(j+1)-1}}{y_{k,2(j+1)}},\quad 1<k<\ell,\nonumber \\
E_\ell&=&\frac{\partial}{\partial y_{\ell,2}},\nonumber\eqa

\bqa H_k=\<\mu\,,\alpha_k^\vee\>\,+\,\sum_{i=1}^\ell
a_{k,i}\sum_{j=1}^{n_i}y_{i,j}\frac{\partial}{\partial y_{i,j}},\eqa

\bqa F_1 &=& \sum_{n=1}^\ell\,y_{1,n}\Big(-\<\mu,\alpha_1^\vee\>+
\sum_{j=1}^{2(n-1)-1}y_{2,j}\frac{\partial}{\partial y_{2,j}}
-2\sum_{j=1}^{n-1}y_{1,j}\frac{\partial}{\partial y_{1,j}}-
y_{1,n}\frac{\partial}{\partial y_{1,n}}\Big),\nonumber\\
F_2 &=&\sum_{n=1}^{2(\ell-1)}\,y_{2,n}\Big(\<\mu,\alpha_2^\vee\>+
2\sum_{j=1}^{[n/2]+1}y_{1,j}\frac{\partial}{\partial y_{1,j}}
-2\sum_{j=1}^{n-1}y_{2,j}\frac{\partial}{\partial y_{2,j}}+\nonumber\\
\nonumber &+&\sum_{j=1}^{2[(n+1)/2]-3}
y_{3,j}\frac{\partial}{\partial y_{3,j}}-
y_{2,n}\frac{\partial}{\partial y_{2,n}}\Big),\eqa

\bqa F_k\,=\,\sum_{n=1}^{2(\ell+1-k)}\,y_{k,n}\Big(
\<\mu,\alpha_k^\vee\>+
2\sum_{j=1}^{2[n/2]+1}y_{k-1,j}\frac{\partial}{\partial y_{k-1,j}}
-2\sum_{j=1}^{n-1}y_{k,j}\frac{\partial}{\partial y_{k,j}}+\\
\nonumber +\sum_{j=1}^{2[(n+1)/2]-3}
y_{k+1,j}\frac{\partial}{\partial y_{k+1,j}}-
y_{k,n}\frac{\partial}{\partial y_{k,n}}\Big), \eqa for $2<k<\ell$
\bqa F_\ell\,=\,(y_{\ell,1}+y_{\ell,2})\left(
-\<\mu,\alpha_\ell^\vee\>+ y_{\ell-1,1}\frac{\partial}{\partial
y_{\ell-1,1}}+
y_{\ell-1,2}\frac{\partial}{\partial y_{\ell-1,2}}\right)+\\
\nonumber +y_{\ell,2}\left( y_{\ell-1,3}\frac{\partial}{\partial
y_{\ell-1,3}}+y_{\ell-1,4}\frac{\partial}{\partial
y_{\ell-1,4}}\right)-\left(y_{\ell,1}^2\frac{\partial}{\partial
y_{\ell,1}}+ 2y_{\ell,1}y_{\ell,2}\frac{\partial}{\partial
y_{\ell,1}}+ y_{\ell,2}^2\frac{\partial}{\partial y_{\ell,2}}
\right),\eqa where $E_i=\pi_{\lambda}(e_i), H_i=\pi_{\lambda}(h_i),
F_i=\pi_{\lambda}(f_i),\,\,\,i=1,\ldots,\ell.$ and $n_1=\ell$,
$n_k=2(\ell+1-k)$ for $1<k\leq\ell$.
\end{prop}
The proof is given in Part II,  Section \ref{gensp2l}. 

For left/right Whittaker vectors in the factorized parametrization
we have the following expressions.

\begin{lem}
Left/right Whittaker vectors in the factorized parametrization are
given by: \bqa \psi_R(y)= \exp\Big\{-\Big(\sum_{n=1}^\ell
y_{1,n}\,+\,\sum_{k=2}^\ell
\sum\limits_{n=1}^{n_k}y_{k,n}\Big)\Big\},\nonumber\eqa \bqa\label{psirc}
\psi_L(y)= \prod_{i=1}^\ell\Big(\prod_{n=1}^\ell y_{1,n}\times
\prod_{k=2}^i\prod_{n=1}^{2(\ell+1-k)}y_{k,n}\times
\prod_{k=i+1}^\ell\prod_{n=1}^{\ell+1-k}y^2_{k,2n-1} \Big)
^{\<\mu,\alpha_i^{\vee}\>}\times \\
\nonumber \times \exp\Big\{ -\Big(\sum_{n=1}^\ell\frac{1}{y_{1,n}}
\Big(1+\frac{y_{2,2(n-1)}}{y_{2,2(n-1)-1}}\Big)^2\prod_{i=n+1}^\ell
\Big(\frac{y_{2,2(i-1)}}{y_{2,2(i-1)-1}}\Big)^2+\\\nonumber 
+\sum_{k=2}^\ell \sum_{n=1}^{n_k/2} \frac{1}{y_{k,2n}}
\Big(1+\frac{y_{k+1,2(n-1)}}{y_{k+1,2(n-1)-1}}\Big)
\prod_{i=n+1}^{n_k/2} \frac{y_{k+1,2(i-1)}}{y_{k+1,2(-1)-1}}
\frac{y_{k,2i-1}}{y_{k,2i}}\Big)\Big\}, \eqa where  $n_1=\ell$ and
$n_k=2(\ell+1-k),$ for $k=2,\ldots,\ell$.
\end{lem}
Proof is given in Part II, Section \ref{melsp2l}.

Using the expressions (\ref{psirc}) for the left/right  Whittaker
vectors  we obtain the integral representation of
$\mathfrak{sp}_{2\ell}$-Whittaker function in terms of factorized
parametrization.

\begin{te}\label{teintfsp2l}
 The eigenfunctions of the $\mathfrak{sp}_{2\ell}$-Toda
chain  admit the integral representation:
\bqa\label{psilc}\Psi_{\lambda_1,\ldots,\lambda_{\ell}}^{C_{\ell}}
(z_1,\ldots,z_{\ell})=
e^{\imath\lambda_1z_1+\ldots+\imath\lambda_{\ell}z_{\ell}}\int_C
\bigwedge_{i=1}^\ell\bigwedge_{k=1}^{n_i}\frac{dy_{i,k}}{y_{i,k}}
\Big(\,\prod\limits_{n=1}^{\ell} y_{1,n}
\prod_{k=2}^\ell\prod_{n=1}^{n_k/2}y^2_{k,2n-1}
\Big)^{\imath\lambda_1}\times \nonumber\\
\nonumber \times \prod_{i=2}^\ell\Big(\prod_{n=1}^\ell y_{1,n}
\prod_{k=2}^i\prod_{n=1}^{n_k}y_{k,n}
\prod_{k=i+1}^\ell\prod_{n=1}^{n_k/2}y^2_{k,2n-1}
\Big)^{\imath(\lambda_i-\lambda_{i-1})}\times \\ \times \exp\Big\{
-\Big(\sum_{n=1}^\ell\frac{1}{y_{1,n}}
\Big(1+\frac{y_{2,2(n-1)}}{y_{2,2(n-1)-1}}\Big)^2\prod_{i=n+1}^\ell
\Big(\frac{y_{2,2(i-1)}}{y_{2,2(i-1)-1}}\Big)^2+\\+
\nonumber\sum_{k=2}^\ell \sum_{n=1}^{n_k/2} \frac{1}{y_{k,2n}}
\Big(1+\frac{y_{k+1,2(n-1)}}{y_{k+1,2(n-1)-1}}\Big)
\prod_{i=n+1}^{n_k/2} \frac{y_{k+1,2(i-1)}}{y_{k+1,2(-1)-1}}
\frac{y_{k,2i-1}}{y_{k,2i}} \,\,+\\ \nonumber
+e^{2z_{1}}\sum_{n=1}^\ell y_{1,n}\,+\,\sum_{k=2}^\ell
e^{z_{k}-z_{k-1}} \sum\limits_{n=1}^{n_k}y_{k,n}\Big)\Big\},\eqa
where  $n_1=\ell$ and $n_k=2(\ell+1-k),$ for $k=2,\ldots,\ell$. The
domain of integration $C \subset N_+$ is a middle-dimensional
non-compact submanifold  such that the integrand decreases
exponentially at the boundaries and infinities. In particular one
can chose $C=\RR_{+}^{\ell^2}$.
\end{te}
The proof is given in Part II, Section \ref{melsp2l}.

\begin{ex} For $\ell=2$  the general formula
(\ref{psilc}) acquires the form
\bqa\label{C2RT}\Psi^{C_2}_{\lambda_1,\lambda_2}(z_{1},z_{2})=
e^{\imath\lambda_1 z_{1}+\imath\lambda_2 z_{2}}
\int\limits_{C}\bigwedge\limits_{i,k=1}^{2} \frac{d
y_{i,k}}{y_{i,k}}(y_{1,1}y_{2,1}^2y_{1,2})^{\imath \lambda_1}
(y_{2,1}y_{1,2}y_{2,2})^{\imath\lambda_2-\imath\lambda_1}\times \\
\nonumber\times
\exp\Big\{-\Big(\,\,\,\frac{1}{y_{1,1}}\Big(\frac{y_{2,2}}{y_{2,1}}\Big)^2+
\frac{1}{y_{1,2}}\Big(\frac{y_{2,2}}{y_{2,1}}+1\Big)^2 + \frac{1}{y_{2,2}}+\\
\nonumber +e^{2z_{1}}(y_{1,1}+y_{1,2})+
e^{z_{2}-z_{1}}(y_{2,1}+y_{2,2})\,\Big\},\eqa  with one can take
$C=\RR^{4}$.
\end{ex}

\subsubsection{ $\mathfrak{sp}_{2\ell}$-Whittaker function:
modified factorized   parametrization}

In this part we introduce a modified  factorized parametrization of
an open part of $N_+\subset Sp(2\ell)$.   We use this
parametrization to construct integral representations for
$\mathfrak{sp}_{2\ell}$-Whittaker functions. Similar to integral
representation of $\mathfrak{so}_{2\ell+1}$-Whittaker functions
considered above these integral representations have  a simple
recursive structure over the rank $\ell$ and can be describe in
purely combinatorial terms using suitable graphs. These
representations can be considered as a generalization of Givental
integral representations to the case of
$\mathfrak{g}=\mathfrak{sp}_{2\ell}$.

We follow the same approach that was used in the description of
modified factorized representation for $\mathfrak{so}_{2\ell+1}$.
There exists a realization of a tautological representation
$\pi_{2\ell}:\mathfrak{sp}_{2\ell}\to End(\mathbb{C}_{2\ell})$ such
that Weyl generators corresponding to  Borel (Cartan) subalgebra of
$\mathfrak{sp}_{2\ell}$ are realized by
 upper triangular (diagonal) matrices.
This defines an embedding $\mathfrak{sp}_{2\ell}\subset
\mathfrak{gl}_{2\ell}$ such that Borel (Cartan) subalgebra maps into
Borel (Cartan) subalgebra (see e.g. \cite{DS}). To define the
corresponding embedding of the groups  consider the following
involution on $GL(2\ell)$: \bqa g\longmapsto g^*:=\dot{w}_0\cdot
(g^{-1})^{t}\cdot\dot{w}_0^{-1},\eqa where $a\to a^{t}$ is induced by the
standard transposition matrices and $\dot{w}_0$ is a lift
of the longest element of the Weyl group of $\mathfrak{gl}_{2\ell}$.
In the matrix form it can be written as
$$ \pi_{2\ell}(\dot{w}_0)=S\cdot J,$$ where
 $S={\rm diag}(1,-1,\ldots,-1,1)$ and $J=\|J_{i,j}\|=\|\delta_{i+j,2\ell+2}\|$.
 The symplectic  group $G=Sp(2\ell)$ then can be defined as a
 following  subgroup of   $GL(2\ell)$  (see i.e. \cite{DS}):
$$ Sp(2\ell)=\{g\in GL(2\ell):g^{*}=g\}.$$

 Let $\epsilon_{i,j}$ stands for an  elementary $(2\ell\times2\ell)$
 matrix with a unit at  $(i,j)$ place and zeros otherwise.
Introduce the following $(2\ell\times2\ell)$ matrices: \bqa
U_n=\sum_{i=1}^\ell
\epsilon_{i,i}+e^{-z_{n-1,1}}\epsilon_{\ell+1,\ell+1}+
\sum_{i=1}^{n-1}e^{z_{n-1,i}}\epsilon_{\ell+1-i,\ell+1+i}\,\,\,\,,\eqa

\bqa \widetilde{U}_n=\sum_{i=1}^\ell
\epsilon_{i,i}+e^{-z_{n-1,1}}\epsilon_{\ell+1,\ell+1}+
\sum_{i=1}^{n-1}e^{z_{n-1,i}}\epsilon_{\ell+1-i,\ell+1+i}+\\
+\nonumber \sum_{i=2}^{n}e^{x_{n,i}}\epsilon_{\ell+i-1,\ell+i}+
\sum_{i=1}^{\ell-n}\epsilon_{\ell+n+i,\ell+n+i}\,\,\,\, ,\eqa \bqa
\widetilde{U}'_n=\sum_{i=1}^\ell
\epsilon_{i,i}+e^{-z_{n-1,1}}\epsilon_{\ell+1,\ell+1}+
\sum_{i=1}^{n-1}e^{z_{n-1,i}}\epsilon_{\ell+1-i,\ell+1+i}+\\
+\nonumber \sum_{i=1}^{n}e^{x_{n,i}}\epsilon_{\ell+i-1,\ell+i}+
\sum_{i=1}^{\ell-n}\epsilon_{\ell+n+i,\ell+n+i}\,\,\,\,  ,\eqa

\bqa V_n=\sum_{i=1}^{\ell-1}
\epsilon_{i,i}+e^{-z_{n,1}}\epsilon_{\ell,\ell}+
\sum_{i=1}^{n}e^{z_{n,i}}\epsilon_{\ell+i,\ell+i}\,\,\,\   .\eqa

\bqa \widetilde{V}_n=\sum_{i=1}^{\ell-1}
\epsilon_{i,i}+e^{-z_{n,1}}\epsilon_{\ell,\ell}+
\sum_{i=1}^{n}e^{z_{n,i}}\epsilon_{\ell+i,\ell+i}+\\ \nonumber
+\sum_{i=2}^{n}e^{x_{n,i}}\epsilon_{\ell+i-1,\ell+i}+
\sum_{i=1}^{\ell-n}\epsilon_{\ell+n+i,\ell+n+i}\,\,\,\, ,\eqa \bqa
\widetilde{V}'_n=\sum_{i=1}^{\ell-1}
\epsilon_{i,i}+e^{-z_{n,1}}\epsilon_{\ell,\ell}+
\sum_{i=1}^{n}e^{z_{n,i}}\epsilon_{\ell+i,\ell+i}+\\ \nonumber
+\sum_{i=1}^{n}e^{x_{n,i}}\epsilon_{\ell+i-1,\ell+i}+
\sum_{i=1}^{\ell-n}\epsilon_{\ell+n+i,\ell+n+i}\,\,\,\, ,\eqa
 We
can define a modified factorized parametrization as follows.

\begin{te}\label{Mparsp2l}
i)  The image of any generic unipotent  element $v^{C_{\ell}}\in
N_+$ in the tautological representation $\pi_{2\ell}:
\mathfrak{sp}_{2\ell}\to End(\mathbb{C}^{2\ell})$ can be presented
in the form \bqa\label{ggc1}
\pi_{2\ell}(v^{C_{\ell}})=\mathfrak{X}_1\mathfrak{X}_{2}\cdots\mathfrak{X}_{\ell},\eqa
where \bqa\label{ggc2}
\mathfrak{X}_1&=&1+e^{x_{11}+z_{11}}\epsilon_{\ell-1,\ell}\,\,\, , \nonumber\\
\mathfrak{X}_n&=&
[\widetilde{U}_nU_n^{-1}]^*\widetilde{U}'_n(U'_n)^{-1}
[(V'_n)^{-1}\widetilde{V}'_n]^*V_n^{-1}\widetilde{V}_n,\qquad
n=2,\ldots,\ell,\eqa and $z_{\ell,k}=0,\,\,k=1,\ldots\ell$ are
assumed.

ii) This defines a parametrization of an open part $N_+^{(0)}$  in
$N_+$.
\end{te}

{\it Proof.} Let  $v^{C_{\ell}}(y)$ be parametrization of an open
part of  $N_+$  according to (\ref{rec1C})-(\ref{rec2C}). Let
$\tilde{X}_i(y)=e^{y e_{i,i+1}}$ be a  one-parametric unipotent
subgroup in $GL(2\ell)$. Then we have
$\tilde{X}_i(y)^*=\tilde{X}_{2\ell+1-i}(y)$.  Embed elementary
unipotent subgroups $X_i(y)$ of $Sp(2\ell)$ into $GL(2\ell)$ as
follows:
$$X_i(y)=\tilde{X}_i(y)^*\,\tilde{X}_i(y).$$ This maps
 an arbitrary regular unipotent element $v^{C_\ell}$
into unipotent subgroup of $GL(2\ell).$ Let us now change the
variables in the following way: \bqa\label{nontwistc}
y_{11}=e^{x_{11}+z_{11}},\qquad y_{1,k}=\Big(e^{z_{k-1,1}+x_{k,1}}+
e^{z_{k,1}+x_{k,1}}\Big),\nonumber\\ \label{1yCn}
y_{k,2r-1}=e^{x_{k+r-1,k}-z_{k+r-2,k-1}},\hspace{2cm}
k=2,\ldots,\ell,\\ \nonumber
y_{k,2r}=e^{x_{k+r-1,k}-z_{k+r-1,k-1}},\hspace{2cm}
r=1,\ldots,\ell+1-k. \eqa Here $z_{\ell,k}=0$ for $k=1,\ldots\ell$ are
assumed.

 By elementary manipulations it is easy to check that after the change
of variables (\ref{1yCn}), the image $\pi_{2\ell}(v^{C_{\ell}})$ of
$v^{C_{\ell}}$ defined by (\ref{rec1C})-(\ref{rec2C})  transforms
into the (\ref{ggc1}) -(\ref{ggc3}). Taking into account that the
change of variables (\ref{1yCn}) is invertible
 we get a parametrization of $N_+^{(0)}\subset N_+$ $\Box$

The modified factorized parametrization of a unipotent group $N_+$
 defines a particular realization of a principal series
representation of $\mathcal{U}(\mathfrak{sp}_{2\ell})$ by
differential operators. It can be   obtained using the change of
variables (\ref{1yCn}) applied to the realization given in
Proposition \ref{realCn}.   We shall use the term Gauss-Givental
representation for this realization of representation of
 $\mathcal{U}(\mathfrak{sp}_{2\ell})$.
\begin{prop}
The following differential operators define a representation
$\pi_{\lambda}$ of $\mathcal{U}(\mathfrak{sp}_{2\ell})$ in $V_{\mu}$
in terms of the modified  factorized  parametrization: \bqa
E_1=\sum_{n=1}^\ell\left( \frac{\partial}{\partial x_{n,1}}-
\frac{\partial}{\partial x_{n+1,1}}
\right)\left(e^{z_{n,1}-x_{n,1}}+\sum_{i=n+1}^\ell\Big(
e^{z_{i-1,1}-x_{i,1}}+e^{z_{i,1}-x_{i,1}}\Big)\right)-\\
\nonumber -\sum_{n=1}^{\ell-1}\frac{\partial}{\partial z_{n,1}}
\sum_{i=n+1}^\ell\Big(
e^{z_{i-1,1}-x_{i,1}}+e^{z_{i,1}-x_{i,1}}\Big)\,\,\, ,\eqa

\bqa E_2=\left(\frac{\partial}{\partial z_{11}}-
\frac{\partial}{\partial x_{11}}\right)\Big(e^{z_{22}-x_{22}}+
\sum_{i=3}^\ell e^{z_{i-1,2}-x_{i,2}}+e^{z_{i,2}-x_{i,2}}\Big)+\\
\nonumber +\sum_{n=2}^\ell\left(\frac{\partial}{\partial x_{n,2}}-
\frac{\partial}{\partial x_{n,3}}\right)\Big(e^{z_{n,2}-x_{n,2}}+
\sum_{i=n+1}^\ell e^{z_{i-1,2}-x_{i,2}}+e^{z_{i,2}-x_{i,2}} \Big)+\\
\nonumber +\sum_{n=2}^\ell\left(\frac{\partial}{\partial z_{n,1}}-
\frac{\partial}{\partial z_{n,2}}\right) \sum_{i=n+1}^\ell\Big(
e^{z_{i-1,2}-x_{i,2}}+e^{z_{i,2}-x_{i,2}}\Big)-\\ \nonumber
-\sum_{n=2}^\ell\frac{\partial}{\partial x_{n,1}}\Big(e^{z_{n-1,1}+
x_{n,1}}\frac{e^{z_{n-1,2}-x_{n,2}}+e^{z_{n,2}-x_{n,2}}}
{e^{z_{n-1,1}+ x_{n,1}}+e^{z_{n,1}+ x_{n,1}}}+ \sum_{i=n+1}^\ell
e^{z_{i-1,2}-x_{i,2}}+e^{z_{i,2}-x_{i,2}}\Big)\,\,\,\, ,\eqa

\bqa  E_k=\left(\frac{\partial}{\partial z_{k-1,k-1}}+
\frac{\partial}{\partial x_{k,k}}\right)\Big(e^{z_{k,k}-x_{k,k}}+
\sum_{i=k+1}^\ell e^{z_{i-1,k}-x_{i,k}}+e^{z_{i,k}-x_{i,k}}\Big)+\\
\nonumber +\sum_{n=k}^\ell\left(\frac{\partial}{\partial z_{n,k-1}}-
\frac{\partial}{\partial z_{n,k}}\right)\sum_{i=n+1}^\ell\Big(
e^{z_{i-1,k}-x_{i,k}}+e^{z_{i,k}-x_{i,k}}\Big)+\eqa\bqa \nonumber
+\sum_{n=k+1}^\ell\left(\frac{\partial}{\partial x_{n,k}}-
\frac{\partial}{\partial x_{n,k+1}}\right)\Big(e^{z_{n,k}-x_{n,k}}+
\sum_{i=n+1}^\ell
e^{z_{i-1,k}-x_{i,k}}+e^{z_{i,k}-x_{i,k}}\Big),\qquad 2<k<\ell,\eqa

\bqa E_\ell=e^{-x_{\ell,\ell}}\left(
\frac{\partial}{\partial z_{k-1,k-1}}+ \frac{\partial}{\partial
x_{k,k}}\right)\,\,\,\, ,\eqa

\bqa H_k=\<\mu,\alpha_k^\vee\>+ \sum_{n=1}^\ell
a_{k,n}\sum_{i=n}^\ell\frac{\partial}{\partial x_{i,n}}\,\,\,\, ,
\eqa

\bqa F_1=-\Big(e^{z_{11}+x_{11}}+\sum_{i=2}^\ell
e^{z_{i-1,1}+x_{i,1}}+e^{z_{i,1}+x_{i,1}}\Big)\left(
\<\mu,\alpha_1^\vee\>+\frac{\partial}{\partial x_{11}}\right)-\\
\nonumber -\sum_{n=2}^\ell\frac{\partial}{\partial x_{n,1}}\Big(
e^{z_{n,1}+x_{n,1}}+\sum_{i=n+1}^\ell
e^{z_{i-1,1}+x_{i,1}}+e^{z_{i,1}+x_{i,1}}\Big)-\\ \nonumber
-\sum_{n=2}^\ell\frac{\partial}{\partial z_{n,1}}
\sum_{i=n+1}^\ell\Big(e^{z_{i-1,1}+x_{i,1}}+
e^{z_{i,1}+x_{i,1}}\Big)\, ,\eqa

\bqa F_2=-\left(\<\mu,\alpha_2^\vee\>-\frac{\partial}{\partial
x_{11}}- \frac{\partial}{\partial z_{11}}\right)\sum_{i=2}^\ell\Big(
e^{x_{i,2}-z_{i-1,1}}+e^{x_{i,2}-z_{i,1}}\Big)+\\
+\nonumber \sum_{n=2}^\ell\left(\frac{\partial}{\partial x_{n,1}}-
\frac{\partial}{\partial x_{n,2}}\right)\Big(e^{x_{n,2}-z_{n,1}}+
\sum_{i=n+1}^\ell e^{x_{i,2}-z_{i-1,1}}+e^{x_{i,2}-z_{i,1}}\Big)+\\
\nonumber +\sum_{n=2}^\ell\left(\frac{\partial}{\partial z_{n,1}}-
\frac{\partial}{\partial z_{n,2}}\right)\sum_{i=n+1}^\ell\Big(
e^{x_{i,2}-z_{i-1,1}}+e^{x_{i,2}-z_{i,1}}\Big)\, ,\eqa

\bqa F_k=-\left(\<\mu,\alpha_k^\vee\>-\frac{\partial}{\partial
x_{k-1,k-1}}- \frac{\partial}{\partial
z_{k-1,k-1}}\right)\sum_{i=k}^\ell\Big(
e^{x_{i,k}-z_{i-1,k-1}}+e^{x_{i,k}-z_{i,k-1}}\Big)+\\
\nonumber +\sum_{n=k}^\ell\left(\frac{\partial}{\partial x_{n,k-1}}-
\frac{\partial}{\partial x_{n,k}}\right)\Big(e^{x_{n,k}-z_{n,k-1}}+
\sum_{i=n+1}^\ell e^{x_{i,k}-z_{i-1,k-1}}+e^{x_{i,k}-z_{i,k-1}}\Big)+\eqa\bqa
\nonumber +\sum_{n=k}^\ell\left(\frac{\partial}{\partial z_{n,k-1}}-
\frac{\partial}{\partial z_{n,k}}\right)\sum_{i=n+1}^\ell\Big(
e^{x_{i,k}-z_{i-1,k-1}}+e^{x_{i,k}-z_{i,k-1}}\Big),\qquad 2<k<\ell,\eqa
\bqa F_\ell\,=\,\Big( e^{x_{\ell,\ell}-z_{\ell-1,\ell-1}}+
e^{x_{\ell,\ell}}\Big) \left(
-\<\mu,\alpha_\ell^\vee\>+\frac{\partial}{\partial x_{k-1,k-1}}+
\frac{\partial}{\partial z_{k-1,k-1}}\right)+\\ \nonumber
+e^{x_{\ell,\ell}}\left(\frac{\partial}{\partial
x_{\ell,\ell-1}}- \frac{\partial}{\partial x_{\ell,\ell}}\right)\,
,\eqa where $z_{\ell,k}=0, k=1,\ldots,\ell$ are assumed and the
derivatives over $x_{i,k},\,\,\, z_{i,k},\,\,\,i < k$ are omitted.
Here we denote $E_i=\pi_{\lambda}(e_i),\,\,\,F_i=\pi_{\lambda}(f_i),
\,\,\,H_i=\pi_{\lambda}(h_i)\,\,\,i=1,\ldots,\ell$.
\end{prop}

 We are going to write down matrix element (\ref{pairing}) explicitly using
Gauss-Givental representation defined above.
 Whittaker vectors $\psi_R$ and $\psi_L$ in this representation
 satisfiy the system of differential equations
\bqa\label{bwitt1} E_i\psi_R(x)= -\psi_R(x),\qquad F_i\psi_L(x)=
-\psi_L(x), \qquad  1\leq i\leq\ell.\eqa Its solution has the
following form.
Using the explicit change of the variables (\ref{1yCn})
we obtain the expressions for Whittaker vectors in modified factorized 
parametrization.

\begin{lem}\label{WhitModC}
 The following expressions for  left/right
Whittaker vectors hold: \bqa\psi_R=\exp\Big\{-\Big(
e^{x_{11}+z_{11}}+\sum_{n=2}^\ell
\Big(e^{z_{n-1,1}-x_{n,1}}+e^{z_{n,1}-x_{n,1}}\Big)\Big)-\\
\nonumber -\sum_{k=2}^\ell\,\,
\sum_{n=k}^\ell\Big(e^{x_{n,k}-z_{n-1,k-1}}+e^{x_{n,k}-z_{n,k-1}}
\Big)\Big\},\eqa
 \bqa\psi_L=e^{\mu_1z_{1,1}}\prod_{n=2}^\ell
\Big(e^{z_{n,1}}+e^{z_{n-1,1}}\Big)^{\mu_n}\times \\
\nonumber\times\prod_{n=1}^\ell
\exp\Big\{\,-\mu_n\Big(\sum_{i=1}^nz_{n,i}-
x_{n,1}-2\sum_{i=2}^{n}x_{n,i}+\sum_{i=1}^{n-1}z_{n-1,i}\Big)\Big\}
\,\,\times \\ \nonumber\times \exp\Big\{-
\sum_{k=1}^\ell\Big(e^{z_{k,k}-x_{k,k}}+\sum_{n=k+1}^\ell
e^{z_{n-1,k}-x_{n,k}}+e^{z_{n,k}-x_{n,k}}\Big) \Big\}, \eqa where
  $z_{\ell,k}=0$ and $\mu_k=\imath\lambda_k-\rho_k,$ $\rho_k=k$ for $k=1,\ldots,\ell$.
\end{lem}

Now we are ready to find an integral representation of the pairing
(\ref{pairing}) for $\mathfrak{g}=\mathfrak{sp}_{2\ell}$.
 To get an explicit expression for the integrand, one
uses the same type of  decomposition of the Cartan element as for
$\mathfrak{gl}_{\ell+1}$ and $\mathfrak{sp}_{2\ell}$ before:
$$e^{-H_z}=\pi_\lambda(\exp(-\sum_{i=1}^{\ell}
\langle\omega_i,z\rangle h_i))=e^{H_L}e^{H_R},$$ where \bqa\label{HRsp2l}
-H_z=H_L+H_R=-\langle\mu,z_{\ell}\rangle-2z_{\ell,1}\sum_{n=1}^\ell
\frac{\partial}{\partial x_{n,1}}+
\sum_{k=1}^{\ell-1}(z_{\ell,i}-z_{\ell,i+1})\sum_{n=k}^\ell
\frac{\partial}{\partial x_{n,k}},\eqa with

\bqa\label{HLC} H_L=\sum_{k=1}^\ell
z_{\ell,k}\Big(\sum_{n=k}^{\ell-1} \frac{\partial}{\partial
z_{n,k}}+ \sum_{n=k}^{\ell} \frac{\partial}{\partial
x_{n,k}}\Big),\eqa \bqa H_R=-\langle\mu,z_{\ell}\rangle-
\sum_{k=1}^{\ell-1} z_{\ell,k}\Big( \sum_{n=k}^{\ell}
\frac{\partial}{\partial x_{n,k}}- \sum_{n=k}^{\ell-1}
\frac{\partial}{\partial z_{n,k}}- \sum_{n=k+1}^{\ell}
\frac{\partial}{\partial x_{n,k+1}}\Big) .\eqa We imply that  $H_L$
acts on the left vector  and $H_R$ acts on the right vector in
(\ref{pairing}). Taking into account  Lemma  \ref{WhitModC} one obtains the
following theorem.

\begin{te}\label{teintmfsp2l}
The eigenfunctions of $\mathfrak{sp}_{2\ell}$-Toda chain
(\ref{pairing}) admit the integral representation
 \bqa\label{waveggc0}
\Psi^{C_\ell}_{\lambda_1,\ldots,\lambda_\ell}(z_{\ell,1},\ldots,z_{\ell,\ell})\,=\,
\int\limits_{C} \bigwedge_{k=1}^{\ell-1}\bigwedge_{i=1}^kdz_{k,i}
\bigwedge_{k=1}^{\ell}\bigwedge_{i=1}^kdx_{k,i}\, e^{{\mathcal
F}^{C_{\ell}}},\eqa where \bqa\label{waveggc1}\nonumber
\mathcal{F}^{C_{\ell}}= \imath\lambda_1z_{1,1}
-\sum\limits_{n=2}^{\ell}\imath\lambda_n \Big(\sum_{i=1}^nz_{n,i}
-x_{n,1}-2\sum_{i=2}^{n}x_{n,i}+\sum_{i=1}^{n-1}z_{n-1,i}
-\ln(e^{z_{n,1}}+e^{z_{n-1,1}})\Big)- \eqa
\bqa-\Big\{\sum_{k=1}^\ell\Big(e^{z_{k,k}-x_{k,k}}+\sum_{n=k+1}^\ell
e^{z_{n-1,k}-x_{n,k}}+e^{z_{n,k}-x_{n,k}}\Big)+ e^{x_{11}+z_{11}}+\\
\nonumber +\sum_{n=2}^\ell\,
\Big(e^{z_{n-1,1}-x_{n,1}}+e^{z_{n,1}-x_{n,1}}\Big)+ \sum_{k=2}^\ell
\sum_{n=k}^\ell\Big(e^{x_{n,k}-z_{n-1,k-1}}+e^{x_{n,k}-z_{n,k-1}}
\Big)\Big\}\,\,\,\,  , \eqa where we set $z_i:=z_{\ell,i},\,\,\,
1\le i\leq \ell $. Here ${C} \subset N_+$ is a middle-dimensional
non-compact submanifold such that the integrand decays exponentially
at the boundaries and at infinities. In particular the domain of
integration can be chosen to be ${C}=\RR^{m}$,  where $m=l(w_0)$.
\end{te}

\begin{ex}  For $\ell=2$ the general expression (\ref{waveggc})
acquires the form \bqa\label{waveggc2}
\Psi^{C_2}_{\lambda_1,\lambda_2}(z_{1},z_{\ell})\,=\,
\int\limits_{C}dx_{11}\wedge x_{11}\wedge dx_{21}\wedge
dx_{22}\wedge dz_{11}\times \nonumber \\\times
\exp\Big\{\imath\lambda_1x_{11}-\imath\lambda_2
\Big(z_{21}+z_{22}-x_{21}-2x_{22}+z_{11}-\log(e^{z_{21}}+e^{z_{11}})\Big)-
\\
-\Big(e^{z_{11}-x_{11}}+e^{z_{22}-x_{22}}+e^{x_{11}+z_{11}}+e^{z_{11}-x_{21}}
+e^{z_{21}-x_{21}}+e^{x_{22}-z_{11}}+e^{x_{22}-z_{21}}\Big)\Big\},
\, \nonumber \eqa where $z_1=z_{2,1}, z_2=z_{2,2}.$ In particular
one can chose $C=\RR^4$.
\end{ex}

There is a simple combinatorial description of the potential
$\mathcal{F}^{C_{\ell}}$ for zero spectrum $\{\la_i=0\}$.  Namely,
it can be presented as the sum over all arrows in the following
diagram.

$$
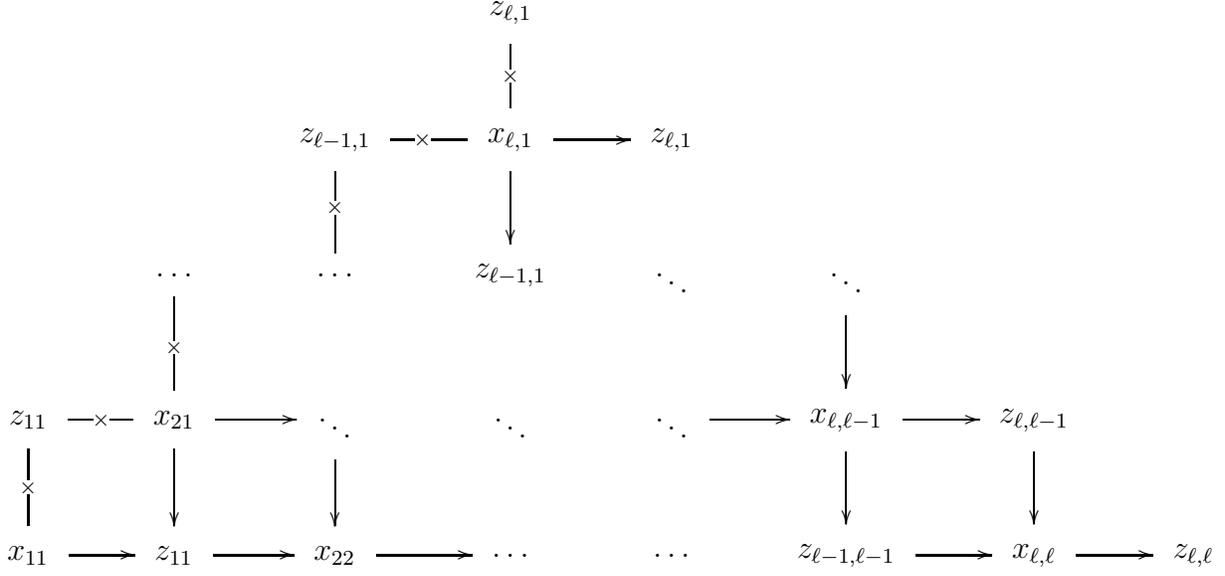
\xymatrix{
 &&& z_{\ell,1}\ar@{-}[d]|{\times} &&&&\\
 && z_{\ell-1,1}\ar@{-}[d]|{\times}\ar@{-}[r]|{\times} &
 x_{\ell,1}\ar[d]\ar[r] & z_{\ell,1} &&&\\
 & \ldots\ar@{-}[d]|{\times} & \ldots & z_{\ell-1,1} & \ddots &
 \ddots\ar[d] &&\\
 z_{11}\ar@{-}[d]|{\times}\ar@{-}[r]|{\times} &
 x_{21}\ar[d]\ar[r] & \ddots\ar[d] & \ddots & \ddots\ar[r] &
 x_{\ell,\ell-1}\ar[d]\ar[r] & z_{\ell,\ell-1}\ar[d] &\\
 x_{11}\ar[r] & z_{11}\ar[r] & x_{22}\ar[r] & \ldots &
 \ldots & z_{\ell-1,\ell-1}\ar[r] & x_{\ell,\ell}\ar[r] &
 z_{\ell,\ell}
}$$ \hspace{1cm} We use the same rule to  assign variables to the
arrows of the diagram as for $A_{\ell}$. In addition we assign to
the symbol $\xymatrix{z\ar@{-}[r]|{\times}& x}$ the exponent $e^{-z-x}$.

Note that the diagram for $C_{\ell}$  can be obtained  by a
factorization of the diagram for $A_{2\ell-1}$. Consider the
following involution \be \label{inv} \iota\,:\quad X\longmapsto
\dot{w}_0^{-1}X^t \dot{w}_0, \ee where  $\dot{w}_0$ is a lift the
longest element of $A_{2\ell-1}$ Weyl group and $X^t$ denotes the
standard transposition. Corresponding action on the modified
factorization parameters is given by \be \label{w0inv}
\dot{w}_0\,:\qquad x_{k,i}\longleftrightarrow -x_{k,k+1-i}.\ee This
defines a factorization  of $A_{2\ell-1}$-diagram  that produces the
diagram for $C_{\ell}$.

One can easily write down $C_{\ell}$-analog of $A_{\ell}$-monomial
relations (\ref{defrelAn}). Let us introduce the  variables \bqa
a_{k,1}=e^{x_{k,1}+z_{k-1,1}},\qquad a_{k,i}=e^{x_{k,i}-z_{k-1,i-1}},\nonumber\\
b_{k,1}=e^{x_{k,1}+z_{k,1}},\qquad
b_{k,i}=e^{z_{k,i}-x_{k,i-1}},\\
c_{k,i}=e^{z_{k,i}-x_{k,i}},\qquad
d_{k,i}=e^{z_{k,i}-x_{k+1,i}}.\nonumber \eqa Then the following
relations hold \bqa\label{torc} c_{k,i}\cdot
b_{k,i}\,=\,d_{k,i}\cdot
a_{k+1,i},\qquad a_{k,i}\cdot d_{k-1,i-1}\,=\,b_{k,i}\cdot c_{k,i-1},\nonumber\\
b_{\ell,1}c_{\ell,1}=e^{2z_{\ell,1}},\qquad\qquad\quad
c_{\ell,i}\cdot b_{\ell,i}=e^{z_{\ell,i}-z_{\ell,i-1}}.\eqa The
above relations can be considered as relations between elementary
paths on the Givental diagram. Using relations for more general
paths that follows from (\ref{torc})  one can define a toric
degeneration of the $C_{\ell}$-flag manifolds thus generalizing
results of \cite{BCFKS}.

\subsubsection{Recursion for $\mathfrak{sp}_{2\ell}$-Whittaker
  functions and  $\mathcal{Q}$-operator \\  for $A^{(2)}_{2\ell-1}$-Toda
  chain }

The integral representation (\ref{waveggc0}), (\ref{waveggc1}) of
$\mathfrak{sp}_{2\ell}$-Whittaker functions possesses a recursive
structure over the rank $\ell$. For any $n=2,\ldots,\ell$ let us
introduce integral operators $Q^{C_n}_{C_{n-1}}$
 with the integral  kernels
\bqa\label{CQC}
Q^{C_n}_{C_{n-1}}(\underline{z}_{n};\,\underline{z}_{n-1};\lambda_n)=\int
\bigwedge_{i=1}^n dx_{n,i}\,\,\,
\Big(e^{z_{n,1}}+e^{z_{n-1,1}}\Big)^{\imath\lambda_n}\times
\nonumber\\ \times
\exp\Big\{\,-\imath\lambda_n\Big(\sum_{i=1}^nz_{n,i}-
x_{n,1}-2\sum_{i=2}^{n}x_{n,i}+\sum_{i=1}^{n-1}z_{n-1,i}\Big)\Big\}
\,\,\times \\ \nonumber \times
Q^{C_n}_{D_{n}}(\underline{z}_{n};\,\underline{x}_{n})\,\,
Q^{D_{n}}_{C_{n-1}}(\underline{x}_{n};\,\underline{z}_{n-1}),\quad\eqa
where  \be Q^{D_n}_{\,\,\,\,C_{n-1}}(\underline{x}_n
;\,\underline{z}_{n-1})=\\
=\exp\Big\{-\Big(\,e^{x_{n,1}+z_{n-1,1}}+\sum_{i=1}^{n-1}\Big(
e^{z_{n-1,i}-x_{n,i}}+ e^{x_{n,i+1}-z_{n-1,i}}\Big)\,\Big)\Big\},\ee
\be Q^{C_n}_{\,\,\,\,D_{n}}(\underline{z}_n ;\,\underline{x}_n)=\\
=\exp\Big\{-\Big(\,e^{x_{n,1}+z_{n,1}}+\sum_{i=1}^{n-1}\Big(e^{z_{n,i}-x_{n,i}}+
e^{x_{n,i+1}-z_{n,i}}\Big)+e^{z_{n,n}-x_{n,n}}\Big)\,\Big\}.\ee For
$n=1$ we define $$Q^{C_1}_{C_0}=\int
dx_{11}e^{\imath\lambda_1x_{11}}\exp\Big\{-\Big(e^{x_{11}+z_{11}}+e^{z_{11}-x_{11}}\Big)\Big\}$$
Using  integral operators $Q^{C_n}_{C_{n-1}},$ the integral
representation for $\mathfrak{sp}_{2\ell}$-Whittaker function can be
written in the recursive form.

\begin{te} The integral representations of  $\mathfrak{sp}_{2\ell}$-Toda chain
eigenfunctions (\ref{waveggc0}) can be written as
\bqa\label{waveiterc}
\Psi^{C_\ell}_{\lambda_1,\ldots,\lambda_\ell}(z_1,\ldots,z_\ell)\,=\,
\int\limits_\mathcal{C}
\bigwedge_{k=1}^{\ell-1}\bigwedge_{i=1}^kdz_{k,i}\,\,
\prod_{n=1}^{\ell}
Q^{C_{n}}_{\,\,C_{n-1}}(\underline{z}_{n};\,\underline{z}_{n};\lambda_{n}),\eqa
or equivalently \bqa\label{iterwa2}
\Psi^{C_\ell}_{\lambda_1,\ldots,\lambda_\ell}(z_{\ell,1},\ldots,z_{\ell,\ell})\,=\,
\int\limits_{\mathcal{C}_{\ell}}
\bigwedge_{i=1}^{\ell-1}dz_{\ell-1,i} \times \,
\\
\nonumber \times
Q^{C_{\ell}}_{\,\,C_{\ell-1}}(\underline{z}_{\ell};\,
\underline{z}_{\ell-1};\lambda_{\ell})
\Psi^{C_{\ell-1}}_{\lambda_1,\ldots,\lambda_{\ell-1}}(z_{\ell-1,1},
\ldots,z_{\ell-1,\ell-1}).\eqa Here  $z_n:=z_{\ell,n},\,\,\,  1\le
n\leq \ell $ and $\mathcal{C} \subset N_+$ is a middle-dimensional
non-compact submanifold such that the integrand decreases
exponentially at the boundaries and at infinities. In particular as
a domain of integration one can chose  $\mathcal{C}=$$\RR^{\ell^2}$.
\end{te}

This recursive form of the integral representation is similar to the
case of $\mathfrak{so}_{2\ell+1}$. Its recursive kernel
$Q^{C_n}_{C_{n-1}}$ is given by a  nontrivial integral in contrast
with $\mathfrak{gl}_{\ell+1}$-case (\ref{QRecA}).  Similar to
$\mathfrak{so}_{2\ell+1}$-Whittaker function  new structure appears
if we consider the Whittaker function for zero spectrum
$\{\la_i=0\}$. As it is clear from (\ref{CQC}) the kernels
$Q^{C_n}_{C_{n-1}}$ at $\la_n=0$ are given by  convolutions of the
kernels $Q^{C_n}_{D_n}(\underline{z}_n, \underline{x}_n)$ and
$Q^{D_{n}}_{C_{n-1}}(\underline{x}_n,\underline{z}_{n-1})$.
 The corresponding integral  operators $Q^{C_n}_{D_n}$,
$Q^{D_{n}}_{C_{n-1}}$ can be regarded as elementary intertwiners
relating  Hamiltonians  of  Toda chains for $C_n$, $D_n$ and $D_n$,
$C_{n-1}$  root systems correspondingly. For example it is easy to
check directly  intertwinig realtions with quadratic Hamiltonians.
Indeed, $D_{\ell}$-Toda chain (for more detailed discussion see
Subsection 2.5.3) has the following  quadratic Hamiltonians
 \bqa \CH_2^{D_{\ell}}(\underline{x}^{(\ell)})=-
\frac{1}{2}\sum_{i=1}^{\ell}\frac{\partial^2}{\partial x_i^2}+
e^{x_1+x_2}+\sum_{i=1}^{\ell-1}e^{x_{i+1}-x_i}. \eqa

\begin{prop}
The integral operators  $Q^{C_n}_{C_{n-1}}$, $Q^{C_n}_{D_{n}}$ and
$Q^{D_n}_{C_{n-1}}$  satisfy the following relations.
\begin{enumerate}
\item Operators $Q^{C_n}_{D_{n}}$ and $Q^{D_{n}}_{C_{n-1}}$ intertwine  quadratic
Hamiltonians of $C$- and $D$-Toda chains:
\bqa\CH_2^{D_{n}}(\underline{x}_n)Q^{D_{n}}_{C_{n-1}}(\underline{x}_n,\,
\underline{z}_{n-1})= Q^{D_n}_{C_{n-1}}
(\underline{x}_n,\,\underline{z}_{n-1})\CH_2^{C_{n-1}}(\underline{x}_{n-1}),
\eqa
\bqa\CH_2^{C_{n}}(\underline{z}_n)Q^{C_{n}}_{D_{n}}(\underline{z}_n,\,
\underline{x}_{n})= Q^{C_n}_{D_{n}}
(\underline{z}_n,\,\underline{x}_{n})\CH_2^{D_{n}}(\underline{x}_{n}).
\eqa

\item The operator $ Q^{C_{n}}_{C_{n-1}}$ at $\la_n=0$
 intertwines the Hamiltonians ${\CH}_2^{C_n}$
and ${\CH}_2^{C_{n-1}}:$
\bqa\CH_2^{C_{n}}(\underline{z}_n)Q^{C_{n}}_{C_{n-1}}(\underline{z}_n,\,
\underline{z}_{n-1})= Q^{C_n}_{C_{n-1}}
(\underline{z}_n,\,\underline{z}_{n-1})\CH_2^{C_{n-1}}(\underline{z}_{n-1}).
\eqa
\end{enumerate}\end{prop}

The integral kernel of $Q^{C_n}_{C_{n-1}}$ can be succinctly encoded
into the following  diagram

$$
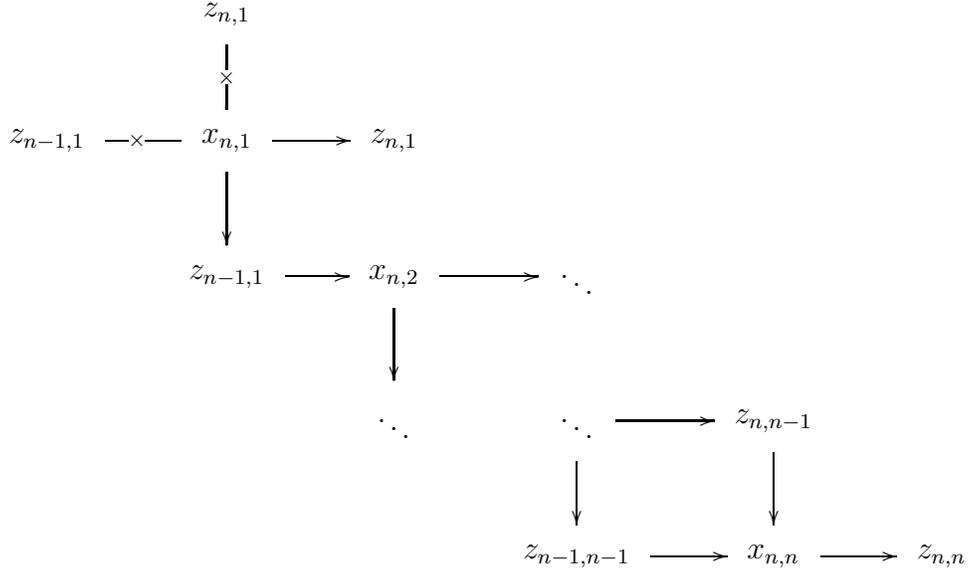
\xymatrix{
 & z_{n,1}\ar@{-}[d]|{\times} &&&&\\
 z_{n-1,1}\ar@{-}[r]|{\times} & x_{n,1}\ar[d]\ar[r] & z_{n,1} &&&\\
 & z_{n-1,1}\ar[r] & x_{n,2}\ar[d]\ar[r] & \ddots &&\\
 && \ddots & \ddots\ar[d]\ar[r] & z_{n,n-1}\ar[d] &\\
 &&& z_{n-1,n-1}\ar[r] & x_{n,n}\ar[r] & z_{n,n}
}
$$

Here the upper (lower) boundary of the oriented diagram corresponds
to the kernels of elementary intertwiner $Q^{C_n}_{D_n}$
($Q^{D_n}_{C_{n-1}}$) and the convolution of the two kernel
 is given  by the integration over variables
$x_{n,1},\ldots, x_{n,n}$ on the diagonal of  the diagram.

Similarly to the cases of $\mathfrak{gl}_{\ell+1}$ and
$\mathfrak{sp}_{2\ell}$ recursion operators $Q^{C_n}_{C_{n-1}}$ can
be considered as degenerations of a Baxter $\mathcal{Q}$-operators
for twisted affine $A^{(2)}_{2\ell-1}$-Toda chain introduced below.
Let us stress that up to now $\mathcal{Q}$-operators for
$A^{(2)}_{2\ell-1}$ were not known. We will not present here a
complete set of the characteristic properties of the introduced
$\mathcal{Q}$-operators and only consider  commutation relations
with quadratic affine Toda chain Hamiltonians. The detailed account
will be given elsewhere.

We start with a description of $A_{2\ell-1}^{(2)}$-Toda chains. The
set of simple roots of the affine root system $A^{(2)}_{2\ell-1}$
can be represented in terms of the orthogonal bases
$\{\epsilon_i\}$, $i=1,\ldots, \ell$ in $\mathbb{R}^{\ell}$ as
follows:
  \be \alpha_1=2\epsilon_1,\qquad \alpha_{i+1}=\epsilon_{i+1}-\epsilon_i,\qquad
1\leq i\leq \ell-1,\qquad
\alpha_{\ell+1}=-\epsilon_{\ell}-\epsilon_{\ell-1}, \ee and
corresponding Dynkin diagram is given  by
$$
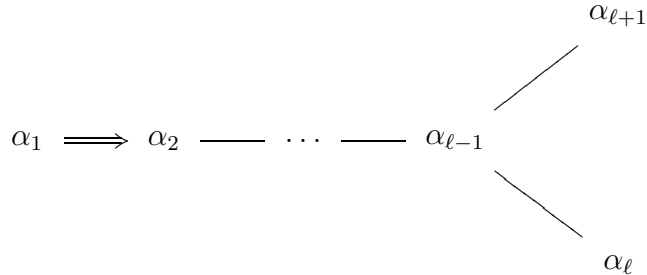
\xymatrix{
 &&&& \alpha_{\ell+1}\\
 \alpha_1\ar@{=>}[r] & \alpha_2\ar@{-}[r] & \ldots &
 \alpha_{\ell-1}\ar@{-}[l]\ar@{-}[ur]\ar@{-}[dr] &\\
 &&&& \alpha_{\ell}
}$$

These root data allows to define affine $A^{(2)}_{2\ell-1}$-Toda
chain with the quadratic Hamiltonian given by
\bqa\CH_2^{A^{(2)}_{2\ell-1}}(\underline{z}^{(\ell)})&=&
-\frac{1}{2}\sum_{i=1}^{\ell} \frac{\partial^2}{\partial
z_i^2}+2e^{2z_1}+ \sum_{i=1}^{\ell-1}e^{z_{i+1}-z_i}+
ge^{-z_{\ell-1}-z_{\ell}}, \eqa where $g$ is an arbitrary parameter.

Define the Baxter $\mathcal{Q}$-operator of
 $A_{2\ell-1}^{(2)}$-Toda chain as an  integral operator with the
 following integral kernel:
\bqa \label{BaxterC}
\mathcal{Q}^{^{A_{2\ell-1}^{(2)}}}(\underline{z}^{(\ell)}
\,,\underline{y}^{(\ell)},\lambda)=\int\bigwedge_{i=1}^{\ell+1}dx_{i}\,\,
\Big(e^{z_1}+e^{y_1}\Big)^{\imath\lambda}
\Big(e^{-z_\ell}+e^{-y_\ell}\Big)^{-2\imath\lambda}\times \\\times
\nonumber \exp\Big\{-\imath\lambda\Big(\sum_{i=1}^\ell z_i- x_1-
2\sum_{i=2}^\ell x_i+\sum_{i=1}^\ell y_i\Big) \Big\}\,\times \\\times
\nonumber Q^{A^{(2)}_{2\ell-1}}_{\,\,\,\,A^{(2)}_{2\ell-1}}
(z_1,\ldots,z_\ell\,;x_1,\ldots,x_{\ell+1})\,\,\,
Q^{\,\,\,\,A^{(2)}_{2\ell-1}}_{A^{(2)}_{2\ell-1}}
(x_1,\ldots,x_{\ell+1}\,;y_1,\ldots,y_\ell), \eqa where
\bqa\label{inttwCn}
Q^{A^{(2)}_{2\ell-1}}_{\,\,\,\,A^{(2)}_{2\ell-1}}
(z_1,\ldots,z_\ell\,;\,x_1,\ldots,x_{\ell+1})=\\ =\nonumber
\exp\Big\{-\Big( e^{z_1+x_1}+\sum_{i=1}^{\ell-1}\Big(
e^{z_i-x_i}+e^{x_{i+1}-z_i} \Big)
+e^{z_\ell-x_\ell}+ge^{-z_\ell-x_{\ell}}\Big)\Big\},\eqa and
\bqa\label{inttwCn}
Q^{A^{(2)}_{2\ell-1}}_{\,\,\,\,A^{(2)}_{2\ell-1}}
(x_1,\ldots,x_{\ell+1}\,;\,y_1,\ldots,y_{\ell})=\\= \nonumber
\exp\Big\{-\Big( e^{y_1+x_1}+\sum_{i=1}^{\ell-1}\Big(
e^{y_i-x_i}+e^{x_{i+1}-y_i}\Big)+e^{y_\ell-x_\ell}+
ge^{-y_\ell-x_{\ell}}\Big)\Big\}.\eqa
Here we use the following notations
$\underline{z}^{(\ell)}=(z_1,\ldots,z_{\ell})$,
$\underline{y}^{(\ell)}=(y_1,\ldots, y_{\ell})$.

The following statement can be verified straightforwardly.
\begin{prop}
The $\mathcal{Q}$-operator (\ref{BaxterC}) commutes with the
quadratic Hamiltonian of  $A_{2\ell-1}^{(2)}$-Toda chain:
\bqa\CH^{A_{2\ell-1}^{(2)}}(\underline{z}^{(\ell)})
Q^{^{A_{2\ell-1}^{(2)}}}(\underline{z}^{(\ell)},\,
\underline{y}^{(\ell)})= Q^{^{A_{2\ell-1}^{(2)}}}
(\underline{z}^{(\ell)},\,\underline{y}^{(\ell)})
\CH^{^{A_{2\ell-1}^{(2)}}}(\underline{y}^{(\ell)}). \eqa
\end{prop}
Now we will demonstrate that the recursion  operator
$Q^{C_{\ell}}_{C_{\ell-1}}$ can be obtained by a degeneration of the
Baxter $\mathcal{Q}$-operator for $A_{2\ell-1}^{(2)}$. 
Consider a slightly modified recursion operator 
$Q^{C_{\ell}}_{C_{\ell-1}\oplus C_{1}}$ with the kernel given by 
\bqa\nonumber
Q^{C_{\ell}}_{C_{\ell-1}\oplus C_{1}}
(\underline{z}^{(\ell)}\,,\,\underline{y}^{(\ell)}\,,
\lambda):=e^{\imath\lambda y_{\ell}}\,\,
Q^{C_{\ell}}_{C_{\ell-1}}
(\underline{z}^{(\ell)}\,,\,\underline{y}^{(\ell-1)}\,,\,\lambda),\eqa
where $\underline{y}^{(\ell-1)}=(y_1,\ldots, y_{\ell-1})$.
Thus defined operator intertwines Hamiltonians of $\mathfrak{sp}_{2\ell}$- and
  $\mathfrak{sp}_{2\ell-2}\oplus \mathfrak{sp}_2$-Toda chains. Thus
  for quadratic Hamiltonians we have
\bqa\label{fusionC}\CH_2^{C_{\ell}}(\underline{z}^{(\ell)})
Q^{C_{\ell}}_{C_{\ell-1}\oplus C_{1}}
(\underline{z}^{(\ell)}\,,\,\underline{y}^{(\ell)},\lambda) =
Q^{C_{\ell}}_{C_{\ell-1}\oplus C_{1}}
(\underline{z}^{(\ell)}\,,\underline{y}^{(\ell)}\,,\lambda)
\Big(\CH_2^{C_{\ell}}(\underline{y}^{(\ell-1)})+
\CH_2^{C_1}(y_{\ell})\Big), \nonumber \eqa where
$\CH_2^{C_1}(y_{\ell})=
-\frac{1}{2}\Big(\partial^2/\partial y_{\ell}^2\Big)$. Obviously the
projection of above equation  on the subspace of functions
$F(\underline{y}^{(\ell)})=\exp(\imath\lambda y_{\ell})
f(\underline{y}^{(\ell-1)})$ recovers the  genuine recursion operator
satisfying:
\bqa\CH_2^{C_{\ell}}(\underline{z}^{(\ell)})
Q^{C_{\ell}}_{C_{\ell-1}}
(\underline{z}^{(\ell)}\,,\,\underline{y}^{(\ell-1)}\,,\,\lambda) =
Q^{C_{\ell}}_{C_{\ell-1} }
(\underline{z}^{(\ell)}\,,\,\underline{y}^{(\ell-1)}\,,\,\lambda)
\Big(\CH_2^{C_\ell}
(\underline{y}^{(\ell-1)})+\frac{1}{2}\lambda^{2})\Big). \eqa Consider
a one-parameter family of the kernels: \bqa
\mathcal{Q}^{^{A_{2\ell-1}^{(2)}}} (\underline{z}^{(\ell)},\,
\underline{y}^{(\ell)},\lambda,\,\varepsilon)= \varepsilon^{\imath\lambda}
e^{\imath\lambda y_{\ell}} \int\bigwedge_{i=1}^{\ell+1}dx_{i}\,\,
\Big(e^{z_1}+e^{y_1}\Big)^{\imath\lambda} \Big(\varepsilon
e^{y_{\ell}-z_\ell}+1\Big)^{-2\imath\lambda}\times
\\\times
\nonumber \exp\Big\{-\imath\lambda\Big(\sum_{i=1}^\ell z_i- x_1-
2\sum_{i=2}^\ell x_i+\sum_{i=1}^{\ell-1} y_i\Big) \Big\}\,\times \\
\nonumber\times  Q^{A^{(2)}_{2\ell-1}}_{\,\,\,\,A^{(2)}_{2\ell-1}}
(z_1,\ldots,z_\ell\,;x_1,\ldots,x_{\ell+1})\,\,\,
Q^{\,\,\,\,A^{(2)}_{2\ell-1}}_{A^{(2)}_{2\ell-1}}
(x_1,\ldots,x_{\ell+1}\,;y_1,\ldots,y_\ell;\,\varepsilon),\eqa where
\bqa\label{inttwCne}
Q^{A^{(2)}_{2\ell-1}}_{\,\,\,\,A^{(2)}_{2\ell-1}}
(x_1,\ldots,x_{\ell+1}\,;y_1,\ldots,y_\ell;\varepsilon)=
\exp\Big\{-\Big( e^{y_1+x_1}+\sum_{i=1}^{\ell-1}\Big(
e^{y_i-x_i}+e^{x_{i+1}-y_i}\Big)+\\
\nonumber +\varepsilon e^{y_\ell-x_\ell}+
\varepsilon^{-1}ge^{-y_\ell-x_{\ell}}\Big)\Big\},\eqa is obtained by
shifting the variable $y_{\ell}=y_{\ell}+\ln\varepsilon$ in
(\ref{BaxterC}). Then the following relation between $\mathcal{Q}$-operator for
$A_{2\ell-1}^{(2)}$-Toda chain and recursive operator for
the $\mathfrak{sp}_{2\ell}$-Whittaker function holds \bqa
Q^{C_{\ell}}_{C_{\ell-1}\oplus C_{1}}
(\underline{z}^{(\ell)}\,,\,\underline{y}^{(\ell)},\la)=
\lim_{\varepsilon\rightarrow 0, \varepsilon^{-1}g\rightarrow 0}
\varepsilon^{-\imath\lambda}
\mathcal{Q}^{^{A_{2\ell-1}^{(2)}}}(
\underline{z}^{(\ell)}\,,\,\underline{y}^{(\ell)},\lambda;\,\varepsilon). \eqa


\newpage
\subsection{ Integral representations of $\mathfrak{so}_{2\ell}$-Toda chain
eigenfunctions }

 In this subsection we provide an analog of the Givental
integral representation of Whittaker functions for
$\mathfrak{so}_{2\ell}$ Lie algebras. As in the previously
considered cases, we start with a derivation of an integral
representation of $\mathfrak{so}_{2\ell}$-Whittaker functions using
the factorized representation. Then we
consider a modification of  the factorized representation leading to
a Givental type integral representation of
$\mathfrak{so}_{2\ell}$-Whittaker functions.

Consider $D_{\ell}$ root system corresponding to
 Lie algebra $\mathfrak{so}_{2\ell}$.
Let $(\epsilon_1,\ldots,\epsilon_{\ell})$ be an orthogonal basis in
$\RR^{\ell}.$ We realize $D_{\ell}$  simple root
 and fundamental weights as the
following vectors in $\mathbb{R}^{\ell}$:
 \bqa\label{RootDn}
\begin{array}{l}
\alpha_1=\epsilon_2-\epsilon_1,\\
\alpha_2=\epsilon_2+\epsilon_1,\\
\alpha_3=\epsilon_3-\epsilon_2,\\
\ldots\\
\alpha_n=\epsilon_\ell-\epsilon_{\ell-1},
\end{array}\hspace{2cm}
\begin{array}{l}
\omega_1=(-\epsilon_1+\epsilon_2+\ldots+\epsilon_\ell)/2,\\
\omega_2=(\epsilon_1+\epsilon_2+\ldots+\epsilon_\ell)/2,\\
\omega_3=\epsilon_3+\ldots+\epsilon_\ell,\\
\ldots\\
\omega_\ell=\epsilon_\ell.
\end{array}
\eqa Coroots $\alpha_i^{\vee}$ can be identified  with the
corresponding roots $\alpha_i$ using the scalar product in
$\mathbb{R}^{\ell}$. One associates with these root data
$\mathfrak{so}_{2\ell}$-Toda chain with a quadratic Hamiltonian
given by \bqa\label{CH2D}
\CH_2^{D_{\ell}}&=&-\frac{1}{2}\sum\limits_{i=1}^{\ell}
\frac{\partial^2}{{\partial x_i}^2}+
e^{x_1+x_2}+\sum\limits_{i=1}^{\ell-1} e^{x_{i+1}-x_{i}}. \eqa One
can complete (\ref{CH2D}) to a full set of $\ell$ mutually commuting
functionally independent Hamiltonians $H^{D_{\ell}}_k$ of the
$\mathfrak{so}_{2\ell}$-Toda chain. We are looking for integral
representations of common eigenfunctions of the full set of the
Hamiltonians. Corresponding eigenfunction problem for the quadratic
Hamiltonian can be written in the following form \be
\CH_2^{D_{\ell}}\,\,\,\Psi^{{D}_{\ell}}_{\la_1,\cdots
,\la_{\ell}} (x_1,\ldots,x_{\ell})=
\frac{1}{2}\sum\limits_{i=1}^{\ell}\lambda_{i}^2\,\,\,
\Psi^{{D}_{\ell}}_{\la_1,\cdots ,\la_{\ell}} (x_1,\ldots,x_{\ell}).
\ee

\subsubsection{ $\mathfrak{so}_{2\ell}$-Whittaker function:
  factorized   parametrization}

The reduced word for the  maximal length element $w_0$ in the Weyl
group of $\mathfrak{so}_{2\ell}$ can be represented in the following
recursive way:
$$I=(i_1,i_2,\ldots,i_m):=
(12,3123,\ldots,(\ell\ldots3123\ldots\ell)),$$
where index $i_k$ corresponds  to an elementary reflection with
respect to the  root $\alpha_k$. Let $N_+\subset G$ be a maximal
unipotent subgroup of $G=SO(2\ell)$. One associates with the reduced
word $I$ the following recursive parametrization of a generic
element  $v^{D_{\ell}}\in N_+$:
 \bqa\label{rec1D}
v^{D_{\ell}}=v^{D_{\ell-1}}\cdot\mathfrak{X}^{D_{\ell}}_{D_{\ell-1}}
,\eqa where
\bqa\label{rec2D}\mathfrak{X}^{D_\ell}_{D_{\ell-1}}=X_{\ell}(y_{\ell,1})\cdots
X_k(y_{k,2(\ell+1-k)-1})\cdots
X_3(y_{3,2\ell-5})X_1(y_{1,\ell-1})\cdot\\
\nonumber X_2(y_{2,\ell-1})X_3(y_{3,2\ell-4})\cdots
X_k(y_{k,2(\ell+1-k)})\cdots X_{\ell}(y_{\ell,2}). \eqa Here
$X_i(y)=e^{ye_i}$ and  $e_i\equiv e_{\alpha_i}$ are simple root
generators. The subset  $N_+^{(0)}$ of the elements allowing
representation (\ref{rec1D}), (\ref{rec2D}) is an open part of
$N_+$. The action of the Lie algebra $\mathfrak{so}_{2\ell}$ on
$N_+$ (\ref{infact}) defines an action on the space of functions on
$N_+^{(0)}$. The explicit description of the action on the space
$V_{\mu}$ of (twisted) functions on $N_+^{(0)}$ is given below.

\begin{prop}\label{GGrepD}
The following differential operators define a realization of
 a principal series representation
$\pi_{\lambda}$ of $\mathcal{U}(\mathfrak{so}_{2\ell})$ in $V_{\mu}$
in terms of factorized parametrization of $N^{(0)}_+$:

 \bqa E_1=\frac{\partial}{\partial y_{1,\ell-1}}+
\sum_{n=1}^{[\ell/2]}\Big(\frac{\partial}{\partial y_{2,\ell-n-1}}-
\frac{\partial}{\partial y_{2,\ell-n}}\Big)\prod_{k=1}^{2n-1}
\Big(\frac{y_{1,\ell-k}}{y_{2,\ell-k}}\Big)^{(-1)^k}
\frac{y_{3,2(\ell-1-k)}}{y_{3,2(\ell-1-k)-1}}+\\ \nonumber
+\sum_{n=2}^{[\ell/2]}\Big(\frac{\partial}{\partial y_{1,\ell-n-1}}-
\frac{\partial}{\partial y_{1,\ell-n}}\Big)\prod_{k=1}^{2(n-1)}
\Big(\frac{y_{1,\ell-k}}{y_{2,\ell-k}}\Big)^{(-1)^k}
\frac{y_{3,2(\ell-1-k)}}{y_{3,2(\ell-1-k)-1}}+\\ \nonumber
+\sum_{n=1}^{[\frac{\ell-1}{2}]}\Big( \frac{\partial}{\partial
y_{3,2(2n-1)-1}}- \frac{\partial}{\partial y_{3,2(2n-1)}}\Big)
\frac{y_{3,2(2n-1)}}{y_{1,2n}}\prod_{k=1}^{\ell-2n-1}
\Big(\frac{y_{1,\ell-k}}{y_{2,\ell-k}}\Big)^{(-1)^k}
\frac{y_{3,2(\ell-1-k)}}{y_{3,2(\ell-1-k)-1}}+\\ \nonumber
+\sum_{n=1}^{[\frac{\ell-2}{2}]}\Big( \frac{\partial}{\partial
y_{3,4n-1}}- \frac{\partial}{\partial y_{3,4n}}\Big)
\frac{y_{3,4n}}{y_{2,2n+1}}\prod_{k=1}^{\ell-2(n+1)}
\Big(\frac{y_{1,\ell-k}}{y_{2,\ell-k}}\Big)^{(-1)^k}
\frac{y_{3,2(\ell-1-k)}}{y_{3,2(\ell-1-k)-1}},\eqa

\bqa E_2= \frac{\partial}{\partial y_{2,\ell-1}}+
\sum_{n=2}^{[\ell/2]}\Big(\frac{\partial}{\partial y_{2,\ell-n-1}}-
\frac{\partial}{\partial y_{2,\ell-n}}\Big)\prod_{k=1}^{2(n-1)}
\Big(\frac{y_{2,\ell-k}}{y_{1,\ell-k}}\Big)^{(-1)^k}
\frac{y_{3,2(\ell-1-k)}}{y_{3,2(\ell-1-k)-1}}+\\+ \nonumber
\sum_{n=1}^{[\ell/2]}\Big(\frac{\partial}{\partial y_{1,\ell-n-1}}-
\frac{\partial}{\partial y_{1,\ell-n}}\Big)\prod_{k=1}^{2n-1}
\Big(\frac{y_{2,\ell-k}}{y_{1,\ell-k}}\Big)^{(-1)^k}
\frac{y_{3,2(\ell-1-k)}}{y_{3,2(\ell-1-k)-1}}+\\ +\nonumber
\sum_{n=1}^{[\frac{\ell-2}{2}]}\Big( \frac{\partial}{\partial
y_{3,4n-1}}- \frac{\partial}{\partial y_{3,4n}}\Big)
\frac{y_{3,4n}}{y_{1,2n+1}}\prod_{k=1}^{\ell-2(n+1)}
\Big(\frac{y_{2,\ell-k}}{y_{1,\ell-k}}\Big)^{(-1)^k}
\frac{y_{3,2(\ell-1-k)}}{y_{3,2(\ell-1-k)-1}}+\\ +\nonumber
\sum_{n=1}^{[\frac{\ell-1}{2}]}\Big( \frac{\partial}{\partial
y_{3,2(2n-1)-1}}- \frac{\partial}{\partial y_{3,2(2n-1)}}\Big)
\frac{y_{3,2(2n-1)}}{y_{1,2n}}\prod_{k=1}^{\ell-2n-1}
\Big(\frac{y_{2,\ell-k}}{y_{1,\ell-k}}\Big)^{(-1)^k}
\frac{y_{3,2(\ell-1-k)}}{y_{3,2(\ell-1-k)-1}},\eqa

\bqa E_k= \frac{\partial}{\partial y_{k,2(\ell+1-k)}}+
\sum_{n=1}^{\ell-k}\Big(\frac{\partial}{\partial y_{k,2n}}-
\frac{\partial}{\partial y_{k,2n+1}}\Big)\prod_{i=1}^{\ell+1-n-k}
\frac{y_{k,2(i+1)-1}}{y_{k,2(k+1)}}
\frac{y_{k+1,2i}}{y_{k+1,2i-1}}+\eqa\bqa +\nonumber
\sum_{n=1}^{\ell-k}\Big(\frac{\partial}{\partial y_{k+1,2n-1}}-
\frac{\partial}{\partial y_{k+1,2n}}\Big)
\frac{y_{k+1,2n}}{y_{k,2(n+1)}}\prod_{i=2}^{\ell+1-n-k}
\frac{y_{k,2(i+1)-1}}{y_{k,2(k+1)}} \frac{y_{k+1,2i}}{y_{k+1,2i-1}},
\qquad 2<k<\ell,\eqa
\bqa E_{\ell}=\frac{\partial}{\partial
y_{\ell,2}},\eqa \bqa
H_i=\<\mu,\,\alpha_i^\vee\>+\sum_{k=1}^\ell
a_{i,k}\,\sum_{j=1}^{n_k}y_{k,j}\frac{\partial}{\partial y_{k,j}},
\eqa \bqa
F_i=-\<\mu,\alpha_i^\vee\>\sum_{n=1}^{\ell-1}y_{i,n}-
\sum_{n=1}^{\ell-1}\Big(y_{i,n}^2\frac{\partial}{\partial y_{i,n}}+
2\sum_{k=n+1}^{\ell-1}y_{i,k}y_{i,n}
\frac{\partial}{\partial\,y_{i,n}}\Big)+\eqa\bqa \nonumber
+\sum_{n=1}^{2(\ell-2)-1}\sum_{k=[n/2]+2}y_{i,k}y_{3,n}
\frac{\partial}{\partial\,y_{3,n}},\qquad i=1,2,\eqa
\bqa
F_k=-\<\mu,\alpha_k^\vee\>\sum_{n=1}^{2(\ell+1-k)}y_{k,n}-
\sum_{n=1}^{2(\ell+1-k)}\Big(y_{k,n}^2
\frac{\partial}{\partial\,y_{k,n}}+
2\sum_{i=n+1}^{2(\ell+1-k)}y_{k,i}y_{k,n}
\frac{\partial}{\partial\,y_{k,n}}\Big)+\\ \nonumber
+\sum_{n=1}^{\ell+2-k}\,\,\sum_{i=2(n-1)}^{2(\ell+1-k)}y_{k,i}
\Big(y_{k-1,2n-1}\frac{\partial}{\partial\,y_{k-1,2n-1}}+
y_{k-1,2n}\frac{\partial}{\partial\,y_{k-1,2n}}\Big)+\eqa\bqa \nonumber
+\sum_{n=1}^{2(\ell-k)-1}\sum_{i=2[n/2]+3}y_{k,i}
y_{k+1,n}\frac{\partial}{\partial\,y_{k+1,n}},\qquad 3\leq k\leq\ell,\eqa
where $n_1=n_2=\ell-1$, $n_k=2(\ell+1-k)$, for $2<k\leq\ell$.
\end{prop}
The proof is given in Part II, Section \ref{genso2l}.

The left/right Whittaker vectors in the factorized parametrization
can be found explicitly.

\begin{lem} The following expressions for the left/right
Whittaker vectors hold:
\bqa\psi_R(y)=\exp\Big\{-\Big(\sum_{n=1}^{\ell-1}y_{1,n}+
\sum_{n=1}^{\ell-1}y_{2,n}+ \sum_{k=3}^\ell
\sum_{n=1}^{2(\ell+1-k)}y_{k,n}\Big)\Big\},\eqa

\bqa\psi_L(y)= \Big(\prod_{n=1}^{\ell/2}y_{1,2n-1}
\prod_{n=1}^{\frac{\ell-1}{2}}y_{2,2n} \prod_{n=3}^\ell\prod_{i=3}^n
y_{i,2(n+1-i)-1}\Big)^{\<\mu,\alpha_1^\vee\>}\times \\\times
\nonumber \Big(\prod_{n=1}^{\ell/2}y_{2,2n-1}
\prod_{n=1}^{\frac{\ell-1}{2}}y_{1,2n} \prod_{n=3}^\ell\prod_{i=3}^n
y_{i,2(n+1-i)-1}\Big)^{\<\mu,\alpha_2^\vee\>}\times \\\times
\nonumber\prod_{k=3}^\ell\Big(\prod_{i=1}^k\prod_{n=1}^{n_i}y_{i,n}
\prod_{i=k+1}^\ell\prod_{n=1}^{n_i/2}y_{i,2n-1}^2
\Big)^{\<\mu,\alpha_k^\vee\>}\times \\ \times \nonumber
\exp\Big\{-\Big( \sum_{n=1}^{\ell-1}\Big[
\frac{1}{y_{1,\ell-1}}\prod_{k=1}^{n-1}\Big(
\frac{y_{1,\ell-k}}{y_{1,\ell-k-1}}\Big)^{p_{k-1}} \Big(
\frac{y_{2,\ell-k}}{y_{2,\ell-k-1}}\Big)^{p_{k}}+\\+ \nonumber
\frac{1}{y_{2,\ell-1}}\prod_{k=1}^{n-1}\Big(
\frac{y_{1,\ell-k}}{y_{1,\ell-k-1}}\Big)^{p_{k}} \Big(
\frac{y_{2,\ell-k}}{y_{2,\ell-k-1}}\Big)^{p_{k+1}} \Big]\times \\
\times \nonumber
\Big(1+\frac{y_{3,2(\ell-n-1)}}{y_{3,2(\ell-n-1)-1}}\Big)
\prod_{k=1}^{n-1}\frac{y_{3,2(\ell-k-1)}}{y_{3,2(\ell-k-1)-1}}+
\sum_{k=3}^\ell\frac{1}{y_{k,2(\ell+1-k)}}\Big)\Big\}, \eqa where
$n_1=n_2=\ell-1$ and $n_k=2(\ell+1-k),\,k>2$, and $p_k=(1-(-1)^k)$
is the parity of $k$.
\end{lem}
The proof is given in Part II, Section \ref{melso2l}.

Using (\ref{psir}) and (\ref{psil}) it is easy to obtain the
integral representations of $\mathfrak{so}_{2\ell}$-Whittaker
function in the factorized parametrization.
\begin{te}\label{teintfso2l}
 The eigenfunctions of  $\mathfrak{so}_{2\ell}$-Toda
chain admit the following integral representation:
\bqa\label{waveggc}\Psi_{\lambda_1,\ldots,\lambda_{\ell}}^{D_{\ell}}
(x_1,\ldots,x_{\ell})= e^{\imath\lambda_1x_1+\ldots+
\imath\lambda_{\ell}x_{\ell}}\int_C\bigwedge_{i=1}^\ell
\bigwedge_{k=1}^{n_i} \frac{dy_{i,k}}{y_{i,k}}\times\\\times
\nonumber \Big(\prod_{n=1}^{\ell/2}y_{1,2n-1}
\prod_{n=1}^{\frac{\ell-1}{2}}y_{2,2n} \prod_{n=3}^\ell\prod_{i=3}^n
y_{i,2(n+1-i)-1}\Big)^{\imath(\lambda_2-\lambda_1)}\times \\\times
\nonumber \Big(\prod_{n=1}^{\ell/2}y_{2,2n-1}
\prod_{n=1}^{\frac{\ell-1}{2}}y_{1,2n} \prod_{n=3}^\ell\prod_{i=3}^n
y_{i,2(n+1-i)-1}\Big)^{\imath(\lambda_1+\lambda_2)}\times \\ \times
\nonumber\prod_{k=3}^\ell\Big(\prod_{i=1}^k\prod_{n=1}^{n_i}y_{i,n}
\prod_{i=k+1}^\ell\prod_{n=1}^{n_i/2}y_{i,2n-1}^2
\Big)^{\imath(\lambda_k-\lambda_{k-1})}\times \eqa
\bqa\nonumber\times \exp\Big\{-\Big( \sum_{n=1}^{\ell-1}\Big[
\frac{1}{y_{1,\ell-1}}\prod_{k=1}^{n-1}\Big(
\frac{y_{1,\ell-k}}{y_{1,\ell-k-1}}\Big)^{p_{k-1}} \Big(
\frac{y_{2,\ell-k}}{y_{2,\ell-k-1}}\Big)^{p_{k}}+\\+ \nonumber
\frac{1}{y_{2,\ell-1}}\prod_{k=1}^{n-1}\Big(
\frac{y_{1,\ell-k}}{y_{1,\ell-k-1}}\Big)^{p_{k}} \Big(
\frac{y_{2,\ell-k}}{y_{2,\ell-k-1}}\Big)^{p_{k+1}} \Big]\times \\
\times \nonumber
\Big(1+\frac{y_{3,2(\ell-n-1)}}{y_{3,2(\ell-n-1)-1}}\Big)
\prod_{k=1}^{n-1}\frac{y_{3,2(\ell-k-1)}}{y_{3,2(\ell-k-1)-1}}+
\sum_{k=3}^\ell\frac{1}{y_{k,2(\ell+1-k)}}+\\+ \nonumber
e^{x_2-x_1}\sum_{n=1}^{\ell-1}y_{1,n}+
e^{x_1+x_2}\sum_{n=1}^{\ell-1}y_{2,n}+ \sum_{k=3}^\ell
e^{x_k-x_{k-1}}\sum_{n=1}^{2(\ell+1-k)}y_{k,n}\Big)\Big\}.\eqa Here
we assume $x_1=x_{\ell,1},\ldots,x_{\ell,\ell}$, $p_k=(1-(-1)^k)/2$,
$n_1=n_2=\ell-1$ and $n_k=2(\ell+1-k),\,k>2$. Domain of integration
$C \subset N_+$ is a middle-dimensional non-compact submanifold such
 that the integrand decreases exponentially  at the possible  boundaries
 and infinities. In particular one can chose $C=\RR_{+}^{m},$ where $m=l(w_0)$.
\end{te}
The proof is given in Part II, Section \ref{melso2l}.

\begin{ex} For $\ell=3$ the general formula (\ref{waveggc}) acquires
the following form
\bqa\Psi^{D_3}_{\lambda_1,\,\lambda_2,\,\lambda_3}(x_{31},x_{32},x_{33})=
\int_C\,\,
\bigwedge_{i=1}^3\bigwedge_{k=1}^{2}\frac{dy_{i,k}}{y_{i,k}} \times \\
\nonumber \times (y_{11}y_{31}y_{22})^{\imath(\lambda_2-\lambda_1)}
(y_{21}y_{31}y_{12})^{\imath(\lambda_2+\lambda_1)}
(y_{31}y_{12}y_{22}y_{32})^{\imath(\lambda_3-\lambda_2)}\times \\
\nonumber \times
\exp\Big\{\frac{1}{y_{12}}\Big(1+\frac{y_{32}}{y_{31}}\Big)+
\frac{1}{y_{12}}\frac{y_{22}}{y_{21}}\frac{y_{32}}{y_{31}}+
\frac{1}{y_{22}}\Big(1+\frac{y_{32}}{y_{31}}\Big)+
\frac{1}{y_{22}}\frac{y_{12}}{y_{11}}\frac{y_{32}}{y_{31}}+
\frac{1}{y_{32}}+\\ \nonumber +e^{x_{32}-x_{31}}(y_{11}+y_{12})+
e^{x_{32}+x_{31}}(y_{21}+y_{22})+e^{x_{33}-x_{32}}(y_{31}+y_{32})\Big\}.\eqa
For the domain of integration one can chose $C=\mathbb{R}^6$. 
\end{ex}

\subsubsection{ $\mathfrak{so}_{2\ell}$-Whittaker function:
modified factorized   parametrization}

In this part we introduce a modified  factorized parametrization of
an open part $N^{(0)}_+$ of maximal unipotent subgroup $N_+\subset
SO(2\ell)$.   We use this parametrization to construct integral
representations for $\mathfrak{so}_{2\ell}$-Whittaker functions.
Similar to other series of classical Lie algebras these integral
representations for $\mathfrak{so}_{2\ell}$-Whittaker functions have
a simple recursive structure over the rank $\ell$ and can be
describe in purely combinatorial terms using suitable graphs. These
representations can be considered as a generalization of Givental
integral representations to $\mathfrak{g}=\mathfrak{so}_{2\ell}$.

We follow the same approach that was used in the description of
modified factorized representation for other classical groups. There
exists a realization of a tautological representation
$\pi_{2\ell}:\mathfrak{so}_{2\ell}\to End(\mathbb{C}^{2\ell})$ such
that Weyl generators corresponding to  Borel (Cartan) subalgebra of
$\mathfrak{so}_{2\ell}$ are realized by
 upper triangular (diagonal) matrices.
This defines an embedding $\mathfrak{so}_{2\ell}\subset
\mathfrak{gl}_{2\ell}$ such that Borel (Cartan) subalgebra maps into
Borel (Cartan) subalgebra (see e.g. \cite{DS}). To define the
corresponding embedding of the groups  consider the following
involution on $GL(2\ell)$: \bqa g\longmapsto g^*:=\dot{w}_0\cdot
(g^{-1})^{t}\cdot\dot{w}_0^{-1},\eqa where $a^{t}$ is induced by the
standard transposition of the  matrix $a$ and $\dot{w}_0$ is a lift
of the longest element of the Weyl group of $\mathfrak{gl}_{2\ell}$.
In the matrix form it can be written as
$$ \pi_{2\ell}(\dot{w}_0)=S\cdot J,$$ where
 $S={\rm diag}(1,-1,\ldots,-1,1)$ and $J=\|J_{i,j}\|=\|\delta_{i+j,2\ell+2}\|$.
 The orthogonal group $G=SO(2\ell)$ then can be defined as a
 following  subgroup of   $GL(2\ell)$  (see i.e. \cite{DS}):
$$ SO(2\ell)=\{g\in GL(2\ell):g^{*}=g\}.$$

Let $\epsilon_{i,j}$ be elementary $(2\ell\times 2\ell)$ matrices
with unites at  $(i,j)$ place and zeroes otherwise. Introduce the
following  matrices

\bqa U_n=\sum_{i=1}^n\epsilon_{\ell-n+i,\ell-n+i}+
e^{-x_{n-1,1}}\epsilon_{\ell+1,\ell+1}+
\sum_{i=1}^{n-1}e^{x_{n-1,i}}\epsilon_{\ell+i+1,\ell+i+1}+\\
+\nonumber e^{z_{n-1,1}}\epsilon_{\ell,\ell+1}+
\sum_{i=1}^{n-1}e^{z_{n-1,i}}\epsilon_{\ell+i,\ell+i+1}+
\sum_{i=1}^{\ell-n}(\epsilon_{i,i}+\epsilon_{2\ell+1-i,2\ell+1-i}),\eqa
\bqa \widetilde{U}'_n=\sum_{i=1}^{n+2}\epsilon_{\ell-n+i,\ell-n+i}+
\sum_{i=2}^{n-1}e^{x_{n-1,i}}\epsilon_{\ell+i+1,\ell+i+1}+
\sum_{i=2}^{n-1}e^{z_{n-1,i}}\epsilon_{\ell+i,\ell+i+1}+\\ \nonumber
+\sum_{i=1}^{\ell-n}(\epsilon_{i,i}+\epsilon_{2\ell+1-i,2\ell+1-i}),\eqa
\bqa \widetilde{U}''_n=\sum_{i=1}^n\epsilon_{\ell-n+i,\ell-n+i}+
e^{-x_{n-1,1}}\epsilon_{\ell+1,\ell+1}+
e^{x_{n-1,1}}\epsilon_{\ell+2,\ell+2}+\\ \nonumber
+e^{z_{n-1,1}}(\epsilon_{\ell,\ell+1}+\epsilon_{\ell+1,\ell+2})+
\sum_{i=3}^n\epsilon_{\ell+i,\ell+i}+
\sum_{i=1}^{\ell-n}(\epsilon_{i,i}+\epsilon_{2\ell+1-i,2\ell+1-i}),\eqa
\bqa
\widetilde{U}_n=\widetilde{U}'_n\widetilde{U}''_n=
\sum_{i=1}^n\epsilon_{\ell-n+i,\ell-n+i}+
e^{-x_{n-1,1}}\epsilon_{\ell+1,\ell+1}+
\sum_{i=1}^{n-1}e^{x_{n-1,i}}\epsilon_{\ell+i+1,\ell+i+1}+\\
+\nonumber e^{z_{n-1,1}}\epsilon_{\ell,\ell+1}+
\sum_{i=1}^{n-1}e^{z_{n-1,i}}\epsilon_{\ell+i,\ell+i+1},\eqa

\bqa V_n=\sum_{i=1}^ne^{x_{n,n+1-i}}\epsilon_{\ell-n+i,\ell-n+i}+
e^{-x_{n,1}}\epsilon_{\ell+1,\ell+1}+\sum_{i=2}^n\epsilon_{\ell+i,\ell+i}+\\
\nonumber
+\sum_{i=1}^{\ell-n}(\epsilon_{i,i}+\epsilon_{2\ell+1-i,2\ell+1-i}),\eqa
\bqa\widetilde{V}'_n=
\sum_{i=1}^{n-1}e^{x_{n,n+1-i}}\epsilon_{\ell-n+i,\ell-n+i}+
\sum_{i=-1}^n\epsilon_{\ell+i,\ell+i}+
+\sum_{i=1}^{n-2}e^{z_{n-1,n-i}}\epsilon_{\ell-n+i,\ell-n+i+1}+\\
\nonumber
+\sum_{i=1}^{\ell-n}(\epsilon_{i,i}+\epsilon_{2\ell+1-i,2\ell+1-i}),\eqa
\bqa\widetilde{V}''_n=e^{x_{n,1}}\epsilon_{\ell,\ell}+
e^{-x_{n,1}}\epsilon_{\ell+1,\ell+1}+
e^{z_{n-1,1}}(\epsilon_{\ell-1,\ell}+\epsilon_{\ell,\ell+1})+\\
\nonumber
+\sum_{i=1}^{\ell-n}(\epsilon_{i,i}+\epsilon_{2\ell+1-i,2\ell+1-i}),\eqa
\bqa\widetilde{V}_n=\widetilde{V}''_n\widetilde{V}'_n=
\sum_{i=1}^ne^{x_{n,n+1-i}}\epsilon_{\ell-n+i,\ell-n+i}+
e^{-x_{n,1}}\epsilon_{\ell+1,\ell+1}+\sum_{i=2}^n\epsilon_{\ell+i,\ell+i}+\\
+\nonumber
\sum_{i=1}^{n-1}e^{z_{n-1,n-i}}\epsilon_{\ell-n+i,\ell-n+i+1}+
e^{z_{n-1,1}}\epsilon_{\ell,\ell+1}+
\sum_{i=1}^{\ell-n}(\epsilon_{i,i}+\epsilon_{2\ell+1-i,2\ell+1-i}).\eqa

\begin{te}\label{Mparso2l}

 i)  The image of a generic unipotent element
  $v^{D_{\ell}}\in N_+$ in the tautological representation
  $\pi_{2\ell}: \mathfrak{so}_{2\ell}\to End(\mathbb{C}^{2\ell})$
can be presented in the form \bqa \label{modfacparD}
v^{D_{\ell}}=\mathfrak{X}_2\mathfrak{X}_{3}\cdots\mathfrak{X}_{\ell},\eqa
with \bqa \mathfrak{X}_2=S_1\widetilde{U}_2U_2^{-1}S_1\cdot
S_3(\widetilde{U}_2U_2^{-1})^*S_3\cdot
S_1(\tilde{V}_2V_2^{-1})^*S_1\cdot S_3\tilde{V}_2V_2^{-1}S_3\nonumber\\
\mathfrak{X}_n= (\widetilde{U}'_n(U'_n)^{-1})^*\cdot
S_{n-1}\widetilde{U}_nU_n^{-1}S_{n-1}\cdot
S_{n+1}(\widetilde{U}''_n(U''_n)^{-1})^*S_{n+1}\cdot\label{modfacparD1}\\
\nonumber \cdot S_{n-1}(\widetilde{V}''_n(V''_n)^{-1})^*S_{n-1}\cdot
S_{n+1}\widetilde{V}_nV_n^{-1}S_{n+1}\cdot
(\widetilde{V}'_n(V'_n)^{-1})^*,\eqa where
$x_{\ell,k}=0,\,\,k=1,\ldots\ell$ is assumed and $S_i$ is defined as
follows: \bqa S_i=\sum_{k=1}^{i-1}\epsilon_{k,k}+
\epsilon_{i,i+1}+\epsilon_{i+1,i}+\sum_{k=i+2}^{2\ell}\epsilon_{k,k}.\eqa

ii) This defines a parametrization of an open part $N_+^{(0)}$ of
$N_+$.
\end{te}

{\it Proof.} Let  $v^{D_{\ell}}(y)$ be a parametrization of  $N_+$
according to (\ref{rec1D})-(\ref{rec2D}). Let $\tilde{X}_i(y)=e^{y
e_{i,i+1}}$ be a  one-parametric unipotent subgroup in $GL(2\ell),$
then $\tilde{X}_i(y)^*=\tilde{X}_{2\ell+1-i}(y).$ Embed elementary
unipotent subgroups $X_i(y)$ of $SO(2\ell)$ into $GL(2\ell)$ as
follows:
$$X_i(y)=\tilde{X}_i(y)^*\cdot\tilde{X}_i(y).$$
 This defines a map of   an arbitrary regular unipotent element
 $v^{D_\ell}$
into unipotent subgroup of $GL(2\ell).$ Let us  change the variables
in the following way: \bqa\label{nontwistd}
y_{1,n}=\Big(e^{z_{n,1}-x_{n,1}}+
e^{z_{n,1}-x_{n+1,1}}\Big),\hspace{2cm}n=1,\ldots,\ell-1,\nonumber \\
y_{2,n}=\Big(e^{z_{n,1}+x_{n,1}}+
e^{z_{n,1}+x_{n+1,1}}\Big),\hspace{2cm}n=1,\ldots,\ell-1,\\
y_{k,2r-1}=e^{z_{k+r-2,k-1}-x_{k+r-2,k-1}},\hspace{2cm}k=3,\ldots,\ell,
\nonumber\\ \nonumber
y_{k,2r}=e^{z_{k+r-2,k-1}-x_{k+r-1,k-1}},\hspace{2cm}
r=1,\ldots,\ell+1-k. \eqa Here the conditions
$x_{\ell,k}=0,\,\,\,k=1,\ldots,\ell$ are implied. By elementary
manipulations it is easy to check that after the change of
variables (\ref{nontwistd}), the image $\pi_{2\ell}(v^{D_{\ell}})$ of $v^{D_{\ell}}$
defined by (\ref{rec1D})-(\ref{rec2D})  transforms into the
(\ref{modfacparD}) -(\ref{modfacparD1}). Taking into account that
the change of variables (\ref{1yDn}) is invertible
 we obtain a parametrization of $N_+^{(0)}\subset N_+$ $\Box$

The modified factorized parametrization of a unipotent group $N_+$
 defines a particular realization of a principal series
representation of $\mathcal{U}(\mathfrak{so}_{2\ell})$ by
differential operators. It can be   obtained using the change of
variables (\ref{1yDn}) applied to a realization given in Proposition
\ref{GGrepD}.   We shall use the term Gauss-Givental representation
for this realization of representation of
 $\mathcal{U}(\mathfrak{so}_{2\ell})$.

\begin{prop}
The following differential operators define a representation
$\pi_{\lambda}$ of $\mathcal{U}(\mathfrak{so}_{2\ell})$ in $V_{\mu}$
in terms of modified factorized parametrization  of $N_+^{(0)}$: \bqa
E_1=e^{x_{22}-z_{11}}\frac{e^{x_{11}}}{e^{x_{11}}+e^{x_{21}}}
\Big(-\frac{\partial}{\partial x_{11}}- \frac{e^{x_{11}}}
{e^{x_{11}}+e^{x_{21}}}\frac{\partial}{\partial z_{11}}\Big)+\\
\nonumber +\sum_{k=2}^{\ell-1}\Big\{\,
\frac{p_ke^{x_{k,1}}+p_{k+1}e^{x_{k+1,1}}}
{e^{x_{k,1}}+e^{x_{k+1,1}}}\Big(e^{x_{k,2}-z_{k,1}}+
e^{x_{k+1,2}-z_{k,1}}\Big)\times \\ \nonumber \times \left(
-\frac{\partial}{\partial x_{11}}-\frac{e^{x_{21}}-e^{x_{11}}}
{e^{x_{11}}+e^{x_{21}}}\frac{\partial}{\partial z_{11}}+
\frac{\partial}{\partial x_{21}}-\frac{\partial}{\partial x_{22}}+
\right.\\ \nonumber \left.
+\sum_{i=3}^k(-1)^{i}\frac{\partial}{\partial x_{i,1}}-
\frac{\partial}{\partial x_{i,2}}+
(-1)^{i-1}\frac{e^{x_{i,1}}-e^{x_{i-1,1}}}
{e^{x_{i,1}}+e^{x_{i-1,1}}}\frac{\partial}{\partial z_{i-1,1}}-
\frac{\partial}{\partial z_{i-1,2}} \right)+\\ \nonumber
+\left(e^{x_{k+1,2}+(-1)^{k}x_{k+1,1}}
\frac{1}{e^{z_{k,1}+x_{k,1}}+e^{z_{k,1}+x_{k+1,1}}}
\frac{p_ke^{x_{k,1}}+p_{k+1}e^{x_{k+1,1}}}
{e^{x_{k,1}}+e^{x_{k+1,1}}}-\right.\\ \nonumber\left.
-e^{x_{k,2}+(-1)^{k-1}x_{k,1}}
\frac{1}{e^{z_{k,1}-x_{k,1}}+e^{z_{k,1}-x_{k+1,1}}}
\frac{p_{k-1}e^{x_{k,1}}+p_{k}e^{x_{k+1,1}}}
{e^{x_{k,1}}+e^{x_{k+1,1}}}\right)\frac{\partial}{\partial z_{k,1}}
-\\ \nonumber -e^{x_{k+1,2}-z_{k,1}}
-\frac{p_ke^{x_{k,1}}+p_{k+1}e^{x_{k+1,1}}}
{e^{x_{k,1}}+e^{x_{k+1,1}}}\frac{\partial}{\partial z_{k,2}}
\Big\},\eqa

\bqa E_2=e^{x_{22}-z_{11}}\frac{e^{x_{21}}}{e^{x_{11}}+e^{x_{21}}}
\Big(\frac{\partial}{\partial x_{11}}+\frac{e^{x_{21}}}
{e^{x_{11}}+e^{x_{21}}}\frac{\partial}{\partial z_{11}}\Big)+\\
\nonumber+ \sum_{k=2}^{\ell-1}\Big\{\,
\frac{p_{k-1}e^{x_{k,1}}+p_{k}e^{x_{k+1,1}}}
{e^{x_{k,1}}+e^{x_{k+1,1}}}\Big(e^{x_{k,2}-z_{k,1}}+
e^{x_{k+1,2}-z_{k,1}}\Big)\times \\ \nonumber \times \left(
\frac{\partial}{\partial x_{11}}+\frac{e^{x_{21}}-e^{x_{11}}}
{e^{x_{11}}+e^{x_{21}}}\frac{\partial}{\partial z_{11}}-
\frac{\partial}{\partial x_{21}}-\frac{\partial}{\partial x_{22}}+
\right.\\ \nonumber \left.
+\sum_{i=3}^k(-1)^{i-1}\frac{\partial}{\partial x_{i,1}}-
\frac{\partial}{\partial x_{i,2}}+
(-1)^{i}\frac{e^{x_{i,1}}-e^{x_{i-1,1}}}
{e^{x_{i,1}}+e^{x_{i-1,1}}}\frac{\partial}{\partial z_{i-1,1}}-
\frac{\partial}{\partial z_{i-1,2}} \right)+\\ \nonumber
+\left(e^{x_{k+1,2}+(-1)^{k-1}x_{k+1,1}}
\frac{1}{e^{z_{k,1}+x_{k,1}}+e^{z_{k,1}+x_{k+1,1}}}
\frac{p_{k-1}e^{x_{k,1}}+p_{k}e^{x_{k+1,1}}}
{e^{x_{k,1}}+e^{x_{k+1,1}}}-\right.\\ \nonumber\left.
-e^{x_{k,2}+(-1)^{k}x_{k,1}}
\frac{1}{e^{z_{k,1}-x_{k,1}}+e^{z_{k,1}-x_{k+1,1}}}
\frac{p_{k}e^{x_{k,1}}+p_{k+1}e^{x_{k+1,1}}}
{e^{x_{k,1}}+e^{x_{k+1,1}}}\right)\frac{\partial}{\partial z_{k,1}}
-\\ \nonumber -e^{x_{k+1,2}-z_{k,1}}
-\frac{p_{k-1}e^{x_{k,1}}+p_{k}e^{x_{k+1,1}}}
{e^{x_{k,1}}+e^{x_{k+1,1}}}\frac{\partial}{\partial z_{k,2}}
\Big\}.\eqa Here $p_k=(1-(-1)^k)/2$ is the parity of $k$.

\bqa E_k=\Big(\frac{\partial}{\partial z_{k-1,k-1}}-
\frac{\partial}{\partial x_{k-1,k-1}}\Big)\Big(
e^{x_{k,k}-z_{k-1,k-1}}+\sum_{n=k}^{\ell-1}e^{x_{n,k}-z_{n,k-1}}+
e^{x_{n+1,k}-z_{n,k}}\Big)+\\ \nonumber
+\sum_{i=k}^{\ell-1}\Big(\frac{\partial}{\partial x_{i,k-1}}-
\frac{\partial}{\partial\,x_{i,k}}\Big)\sum_{n=i}^{\ell-1}\Big(
e^{x_{n,k}-z_{n,k-1}}+e^{x_{n+1,k}-z_{n,k-1}}\Big)+\\ \nonumber
+\sum_{i=k}^{\ell-1}\Big(\frac{\partial}{\partial z_{i,k-1}}-
\frac{\partial}{\partial\,z_{i,k}}\Big)\Big(
e^{x_{i+1,k}-z_{i,k-1}}+ \sum_{n=i+1}^{\ell-1}
e^{x_{n,k}-z_{n,k-1}}+e^{x_{n+1,k}-z_{n,k-1}}\Big),\eqa where
$2<k<\ell$ and \bqa
E_{\ell}=e^{-z_{\ell-1,\ell-1}}\Big(
\frac{\partial}{\partial\,z_{\ell-1,\ell-1}}+
\frac{\partial}{\partial\,x_{\ell-1,\ell-1}}\Big),\eqa \bqa
H_1=\<\mu,\alpha_1^\vee\>)+ 2\Big(-\frac{\partial}{\partial
x_{\ell-1,1}}-\frac{\partial}{\partial x_{11}}+\\
+\nonumber \sum_{n=1}^{\ell-1}
\frac{p_ne^{x_{n,1}}+p_{n+1}e^{x_{n+1,1}}}
{e^{x_{n,1}}+e^{x_{n+1,1}}}\frac{\partial}{\partial z_{n,1}}\Big)-
\sum_{k=2}^{\ell-1}\frac{\partial}{\partial z_{k,2}},\eqa

\bqa H_2=\<\mu,\alpha_2^\vee\>+ 2\Big(\frac{\partial}{\partial
x_{\ell-1,1}}+\frac{\partial}{\partial x_{11}}+\\
\nonumber+ \sum_{n=1}^{\ell-1}
\frac{p_{n-1}e^{x_{n,1}}+p_ne^{x_{n+1,1}}}
{e^{x_{n,1}}+e^{x_{n+1,1}}}\frac{\partial}{\partial z_{n,1}}\Big)-
\sum_{k=2}^{\ell-1}\frac{\partial}{\partial z_{k,2}}, \eqa where
$p_k=(1-(-1)^k)/2$. \bqa
H_i=\<\mu,\alpha_i^\vee\>+\sum_{k=2}^\ell\,a_{i,k}\,
\sum_{j=k-1}^{\ell-1}\frac{\partial}{\partial z_{j,k-1}},\qquad
2<i\leq\ell,\eqa
 \bqa
F_1=-\sum_{n=1}^{\ell-1}
\Big(e^{z_{n,1}-x_{n+1,1}}+e^{z_{n,1}-x_{n+1,1}}\Big)
\Big[\<\mu,\alpha_1^\vee\>-\frac{\partial}{\partial x_{11}}+\\
+\nonumber\frac{e^{x_{11}}}{e^{x_{11}}+e^{x_{21}}}
\frac{\partial}{\partial z_{11}} \Big]-e^{x_{\ell,1}-x_{\ell,2}}
\sum_{k=2}^{\ell-1}\Big[-\frac{\partial}{\partial x_{k,1}}+
\frac{\partial}{\partial x_{k,2}}+\\ \nonumber
+\frac{e^{x_{k,1}}}{e^{x_{k,1}}+e^{x_{k+1,1}}}
\frac{\partial}{\partial z_{k,1}}+
\frac{e^{x_{k,1}}}{e^{x_{k-1,1}}+e^{x_{k,1}}}
\frac{\partial}{\partial z_{k-1,1}}\Big] \sum_{n=k}^{\ell-1}
\Big(e^{z_{n,1}-x_{n+1,1}}+e^{z_{n,1}-x_{n+1,1}}\Big),\eqa

\bqa F_2=-\sum_{n=1}^{\ell-1}
\Big(e^{z_{n,1}+x_{n+1,1}}+e^{z_{n,1}+x_{n+1,1}}\Big)
\Big[\<\mu,\alpha_2^\vee\>+\frac{\partial}{\partial x_{11}}+\\
+\nonumber\frac{e^{x_{21}}}{e^{x_{11}}+e^{x_{21}}}
\frac{\partial}{\partial z_{11}} \Big]-
\sum_{k=2}^{\ell-1}\Big[\frac{\partial}{\partial x_{k,1}}+
\frac{\partial}{\partial x_{k,2}}+\\ \nonumber
+\frac{e^{x_{k+1,1}}}{e^{x_{k,1}}+e^{x_{k+1,1}}}
\frac{\partial}{\partial z_{k,1}}+
\frac{e^{x_{k-1,1}}}{e^{x_{k-1,1}}+e^{x_{k,1}}}
\frac{\partial}{\partial z_{k-1,1}}\Big] \sum_{n=k}^{\ell-1}
\Big(e^{z_{n,1}+x_{n+1,1}}+e^{z_{n,1}+x_{n+1,1}}\Big),\eqa

\bqa F_k=\Big(-\<\mu,\alpha_k\>+\\
+\nonumber \frac{\partial}{\partial x_{k-1,k-1}}+
\frac{\partial}{\partial z_{k-2,k-2}}\Big) \sum_{n=k-1}^{\ell-1}
\Big(e^{z_{n,k-1}-x_{n,k-1}}+ e^{z_{n,k-1}-x_{n+1,k-1}}\Big)-\\
-\nonumber \sum_{i=k-1}^{\ell-1}\left(
\frac{\partial}{\partial z_{i,k-1}}-\frac{\partial}{\partial
z_{i,k-2}}\right)\Big(e^{z_{i,k-1}-x_{i+1,k-1}}+\\ \nonumber
+\sum_{j=i+1}^{\ell-1}e^{z_{j,k-1}-x_{j,k-1}}+
e^{z_{j,k-1}-x_{j+1,k-1}}\Big)-\\ \nonumber
-\sum_{i=k}^{\ell-1}\left(
\frac{\partial}{\partial x_{i,k-1}}-\frac{\partial}{\partial
x_{i,k-2}}\right)\sum_{j=i+1}^{\ell-1}\Big(e^{z_{j,k-1}-x_{j,k-1}}+
e^{z_{j,k-1}-x_{j+1,k-1}}\Big),\eqa where $3\leq k\leq\ell$ and 
$x_{\ell,k}=0$ is assumed.
\end{prop}

We are going to write down the matrix element (\ref{pairing}) for
$\mathfrak{g}=\mathfrak{so}_{2\ell}$ explicitly using
Gauss-Givental representation defined above. Whittaker vectors
$\psi_R$ and $\psi_L$ in this representation should  satisfiy the
system of differential equations \bqa\label{bwitt1} E_i\psi_R(x)=
-\psi_R(x),\qquad F_i\psi_L(x)= -\psi_L(x), \qquad 1\leq
i\leq\ell.\eqa Its solution has the following form.

\begin{lem}\label{WmrD}
The following expressions for the  left/right Whittaker vectors hold:
\bqa \psi_R=\exp\Big\{-\sum_{n=1}^{\ell-1}\,\, \Big(
e^{z_{n,1}-x_{n,1}}+e^{z_{n,1}-x_{n+1,1}}+
e^{z_{n,1}+x_{n,1}}+e^{z_{n,1}+x_{n+1,1}}\Big)\,-\\- \nonumber
\sum_{k=3}^{\ell}\,\,\sum_{n=1}^{\ell+1-k}
\Big(e^{z_{k+n-2,k-1}-x_{k+n-2,k-1}}+e^{z_{k+n-2,k-1}-x_{k+n-1,k-1}}
\Big)\Big\},\eqa \bqa\psi_L=e^{2\mu_1x_{1,1}}
\prod_{n=2}^\ell\Big(e^{x_{n,1}}+e^{x_{n-1,1}}\Big)^{2\mu_n}\times\\
\nonumber \times\prod_{n=1}^\ell\exp\Big\{
-\mu_n\Big(\sum_{i=1}^nx_{n,i}-
2\sum_{i=1}^{n-1}z_{n-1,i}+\sum_{i=1}^{n-1}x_{n-1,i}\Big)
\Big\}\times
\\\times \nonumber
\exp\Big\{-\sum_{k=1}^{\ell-1}\Big(e^{x_{k+1,k+1}-z_{k,k}}+
\sum_{i=k+1}^{\ell-1}e^{x_{i,k+1}-z_{i,k}}+e^{x_{i+1,k+1}-z_{i,k}}
\Big)\Big\},\eqa where we set $x_{\ell,k}=0,\,\,\,\,
k=1,\ldots,\ell$ and $\mu_n=\imath\lambda_n-\rho_n,\,\,\,\ \rho_1=0$
and $\rho_n=n-1$ for $1<n\leq\ell$. $(\sum_i^j=0$ when $j<i).$
\end{lem}

Now we are ready to find the integral representation of the pairing
(\ref{pairing}) for $\mathfrak{g}=\mathfrak{so}_{2\ell}$.
 To get an explicit expression for the integrand, one
uses the same type of  decomposition of the Cartan element as for
other classical groups  in the previous subsections:
$$e^{-H_x}=\pi_\lambda(\exp(-\sum_{i=1}^{\ell}
\langle\omega_i,x\rangle h_i))=e^{H_L}e^{H_R},$$ where \bqa\label{HRso2l}
-H_x=H_L+H_R= \sum_{i=1}^\ell\mu_ix_{\ell,i}+ \sum_{k=3}^\ell
(x_{\ell,k}-x_{\ell,k-1})\sum_{i=k-1}^{\ell-1}
\frac{\partial}{\partial z_{i,k-1}}+
x_{\ell,2}\sum_{i=1}^{\ell-1}\frac{\partial}{\partial z_{i,1}}+\\+
\nonumber x_{\ell,1}\Big(\frac{\partial}{\partial x_{\ell-1,1}}+
\frac{\partial}{\partial x_{1,1}}-\sum_{k=1}^{\ell-1}\,(-1)^k\,
\frac{e^{x_{k+1,1}}-e^{x_{k,1}}}{e^{x_{k,1}}+e^{x_{k+1,1}}}
\frac{\partial}{\partial z_{k,1}}\,\,\Big),\eqa
 with
  \bqa
H_L=\sum_{k=1}^\ell\,x_{\ell,k}\Big(
\sum_{i=k}^{\ell-1}\frac{\partial}{\partial x_{i,k}}+
\sum_{i=k-1}^{\ell-1}\frac{\partial}{\partial z_{i,k-1}}\Big),\eqa
\bqa H_R=-H_x-H_L .\eqa
 We imply that $H_L$ acts on the left vector  and $H_R$ acts
on the right vector in (\ref{pairing}). Taking into account
 Proposition \ref{WmrD} one obtains the
following theorem.

\begin{te}\label{teintmfso2l}
The eigenfunctions of $\mathfrak{so}_{2\ell}$-Toda chain
(\ref{pairing}) admit the integral representation:
 \bqa\label{waveggd}
\Psi^{D_\ell}_{\lambda_1,\ldots,\lambda_\ell}
(x_{\ell,1},\ldots,x_{\ell,\ell})\,=\, \int\limits_{C}
\bigwedge_{k=1}^{\ell-1}\bigwedge_{i=1}^kdx_{k,i} \wedge
dz_{k,i}\quad \nonumber e^{{\mathcal F}^{D_{\ell}}},\eqa where

\bqa \mathcal{F}^{D_{\ell}}= \imath\lambda_1x_{1,1}-\sum_{n=2}^\ell
\imath\lambda_n\Big(\sum_{i=1}^nx_{n,i}-\\- \nonumber
2\sum_{i=1}^{n-1}z_{n-1,i}+
\sum_{i=1}^{n-1}x_{n-1,i}-2\ln\,(e^{x_{n,1}}+e^{x_{n-1,1}})\Big)-
\\ \nonumber -\sum_{k=1}^{\ell-1}\Big(e^{x_{k+1,k+1}-z_{k,k}}+
\sum_{i=k+1}^{\ell-1}e^{x_{i,k+1}-z_{i,k}}+e^{x_{i+1,k+1}-z_{i,k}}
\Big)-\\ \nonumber -\sum_{n=1}^{\ell-1}\,\,\Big(
e^{z_{n,1}-x_{n,1}}+e^{z_{n,1}-x_{n+1,1}}+
e^{z_{n,1}+x_{n,1}}+e^{z_{n,1}+x_{n+1,1}}\Big)\,-\\- \nonumber
\sum_{k=3}^{\ell}\,\,\sum_{n=1}^{\ell+1-k}
\Big(e^{z_{k+n-2,k-1}-x_{k+n-2,k-1}}+e^{z_{k+n-2,k-1}-x_{k+n-1,k-1}}\Big)
,\eqa where  $x_i:=x_{\ell,i},\,\,\, 1\le i\leq \ell$ and
${C} \subset N_+$ is a middle-dimensional non-compact submanifold
such that the integrand decreases  exponentially at the boundaries
and at infinities. In particular one can take ${C}=\RR^{m}$,
$m=l(w_0)$ as a domain of integration.

\end{te}

\begin{ex} For $\ell=2$ the general expression (\ref{waveggd}) acquires
the following form
\bqa\Psi^{D_2}_{\lambda_1,\,\lambda_2}(x_{21},x_{22})=\int_Cd
x_{11}dz_{11}\,\,
\Big(\,e^{-x_{22}}+e^{-x_{11}}\,\Big)^{2\imath\lambda_2}\times \\
\nonumber \times \exp\Big\{\imath\lambda_2(x_{21}+x_{22}-2z_{11}+x_{11})-
\imath\lambda_1x_{11}\Big\}\times \\ \nonumber
\times \exp\Big\{-\Big(e^{z_{11}-x_{21}}+e^{x_{22}-z_{11}}+e^{-x_{22}-z_{11}}
+e^{x_{11}-z_{11}}+e^{-x_{11}-z_{11}}\Big)\Big\}. \eqa One can chose
$C=\mathbb{R}^2$ as an integration domain.
\end{ex}

There is a simple combinatorial description of the potential
$\mathcal{F}^{D_{\ell}}$ for zero spectrum $\{\la_i=0\}$.
 Namely, it can be presented as a sum over
arrows in the following diagram.

$$
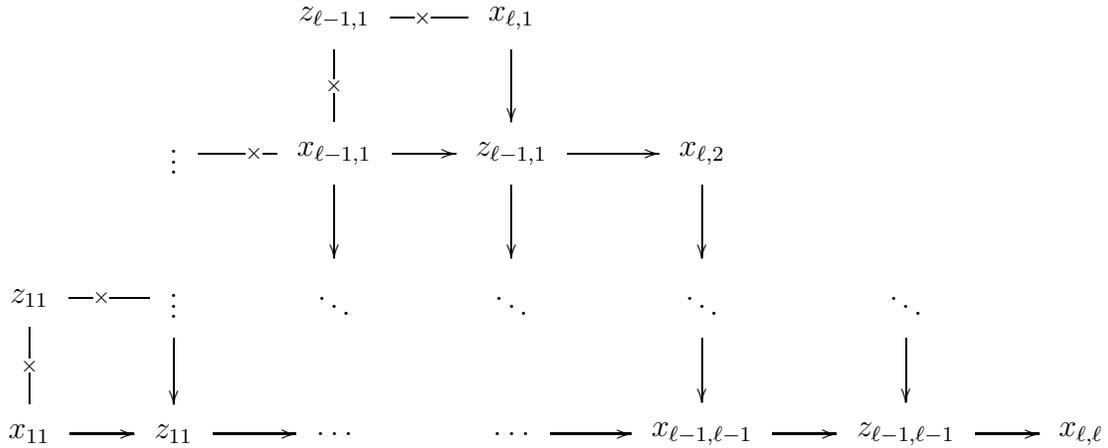
\xymatrix{
 && z_{\ell-1,1}\ar@{-}[d]|{\times}\ar@{-}[r]|{\times} &
 x_{\ell,1}\ar[d] &&&\\
 & \vdots\ar@{-}[r]|{\times} & x_{\ell-1,1}\ar[d]\ar[r] &
 z_{\ell-1,1}\ar[d]\ar[r] & x_{\ell,2}\ar[d] &&\\
 z_{11}\ar@{-}[d]|{\times}\ar@{-}[r]|{\times} & \vdots\ar[d] &
 \ddots & \ddots & \ddots\ar[d] & \ddots\ar[d] &\\
 x_{11}\ar[r] & z_{11}\ar[r] & \ldots & \ldots\ar[r] &
 x_{\ell-1,\ell-1}\ar[r] & z_{\ell-1,\ell-1}\ar[r] & x_{\ell,\ell}
}
$$
\hspace{1cm}

Note that the diagram for $D_{\ell}$  can be obtained  by a
factorization of the diagram for $A_{2\ell-1}$. Consider the
following involution: \be \label{inv} \iota\,:\quad X\longmapsto
\dot{w}_0^{-1}X^t \dot{w}_0, \ee where $\dot{w}_0$ is a lift the
longest element of $A_{2\ell-1}$ Weyl group and $X^t$ denotes the
standard transposition. Corresponding action on the modified
factorization parameters is given by
 \be \label{w0inv}
\dot{w}_0\,:\qquad x_{k,i}\longleftrightarrow -x_{k,k+1-i} \ee This
defines a factorization  of $A_{2\ell-1}$-diagram  that produce the
diagram for $D_{\ell}$. Note that diagram for $D_{\ell}$ can be also
obtained by erasing the last row of vertexes and arrows on the right
slope from the diagram for $C_{\ell}$

An analog of the monomial relations (\ref{defrelAn}) is as follows.
Introduce   variables $a_{i,k}$, $b_{i,k}$, $c_{i,k}$, $d_{i,k}$
associated with the arrows of the diagram \be
a_{k,1}=e^{x_{k,1}+z_{k-1,1}},\qquad
a_{k,i}=e^{z_{k-1,i-1}-x_{k,i}},\\
b_{k,1}=e^{x_{k,1}+z_{k,1}},\qquad
b_{k,i}=e^{x_{k,i}-z_{k,i-1}},\\
c_{k,i}=e^{z_{k,i}-x_{k,i}},\qquad d_{k,i}=e^{z_{k,i}-x_{k+1,i}}.\ee
Then the following relations hold \bqa
b_{k,i}c_{k,i}=a_{k+1,i}d_{k,i},\qquad
a_{k+1,i+1}d_{k,i}=b_{k+1,i+1}c_{k+1,i},\\
a_{\ell,1}a_{\ell,2}=e^{x_{\ell,1}+x_{\ell,2}},\qquad\qquad
a_{\ell,i}d_{\ell-1,i-1}=e^{x_{\ell,i}-x_{\ell,i-1}}.\eqa

\subsubsection{Recursion for $\mathfrak{so}_{2\ell}$-Whittaker
  functions and  $\mathcal{Q}$-operator for $D^{(1)}_{\ell}$-Toda
  chain }

The integral representation (\ref{waveggd}) of
$\mathfrak{so}_{2\ell}$-Whittaker functions possesses a  recursive
structure over the rank $\ell$. For any $n=2,\ldots,\ell$ let us
introduce   integral operators $Q^{D_n}_{D_{n-1}}$ with the kernels
$Q^{D_n}_{D_{n-1}}(\underline{x}_{n};\,\underline{x}_{n-1};\lambda_n)$
defined as follows \be\label{DQD}
Q^{D_n}_{D_{n-1}}(\underline{x}_{n};\,\underline{x}_{n-1};\lambda_n)=
\int \bigwedge_{i=1}^{n-1} dz_{n,i}\,\,
\Big(e^{x_{n-1,1}}+e^{x_{n,1}}\Big)^{2\imath\lambda_n}\times \\
\nonumber \times \exp\Big\{\,-\imath\lambda_n\Big(\sum_{i=1}^nx_{n,i}-
2\sum_{i=1}^{n-1}z_{n-1,i}+\sum_{i=1}^{n-1}x_{n-1,i}\Big)\Big\}
\,\,\times  \\ \nonumber
\times Q^{D_n}_{C_{n-1}}(\underline{x}_{n};\,\underline{z}_{n-1})\,\,\,
Q^{C_{n-1}}_{D_{n-1}}(\underline{z}_{n-1};\,\underline{x}_{n-1}),\quad\ee
where  \be Q^{D_n}_{\,\,\,\,C_{n-1}}(\underline{x}_n
;\,\underline{z}_{n-1})=
\exp\Big\{-\Big(e^{x_{n,1}+z_{n-1,1}}+\sum_{i=1}^{n-1}\Big(e^{z_{n-1,i}-x_{n,i}}+
e^{x_{n,i+1}-z_{n-1,i}}\Big)\,\Big)\Big\},\ee \be
Q^{C_{n-1}}_{\,\,\,\,D_{n-1}}(\underline{z}_{n-1}
;\,\underline{x}_{n-1})=\exp\Big\{-\Big(e^{x_{n-1,1}+z_{n-1,1}}+\\
\nonumber+\sum_{i=1}^{n-2}\Big(e^{z_{n-1,i}-x_{n-1,i}}+
e^{x_{n-1,i+1}-z_{n-1,i}}\Big)+e^{z_{n-1,n-1}-x_{n-1,n-1}}\,\Big)\Big\},
\ee and for $n=1$ we define
$$Q^{D_1}_{D_{0}}(x_{1,1};\lambda_1)=e^{\imath\lambda_1x_{1,1}}.$$
Using  $Q^{D_n}_{D_{n-1}}$, $n=1,\ldots,\ell$ the integral
representation (\ref{waveggd}) can be written   in the recursive
form.

\begin{te} The eigenfunction for $\mathfrak{so}_{2\ell}$-Toda chain can be written
in the following recursive form: \bqa\label{waveiterd}
\Psi^{D_\ell}_{\lambda_1,\ldots,\lambda_\ell}(x_1,\ldots,x_\ell)\,=\,
\int\limits_\mathcal{C}
\bigwedge_{k=1}^{\ell-1}\bigwedge_{i=1}^kdx_{k,i}\prod_{k=1}^{\ell}
Q^{D_{k}}_{\,\,D_k-1}(\underline{x}_{k};\,\underline{x}_{k-1};\lambda_{k}),\eqa
or equivalently \bqa\label{iterwaved}
\Psi^{D_\ell}_{\lambda_1,\ldots,\lambda_\ell}(x_{\ell,1},\ldots,x_{\ell,\ell})\,=\,
\int\limits_{\mathcal{C}_{\ell}}
\bigwedge_{i=1}^{\ell-1}dx_{\ell-1,i} \,
Q^{D_{\ell}}_{\,\,D_{\ell-1}}(\underline{x}_{\ell};\,\underline{x}_{\ell-1};
\lambda_{\ell})
\Psi^{D_{\ell-1}}_{\lambda_1,\ldots,\lambda_{\ell-1}}(x_{\ell-1,1},\ldots,
x_{\ell-1,\ell-1}),\nonumber \eqa where we assume
$x_n:=x_{\ell,n},\,\,\,  1\leq n\leq \ell $. Here $\mathcal{C}
\subset N_+$ is a middle-dimensional non-compact submanifold such
that the integrand decreases  exponentially at the boundaries and at
infinities. In particular the domain of integration can be chosen to
be $\cal{C}=$$\RR^{m},$ where $m=l(w_0).$
\end{te}

As for other classical Lie algebras, different from
$\mathfrak{gl}_{\ell+1}$, the specailization to zero spectrum
$\{\la_n=0\}$  revieals a more refined recursive structure. In this case the kernel
of the operator  $Q^{D_n}_{D_{n-1}}$ is reduced to a convolution of
two kernels
$Q^{D_{n}}_{\,\,\,\,C_{n-1}}(\underline{x}_{n};\,\underline{z}_{n-1})$
 and
 $Q^{C_{n-1}}_{\,\,\,\,D_{n-1}}(\underline{z}_{n-1};\,\underline{x}_{n-1})$.
The corresponding integral operators  $Q^{D_n}_{C_{n-1}}$,
$Q^{C_{n-1}}_{D_{n-1}}$ can be regarded as elementary intertwiners
relating Toda chains for $D_n$, $C_{n-1}$ and $C_{n-1}$, $D_{n-1}$
root systems. Thus for quadratic Hamiltonians one can directly check
the following relations

\begin{lem} The operators $Q^{D_n}_{D_{n-1}}$, $Q^{D_n}_{C_{n-1}}$ and
$Q^{C_{n-1}}_{D_{n-1}}$ satisfy the following intertwining relations
with quadratic Toda Hamiltonians.
\begin{enumerate}
\item Operators $Q^{D_n}_{C_{n-1}}$ and $Q^{C_{n-1}}_{D_{n-1}}$ intertwine quadratic
Hamiltonians of $C$- and $D$-Toda chains:
\bqa\CH_2^{D_{n}}(\underline{x}_n)Q^{D_{n}}_{C_{n-1}}(\underline{x}_n,\,
\underline{z}_{n-1})= Q^{D_n}_{C_{n-1}}
(\underline{x}_n,\,\underline{z}_{n-1})\CH_2^{C_{n-1}}(\underline{x}_{n-1}),
\eqa
\bqa\CH_2^{C_{n}}(\underline{z}_n)Q^{C_{n}}_{D_{n}}(\underline{z}_n,\,
\underline{x}_{n})= Q^{C_n}_{D_{n}}
(\underline{z}_n,\,\underline{x}_{n})\CH_2^{D_{n}}(\underline{x}_{n}).
\eqa

\item Operator $Q^{D_{n}}_{D_{n-1}}$ for $\la_n=0$ intertwines  Hamiltonians ${\CH}_2^{D_n}$
and ${\CH}_2^{D_{n-1}}:$
\bqa\CH_2^{D_{n}}(\underline{x}_n)Q^{D_{n}}_{D_{n-1}}(\underline{x}_n;\,
\underline{x}_{n-1};\la_n=0)= Q^{D_n}_{D_{n-1}}
(\underline{x}_n;\,\underline{x}_{n-1};\la_n=0)\CH_2^{D_{n-1}}(\underline{x}_{n-1}),
\eqa where
 \bqa \CH_2^{C_{n}}=-
\frac{1}{2}\sum_{i=1}^{n}\frac{\partial^2}{\partial z_i^2}+
2e^{2z_1}+\sum_{i=1}^{n-2}e^{z_{i+1}-z_i}, \eqa
 \bqa \CH_2^{D_{n}}=-
\frac{1}{2}\sum_{i=1}^{n}\frac{\partial^2}{\partial x_i^2}+
e^{x_1+x_2}+\sum_{i=1}^{n-1}e^{x_{i+1}-x_i}. \eqa

\end{enumerate}\end{lem}

The integral  kernel of the operator
$Q^{D_{n}}_{\,\,\,\,D_{n-1}}(\underline{x}_{n};\,\underline{x}_{n-1})$
at $\la_n=0$ can be succinctly encoded into the following diagram

\bqa
\label{Dsubdiag}
\xymatrix{
 z_{n-1,1}\ar@{-}[d]|{\times}\ar@{-}[r]|{\times} &
 x_{n,1}\ar[d] &&&&\\
 x_{n-1,1}\ar[r] & z_{n-1,1}\ar[d]\ar[r] & x_{n,2}\ar[d] &&&\\
 & x_{n-1,2}\ar[r] & z_{n-1,2}\ar[d]\ar[r] & \ddots &&\\
 && \ddots & \ddots\ar[d]\ar[r] & x_{n,n-1}\ar[d] &\\
 &&& x_{n-1,n-1}\ar[r] & z_{n-1,n-1}\ar[r] & x_{n,n}
}\eqa

 Here the upper and lower descending paths of the oriented diagram
 correspond
  to the kernels of elementary intertwiners
$Q^{D_n}_{C_{n-1}}$ and $Q^{C_n}_{D_{n-1}}$ respectively. The
convolution of the kernels $Q_{C_n}^{D_n}$ and
$Q_{D_{n-1}}^{C_{n-1}}$ corresponds to the integration over the
variables $z_{n-1,i}$ associated with the inner vertexes of the
sub-diagram (\ref{Dsubdiag}.)

Similarly to the cases of other classical series of Lie algebras
recursion operators $Q^{D_n}_{D_{n-1}}$ can be considered as
degenerations of  Baxter $\mathcal{Q}$-operators for affine
$D^{(1)}_{\ell}$-Toda chains. Let us recall the root data for
$D^{(1)}_{\ell}$.   Simple roots of  the affine root system
$D^{(1)}_{\ell}$ can be represented as vectors in
$\mathbb{R}^{\ell}$ in the following way
 \be \alpha_1=\epsilon_1+\epsilon_2,\,\,\,\,\,\,\,\,\,
\alpha_{i+1}=\epsilon_{i+1}-\epsilon_i,\,\,\, \,1\leq i\leq
\ell-1,\\ \nonumber
\alpha_{\ell+1}=-\epsilon_\ell-\epsilon_{\ell-1},
 \ee
and Dynkin diagram is given  by

$$
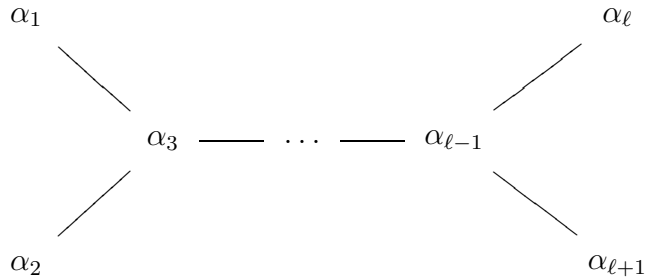
\xymatrix{
 \alpha_1\ar@{-}[dr] &&&& \alpha_{\ell}\\
 & \alpha_3\ar@{-}[r] & \ldots\ar@{-}[r] &
 \alpha_{\ell-1}\ar@{-}[ur]\ar@{-}[dr] &\\
 \alpha_2\ar@{-}[ur] &&&& \alpha_{\ell+1}
}
$$

The corresponding
 $D_{\ell}^{(1)}$-Toda chain quadratic  Hamiltonian is defined by
\bqa\CH_2^{D^{(1)}_{\ell}}&=&-\frac{1}{2}\sum_{i=1}^{\ell}
\frac{\partial^2}{\partial x_i^2}+e^{x_1+x_2}+
\sum_{i=1}^{\ell-2}e^{x_{i+1}-x_i}+ge^{x_{\ell}-x_{\ell-1}}+
ge^{-x_{\ell}-x_{\ell-1}}.\eqa Define the Baxter
$\mathcal{Q}$-operator of $D^{(1)}_{\ell}$-Toda chain as an integral
operator with the following integral kernel
 \bqa\label{QDDEF}\mathcal{Q}^{^{D_{\ell}^{(1)}}}(
\underline{x}^{(\ell)},\underline{y}^{(\ell)},\lambda)
=\int\bigwedge_{i=1}^{\ell-1}dz_{i}\,\,
\Big(e^{x_1}+e^{y_1}\Big)^{2\imath\lambda}
\Big(e^{-x_\ell}+e^{-y_\ell}\Big)^{-2\imath\lambda}\times
\nonumber\\ \times \exp\Big\{-\imath\lambda\Big(\sum_{i=1}^\ell x_i-
2\sum_{i=1}^{\ell-1}z_i+\sum_{i=1}^\ell y_i\Big) \Big\}\,\times \\
\nonumber \times Q_{C^{(1)}_{\ell-1}}^{\,\,\,\,D^{(1)}_{\ell}}
(x_1,\ldots,x_\ell\,;z_1,\ldots,z_{\ell-1})\,\,\,\,\,
Q_{\,\,\,\,D^{(1)}_{\ell}}^{C^{(1)}_{\ell-1}}
(z_1,\ldots,z_{\ell-1}\,;y_1,\ldots,y_\ell) ,\eqa where
\bqa\label{inttwCn} Q_{C^{(1)}_{\ell-1}}^{\,\,\,\,D^{(1)}_{\ell}}
(x_1,\ldots,x_\ell\,;\,z_1,\ldots,z_{\ell-1})=\\ \nonumber
=\exp\Big\{ e^{z_1+x_1}+\sum_{i=1}^{\ell-1}\Big(
e^{z_i-x_i}+e^{x_{i+1}-z_i} \Big)+ge^{-x_{\ell}-z_{\ell-1}}\Big\},
\eqa and \be Q^{D^{(1)}_{\ell}}
_{\,\,\,\,C^{(1)}_{\ell-1}}(x_1,\ldots ,x_{\ell};\,z_1,\ldots
,z_{\ell})=Q_{\,\,\,\,D^{(1)}_{\ell}}^{C^{(1)}_{\ell-1}}(z_1,\ldots
,z_{\ell};x_1,\ldots ,x_{\ell}).\ee
Here we use the following notations
$\underline{x}^{(\ell)}=(x_1,\ldots ,x_{\ell})$, 
$\underline{y}^{(\ell)}=(y_1,\ldots ,y_{\ell})$.

\begin{prop}
The $\mathcal{Q}$-operator (\ref{QDDEF})  commutes with quadratic
Hamiltonian of the $D_\ell^{(1)}$-Toda chain:
\bqa\CH^{D_{\ell}^{(1)}}(\underline{x}^{(\ell)})\mathcal{Q}^{^{D_{\ell}^{(1)}}}
(\underline{x}^{(\ell)},\underline{y}^{(\ell)})= \mathcal{Q}^{^{D_{\ell}^{(1)}}}
(\underline{x}^{(\ell)},\,\underline{y}^{(\ell)})
\CH^{^{D_{\ell}^{(1)}}}(\underline{y}^{(\ell)}). \eqa
\end{prop}
Now we will demonstrate that recursion operator
$Q^{D_{\ell}}_{D_{\ell-1}}$ can be considered as  a degeneration of
Baxter $\mathcal{Q}$-operators for $D_{\ell}^{(1)}$. Let us
introduce a slightly modified recursion  operator with the kernel:
$Q^{D_{\ell}}_{D_{\ell-1}\oplus D_{1}}$: \bqa
Q^{D_{\ell}}_{D_{\ell-1}\oplus D_{1}}
(\underline{x}^{(\ell)}\,,\,\underline{y}^{(\ell)},\,\lambda)
=e^{\imath\lambda y_{\ell}} Q^{D_{\ell}}_{D_{\ell-1}}
(\underline{x}^{(\ell)}\,,\,\underline{y}^{(\ell-1)}\,,\,\lambda),\eqa
where we use the notations $\underline{y}^{(\ell-1)}=(y_1,\ldots ,y_{\ell-1})$.
This operator intertwines Hamiltonians of $\mathfrak{so}_{2\ell}$-
and   $\mathfrak{so}_{2\ell-2}\oplus \mathfrak{so}_2$-Toda chains.
 For instance we have   for quadratic Hamiltonians
\bqa\nonumber \CH_2^{D_{\ell}}(\underline{x}^{(\ell)})
Q^{D_{\ell}}_{D_{\ell-1}\oplus
D_{1}}(\underline{x}^{(\ell)}\,,\,\underline{y}^{(\ell)},\la)
= Q^{D_{\ell}}_{D_{\ell-1}\oplus D_{1}}
(\underline{x}^{(\ell)}\,,\,\underline{y}^{(\ell)},\la)
\Big(\CH_2^{D_{\ell-1}}(\underline{y}^{(\ell-1)})+
\CH_2^{D_1}(y_\ell)\Big), \eqa where
$\CH_2^{D_1}(y_{\ell})=-\frac{1}{2}\Big(\partial^2/\partial
y_{\ell}^2\Big).$ Obviously the projection of the  above relation on the
subspace of functions
$F(\underline{y}^{(\ell)})=\exp(\imath\lambda y_{\ell})
f(\underline{y}^{(\ell-1)})$ recovers the genuine recursion operator
satisfying:
\bqa\CH_2^{D_{\ell}}(\underline{x}^{(\ell)}) Q^{D_{\ell}}_{D_{\ell-1}}
(\underline{x}^{(\ell)}\,,\,\underline{y}^{(\ell-1)}\,,\,\lambda)=
Q^{D_{\ell}}_{D_{\ell-1} }(\underline{x}^{(\ell)}\,,
\,\underline{y}^{(\ell-1)}\,,\,\lambda)\Big(\CH_2^{D_{\ell-1}}
(\underline{y}^{(\ell-1)})+\frac{1}{2}\lambda^{2}\Big). \eqa Let us
introduce a one-parameter family of the operators with the kernels \bqa
\mathcal{Q}^{^{D_{\ell}^{(1)}}}(\underline{x}^{(\ell)},\,
\underline{y}^{(\ell)},\lambda;\varepsilon):=
\varepsilon^{\imath\lambda}e^{\imath\lambda y_{\ell}}
\int\bigwedge_{i=1}^{\ell-1}dz_{i}\,\,
\Big(e^{x_1}+e^{y_1}\Big)^{2\imath\lambda}
\Big(\varepsilon e^{-x_\ell+y_{\ell}}+1\Big)^{-2\imath\lambda}\times \\
\nonumber\times  \exp\Big\{-\imath\lambda\Big(\sum_{i=1}^\ell x_i-
2\sum_{i=1}^{\ell-1}z_i+\sum_{i=1}^{\ell-1} y_i\Big) \Big\}\,\times \\
\nonumber \times Q_{C^{(1)}_{\ell-1}}^{\,\,\,\,D^{(1)}_{\ell}}
(x_1,\ldots,x_\ell\,;z_1,\ldots,z_{\ell-1})\,\,\,\,\,
Q_{\,\,\,\,D^{(1)}_{\ell}}^{C^{(1)}_{\ell-1}}
(z_1,\ldots,z_{\ell-1}\,;y_1,\ldots,y_\ell;\varepsilon) ,\eqa where
\bqa\label{inttwCn} Q_{C^{(1)}_{\ell-1}}^{\,\,\,\,D^{(1)}_{\ell}}
(x_1,\ldots,x_\ell\,;\,z_1,\ldots,z_{\ell-1})=\\ \nonumber
=\exp\Big\{ e^{z_1+x_1}+\sum_{i=1}^{\ell-1}\Big(
e^{z_i-x_i}+e^{x_{i+1}-z_i} \Big)+ge^{-x_{\ell}-z_{\ell-1}}\Big\},
\eqa and \bqa
 Q^{C^{(1)}_{\ell-1}}_{\,\,\,\,D^{(1)}_{\ell}}
(z_1,\ldots,z_{\ell-1}\,,\,y_1,\ldots,y_\ell;\varepsilon)=
\nonumber\\= \exp\Big\{-\Big( e^{z_1+y_1}+\sum_{i=1}^{\ell-2}\Big(
e^{z_i-y_i}+e^{y_{i+1}-z_i}+e^{z_{\ell-1}-y_{\ell-1}}\Big)+\\
+\varepsilon e^{y_{\ell}-z_{\ell-1}}
+\varepsilon^{-1}ge^{-x_{\ell}-z_{\ell-1}}\Big)\Big\}.\nonumber \eqa
These operators are obtained  by a shift of the variable
$y_{\ell}=y_{\ell}+\ln\varepsilon$ in (\ref{DQD}).
 Then the following relation
between $\mathcal{Q}$-operator for $D_{\ell}^{(1)}$-Toda chain and
(modified) recursion operator for $\mathfrak{so}_{2\ell}$-Whittaker
function holds \bqa Q^{D_{\ell}}_{D_{\ell-1}\oplus
D_{1}}(\underline{x}^{(\ell)}\,,\,\,
\underline{y}^{(\ell)},\,\lambda)
=\lim_{\varepsilon\rightarrow 0  \,\,\varepsilon^{-1}g\rightarrow
0}\varepsilon^{-\imath\lambda}
{Q}^{^{D_{\ell}^{(1)}}}(\underline{x}^{(\ell)}\,,\,
\underline{y}^{(\ell)}\,,\,\lambda\,;\,\varepsilon).\eqa

\newpage
\section{Part II. Proofs}

Let $G$ be a complex connected simply-connected semisimple Lie group
of finite rank $\ell$, $\mathfrak{g}={\rm Lie}(G)$ be  the
corresponding semisimple Lie algebra with the Chevalley generators
$f_i,\,h_i,\,e_i$. Let us fix a Borel subgroup $B_+$
and let $T$ be the maximal torus  $T\subset B_+$. This defines
a pair $N_+,\,N_-$  of opposite unipotent subgroups in $G$,
$N_+\subset B_+$. Let $\Gamma$ be the set of vertices of Dynkin graph of
$\mathfrak{g}$, $\{\alpha_i,\,i\in\Gamma\}$ be the set of simple roots,
$\{\gamma_k, k=1,\ldots,\frac{1}{2} (\dim\,\mathfrak{g}/\mathfrak{h})\}$
 be the set of all positive roots
 and $\{\alpha_i^\vee,\,i\in\Gamma\}$ be the set of
simple coroots. For every  $i\in\Gamma$  there is a group
homomorphism
\bqa\label{oneparsub}
\varphi_i\,:\quad SL_2\longrightarrow G,
\eqa
defined as follows. Introduce a set of one-parameter subgroups
$e^{te_i}=X_i(t)\subset N_+$, $e^{tf_i}=Y_i(t)\subset N_-$ and
$e^{th_i}=\alpha_i^\vee(t)\subset T$. Homomorphisms (\ref{oneparsub})
are defined as \bqa\label{homomorone} \varphi_i(e^{te})=e^{te_i},\qquad
\varphi_i(e^{tf})=e^{tf_i},\qquad
\varphi_i(e^{th})=\alpha_i^\vee(t), \eqa where $e,f,h$ are standard
generators of $\mathfrak{sl}_2$. Let us fix the
lifts $\dot{s}_i\subset G$, $\dot{s}\subset SL(2)$ of the generators $s_i$
of the Weyl group of $G$ and the generator of the Weyl group of
$SL(2)$
\bqa \label{liftWeyl}
\dot{s}=e^{e}e^{-f}e^{e},\qquad  \dot{s_i}=e^{e_i}e^{-f_i}e^{e_i}.\eqa
Thus defined lifts of  Weyl group generators
are obviously compatible $\varphi_i(\dot{s})=\dot{s}_i$ 
with  homomorphisms (\ref{oneparsub}).
 We have the following relations
\bqa\label{sl2conj} \dot{s}^{-1}\,f\,\dot{s}\,=\,-e,\qquad \quad
\dot{s}_i^{-1}\,f_i\,\dot{s}_i\,=\,-e_i. \eqa

The action $w_0(\alpha_i)=-\alpha_{i^*}$
of the maximal length element $w_0$ of the Weyl group
on  simple roots defines  an involution
$i\mapsto i^*$. The corresponding action of $\dot{w}_0$ 
is given by 
\bqa\label{Star}\dot{w}_0^{-1}\,f_i\,\dot{w}_0\,=\,-e_{i^*}.\eqa

\begin{rem}
 For classical Lie groups one has
  $i^*=\ell+1-i$ for $G=SL(\ell+1)$, $i^*=i$ for
$G=SO(2\ell+1)$ and for $G=Sp(2\ell)$. In the case $G=SO(2\ell)$ (for
  $\ell\geq2$) the action   of the involution $*$ is as follows:
\bqa*\,:\qquad 1\longmapsto \left\{\begin{array}{lcr}
1\,, & & \ell\,\,\mbox{even}\\
2\,, & & \ell\,\,\mbox{odd}
\end{array}\right.
\hspace{1.5cm} 2\longmapsto \left\{\begin{array}{lcr}
2\,, & & \ell\,\,\mbox{even}\\
1\,, & & \ell\,\,\mbox{odd}
\end{array}\right.\eqa
\bqa \nonumber k^*=k\,,\qquad2<k\leq\ell, \eqa where the
enumeration of roots of $SO(2\ell)$ is given by (\ref{RootDn}).
\end{rem}

In the following we will be considering  matrix elements of
finite-dimensional representations $V_{\omega_i}$ of $\mathfrak{g}$
corresponding to the fundamental weights $\omega_i$, $i\in \Gamma$.
Let $\xi_{\omega_i}^+$and $\xi_{\omega_i}^-$ be  highest and lowest
vectors in $V_{\omega_i}$ such that $\langle \xi_{\omega_i}^-|
\xi_{\omega_i}^+\rangle=1$. For the lift (\ref{liftWeyl}) of the
elements of the Weyl group we have (see e.g. \cite{K} Lemma 3.8,
\cite{FZ} eq. (2.29))  \bqa\label{actionof lift}
\dot{w}_0^{-1}\xi_{\omega_i}^+ =\xi_{\omega_i}^-,\qquad \quad
\dot{s}_i^{-1}\xi_{\omega_i}^+ =f_i\xi^+_{\omega_i}.\eqa

Consider the following parametrization of a generic
group element $g\in G$
\bqa\label{GaussParam}
g=g^{(-)}\,g^{(0)}\,g^{(+)}=\exp(\sum_{\alpha\in
  \Delta_+} u_{-\alpha}f_{\alpha})\,\,\,
\exp(\sum_{i=1}^{\ell}u_ih_i)\,\,\,\exp(\sum_{\alpha\in
  \Delta_+} u_{\alpha}e_{\alpha}).\eqa
For coordinates $u_i$ corresponding to Cartan generators $h_i$
and for coordinates $u_{\pm\alpha_i}$ corresponding to simple root
generators  $e_{\alpha_i}$, $f_{\alpha_i}$  there exits   simple
expressions in terms of matrix elements of  fundamental
representations $V_{\omega_i}$:
\bqa\label{minorsparam}
u_{\alpha_i}(g)=\frac{\langle\xi_{\omega_i}^-|\pi_i(g)\pi_i(f_i)|\xi_{\omega_i}^+\rangle}
{\langle\xi_{\omega_i}^-|\pi_i(g)|\xi_{\omega_i}^+\rangle},\hspace{8 mm}
u_{-\alpha_i}(g)=\frac{\langle\xi_{\omega_i}^-|\pi_i(g)\pi_i(e_i)|
\xi_{\omega_i}^+\rangle}
{\langle\xi_{\omega_i}^-|\pi_i(g)|\xi_{\omega_i}^+\rangle},\eqa
\bqa
 u_i(g)=\langle\xi_{\omega_i}^-|\pi_i(g)|\xi_{\omega_i}^+\rangle,\nonumber
 \eqa
where $\pi_i\equiv \pi_{\omega_i}$ is a fundamental representation in $V_{\omega_i}$. 
Define generalized twisted minors as 
\bqa\label{twistedmin}
\Delta_{\omega_i,\dot{w}}(g)=\<\xi_{\omega_i}^-|\pi_{\omega_i}(g)
\pi_{\omega_i}(\dot{w})|\xi_{\omega_i}^+\>,\qquad g\in G.
\eqa
Then coordinate $u_i$ and $u_{\alpha_i}$ of
a twisted unipotent element 
$v\dot{w}_0^{-1}\in G$ (where $v\in N_+$) can be expressed in terms of twisted minors
(\ref{twistedmin}) as follows
\bqa
e^{u_i(v\dot{w}_0^{-1})}=\Delta_{\omega_i,\dot{w}_0^{-1}}(v),\nonumber
\eqa
\bqa
u_{\alpha_i}(v\dot{w}_0^{-1})=
\frac{\langle\xi_i^{-}|\pi_i(v\dot{w}_0^{-1})\pi_i(f_i)|\xi_i^+\rangle}
{\langle\xi_i^{-}|\pi_i(v\dot{w}_0^{-1})|\xi_i^+\rangle}=
\frac{\langle\xi_i^{-}|
\pi_i(v)\pi_i(\dot{w}_0^{-1})\pi_i(\dot{s}_i^{-1})|
\xi_{\omega_i}^+\rangle}
{\langle\xi_i^{-}|\pi_i(v\dot{w}_0^{-1})|\xi_{\omega_i}^+\rangle}=\\
= -\frac{\langle\xi_{\omega_i}^{-}|
\pi_i(v)\pi_i(e_{i^*})\pi_i(\dot{w}_0^{-1})|
\xi_{\omega_i}^+\rangle}
{\langle\xi_i^{-}|\pi_i(v\dot{w}_0^{-1})|\xi_i^+\rangle}=
\frac{\Delta_{\omega_i,\dot{w}_0^{-1}\dot{s}_i^{-1}}(v)}{\Delta_
{\omega_i,\dot{w}_0^{-1}}(v)}.\nonumber\eqa
In the following we will use the shorthand notations
\bqa\label{primes}\Delta_{i}'(v):&=&\langle\xi_{\omega_i}^-|\pi_i(v
e_{i*}\dot{w}_0^{-1})\xi_{\omega_i}^+\rangle=-
\Delta_{\omega_i,\dot{w}_0^{-1}\dot{s}_i^{-1}}(v)\,\,\,, \\
\Delta_i(v):&=&\Delta_{\omega_i,\dot{w}_0^{-1}}(v).\nonumber\eqa

\subsection{ Measure on $N_+$: Proof of Lemma \ref{lemmeasure} }
\label{measureN}

In this part we derive an explicit expression 
(\ref{meas}) for a measure $d\mu_{N_+}(x)$ on a unipotent subgroup $N_+\subset G$
of any classical Lie group  using a factorized
 parametrization  (\ref{factparam})  of $N_+$.
Recall that for a reduced word $I_\ell=(i_1,\ldots,i_{m_\ell})$ of
$w_0$ there is a birational isomorphism
$\mathbb{C}^{m_\ell}\rightarrow N_+$. Particularly, given an
unipotent element $v\in N_+$ the following factorized representation holds.
\bqa\label{FactorParametr}
v(t)=X_{i_1}(t_1)X_{i_2}(t_2)\cdot\ldots\cdot X_{i_m}(t_{m_\ell}),\eqa
where $X_{i}(t)=e^{te_i}$. The variables $t_i$ are called factorization parameters of $v$.

\begin{prop}\label{measureprop}
Let $v(t)\in N^{(0)}_+$ be a factorized parametrization (\ref{FactorParametr})
corresponding to a reduced word $I=(i_1,\ldots,i_{m_\ell})$. Then
\bqa\label{ing} d\mu_{N_+}(v(t))\,=\,\prod_{k=1}^\ell
\prod_{i=1}^{m_\ell} (t_i)^{\langle \omega_k,\gamma_i\rangle}
\cdot\bigwedge_{i=1}^{m_\ell}\frac{dt_i}{t_i},\eqa is a restriction
of the right-invariant measure $d\mu_{N_+}$  to $N^{(0)}_+$, that is
\be\label{invar} \,\,\,\, d\mu_{N_+} (v(t))=d\mu_{N_+} (v(t)\cdot
X_j(\tau)), \,\,\,\,j=1,\ldots,\ell. \ee
\end{prop}

{\it Proof.} To prove the Proposition consider a dependence on a
choice of a reduced word $I=(i_1,\ldots,i_m)$ explicitly. Let
$t^I=(t^I_1,\ldots,t^I_{m_\ell})$ be factorization parameters
corresponding to a reduced word $I$. 
According to \cite{BZ} (Theorem 4.3)
one has the following expressions for matrix elements  \bqa
\Delta_k(t^I):=\Delta_k(x(t^I)w_0^{-1})\,=\,\prod_{i=1}^{m_\ell}
(t^{I}_i)^{\langle \omega_k,\gamma_i\rangle} \eqa
 Two parameterizations
$x(t^{I})$ and $x(t^{I'})$ of $N_+^{(0)}$ corresponding to reduced
words $I$ and $I'$ are related by a birational transformation.
\begin{lem} For any reduced decompositions of $w_0$
corresponding to  reduced words $I$ and $I'$ the following relations
hold
\begin{enumerate}
\item\bqa \Delta_k(t^I)=\Delta_k(t^{I'}),\qquad 1\leqq
k\leqq\ell.\label{deltaeq}\eqa

\item\bqa \bigwedge_{j=1}^{m_\ell} \,\frac{dt^{I}_j}
{t^{I}_j}\,=\,\bigwedge_{j=1}^{m_\ell}\, \frac{dt^{I'}_j}{t^{I'}_j}.
\label{deltaeq1}\eqa
\end{enumerate}
\end{lem}

\emph{Proof of Lemma. } It is shown in \cite{Lu}  that
 birational transformations $R^{I'}_{I}$ of $N_+$ corresponding to
 any two reduced words $I$ and $I'$ can be represented as a composition of elementary
 transformations (so-called 3- and 4-moves). Therefore  to prove (\ref{deltaeq}),
(\ref{deltaeq1}) one should check these identities for the
elementary moves only  in the following In the case of classical Lie
groups it is enough to consider the following two birational
transformations $R^{I'}_I:t^I\to t^{I'}$
\begin{enumerate}
\item $X_i(t_1)X_j(t_2)X_i(t_3)=X_j(t'_1)X_i(t'_2)X_j(t'_3)$ for
$a_{ij}=a_{ji}-1$,
\item $X_j(t_1)X_i(t_2)X_j(t_3)X_i(t_4)=
X_i(t'_1)X_j(t'_2)X_i(t'_3)X_j(t'_4)$ for $a_{ij}=-1$ and
$a_{ji}=-2$, 
\end{enumerate} where we denote $t=t^I$ and $t'=t^{I'}$.

The proof of the  identity (\ref{deltaeq}) for elementary 3- and
4-moves follows straightforwardly from the results in \cite{BZ}.
Thus we consider only the proof of  (\ref{deltaeq1}) below.

1) In the case $a_{ij}=a_{ji}=-1$ we should consider the birational
transformation between the parametrizations associated with reduced
words $I=(\ldots iji\ldots )$ and $I'=(\ldots jij\ldots )$. We have
the following relation between parameters
$$v=X_i(t_1)X_j(t_2)X_i(t_3)\,=\,X_j(t_1')X_i(t_2')X_j(t_3'),$$
where $$t'_1=\frac{t_2t_3}{t_1+t_3},\qquad t'_2=t_1+t_3,\qquad
t'_3=\frac{t_1t_2}{t_1+t_3}\,.$$ Direct check gives \be d\log
t'_1\wedge d\log t'_2\wedge d\log t'_3\,=\,d\log t_1\wedge d\log
t_2\wedge d\log t_3\,.\ee

2) In the case $a_{ij}=-1,\,a_{ji}=-2$ we should consider the
birational transformation between the parameterizations associated
with reduced words $I=(\ldots jiji\ldots )$ and $I'=(\ldots
ijij\ldots)$. Thus we have the following relation between parameters
$$X=X_j(t_1)X_i(t_2)X_j(t_3)X_i(t_4)\,=\,
X_i(t'_1)X_j(t'_2)X_i(t'_3)X_j(t'_4), $$ with
 \bqa t'_1=\frac{t_2t_3^2t_4}{t_1^2t_2+(t_1+t_3)^2t_4}\,,\qquad
t'_2=\frac{t_1^2t_2+(t_1+t_3)^2t_4}{t_1t_2+(t_1+t_3)t_4}\,,\\
\nonumber t'_3=\frac{\Big(t_1t_2+(t_1+t_3)t_4\Big)^2}
{t_1^2t_2+(t_1+t_3)^2t_4}\,,\qquad
t'_4=\frac{t_1t_2t_3}{t_1t_2+(t_1+t_3)t_4}.\eqa

One can readily verify the following identity:  \bqa d\log
t'_1\wedge d\log t'_2\wedge d\log t'_3\wedge d\log t'_4\,=
 \nonumber d\log t_1\wedge d\log t_2\wedge d\log
t_3\wedge d\log t_4.\eqa This completes the proof of the Lemma.

Now we can complete the proof of the Proposition \ref{measureprop}.
To establish the right-invariance of measure $d\mu_{N_+}(v)$ we use
(\ref{deltaeq}), (\ref{deltaeq1}). For any simple root $\alpha_i$
one can find a reduced word $I(\alpha_i)=(j_1,\ldots,j_{m})$ with
$m=m_\ell$ such that $j_m=i$. Then identities (\ref{deltaeq}),
(\ref{deltaeq1}) imply that
$$d\mu_{N_+}(v(t^{I(\alpha_i)}))=d\mu_{N_+}(v(t^{I_\ell})).$$
In this way we obtain \bqa v(t^{I(\alpha_i)})\cdot
X_i(\tau)=X_{j_1}(t_1)\cdot\ldots\cdot
X_{j_{m-1}}(t_{m-1})X_i(t_m+\tau) \eqa
By construction the factorization parameter $t_m$ enters only in the (monomial)
expression for $\Delta_j(v(t))$ as a homogeneous factor of degree
one. In this way, the factorization parameter $t_m$ appears in the
measure $d\mu_{N_+}$ only in the $\alpha_j$-component $\Delta_j(v(t))d\ln t_m$,
and hence, the measure $d\mu_{N_+}$ is invariant under the shift 
$t_m\to t_m+\tau$. 
Thus the measure is right-invariant with respect to the action of
$X_j(\tau)$ for any $j=1,\ldots,\ell$, and eventually it is
right-invariant with respect to the whole $N_+$. 
  This  completes the proof of Lemma \ref{lemmeasure}.

\subsection{Whittaker vectors  for classical Lie groups:\\ Proof of 
Lemma \ref{Whittakervectlem} and Proposition \ref{propwhitgen}}
\label{Witvec}

In this subsection we derive  expressions for left and right
$\mathfrak{g}$-Whittaker vectors in terms of the matrix elements of
finite-dimensional representations $\mathfrak{g}$.
 The Whittaker vectors satisfy the following equations
\be  e_i\psi_R=-\psi_{R}\,,\qquad f_i\psi_L=-\psi_{L}\,,
 \hspace{1cm} i=1,\ldots,\ell.\ee
Integrating actions of the
nilpotent Lie subalgebras $\mathfrak {n}_{\pm}\subset \mathfrak{g}$
 to actions of the nilpotent Lie subgroups $N_{\pm}\subset G$, 
 equations on $\mathfrak{g}$-Whittaker
can be written in terms of one-parameter subgroups
$X_{i}(t)\subset N_+$, $Y_{i}(t)\subset N_-$ as follows
$$\pi_{\lambda}(X_{{i}}(t))\psi_R(v)=e^{-t}\psi_R(v),\,\,\,\,
\pi_{\lambda}(Y_{{i}}(t))\psi_L(v)=e^{-t}\psi_L(v),\,\,\,\,i=1,\ldots,\ell,\,\,\,\,
v\in N_+.$$
Equivalently one has for any $z_{\pm}\in N_{\pm}$
\bqa\label{psi_r}\pi_{\lambda}(z_+)\psi_R(v)=\exp\Big\{-\sum\limits_{i=1}^{\ell}(z_+)_{i}
\Big\}\psi_R(v), \quad
\pi_{\lambda}(z_-)\psi_L(v)=\exp\Big\{-\sum\limits_{i=1}^{\ell}(z_-)_{i}
\Big\}\psi_L(v),\eqa where $(z_{\pm})_i:=u_{\pm\alpha_i}(z_{\pm})$.
Construction of the right Whittaker vector is pretty
straightforward. Note that we have a simple identity
$$
u_{\alpha_i}(v_1\,v_2)=u_{\alpha_i}(v_1)+u_{\alpha_i}(v_2),
\qquad v_1,v_2\in N_+.
$$
Then from (\ref{psi_r}) we infer that the right Whittaker
vector is given by a multiplicative character of the  maximal unipotent
subgroup $N_+$
 \bqa
\psi_R(v)=\exp\Big\{-\sum\limits_{i=1}^{\ell}v_{i}\Big\}=
\exp\Big\{-\frac{\Delta_{\omega_i,\dot{s}_i^{-1}}(v)}{\Delta_{\omega_i,1}(v)}\Big\},
\qquad v\in N_+.\eqa
where $v_i:=u_{\alpha_i}(v)$ 
and we use (\ref{minorsparam}) to express $v_i$ in terms of matrix elements.

To construct the left Whittaker vector in terms of matrix elements  we
use  an inner automorphism of $G$, acting on $z\in G$
 as $z^{\tau}=\dot{w}_0^{-1}z\dot{w}_0$. Taking into account that
$\dot{w}_0^{-1}X_{i^*}(-t)\dot{w}_0=Y_{i}(t)$ we have
$\dot{w}_0^{-1}N_+\dot{w}_0=N_-$. Now the equation for the left
Whittaker
$$\pi_{\lambda}(Y_{{i}}(t))\psi_L(v)=e^{-t}\psi_L(v),\,\,\,\,i=1,\ldots,\ell$$
can be written in the following  form
\bqa\label{psi_l}\pi_{\lambda}(z^{\tau})\psi_{L}(v)=
\exp\Big\{-\sum\limits_{i=1}^{\ell}z_{i} \Big\}\psi_{L}(v),\qquad
z\in N_+,\eqa
The left Whittaker vector can be obtained by the twist of the right
vector 
$$
\psi_{L}(v)=\psi_R(v\dot{w}_0^{-1})
$$
where the function $\psi_R$ is considered as a $B_-$-equivariant
function on $G$ (see (\ref{twistf})  for the precise definition). 
Using Gauss decomposition and the parametrization 
(\ref{GaussParam}), (\ref{minorsparam})
we get for  the left Whittaker vector 
$$\psi_L(v)=e^{\<\imath\lambda-\rho,\sum_{i=1}^{\ell}u_i(v\dot{w}_0^{-1})h_i\>}\,\,\,
e^{\sum_iu_{\alpha_i}(v\dot{w}_0^{-1})}.$$ 
In terms of the matrix elements of finite-dimensional representations 
we have the following representation 
\be \psi_L(v)= \prod_{i=1}^{\ell}\,
\Delta_{\omega_i,\dot{w}_0^{-1}}(v)^{\<\imath\lambda-\rho,\alpha_i^{\vee}\>}
\cdot\exp\Big\{\frac{\Delta_{\omega_i,\dot{w}_0^{-1}\dot{s}_i^{-1}}(v)}{\Delta_
{\omega_i,\dot{w}_0^{-1}}(v)}\Big\}=\\
=\prod_{i=1}^{\ell}\,\Delta_{i}(v)^{\<\imath\lambda-\rho,\alpha_i^{\vee}\>}
\,\,\exp\Big\{-\frac{\Delta'_{i}(v)}{\Delta_{i}(v)}\Big\}. \ee
This completes the proof of Lemma \ref{Whittakervectlem}.
The proof of the Proposition \ref{propwhitgen}
is then obtained by combining the expressions for right Whittaker
vector and  left Whittaker vector twisted 
by the action of Cartan generator $\exp h_x=\exp -(\sum_{i=1}^\ell \<\omega_i,x\>h_i)$.

\subsection{Explicit evaluation of matrix elements}

To construct integral representations of Whittaker functions one
should express various matrix elements entering the integral
formulas (\ref{int})  using  factorized and modified
factorized parametrizations of group elements. This can be done
rather straightforwardly using results of \cite{BZ}, \cite{BF}.
Below we shall use a recursive structure of reduced word $I$
 corresponding to a maximal length element $w_0$ of Weyl group
of classical Lie algebras.  This recursive structure
translates into recursive formulas for the relevant ratios of
matrix elements.  Resolving  recursive equations we find
explicit expressions of $\psi_L $ and $\psi_R$ in a (modified)
factorized parametrization. This provides
corresponding integral representations
for Whittaker functions of classical Lie groups.
In the case of the modified factorized  parametrization
 we obtain a generalization of Givental integral representation
for $\mathfrak{g}=\mathfrak{gl}_{\ell+1}$.

\subsubsection{Expressions for  $\mathfrak{gl}_{\ell+1}$-matrix elements:\\
Proofs of Theorem \ref{teintfgl} and Theorem \ref{teintmfgl}}
\label{melgl}

In this subsection we introduce  expressions for  matrix elements
relevant for the construction of integral representations of
$\mathfrak{gl}_{\ell+1}$-Whittaker functions using 
 factorized parametrization of an open part of
$N_+\subset GL(\ell+1)$. This  provides a proof of the integral
representations of $\mathfrak{gl}_{\ell+1}$-Whittaker  functions
presented in Part I.

The eigenfunctions of $\mathfrak{gl}_{\ell+1}$ and
$\mathfrak{sl}_{\ell+1}$ Toda chains differ by a simple factor
(\ref{psisl}), and the Whittaker vectors $\psi_L,\psi_R$ are the
same for both Lie algebras. Thus we use the $\mathfrak{sl}_{\ell+1}$
root data for calculations of the matrix elements
$\Delta_{\omega_i,\dot{w}_0^{-1}}(v),$
$\Delta_{\omega_i,\dot{w}_0^{-1}\dot{s}_i^{-1}}(v),$
$i=1,\ldots,\ell$ in the fundamental representations of
$\mathfrak{sl}_{\ell+1}$ and set in addition
$\Delta_{\omega_{\ell+1,\dot{w}_0^{-1}}}(v)=1.$
The $\mathfrak{sl}_{\ell+1}$ root data given by (\ref{weightsSL}).
Reduced decomposition $w_0=s_{i_1}s_{i_2}\cdots s_{i_m}$ of the
maximal length element $w_0\in W$  corresponding to a reduced word
$I_{\ell}=(i_1,\ldots,i_m)$ with $m=m_{\ell}=\ell(\ell+1)/2$,
provides a total ordering of positive co-roots   by
$R^{\vee}_+=\{\gamma_k=s_{i_1}\cdots s_{i_{k-1}}\alpha^{\vee}_k\}$
of $\mathfrak{sl}_{\ell+1}$.
   We consider a decomposition of $w_0$ described by the following   reduced word
$$I_\ell=(1,21,\ldots,(\ell\ldots21)).$$
Corresponding ordering of positive co-roots is given by:
\bqa\label{orderCn} \gamma^\vee_1=\alpha_1^\vee,\quad
\begin{array}{l}
\gamma^\vee_2=\alpha_1^\vee+\alpha_2^\vee,\\
\gamma^\vee_3=\alpha_2^\vee,\end{array} \quad\ldots\qquad
\begin{array}{l}
\gamma_{m_{\ell-1}+1}^\vee=\alpha_1^\vee+\ldots+\alpha_\ell^\vee,\\
\vdots\\
\gamma^\vee_{m_\ell}=\alpha_\ell^\vee.
\end{array}
\eqa Recursive parametrization of an open part $N_+^{(0)}$ of
$N_+$ corresponding to a reduced word $I_{\ell}$ is as follows.
Given $v^{A_{\ell}}\in N_+^{(0)}$ we have
 \bqa\label{FactParamAn} v^{A_\ell}(y)=
\mathfrak{X}_1(y)\mathfrak{X}_2(y)\cdot\ldots\cdot\mathfrak{X}_\ell(y)\,,\eqa
where $$\nonumber
\mathfrak{X}_k(y)= X_k(y_{k,n_{k,k}})\cdot\ldots\cdot
X_2(y_{2,n_{k,2}})X_1(y_{1,n_{k,1}}),$$ and
$\mathfrak{X}_1=X_1(y_{11})$. Here we adopt the following notations.
Let $|I_\ell|=m_\ell$ be the length of $w_0$. For the root system of
type  $A_\ell$ one has $m_\ell=\ell(\ell+1)/2$.
 Then for any $k\in\{1,\ldots,\ell\}$ consider a subword
$$I_k=(i_1,\ldots,i_k)\subset
I_\ell=(i_1,\ldots,i_k,\,i_{k+1},\ldots,i_\ell),$$ with
$|I_k|=m_k=k(k+1)/2$. Let $A_k$ be a corresponding root subsystem in
$R_+$ and $v^{C_k}=\mathfrak{X}_1\cdots\mathfrak{X}_k$ be a
factorized parametrization of the corresponding subgroup.
Factorization parameters  for $v^{A_k}(y)$ can be naturally
enumerated as $\{y_{i,n}\}$ with $1\leq i\leq k$, $1\leq n\leq
n_{k,i}$ and \bqa n_{k,i}=k+1-i\,,\hspace{2cm} 1\leq i\leq\ell.\eqa

We are interested in explicit expressions for the following matrix
elements in terms of factorization parameters $\{y_{i,n}\}$
$$\Delta_i(v):=\langle\xi_{\omega_i}^-|\pi_i(v
\dot{w}_0^{-1})\,|\xi_{\omega_i}^+\rangle,  \qquad \Delta_{i}'(v):=
\langle\xi_{\omega_i}^-|\pi_i(v
e_{i*}\dot{w}_0^{-1})\,|\xi_{\omega_i}^+\rangle,$$ where
$\pi_i=\pi_{\omega_i}$ is a fundamental representation with the
highest weight $\omega_i$, $\xi_{\omega_i}^+$ and $\xi_{\omega_i}^-$
are the highest and lowest weight vectors in the representation
$\pi_i$  such that
$\langle\xi_{\omega_i}^+|\,\xi_{\omega_i}^-\rangle=1$. Note that for
Lie algebra $\mathfrak{sl}_{\ell+1}$ according to (\ref{Star}) we
have $i^*=\ell+1-i$. The proof of the following statement
is obtained by an iterative evaluation  of the matrix elements 
taking into account Serre relations and defining ideals
of the fundamental representations and using the technique of \cite{BZ}.

\begin{lem}
\bqa\nonumber\label{DeltaAn}
\Delta_i(v)^{A_{\ell}}=\Big(\prod_{k=1}^i
y_{\ell+1-k,i}\Big)\Delta_i(v)^{A_{\ell}},\hspace{2cm}
i=1,\ldots,\ell,\eqa 
\bqa\label{2rRatioAn}
\Big(\frac{\Delta'_i(v)}{\Delta_i(v)}\Big)^{A_{\ell}}=
e^{x_{\ell,\ell-i+1}-x_{\ell+1,\ell-i+1}}+
\frac{e^{x_{\ell,\ell-i+1}-x_{\ell,\ell-i}}}
{e^{x_{\ell+1,\ell-i+1}-x_{\ell+1,\ell-i}}}\Big(
\frac{\Delta'_i(v)}{\Delta_i(v)}\Big)^{A_{\ell-1}},\quad i=1,\ldots \ell-1,\eqa \bqa
\Big(\frac{\Delta'_{\ell}(v)}{\Delta_{\ell}(v)}\Big)^{A_{\ell}}&=&
\frac{1}{y_{1,\ell}},\nonumber\eqa for $k=2,\ldots,\ell-1$.
\end{lem}
The matrix elements then can be found by resolving the recursive
relations  (\ref{2rRatioAn}).

\begin{lem}
Let $v$ be defined by (\ref{rec1A}) and (\ref{rec2A}). The following
relations for  matrix elements of $v$ in terms of the variables
$y_{i,k}$ hold:
\bqa \nonumber
\Big(\Delta_{\omega_i,\dot{s}_i}(v)\Big)^{A_{\ell}}&=&
\sum\limits_{n=i}^{\ell}y_{i,n},\hspace{2cm}i=1,\ldots,\ell,\\
\Big(\Delta_{\omega_i,\dot{w}_0^{-1}}(v)\Big)^{A_{\ell}} &=&
\prod_{k=i}^\ell\prod_{n=1}^iy_{k+1-n,n},\\
\nonumber \Big(\frac{\Delta'_k(v)}
{\Delta_k(v)}\Big)^{A_{\ell}}&=&\frac{1}{y_{\ell+1-k,k}}\Big(1+
\sum_{n=1}^{\ell-k}\prod_{i=1}^n
\frac{y_{\ell+1-k-i,k+1}}{y_{\ell+1-k-i,k}}\Big),\\ \nonumber
\Big(\frac{\Delta_{\ell}'(v)}{\Delta_{\ell}(v)}\Big)^{A_{\ell}}
&=&\frac{1}{y_{1,\ell}},\qquad k=2,\ldots,\ell-1. \eqa 
 \end{lem}
Combining these expressions with the expression (\ref{meas})  for
the invariant measure on $N_+$
 and substituting into (\ref{psir}), (\ref{psil}) and
(\ref{int}) one completes the proof of  Theorem \ref{teintfgl}.

Now consider an integral representation for
$\mathfrak{gl}_{\ell+1}$-Whittaker function in a modified factorized
parametrization (\ref{anzacgl}). We start with an  analog
of the recursive relations (\ref{2rRatioAn}) for  matrix elements in the modified
factorized parametrization. To simplify the formulation of the
recursive relations it turns out to be useful to 
consider a twisted version 
\bqa\label{GGGL}
 y_{i,n}=e^{x_{\ell+1,i}-x_{\ell+1,i+1}}e^{x_{n+i,i+1}-x_{n+i-1,i}},\qquad 
\eqa
of the modified parametrization (\ref{anzacgl})  by taking
into account  the action of the part $H_R$ of the Cartan generators
(\ref{HRgl}). 
 The simple change of variables (\ref{GGGL}) applied to
 (\ref{2rRatioAn}) gives the following.

\begin{lem}
\begin{enumerate}
\item In the modified factorized parametrization (\ref{GGGL}) recursive relations
 (\ref{2rRatioAn}) are given by
\bqa\label{PartRecursionAn}
\Big(\frac{\Delta'_i(v)}{\Delta_i(v)}\Big)^{A_{\ell}}=
e^{x_{\ell,\ell-i+1}-x_{\ell+1,\ell-i+1}}+
\frac{e^{x_{\ell,\ell-i+1}-
x_{\ell,\ell-i}}}{e^{x_{\ell+1,\ell-i+1}-x_{\ell+1,\ell-i}}}\Big(
\frac{\Delta'_i(v)}{\Delta_i(v)}\Big)^{A_{\ell-1}}.\eqa

\item Solution of the recursive equations  read as
\bqa\label{3rRatioAn} \Big(\frac{\Delta'_i(v)}{\Delta_i(v)}\Big)^{A_\ell}=
\sum_{n=1}^{\ell+1-i}e^{x_{n+i-1,n}-x_{n+i,n}},\qquad i=1,\ldots,\ell,\eqa 

\bqa\label{2DeltaAn}
\Big(\Delta_{\omega_i,\dot{w}_0^{-1}}(v)\Big)^{A_\ell}&=&
\exp\Big\{\sum_{n=1}^{i}(x_{\ell+1,n}-x_{i,n})\Big\},\hspace{2cm}
1\leq i\leq\ell,\\ \nonumber
\Big(\frac{\Delta_i(v)}{\Delta_{i+1}(v)}\Big)^{A_\ell}&=&e^{-x_{\ell+1,i+1}}
\exp\Big\{\sum_{k=1}^{i+1}x_{i+1,k}-\sum_{k=1}^ix_{i,k}\Big\}, \eqa
where $\Delta_{\ell+1}(v)=1$ is assumed.
\end{enumerate}
\end{lem}
Now substitute (\ref{3rRatioAn}), (\ref{2DeltaAn})
into (\ref{psir}), (\ref{psil}) we obtain Whittaker vectors in the  
parametrization (\ref{GGGL}). Taking $\{x_{\ell,k}=0\}$ we recover the 
expressions for Whittaker vectors given in Lemma \ref{Witmfgl}. To
prove the Theorem \ref{teintmfgl}
 one remains to take into account the  measure
$d\mu_{N_+}$ in the modified factorized parametrization.  
This completes the proofs of  the Theorem \ref{teintmfgl}.

\subsubsection{Expressions for
$\mathfrak{so}_{2\ell+1}$-matrix elements:\\
 Proofs of Theorem \ref{tewfso2l1} and Theorem \ref{teintmfso2l1}}
\label{melso2l1}

In this subsection we introduce expressions for  matrix elements
relevant for the construction of integral representations of
$\mathfrak{so}_{2\ell+1}$-Whittaker functions using 
 factorized parametrization of an open
part of $N_+\subset SO(2\ell+1)$. This  provides a proof of the
integral representations of $\mathfrak{so}_{2\ell+1}$-Whittaker
functions presented in Part I.

We are using the root data given by (\ref{RootBn}). Reduced
decomposition $w_0=s_{i_1}s_{i_2}\cdots s_{i_{m_\ell}}$ of the maximal
length element $w_0\in W$  corresponding to a reduced word
$I_{\ell}=(i_1,\ldots,i_{m_{\ell}})$ with $m_{\ell}=\ell^2$ provides a total
ordering of positive coroots   by
$R^{\vee}_+=\{\gamma_k=s_{i_1}\cdots s_{i_{k-1}}\alpha^{\vee}_k\}$
of $\mathfrak{so}_{2\ell+1}$.
   We consider a decomposition of $w_0$ described by the following
   reduced word
$$I_\ell=(1,212,\ldots,(\ell\ldots212\ldots\ell)).$$
Corresponding ordering of positive co-roots is given by:
\bqa\label{orderBn} \gamma^\vee_1=\alpha_1^\vee,\hspace{2cm}
\begin{array}{l}
\gamma^\vee_2=\alpha_1^\vee+\alpha_2^\vee,\\
\gamma^\vee_3=\alpha_1^\vee+2\alpha_2^\vee,\\
\gamma^\vee_4=\alpha_2^\vee,\end{array} \qquad\ldots\\
\nonumber
\begin{array}{l}
\gamma^\vee_{(\ell-1)^2+1}=\alpha_1^\vee+2(\alpha_2^\vee+\ldots+
\alpha_{\ell-1}^\vee)+\alpha_\ell^\vee,\\
\gamma^\vee_{(\ell-1)^2+2}=\alpha_1^\vee+2(\alpha_2^\vee+\ldots+
\alpha_{\ell-2}^\vee)+\alpha_{\ell-1}^\vee+\alpha_\ell^\vee,\\
\vdots\\
\gamma^\vee_{\ell(\ell-1)}=\alpha_1^\vee+\alpha_2^\vee+\ldots+
\alpha_\ell^\vee,\\
\gamma^\vee_{\ell(\ell-1)+1}=\alpha_1^\vee+
2(\alpha_2^\vee+\ldots+\alpha_\ell^\vee),\\
\gamma^\vee_{\ell(\ell-1)+2}=\alpha_2^\vee+\ldots+
\alpha_\ell^\vee,\\
\vdots\\
\gamma^\vee_{\ell^2}=\alpha_\ell^\vee.
\end{array} \eqa
Recursive parametrization of an open part $N_+^{(0)}$ of $N_+$
corresponding to a reduced word $I_{\ell}$ is as follows. Given
$v^{B_{\ell}}\in N_+^{(0)}$ we have
 \bqa\label{FactParamBn} v^{B_\ell}(y)=
\mathfrak{X}_1(y)\mathfrak{X}_2(y)\cdot\ldots\cdot\mathfrak{X}_\ell(y)\,,\eqa
where 
$$
\mathfrak{X}_k(y)=X_k(y_{k,n_{k,k}-1})\cdot\ldots\cdot
X_2(y_{2,n_{k,2}-1})X_1(y_{1,n_{k,1}})X_2(y_{2,n_{k,2}})\cdot\ldots\cdot
X_k(y_{k,n_{k,k}}),\nonumber
$$
and $\mathfrak{X}_1=X_1(y_{11})$.
Here we adopt the following notations.  Let $|I_\ell|=m_\ell$ be the
length of $w_0$. For the root system of type  $B_\ell$ one has
$m_\ell=\ell^2$.
 Then for any $k\in\{1,\ldots,\ell\}$ consider a sub-word
$$I_k=(i_1,\ldots,i_k)\subset
I_\ell=(i_1,\ldots,i_k,\,i_{k+1},\ldots,i_\ell),$$ with
$|I_k|=m_k=k^2$. Let $B_k$ be a corresponding root subsystem in
$R_+$ and $v^{B_k}=\mathfrak{X}_1\cdots\mathfrak{X}_k$ be a
factorized parametrization of the corresponding subgroup.
Factorization parameters  for $v^{B_k}(y)$ can be naturally
enumerated as $\{y_{i,n}\}$ with $1\leq i\leq k$, $1\leq n\leq
n_{k,i}$ and \bqa n_{k,1}=k\,,\qquad n_{k,i}=2(k+1-i),\quad
1<i\leq\ell.\eqa We also use the notation $n_i:=n_{\ell,i}$.

We are interested in explicit expressions for the following matrix
elements in terms of factorization parameters $\{y_{i,n}\}$
$$\Delta_i(v):=\langle\xi_{\omega_i}^-|\pi_i(v
\dot{w}_0^{-1})\,|\xi_{\omega_i}^+\rangle,  \qquad \Delta_{i}'(v):=
\langle\xi_{\omega_i}^-|\pi_i(v
e_{i*}\dot{w}_0^{-1})\,|\xi_{\omega_i}^+\rangle,$$ where
$\pi_i=\pi_{\omega_i}$ is a fundamental representation with the
highest weight $\omega_i$, $\xi_{\omega_i}^+$ and $\xi_{\omega_i}^-$
are the highest and lowest weight vectors in the representation
$\pi_i$  such that
$\langle\xi_{\omega_i}^+|\,\xi_{\omega_i}^-\rangle=1$. Note that for
Lie algebra $\mathfrak{so}_{2\ell+1}$  we have $i\to i^*$ for the
involution defined by (\ref{Star}). The proof of the following statement
is obtained by an iterative evaluation  of the matrix elements 
taking into account Serre relations and defining ideals
of the fundamental representations and using the technique of \cite{BZ}.

\begin{lem}\label{iterB}
Let $v:=v^{B_{\ell}}$ be  defined by (\ref{FactParamBn}). The
following recursive  equations  hold:

\bqa\nonumber \Delta_1(v)^{B_\ell}&=& \Big(\prod_{k=1}^\ell
y_{k,n_k-1}\Big)\cdot \Delta_1(v^{B_{\ell-1}}),\\ \nonumber
\Delta_i(v)^{B_\ell}&=&\Big(\,y_{1,\ell}^2 \prod_{k=2}^i
y_{k,n_k-1}y_{k,n_k}\prod_{k=i+1}^\ell y_{k,n_k-1}^2\Big)\cdot
\Delta_i(v^{B_{\ell-1}}), \qquad 1<i<\ell,\eqa

\bqa\nonumber \Big(\frac{\Delta'_{1}(v)}
{\Delta_{1}(v)}\Big)^{B_{\ell}}&=& \frac{1}{y_{1,{\ell}}}\Big(1+
\frac{y_{2,2(\ell-1)}}{y_{2,2(\ell-1)-1}}\Big)+
\frac{y_{2,2(\ell-1)}}{y_{2,2(\ell-1)-1}} \Big(\frac{\Delta'_{1}(v)}
{\Delta_{1}(v)}\Big)^{B_{\ell-1}},\\
\label{2rRatioBn} \Big(\frac{\Delta'_{k}(v)}
{\Delta_{k}(v)}\Big)^{B_{\ell}}&=& \frac{1}{y_{k,2(\ell+1-k)}}\Big(1+
\frac{y_{k+1,2(\ell-k)}}{y_{k+1,2(\ell-k)-1}}\Big)+\\
\nonumber&& +\frac{y_{k,2(\ell+1-k)-1}}{y_{k,2(\ell+1-k)}}
\frac{y_{k+1,2(\ell-k)}}{y_{k+1,2(\ell-k)-1}}
\Big(\frac{\Delta'_{k}(v)}
{\Delta_{k}(v)}\Big)^{B_{\ell-1}},\\
\Big(\frac{\Delta'_{{\ell}}(v)} {\Delta_{{\ell}}(v)}\Big)^{B_\ell}&=&
\frac{1}{y_{\ell,2}},\qquad k=2,\ldots,\ell-1.\nonumber\eqa
\end{lem}
Now matrix elements can be found by resolving recursive relations
given above. 
\begin{lem}
Let $v$ is defined by (\ref{rec1B}) and (\ref{rec2B}.) The following
expressions of  matrix elements of $v$ in terms of the variables
$y_{i,k}$ hold:
\bqa\label{pmin1}
\Big( \Delta_{\omega_i,s_i}(v)\Big)^{B_{\ell}}=
\sum\limits_{n=i}^{\ell}y_{i,n},\eqa
\bqa\label{pmin2}\Big(\Delta_{1}(v)\Big)^{B_{\ell}}
=\prod_{n=1}^{\ell}y_{1,n}\times
\prod_{k=2}^{\ell}\prod_{n=1}^{\ell+1-k}y_{k,2n-1},\eqa
\bqa\label{pmin3}\Big(\Delta_{k}(v)\Big)^{B_{\ell}}=\prod_{n=2}^{\ell}y_{1,n}^2\times
\prod_{i=k+1}^{\ell}\prod_{n=1}^{\ell+1-i}y_{i,2n-1}^2\times
\prod_{i=2}^k\prod_{n=1}^{\ell+1-i}y_{i,2n-1}y_{i,2n},\\\nonumber
i=1,\ldots,\ell,\,\,\,\,\,k=2,\ldots,\ell.\eqa \bqa 
\Big(\frac{\Delta'_1(v)}{\Delta_1(v)}\Big)^{B_{\ell}}=\sum_{n=1}^\ell\frac{1}{y_{1,n}}
\Big(1+\frac{y_{2,2(n-1)}}{y_{2,2(n-1)-1}}\Big)\prod_{i=n+1}^\ell
\frac{y_{2,2(i-1)}}{y_{2,2(i-1)-1}},\\
\Big( \frac{\Delta'_k(v)}{\Delta_k(v)}\Big)^{B_{\ell}}=
\sum_{n=1}^{n_k/2} \frac{1}{y_{k,2n}}
\Big(1+\frac{y_{k+1,2(n-1)}}{y_{k+1,2(n-1)-1}}\Big)
\prod_{i=n+1}^{n_k/2} \frac{y_{k+1,2(i-1)}}{y_{k+1,2(-1)-1}}
\frac{y_{k,2i-1}}{y_{k,2i}},\\
k=2,\ldots,\ell,\eqa where  $n_1=\ell$ and $n_k=2(\ell+1-k)$.
\end{lem}

Now consider an integral representation for
$\mathfrak{so}_{2\ell+1}$-Whittaker function in a modified factorized
parametrization (\ref{1yBn}), (\ref{2yBn}). We start with an  analog
of the recursive  relations (\ref{2rRatioBn})
for the matrix elements in the modified
factorized parametrization. To simplify the formulation of the
recursive relations it turns out to be useful to 
consider a twisted version 
\bqa\nonumber  y_{1,1}=e^{-x_{\ell,1}}e^{x_{11}-z_{11}},\qquad
y_{1,k}=e^{-x_{\ell,1}}\Big(e^{x_{k-1,1}-z_{k,1}}+ e^{x_{k,1}-z_{k,1}}\Big),\\
\label{1yBntw}
 y_{k,2r-1}=e^{x_{\ell,k-1}-x_{\ell,k}}e^{z_{k+r-1,k}-x_{k+r-2,k-1}},\\
\nonumber
 y_{k,2r}=e^{x_{\ell,k-1}-x_{\ell,k}}e^{z_{k+r-1,k}-x_{k+r-1,k-1}} \eqa for
$k=2,\ldots,\ell$ and $r=1,\ldots,\ell+1-k.$ of the modified parametrization
(\ref{1yBn})  by taking 
into account  the action of the part $H_R$ of the Cartan generators
(\ref{HRso2l1}).

\begin{lem} Choose  an unipotent element $v\in N_+$.
The  following expressions for  the matrix elements of $v$ in
variables $x_{k,i}, z_{k,i}$ defined by (\ref{1yBntw})
 hold: \begin{enumerate}
\item\bqa\nonumber 
\frac{\Delta_k(v)}{\Delta_{k+1}(v)}&=&\exp\Big\{
-\sum_{i=1}^{k}x_{k,i}-2z_{k,1}+2\sum_{i=2}^{k}z_{k,i}-
\sum_{i=1}^{k-1}x_{k-1,i}\Big\}\Big(e^{x_{k-1,1}}+e^{x_{k,1}}\Big)^2,\\
\frac{\Delta_1^2(v)}{\Delta_2(v)}&=&e^{x_{11}-2z_{11}},\qquad
k=2,\ldots,\ell,\eqa
 and  $\Delta_{\ell+1}(v)=1$ is assumed.

\item\bqa\label{3rRatioBn}
\frac{\Delta_1'(v)}{\Delta_1(v)}&=&\sum_{k=1}^\ell e^{z_{k,1}}\\
\nonumber \frac{\Delta_k'(v)}{\Delta_k(v)}&=&e^{x_{k,k}-z_{k,k}}+
\sum_{n=k+1}^\ell\Big(e^{x_{n-1,k}-z_{n,k}}+
e^{x_{n,k}-z_{n,k}}\Big),\qquad k=2,\ldots,\ell.\eqa Here we let
$x_{\ell,k}=0,\,\,\,k=1,\ldots,\ell$.
 We assume, that the terms like
$e^{z_{\ell+1,i}}$ in (\ref{3rRatioBn}) are deleted and  as usual we
suppose that  $\sum_{n=i}^{j}=0$ whenever $i>j.$
\end{enumerate}\end{lem}

Now substitute (\ref{3rRatioAn}), (\ref{2DeltaAn})
into (\ref{psir}), (\ref{psil}) we obtain Whittaker vectors in the  
parametrization  (\ref{1yBntw}). Taking $\{x_{\ell,k}=0\}$ we recover the 
expressions for Whittaker vectors given in Lemma \ref{BWvector}. To
prove the Theorem \ref{teintmfso2l1} 
 one remains to take into account  the measure
$d\mu_{N_+}$ in the modified factorized parametrization.  
This completes the proofs of  the Theorem \ref{teintmfso2l1}.

\subsubsection{Expressions for  $\mathfrak{sp}_{2\ell}$-matrix elements: \\
Proofs  of  Theorem \ref{teintfsp2l} and Theorem \ref{teintmfsp2l}}
\label{melsp2l}

In this subsection we introduce  expressions for  matrix
elements relevant for the construction of integral representations
of $\mathfrak{sp}_{2\ell}$-Whittaker functions using 
the factorized parametrization of an open part of  $N_+\subset
Sp(2\ell)$. This  provides proof of the integral representations
of $\mathfrak{sp}_{2\ell}$-Whittaker  functions presented in Part I.

We are using the root data for $\mathfrak{g}=\mathfrak{sp}_{2\ell}$
given by (\ref{rootc}). Reduced 
decomposition $w_0=s_{i_1}s_{i_2}\cdots s_{i_m}$ of the maximal
length element $w_0\in W$  corresponding to a reduced word
$I_{\ell}=(i_1,\ldots,i_m)$ provides a total ordering of positive
coroots   by $R^{\vee}_+=\{\gamma_k=s_{i_1}\cdots
s_{i_{k-1}}\alpha^{\vee}_k\}$.
   We consider a decomposition of $w_0$ described by the following   reduced word
$$I_\ell=(1,212,\ldots,(\ell\ldots212\ldots\ell)).$$
Corresponding ordering of positive coroots is given by: \bqa\label{orderCn}
\gamma^\vee_1=\alpha_1^\vee\quad
\begin{array}{l}
\gamma^\vee_2=2\alpha_1^\vee+\alpha_2^\vee,\\
\gamma^\vee_3=\alpha_1^\vee+\alpha_2^\vee,\\
\gamma^\vee_4=\alpha_2^\vee,\end{array} \quad\ldots\qquad
\begin{array}{l}
\gamma^\vee_{(\ell-1)^2+1}=
2\alpha_1^\vee+\ldots+2\alpha_{\ell-1}^\vee+\alpha_\ell^\vee,\\
\qquad\vdots\\
\gamma^\vee_{\ell(\ell-1)}=
2\alpha_1^\vee+\alpha_2^\vee\ldots+\alpha_\ell^\vee,\\
\gamma^\vee_{(\ell-1)^2+\ell}=
\alpha_1^\vee+\ldots+\alpha_\ell^\vee,\\
\qquad\vdots\\
\gamma^\vee_{\ell^2}=\alpha_\ell^\vee.
\end{array}
\eqa Recursive parametrization of an open part 
$N_+^{(0)}$ of $N_+$ defined by the reduced word $I_{\ell}$ is
as follows. Given $v^{C_{\ell}}\in N_+^{(0)}$ we have
 \bqa\label{FactParamCn} v^{C_\ell}(y)=
\mathfrak{X}_1(y)\mathfrak{X}_2(y)\cdot\ldots\cdot\mathfrak{X}_\ell(y)\,,\eqa
where  \bqa
\mathfrak{X}_k(y)= 
X_k(y_{k,n_{k,k}-1})\cdot\ldots\cdot
X_2(y_{2,n_{k,2}-1})X_1(y_{1,n_{k,1}})X_2(y_{2,n_{k,2}})\cdot\ldots\cdot
X_k(y_{k,n_{k,k}}),\nonumber\eqa and $\mathfrak{X}_1=X_1(y_{11})$.
Here we adopt the following notations.  Let $|I_\ell|=m_\ell$ be the
length of $w_0$. For the root system of type  $C_\ell$ one has
$m_\ell=\ell^2$.
 Then for any $k\in\{1,\ldots,\ell\}$ consider a subword
$$I_k=(i_1,\ldots,i_k)\subset
I_\ell=(i_1,\ldots,i_k,\,i_{k+1},\ldots,i_\ell)$$ with
$|I_k|=m_k=k^2$. Let $C_k$ be a corresponding root subsystem in
$R_+$ and $v^{C_k}=\mathfrak{X}_1\cdots\mathfrak{X}_k$ be a
factorized parametrization of the corresponding subgroup.
Factorization parameters  for $v^{C_k}(y)$ can be naturally
enumerated as $\{y_{i,n}\}$ with $1\leq i\leq k$, $1\leq n\leq
n_{k,i}$ and \bqa n_{k,1}=k\,,\qquad n_{k,i}=2(k+1-i),\quad
1<i\leq\ell.\eqa Denote also $n_i:=n_{\ell,i}$.

We are interested in explicit expressions for the following matrix
elements in terms of factorization parameters $\{y_{i,n}\}$
$$\Delta_i(v):=\langle\xi_{\omega_i}^-|\pi_i(v
\dot{w}_0^{-1})\,|\xi_{\omega_i}^+\rangle,  \qquad \Delta_{i}'(v):=
\langle\xi_{\omega_i}^-|\pi_i(v
e_{i*}\dot{w}_0^{-1})\,|\xi_{\omega_i}^+\rangle,$$ where
$\pi_i=\pi_{\omega_i}$ is a fundamental representation with the
highest weight $\omega_i$, $\xi_{\omega_i}^+$ and $\xi_{\omega_i}^-$
are the highest and lowest weight vectors in the representation
$\pi_i$  such that
$\langle\xi_{\omega_i}^+|\,\xi_{\omega_i}^-\rangle=1$.
Note that for Lie algebra $\mathfrak{sp}_{2\ell}$
the involution $i\to i^*$  defined by (\ref{Star}) is trivial $i^*=i$.
The proof of the following statement
is obtained by an iterative evaluation  of the matrix elements 
taking into account Serre relations and defining ideals
of the fundamental representations and using the technique of \cite{BZ}.

\begin{lem}
\bqa\nonumber\label{DeltaCn}
\Delta_i(v)^{C_{\ell}}=\Big(y_{1,n_1}
\prod_{k=2}^i y_{k,n_k-1}\,y_{k,n_k} \prod_{k=i+1}^\ell
y^2_{k,n_k-1}\Big)\, \Delta_i(v)^{C_{\ell-1}},\qquad
i=1,\ldots\ell,\eqa
\bqa\nonumber \Big(\frac{\Delta'_1(v)}
{\Delta_1(v)}\Big)^{C_{\ell}}&=& \frac{1}{y_{1,\ell}}\Big(1+
\frac{y_{2,2(\ell-1)}}{y_{2,2(\ell-1)-1}}\Big)^2+
\Big(\frac{y_{2,2(\ell-1)}}{y_{2,2(\ell-1)-1}}\Big)^2
\Big(\frac{\Delta'_1(v)}{\Delta_1(v)}\Big)^{C_{\ell-1}},\\
\label{2rRatioCn} \Big(\frac{\Delta'_k(v)}{\Delta_k(v)}
\Big)^{C_{\ell}}&=& \frac{1}{y_{k,2(\ell+1-k)}}\Big(1+
\frac{y_{k+1,2(\ell-k)}}{y_{k+1,2(\ell-k)-1}}\Big)+\\
\nonumber &&+\frac{y_{k,2(\ell+1-k)-1}}{y_{k,2(\ell+1-k)}}
\frac{y_{k+1,2(\ell-k)}}{y_{k+1,2(\ell-k)-1}}
\Big(\frac{\Delta'_k(v)}{\Delta_k(v)}\Big)^{C_{\ell-1}},\\
\Big(\frac{\Delta'_{\ell}(v)}{\Delta_{\ell}(v)}\Big)^{C_{\ell}}&=&
\frac{1}{y_{\ell,2}},\nonumber\eqa for $k=2,\ldots,\ell-1$.
\end{lem}
Now  matrix elements  can be found by resolving recursive
relations  (\ref{2rRatioCn}).

\begin{lem}
Let $v$ be defined by (\ref{rec1C}) and (\ref{rec2C}). The following
relations for  matrix elements of $v$ in terms of the variables
$y_{i,k}$ hold:
\bqa \nonumber \Big(\Delta_{\omega_i,\dot{s}_i}(v)\Big)^{C_{\ell}}&=&
\sum\limits_{n=i}^{\ell}y_{i,n},\,\,\,\,\,i=1,\ldots,\ell,\\
\Big(\Delta_{\omega_i,\dot{w}_0^{-1}}(v)\Big)^{C_{\ell}}
&=&\prod_{n=1}^\ell y_{1,n}\times
\prod_{k=2}^i\prod_{n=1}^{2(\ell+1-k)}y_{k,n}\times
\prod_{k=i+1}^\ell\prod_{n=1}^{\ell+1-k}y^2_{k,2n-1}, \eqa
\bqa\nonumber
\Big(\frac{\Delta_1'(v)}
{\Delta_1(v)}\Big)^{C_{\ell}}&=&\sum_{n=1}^\ell\frac{1}{y_{1,n}}
\Big(1+\frac{y_{2,2(n-1)}}{y_{2,2(n-1)-1}}\Big)^2\prod_{i=n+1}^\ell
\Big(\frac{y_{2,2(i-1)}}{y_{2,2(i-1)-1}}\Big)^2,\\\nonumber
\Big(\frac{\Delta'_k(v)} {\Delta_k(v)}\Big)^{C_{\ell}}&=&\sum_{n=1}^{\ell+1-k}
\frac{1}{y_{k,2n}}
\Big(1+\frac{y_{k+1,2(n-1)}}{y_{k+1,2(n-1)-1}}\Big)
\prod_{i=n+1}^{\ell+1-k} \frac{y_{k+1,2(i-1)}}{y_{k+1,2(-1)-1}}
\frac{y_{k,2i-1}}{y_{k,2i}},\\\nonumber
\Big( \frac{\Delta_{\ell}'(v)}{\Delta_{\ell}(v)}\Big)^{C_{\ell}}
&=&\frac{1}{y_{\ell,2}},\eqa for
$k=2,\ldots,\ell-1.$
 \end{lem}
Combining these expressions with the expression
(\ref{meas})  for the invariant measure on $N_+$
 and substituting into (\ref{psir}), (\ref{psil}),
(\ref{int}) one completes the proof of Theorem \ref{teintfsp2l}.

Now consider an integral representation for
$\mathfrak{sp}_{2\ell}$-Whittaker function in a modified factorized
parametrization (\ref{ggc1})-(\ref{ggc2}). We start with an  analog
of the recursive  relations (\ref{2rRatioCn}) for the matrix elements in the modified
factorized parametrization. To simplify the formulation of the
recursive relations it turns out to be useful to 
consider a twisted version of the modified parametrization
(\ref{nontwistc})  by taking
into account  the action of the part $H_R$ of the Cartan generators
(\ref{HRsp2l}). Thus we consider the following change of the
variables:
 \bqa\label{twistc} y_{11}=e^{-2z_{\ell,1}}e^{x_{11}+z_{11}},\qquad
y_{1,k}=e^{-2z_{\ell,1}}\Big(e^{z_{k-1,1}+x_{k,1}}+
e^{z_{k,1}+x_{k,1}}\Big),\nonumber\eqa\bqa \nonumber
y_{k,2r-1}=e^{z_{\ell,k-1}-z_{\ell,k}}e^{x_{k+r-1,k}-z_{k+r-2,k-1}},\eqa
\bqa \nonumber
y_{k,2r}=e^{z_{\ell,k-1}-z_{\ell,k}}e^{x_{k+r-1,k}-z_{k+r-1,k-1}},
\qquad r=1,\ldots,\ell+1-k.\eqa Here $k=2,\ldots,\ell.$ 

\begin{lem}
\begin{enumerate}
\item  In a modified factorized parametrization  recursive relations
 (\ref{2rRatioCn}) are given by
\bqa\label{PartRecursionCn}
\Big(\frac{\Delta'_k}{\Delta_k}\Big)^{C_n}\, =\,
e^{z_{n-1,k}-x_{n,k}}+e^{z_{n,k}-x_{n,k}}+
\frac{e^{\langle\alpha_k\,,\,\underline{z}_{n-1}\rangle}}
{e^{\langle\alpha_k\,,\,\underline{z}_{n}\rangle}}\,
\Big(\frac{\Delta'_k}{\Delta_k}\Big)^{C_{n-1}},\qquad 1\leq k<n<\ell,
\nonumber\eqa
with the solution
\bqa\label{3rRatioCn} \Big(\frac{\Delta'_k(v)}{\Delta_k(v)}\Big)^{C_n}=
e^{z_{k,k}-x_{k,k}}+\sum_{n=k+1}^\ell\Big(
e^{z_{n-1,k}-x_{n,k}}+e^{z_{n,k}-x_{n,k}}\Big),
\qquad k=1,\ldots,\ell,\eqa
 where  $\underline{z}_n=(z_{n,1},\ldots,z_{n,n})$ and  we
define:
$\langle\alpha_k\,,\,\underline{z}_{n}\rangle=z_{n,k+1}-z_{n,k}$,
$\langle\alpha_k\,,\,\underline{z}_{n-1}\rangle=z_{n-1,k+1}-z_{n-1,k}$.

\item The  following expressions for  $\Delta_k(v)$ in terms of
variables $x_{k,i}, z_{k,i}$ hold:
\bqa\label{spectralc}
\Big(\frac{\Delta_1(v)}{\Delta_2(v)}\Big)^{C_n}=e^{-z_{\ell,1}}e^{x_{11}},\eqa
\bqa\nonumber \Big(\frac{\Delta_k(v)}{\Delta_{k+1}(v)}\Big)^{C_n}=
e^{-z_{\ell,k}}\Big(e^{z_{k,1}}+e^{z_{k-1,1}}\Big)\exp\Big\{
-\sum_{i=1}^{k}z_{k,i}+x_{k,1}+2\sum_{i=2}^{k}x_{k,i}-
\sum_{i=1}^{k-1}z_{k-1,i}\Big\},\eqa where $k=2,\ldots,\ell$ and
$\Delta_{\ell+1}=1$ is assumed.
\end{enumerate}
\end{lem}

 Now substitute (\ref{3rRatioCn}),(\ref{spectralc})
into (\ref{psir}), (\ref{psil}) we obtain the left/right Whittaker
vectors  in a twisted parametrization (\ref{twistc}).
Taking $\{z_{\ell,k}=0\}$ we recover the formulas
for Whittaker vectors given in Lemma \ref{WhitModC}. To prove the  Theorem
\ref{teintmfsp2l} one remains to take into account the measure
$d\mu_{N_+}$ in the modified factorized parametrization.
This completes the proofs of  the Theorem \ref{teintmfsp2l}.


\subsubsection{Expressions for $\mathfrak{so}_{2\ell}$-matrix elements: \\
Proofs of  Theorem \ref{teintfso2l} and Theorem \ref{teintmfso2l}}
\label{melso2l}

In this subsection we introduce  expressions for  matrix elements
relevant for the construction of integral representations of
$\mathfrak{so}_{2\ell}$-Whittaker functions using 
the  factorized parametrization of an open part of
$N_+\subset SO(2\ell)$. This  provides a proof of the integral
representations of $\mathfrak{so}_{2\ell}$-Whittaker  functions
presented in Part I.

We are using the root data given by (\ref{RootDn}). Reduced
decomposition $w_0=s_{i_1}s_{i_2}\cdots s_{i_{m_{\ell}}}$ of the maximal
length element $w_0\in W$  corresponding to a reduced word
$I_{\ell}=(i_1,\ldots,i_{m_{\ell}})$ with $m_\ell=\ell(\ell-1)$ provides a
total ordering of positive co-roots   by
$R^{\vee}_+=\{\gamma_k=s_{i_1}\cdots s_{i_{k-1}}\alpha^{\vee}_k\}$
of $\mathfrak{sp}_{2\ell}$.
   We consider a decomposition of $w_0$ described by the following
    reduced word
$$I_\ell=(12,3123,\ldots,(\ell\ldots3123\ldots\ell)).$$
Corresponding ordering of positive coroots is given by:

 \bqa\label{orderDn}
\begin{array}{l}
\gamma^\vee_1=\alpha_1^\vee,\\
\gamma^\vee_2=\alpha_2^\vee,
\end{array}
\hspace{2cm}
\begin{array}{l}
\gamma^\vee_3=\alpha_1^\vee+\alpha_2^\vee+\alpha_3^\vee,\\
\gamma^\vee_4=\alpha_2^\vee+\alpha_3^\vee,\\
\gamma^\vee_5=\alpha_2^\vee+\alpha_3^\vee,\\
\gamma^\vee_6=\alpha_3^\vee,
\end{array} \qquad\ldots\\ \nonumber
\begin{array}{l}
\gamma^\vee_{m_{\ell-1}+1}=\alpha_1^\vee+\alpha_2^\vee+
2(\alpha_3^\vee+\ldots+\alpha_{\ell-1}^\vee)+\alpha_\ell^\vee,\\
\gamma^\vee_{m_{\ell-1}+2}=\alpha_1^\vee+\alpha_2^\vee+
2(\alpha_3^\vee+\ldots+\alpha_{\ell-2}^\vee)+\alpha_{\ell-1}^\vee+
\alpha_\ell^\vee,\\
\vdots\\
\gamma^\vee_{m_{\ell-1}+\ell-2}=\alpha_1^\vee+\alpha_2^\vee+
\alpha_3^\vee+\ldots+\alpha_\ell^\vee,\\
\gamma^\vee_{m_{\ell-1}+\ell-1}=p_{\ell-1}\alpha_1^\vee+
p_{\ell}\alpha_2^\vee+\alpha_3^\vee+\ldots+\alpha_\ell^\vee,\\
\gamma^\vee_{m_{\ell-1}+\ell}\,\,=p_{\ell}\alpha_1^\vee+
p_{\ell+1}\alpha_2^\vee+\alpha_3^\vee+\ldots+\alpha_\ell^\vee,\\
\gamma^\vee_{m_{\ell-1}+\ell+1}=
\alpha_3^\vee+\ldots+\alpha_\ell^\vee,\\
\vdots\\
\gamma^\vee_{m_{\ell}}=\alpha_\ell^\vee.
\end{array}
\eqa Recursive parametrization of an open part $N_+^{(0)}$ of
$N_+$ corresponding to a reduced word $I_{\ell}$ is as follows.
Given $v^{D_{\ell}}\in N_+^{(0)}$ we have
 \bqa\label{FactParamDn} v^{D_\ell}(y)=
\mathfrak{X}_2(y)\mathfrak{X}_2(y)\cdot\ldots\cdot
\mathfrak{X}_\ell(y)\,,\eqa where
$$\mathfrak{X}_k(y)=
X_k(y_{k,n_{k,k}-1})\cdot\ldots\cdot
X_3(y_{2,n_{k,3}-1})X_1(y_{1,n_{k,1}})X_2(y_{1,n_{k,2}})
X_3(y_{2,n_{k,3}})\cdot\ldots\cdot X_k(y_{k,n_{k,k}}),\nonumber$$
and $\mathfrak{X}_2=X_1(y_{11})X_2(y_{21})$. Here we adopt the
following notations.  Let $|I_\ell|=m_\ell$ be the length of $w_0$.
For the root system of type  $D_\ell$ one has $m_\ell=\ell(\ell-1)$.
 Then for any $k\in\{1,\ldots,\ell\}$ consider a subword
$$I_k=(i_1,\ldots,i_k)\subset
I_\ell=(i_1,\ldots,i_k,\,i_{k+1},\ldots,i_\ell),$$ with
$|I_k|=m_k=k(k-1)$. Let $D_k$ be a corresponding root subsystem in
$R_+$ and $v^{D_k}=\mathfrak{X}_2\cdots\mathfrak{X}_k$ be a
factorized parametrization of the corresponding subgroup.
Factorization parameters  for $v^{D_k}(y)$ can be naturally
enumerated as $\{y_{i,n}\}$ with $1\leq i\leq k$, $1\leq n\leq
n_{k,i}$ and \bqa n_{k,1}=n_{k,2}=k-1\,,\qquad
n_{k,i}=2(k+1-i),\quad 2<i\leq\ell\eqa We also denote  $n_i=n_{\ell,i}$.

We are interested in explicit expressions for the following matrix
elements in terms of factorization parameters $\{y_{i,n}\}$
$$\Delta_i(v):=\langle\xi_{\omega_i}^-|\pi_i(v
\dot{w}_0^{-1})\,|\xi_{\omega_i}^+\rangle,  \qquad \Delta_{i}'(v):=
\langle\xi_{\omega_i}^-|\pi_i(v
e_{i*}\dot{w}_0^{-1})\,|\xi_{\omega_i}^+\rangle,$$ where
$\pi_i=\pi_{\omega_i}$ is a fundamental representation with the
highest weight $\omega_i$, $\xi_{\omega_i}^+$ and $\xi_{\omega_i}^-$
are the highest and lowest weight vectors in the representation
$\pi_i$  such that
$\langle\xi_{\omega_i}^+|\,\xi_{\omega_i}^-\rangle=1$. Note that for
Lie algebra $\mathfrak{so}_{2\ell}$ the we have $i\to i^*$ for the
involution defined by (\ref{Star}).
The proof of the following statement
is obtained by an iterative evaluation  of the matrix elements 
taking into account Serre relations and defining ideals
of the fundamental representations and using the technique of \cite{BZ}.

\begin{lem} The following recursive relations hold.

\bqa
\Delta_1(v)^{D_\ell}&=&\Big(\,(y_{1,n_1})^{p_{\ell-1}}
(y_{2,n_2})^{p_{\ell}}\prod_{k=3}^\ell y_{k,n_k-1}\Big)
\Delta_1(v)^{D_\ell-1}, \nonumber\\\nonumber 
\Delta_2(v)^{D_\ell}&=&\Big(\,(y_{1,n_1})^{p_{\ell}}
(y_{2,n_2})^{p_{\ell+1}}\prod_{k=3}^\ell y_{k,n_k-1}\Big)
\Delta_2(v)^{D_\ell-1},\\ \nonumber
\Delta_i(v)^{D_\ell}&=&\Big(\,y_{1,n_1}y_{2,n_2}\prod_{k=3}^i
y_{k,n_k-1}y_{k,n_k}\prod_{k=i+1}^\ell y_{k,n_k-1}^2\Big)
\Delta_i(v)^{D_{\ell-1}},\quad 2<i<\ell,\eqa
\bqa
\Big(\frac{\Delta_1'}{\Delta_1}\Big)^{D_{2r}}&=&\frac{1}{y_{1,2r-1}}
\Big(1+\frac{y_{3,2(2r-2)}}{y_{3,2(2r-2)-1}}\Big)+
\frac{y_{2,2r-1}}{y_{1,2r-1}}\frac{y_{3,2(2r-2)}}{y_{3,2(2r-2)-1}}
\Big(\frac{\Delta_1'}{\Delta_1}\Big)^{D_{2r-1}},\nonumber\\\label{1rRatioDn}
\Big(\frac{\Delta_1'}{\Delta_1}\Big)^{D_{2r+1}}&=&\frac{1}{y_{1,2r}}
\Big(1+\frac{y_{3,2(2r-1)}}{y_{3,2(2r-1)-1}}\Big)+
\frac{y_{1,2r}}{y_{2,2r}}\frac{y_{3,2(2r-1)}}{y_{3,2(2r-1)-1}}
\Big(\frac{\Delta_1'}{\Delta_1}\Big)^{D_{2r}},\eqa

\bqa
\Big(\frac{\Delta_2'}{\Delta_2}\Big)^{D_{2r}}&=&\frac{1}{y_{1,2r-1}}
\Big(1+\frac{y_{3,4(r-1)}}{y_{3,4(r-1)-1}}\Big)+
\frac{y_{1,2r-1}}{y_{1,2r-1}}\frac{y_{3,4(r-1)}}{y_{3,4(r-1)-1}}
\Big(\frac{\Delta_2'}{\Delta_2}\Big)^{D_{2r-1}},\nonumber \\ \nonumber
\Big(\frac{\Delta_2'}{\Delta_2}\Big)^{D_{2r+1}}&=&\frac{1}{y_{1,2r-1}}
\Big(1+\frac{y_{3,2(2r-1)}}{y_{3,2(2r-1)-1}}\Big)+
\frac{y_{2,2r}}{y_{1,2r}}\frac{y_{3,2(2r-1)}}{y_{3,2(2r-1)-1}}
\Big(\frac{\Delta_2'}{\Delta_2}\Big)^{D_{2r}}.\eqa
\end{lem}

The matrix elements can be evaluated by resolving recursion
equations and the results of calculation are presented in the
following lemma.

\begin{lem}
Let $v$ be  defined by (\ref{FactParamDn}). The following expressions
of matrix elements of $v$ in terms of the variables $y_{i,k}$ hold:
\bqa \Delta_1(v)=\prod_{n=1}^{\ell/2}y_{1,2n-1}
\prod_{n=1}^{\frac{\ell-1}{2}}y_{2,2n} \prod_{n=3}^\ell\prod_{i=3}^n
y_{i,2(n+1-i)-1},\nonumber\\
 \Delta_2(v)=\prod_{n=1}^{\ell/2}y_{2,2n-1}
\prod_{n=1}^{\frac{\ell-1}{2}}y_{1,2n} \prod_{n=3}^\ell\prod_{i=3}^n
y_{i,2(n+1-i)-1},\\\nonumber
\Delta_k(v)=\prod_{i=1}^k\prod_{n=1}^{n_i}y_{i,n}
\prod_{i=k+1}^\ell\prod_{n=1}^{n_i/2}y_{i,2n-1}^2,\eqa where
 $n_1=n_2=\ell-1$ and $n_k=2(\ell+1-k),\,2<k\leq\ell$.
\bqa\Big(\frac{\Delta_1'}{\Delta_1}+
\frac{\Delta_2'}{\Delta_2}\Big)^{D_\ell}=\sum_{n=1}^{\ell-1}\Big\{
\frac{1}{y_{1,\ell-1}}\prod_{k=1}^{n-1}\Big(
\frac{y_{1,\ell-k}}{y_{1,\ell-k-1}}\Big)^{p_{k-1}} \Big(
\frac{y_{2,\ell-k}}{y_{2,\ell-k-1}}\Big)^{p_{k}}+\\
\nonumber \frac{1}{y_{2,\ell-1}}\prod_{k=1}^{n-1}\Big(
\frac{y_{1,\ell-k}}{y_{1,\ell-k-1}}\Big)^{p_{k}} \Big(
\frac{y_{2,\ell-k}}{y_{2,\ell-k-1}}\Big)^{p_{k+1}}
\Big\}\Big(1+\frac{y_{3,2(\ell-n-1)}}{y_{3,2(\ell-n-1)-1}}\Big)
\prod_{k=1}^{n-1}\frac{y_{3,2(\ell-k-1)}}{y_{3,2(\ell-k-1)-1}}\eqa
where $p_s=(1+(-1)^s)/2$.
\bqa\label{2rRatioDn} \Big(\frac{\Delta'_k (v)}
{\Delta_k(v)}\Big)^{D_\ell}= \frac{1}{y_{k,2(\ell+1-k)}},\eqa for
$k=3,\ldots,\ell$.
\end{lem}

Now consider an integral representation for
$\mathfrak{so}_{2\ell}$-Whittaker function in a modified factorized
parametrization (\ref{nontwistd}). We start with an  analog
of the recursive relations (\ref{1rRatioDn}) for the matrix elements in the modified
factorized parametrization. To simplify the formulation of the
recursive relations it turns out to be useful to 
consider a twisted version of the modified parametrization
(\ref{nontwistd})  by taking
into account  the action of the part $H_R$ of the Cartan generators
(\ref{HRso2l}). Thus we consider the following change of the
variables:
\bqa
y_{1,n}&=&e^{x_{\ell,1}-x_{\ell,2}}\Big(e^{z_{n,1}-x_{n,1}}+
e^{z_{n,1}-x_{n+1,1}}\Big),\quad n=1,\ldots \ell-1,\nonumber \\
y_{2,n}&=&e^{-x_{\ell,1}-x_{\ell,2}}\Big(e^{z_{n,1}+x_{n,1}}+
e^{z_{n,1}+x_{n+1,1}}\Big), \quad n=1,\ldots \ell-1,\label{1yDn}\\
y_{k,2r-1}&=&e^{x_{\ell,k-1}-x_{\ell,k}}e^{z_{k+r-2,k-1}-x_{k+r-2,k-1}},
\nonumber\\ \nonumber
y_{k,2r}&=&e^{x_{\ell,k-1}-x_{\ell,k}}e^{z_{k+r-2,k-1}-x_{k+r-1,k-1}},
\eqa for $k=3,\ldots,\ell$ and $r=1,\ldots,\ell+1-k$.

\begin{lem}  The following recursive
relations in the variables defined by (\ref{1yDn}) hold.
\begin{enumerate}
\item\bqa\Big(\frac{\Delta'_1}{\Delta_1}\Big)^{D_n}+
\Big(\frac{\Delta'_2}{\Delta_2}\Big)^{D_n}=
e^{x_{n-1,2}-z_{n-1,1}}+e^{x_{n,2}-z_{n-1,1}}+\\ \nonumber
\frac{e^{x_{n-1}(\alpha_1)}}{e^{x_{n}(\alpha_1)}}
\Big(\frac{\Delta'_1}{\Delta_1}\Big)^{D_{n-1}}+
\frac{e^{x_{n-1}(\alpha_2)}}{e^{x_{n}(\alpha_2)}}
\Big(\frac{\Delta'_2}{\Delta_2}\Big)^{D_{n-1}}\hspace{1cm}n=2r-1\eqa
\item\bqa\Big(\frac{\Delta'_1}{\Delta_1}\Big)^{D_n}+
\Big(\frac{\Delta'_2}{\Delta_2}\Big)^{D_n}=
e^{x_{n-1,2}-z_{n-1,1}}+e^{x_{n,2}-z_{n-1,1}}+\\ \nonumber
\frac{e^{x_{n-1}(\alpha_2)}}{e^{x_{n}(\alpha_2)}}
\Big(\frac{\Delta'_1}{\Delta_1}\Big)^{D_{n-1}}+
\frac{e^{x_{n-1}(\alpha_1)}}{e^{x_{n}(\alpha_1)}}
\Big(\frac{\Delta'_2}{\Delta_2}\Big)^{D_{n-1}}\hspace{1cm}n=2r.\eqa
\end{enumerate}
\end{lem}
Resolving recursive equations one can easily obtains the following
result. 
\begin{lem} Given an unipotent element $v\in N_+$,
the  following expressions for  the matrix elements of $v$ in
 a modified parametrization hold:
\bqa\label{3rRatioDn}\nonumber
\frac{\Delta'_1}{\Delta_1}+\frac{\Delta'_2}{\Delta_2}\,&=&\,
\sum_{k=1}^{n-1}(e^{z_{k,k}-x_{k,k-1}}+e^{z_{k,k}-x_{k+1,k}}),\\\nonumber
\frac{\Delta'_k}{\Delta_k}\,&=&\,e^{z_{n-1,n+1-k}-x_{n,n+1-k}},\qquad 
k=3,\ldots,n,\eqa
\bqa\nonumber 
\frac{\Delta_2}{\Delta_1}&=&e^{-x_{\ell,1}}e^{x_{11}},\\
\frac{\Delta_1\Delta_2}{\Delta_3}&=&e^{-x_{\ell,2}}\exp\Big\{
-(x_{21}+x_{22})+2z_{11}-x_{11}\Big\}
\Big(e^{x_{11}}+e^{x_{21}}\Big)^2,\\\nonumber 
\frac{\Delta_k}{\Delta_{k+1}}&=&e^{-x_{\ell,k}}\exp\Big\{
-\sum_{i=1}^kx_{k,i}+\\ &&\nonumber+2\sum_{i=1}^{k-1}z_{k-1,i}\,-
\sum_{i=1}^{k-1}x_{k-1,i}\Big\}
\Big(e^{x_{k-1,1}}+e^{x_{k,1}}\Big)^2,\qquad k=3,\ldots,n,\eqa and
$\Delta_{n+1}=1$ is assumed.
\end{lem}

Now substitute (\ref{3rRatioDn}) 
into (\ref{psir}), (\ref{psil}) we obtain Whittaker vectors in a 
parametrization (\ref{1yDn}. Taking $\{x_{\ell,k}=0\}$ we recover the 
expressions for Whittaker vectors given in Lemma \ref{WmrD}. To
prove the Theorem \ref{teintmfso2l}
  one remains to take into account the  measure
$d\mu_{N_+}$ in the modified factorized parametrization. 
This completes the proofs of  the Theorem \ref{teintmfso2l}.

\subsection{Realization of $\mathcal{U}(\mathfrak{g})$ by differential operators }

In this part we  prove  formulae for realization of classical  Lie
algebra generators by differential operators acting in the  space of
(twisted) functions on $N_+$ supplied with the  factorize
parametrization. The analogous formulae for realization of Lie
algebra generators in the modified factorized parametrization
(Gauss-Givental representation) can be  straightforwardly obtained
by a simple change of the variables discussed in Part I and will not
be considered in this section.

Let us outline the general strategy for the derivation  of the
realization of Lie algebra by differential operators used below.
 Let $V_{\mu}$ be  a space of equivariant functions
on $B_-\backslash G$ \bqa \label{Bequivar} f(bg)=\chi_\mu(b)\,\,
f(g)\,,\hspace{2cm} b\in B_-,\eqa where $\chi_{\mu}$   is a
character of the Borel subgroup $B_-\subset G$. Principle series
representation of $\mathcal{U}(\mathfrak{g})$  in $V_{\mu}$ is
defined as  \bqa (Xf)(v)=\frac{\partial}{\partial \varepsilon}
f(ve^{\varepsilon
  X})|_{\varepsilon\rightarrow0}\,, \hspace{2cm} X\in \mathfrak{g}.\eqa
Let  $I=(i_1,\ldots,i_m)$  be a reduced word corresponding to the
reduced decomposition $w_0=s_{i_1}\cdots s_{i_m}$ of the longest
Weyl group element $w_0$. For classical Lie algebras one can chose
$I$ having recursive structure with respect to the rank $\ell$ of
the Lie algebra.   Consider corresponding recursive factorized
parametrization of  unipotent elements of a classical Lie group $G$,
$\mathfrak{g}={\rm Lie}(G)$: \bqa \label{recstruct}
v^{(\ell)}=\mathfrak{X}^{2}_{1}\cdots \mathfrak{X}^{\ell}_{\ell-1}=
v^{(\ell-1)}\cdot\mathfrak{X}^{\ell}_{\ell-1}.\eqa 
We will derive explicit formulae defining representations of
$\mathcal{U}(\mathfrak{g})$ in $V_{\mu}$ in two steps. At the first
step we use recursive structure (\ref{recstruct}) to construct
recursive relations between classical Lie algebra generators for Lie
algebras of adjacent ranks. At the second step we resolve recursion
relations to get explicit formulae for Lie generators of all
classical Lie algebra.

We start with a list of relevant relations  between one-parameter
subgroups in $G$ (see e.g. \cite{Lu}, \cite{BZ}).  Let
$e_i,\,h_i,\,f_i$ be a Chevalley basis of $\mathfrak{g}$, and let
$A=\|a_{ij}\|$ be the Cartan matrix. Let us introduce  
one-parameter subgroups:
\bqa\label{oneparamsub} X_i(y)=e^{ye_i}\,,\qquad
\alpha_i^\vee(y)=e^{yh_i}\,,\qquad Y_i(y)=e^{yf_i}.\eqa
Then the following relations hold:
 \bqa\label{RelationXH}
X_i(y)\alpha_j^\vee(1+\varepsilon)\,=\,
\alpha_j^\vee(1+\varepsilon)X_i(y-a_{ji}\varepsilon y)\quad  mod
(\varepsilon^2), \eqa \bqa\label{RelationXY}
X_i(y)Y_i(\varepsilon)\,=\,
Y_i(\varepsilon)\alpha_i^\vee(1+\varepsilon y)X_i(y-\varepsilon y^2)
\quad mod(\varepsilon^2). \eqa For  $a_{ij}=a_{ji}=-1$ we have
\bqa\label{RelationXXX} X_i(y_1)X_j(y_2)X_i(\varepsilon)\,=\,
X_j\Big(\varepsilon\frac{y_2}{y_1}\Big)\,
X_i\Big(y_1+\varepsilon\Big)\,
X_j\Big(y_2-\varepsilon\frac{y_2}{y_1}\Big)\quad
mod(\varepsilon^2).\eqa For $a_{ij}=-2$ and $a_{ji}=-1$ we have
 \bqa\label{RelationXXXXB}
X_j(y_1)X_i(y_2)X_j(y_3)X_i(\varepsilon)\,=\eqa \bqa =
X_i\Big(\varepsilon\frac{y_3}{y_1}\Big)\,
X_j\Big(y_1+2\varepsilon\frac{y_3}{y_2}\Big)\,
X_i\Big(y_2+\varepsilon-\varepsilon\frac{y_3}{y_1}\Big)\,
X_j\Big(y_3-2\varepsilon\frac{y_3}{y_2}\Big) \quad
mod(\varepsilon^2).\nonumber \eqa For $a_{ji}=-2$ and $a_{ij}=-1$ we
have
 \bqa\label{RelationXXXXC}
 X_j(y_1)X_i(y_2)X_j(y_3)X_i(\varepsilon)\,= \eqa
\bqa = X_i\Big(\varepsilon\frac{y_3^2}{y_1^2}\Big)\,
X_j\Big(y_1+\varepsilon\frac{y_1y_3+y_3^2}{y_1y_2}\Big)\,
X_i\Big(y_2+\varepsilon\frac{y_1^2-y_3^2}{y_1^2}\Big)\,
X_j\Big(y_3-\varepsilon\frac{y_1y_3+y_3^2}{y_1y_2}\Big)\,\,
mod(\varepsilon^2).\nonumber\eqa The derivation of the recursive
relation for the generators of Lie algebra is as follows. Consider the
right action of one-parameter subgroups  (\ref{oneparamsub})
on the recursive factorized representation (\ref{recstruct}) of an element
$v^{\ell}\subset N_+$. One uses the
relations (\ref{RelationXH})-(\ref{RelationXXXXC}) to move the
generators one step to the left. For example, in the case of the
one-parameter subgroup generated by $e_i$ we have: \bqa v^{(\ell)}\,
X_i(\varepsilon)=v^{(\ell-1)}\,\mathfrak{X}^{\ell}_{\ell-1}(y)\,
X_i(\varepsilon)=v^{(\ell-1)}\,\, X_i(c_i(y)\varepsilon)\,\,
\mathfrak{X}^{\ell}_{\ell-1}(y'(y))\quad mod(\varepsilon^2).\eqa
This leads to recursive relations expressing generators of rank
$\ell$ classical Lie algebra in terms of the generators of rank
$(\ell-1)$ classical Lie algebra and the differential operators over
$y_{i,n}$ parametrizing $\mathfrak{X}^{\ell}_{\ell-1}$.  At the final step of the reduction
we use (\ref{Bequivar}). In the following subsections we provide
recursive relations and resolved formulae for generators of all
classical Lie algebras without further comments.

\subsubsection{Generators of $\mathfrak{gl}_{\ell+1}$:  Proof of
  Proposition \ref{propgengl} }
\label{gengl}

Let
$E_{i,i+1}^{(\ell+1)}\,,E_{i,i}^{(\ell+1)}\,,E_{i+1,i}^{(\ell+1)}$
be Chevalley  generators of $\mathfrak{gl}_{\ell+1}$.  Below we
present recursive relations and resolved expressions for these
generators.

Recursive relations are given by:  \be E^{(\ell+1)}_{i,i+1}=
\left(\frac{\partial}{\partial y_{i,\ell+1-i}}+
\frac{y_{i-1,\ell+2-i}}{y_{i,\ell+1-i}}\left(E^{(\ell)}_{i-1,i}-
\frac{\partial}{\partial y_{i-1,\ell+2-i}}\right)\right),
\nonumber\ee \be \hspace{-1cm}
E^{(\ell+1)}_{i,i}\!=\left(\mu_i^{(\ell+1)}-\mu_i^{(\ell)}+
E^{(\ell)}_{i,i}+y_{i-1,\ell+2-i}\frac{\partial}{\partial
y_{i-1,\ell+2-i}}-y_{i,\ell+1-i}\frac{\partial}{\partial
y_{i,\ell+1-i}}\right)\!\!,\,\,\,\ i\neq\ell+1,\nonumber \ee \be
E^{(\ell+1)}_{\ell+1,\ell+1}=(\mu_{\ell+1}^{(\ell+1)}+
y_{\ell,1}\frac{\pr}{\pr y_{\ell,1}}), \ee \be E^{(\ell+1)}_{i+1,i}=
\left(E^{(\ell)}_{i+1,i}+y_{i,\ell+1-i}\left(
E^{(\ell)}_{i,i}-E^{(\ell)}_{i+1,i+1}\right)\right.-\nonumber\\-
y_{i,\ell+1-i}\left.\left(y_{i,\ell+1-i}\frac{\partial}{\partial
y_{i,\ell+1-i}}+y_{i+1,\ell-i}\frac{\partial}{\partial
y_{i+1,\ell-i}}\right)\right). \nonumber\ee

Resolving the recursion  we obtain: \be \hspace{-0.6cm}
E^{(\ell+1)}_{i,i+1}=\sum_{k=0}^{i-1}\prod_{s=0}^k
\frac{y_{i-s,\ell+2-i}}{y_{i+1-s,\ell+1-i}} \frac{\partial}{\partial
y_{i-k,\ell+1-i}}-
\prod_{s=0}^k\frac{y_{i-(s+1),\ell+2-i}}{y_{i-s,\ell+1-i}}
\frac{\partial}{\partial y_{i-(k+1),\ell+2-i}},\nonumber\ee \bqa
\label{cart} E^{(\ell+1)}_{i,i}=
\mu^{(\ell+1)}_i-\sum_{l=1}^{\ell+1-i}
y_{i,l}\frac{\partial}{\partial y_{i,l}}+ \sum_{l=1}^{\ell+2-i}
y_{i-1,l}\frac{\partial}{\partial y_{i-1,l}},\eqa \be
E^{(\ell+1)}_{i+1,i}=\sum_{k=1}^{\ell}\left[(\mu^{(\ell+1)}_{i+1}-\mu^{(\ell+1)}_i)
y_{i,k+1-i}-y_{i,k+1-i}\Big(y_{i,k+1-i}\frac{\partial}{\partial
y_{i,k+1-i}}-
y_{i+1,k-i}\frac{\partial}{\partial y_{i+1,k-i}}\Big)+\right.\\
\nonumber+\left. y_{i,k+1-i}\sum_{s=1}^{k-1}\Big(
y_{i-1,s+2-i}\frac{\partial}{\partial y_{i-1,s+2-i}}-
2y_{i,s+1-i}\frac{\partial}{\partial y_{i,s+1-i}}+
y_{i+1,s-i}\frac{\partial}{\partial y_{i+1,s-i}}\Big)\right].\ee
This completes the proof of Proposition \ref{propgengl}.

\subsubsection{Generators of  $\mathfrak{so}_{2\ell+1}$:  Proof of
  Proposition \ref{Bfactrep} }
\label{genso2l1}

Let  $e_i^{(\ell)}\,,h_i^{(\ell)}\,,f_i^{(\ell)}$ be Chevalley
generators of $\mathfrak{so}_{2\ell+1}$. Below we present recursive
relations and resolved expressions for these generators.

Recursive relations are given by: \bqa
e_1^{(\ell)}=\frac{\partial}{\partial y_{1,\ell}}\,+\,
\frac{y_{2,2(\ell-1)}}{y_{2,2(\ell-1)-1}}\left[
e_1^{(\ell-1)}-\frac{\partial}{\partial y_{1,\ell}}\right]\,+
 2\,\frac{y_{2,2(\ell-1)}}{y_{1,\ell}}
\left(\frac{\partial}{\partial y_{2,2(\ell-1)-1}}-
\frac{\partial}{\partial y_{2,2(\ell-1)}}\right)\nonumber,\eqa

\bqa e_k^{(\ell)}=\frac{\partial}{\partial y_{k,2(\ell+1-k)}}\,+\,
\frac{y_{k,2(\ell+1-k)-1}}{y_{k,2(\ell+1-k)}}
\frac{y_{k+1,2(\ell-k)}}{y_{k+1,2(\ell-k)-1}}\left[e_k^{(\ell-1)}-
\frac{\partial}{\partial y_{k,2(\ell+1-k)-1}}\right]\,+\nonumber
\eqa \bqa \nonumber
+\frac{y_{k+1,2(\ell-k)}}{y_{k,2(\ell+1-k)}}\left(
\frac{\partial}{\partial y_{k+1,2(\ell-k)-1}}-
\frac{\partial}{\partial y_{k+1,2(\ell-k)}}\right),\qquad 1<k<\ell,
\eqa \bqa
f_1^{(\ell)}\,=\,f_1^{(\ell-1)}\,-\,y_{1,\ell}h_1^{(\ell-1)}\,+\,
2y_{1,\ell}y_{2,2(\ell-1)-1}\frac{\partial}{\partial
y_{2,2(\ell-1)-1}}- y_{1,\ell}^2\frac{\partial}{\partial
y_{1,\ell}},\eqa

\bqa f_2^{(\ell)}\,=\,f_2^{(\ell-1)}\,-\,
\Big(y_{2,2(\ell-1)-1}+y_{2,2(\ell-1)}\Big)h_2^{(\ell-1)}\,+\nonumber\eqa
\bqa\nonumber +\Big(y_{2,2(\ell-1)-1}+
y_{2,2(\ell-1)}\Big)y_{3,2(\ell-2)-1} \frac{\partial}{\partial
y_{3,2(\ell-2)-1}}+
2y_{2,2(\ell-1)}y_{1,\ell}\frac{\partial}{\partial y_{1,\ell}}-\\
\nonumber -\left(y_{2,2(\ell-1)-1}^2\frac{\partial}{\partial
y_{2,2(\ell-1)-1}}+
2y_{2,2(\ell-1)}y_{2,2(\ell-1)-1}\frac{\partial}{\partial
y_{2,2(\ell-1)-1}}+y_{2,2(\ell-1)}^2\frac{\partial}{\partial
y_{2,2(\ell-1)}} \right),\eqa

\bqa f_k^{(\ell)}\,=\,f_k^{(\ell-1)}\,-\,
\Big(y_{k,2(\ell+1-k)-1}+y_{k,2(\ell+1-k)}\Big)h_k^{(\ell-1)}\,+\nonumber
\eqa\bqa \nonumber+
\Big(y_{k,2(\ell+1-k)-1}+y_{k,2(\ell+1-k)}\Big)y_{k+1,2(\ell-k)-1}
\frac{\partial}{\partial y_{k+1,2(\ell-k)-1}}+\\ +\nonumber
y_{k,2(\ell+1-k)}\left(y_{k-1,2(\ell+2-k)-1}\frac{\partial}{\partial
y_{k-1,2(\ell+2-k)-1}}+ y_{k-1,2(\ell+2-k)}\frac{\partial}{\partial
y_{k-1,2(\ell+2-k)}}\right)-\\- \nonumber
\left(y_{k,2(\ell+1-k)-1}^2 \frac{\partial}{\partial
y_{k,2(\ell+1-k)-1}}+ 2y_{k,2(\ell+1-k)-1}y_{k,2(\ell+1-k)}
\frac{\partial}{\partial y_{k,2(\ell+1-k)-1}}+\right.\\
\nonumber\left. +y_{k,2(\ell+1-k)}^2 \frac{\partial}{\partial
y_{k,2(\ell+1-k)}} \right).\eqa We have for $h_i$: \bqa
h_k^{(\ell)}=\<\mu^{(\ell)}\,,\alpha_k^\vee\>\,+\,\sum_{i=1}^\ell
a_{k,i}\sum_{j=1}^{n_i}y_{i,j}\frac{\partial}{\partial y_{i,j}},
\eqa where  $n_1=\ell$, $n_k=2(\ell+1-k)$ for $1<k\leq\ell$.

Resolving the recursion  one obtains: \bqa
e_1^{(\ell)}=\frac{\partial}{\partial y_{1,\ell}}\,+\,
\sum_{n=1}^{\ell-1}\left(\frac{\partial}{\partial y_{1,n}}-
\frac{\partial}{\partial y_{1,n+1}}\right) \prod_{i=n}^{\ell-1}
\frac{y_{2,2i}}{y_{2,2i-1}}\,+\nonumber\\
\nonumber +2\left(\frac{\partial}{\partial y_{2,2n-1}}-
\frac{\partial}{\partial y_{2,2n}}\right)
\frac{y_{2,2(n-1)}}{y_{1,n}}
\prod_{i=n+1}^{\ell-1}\frac{y_{2,2i}}{y_{2,2i-1}}, \eqa \bqa
e_k^{(\ell)}=\frac{\partial}{\partial y_{k,2(\ell+1-k)}}\,+\,
\sum_{n=1}^{n-k}\left(\frac{\partial}{\partial y_{k,2n}}-
\frac{\partial}{\partial y_{k,2n+1}}\right)
\prod_{i=n}^{\ell-k}\frac{y_{k+1,2i}}{y_{k+1,2i-1}}
\frac{y_{k,2(i+1)-1}}{y_{k,2(i+1)}}\,+\\\nonumber +\left(
\frac{\partial}{\partial y_{k+1,2n-1}}- \frac{\partial}{\partial
y_{k+1,2n}}\right) \frac{y_{k+1,2n}}{y_{k,2(n+1)}}
\prod_{i=n+1}^{\ell-k}\frac{y_{k+1,2i}}{y_{k+1,2i-1}}
\frac{y_{k,2(i+1)-1}}{y_{k,2(i+1)}},\eqa \bqa
e_\ell^{(\ell)}=\frac{\partial}{\partial y_{\ell,2}},\eqa
 \bqa
f_1^{(\ell)}\,=\,
\sum_{n=1}^\ell\,y_{1,n}\Big(\<\mu,\alpha_1^\vee\>+
\sum_{j=1}^{2(n-1)-1}2y_{2,j}\frac{\partial}{\partial y_{2,j}}
-2\sum_{j=1}^{n-1}y_{1,j}\frac{\partial}{\partial y_{1,j}}-
y_{1,n}\frac{\partial}{\partial y_{1,n}}\Big),\eqa \bqa
f_2^{(\ell)}\,=\,\sum_{n=1}^{2(\ell-1)}\,y_{2,n}\Big(\<\mu,\alpha_2^\vee\>+
2\sum_{j=1}^{[n/2]+1}y_{1,j}\frac{\partial}{\partial y_{1,j}}
-2\sum_{j=1}^{n-1}y_{2,j}\frac{\partial}{\partial y_{2,j}}+\nonumber\\
\nonumber +\sum_{j=1}^{2[(n+1)/2]-3} y_{3,j}
\frac{\partial}{\partial y_{3,j}}-y_{2,n}\frac{\partial}{\partial
y_{2,n}}\Big),\eqa \bqa
f_k^{(\ell)}\,=\,\sum_{n=1}^{2(\ell+1-k)}\,y_{k,n}\Big(\<\mu,\alpha_k^\vee\>+
\sum_{j=1}^{2[n/2]+1}y_{k-1,j}\frac{\partial}{\partial y_{k-1,j}}
-2\sum_{j=1}^{n-1}y_{k,j}\frac{\partial}{\partial
y_{k,j}}+\nonumber\eqa\bqa \nonumber +\sum_{j=1}^{2[(n+1)/2]-3}
y_{k+1,j}\frac{\partial}{\partial y_{k+1,j}}-
y_{k,n}\frac{\partial}{\partial y_{k,n}}\Big), \qquad 2<k<\ell,\eqa

\bqa f_\ell^{(\ell)}\,=\,(y_{\ell,1}+y_{\ell,2})\left(
\<\mu,\alpha_\ell^\vee\>+ y_{\ell-1,1}\frac{\partial}{\partial
y_{\ell-1,1}}+ y_{\ell-1,2}\frac{\partial}{\partial y_{\ell-1,2}}\right)+\nonumber\\
\nonumber +y_{\ell,2}\left( y_{\ell-1,3}\frac{\partial}{\partial
y_{\ell-1,3}}+ y_{\ell-1,4}\frac{\partial}{\partial
y_{\ell-1,4}}\right)-\left(y_{\ell,1}^2\frac{\partial}{\partial
y_{\ell,1}}+ 2y_{\ell,1}y_{\ell,2}\frac{\partial}{\partial
y_{\ell,1}}+ y_{\ell,2}^2\frac{\partial}{\partial
y_{\ell,2}}\right). \eqa This completes the proof of 
 Proposition \ref{Bfactrep}.


\subsubsection{Generators of  $\mathfrak{sp}_{2\ell}$: Proof of 
Proposition \ref{realCn} }
\label{gensp2l}

Let $e_i^{(\ell)}\,,h_i^{(\ell)}\,,f_i^{(\ell)}$ be Chevalley
generators of $\mathfrak{sp}_{2\ell}$.  Below we present recursive
relations and resolved expressions for these generators.

Recursive relations are given by: \bqa
e_1^{(\ell)}=\frac{\partial}{\partial y_{1,\ell}}\,+\,
\Big(\frac{y_{2,2(\ell-1)}}{y_{2,2(\ell-1)-1}}\Big)^2\left[
e_1^{(\ell-1)}-\frac{\partial}{\partial y_{1,\ell}}\right]\,+\nonumber\\
+\nonumber \frac{y_{2,2(\ell-1)-1}y_{2,2(\ell-1)}+
y_{2,2(\ell-1)}^2}{y_{1,\ell}y_{2,2(\ell-1)}}
\left(\frac{\partial}{\partial y_{2,2(\ell-1)-1}}-
\frac{\partial}{\partial y_{2,2(\ell-1)}}\right),\eqa \bqa
e_k^{(\ell)}=\frac{\partial}{\partial y_{k,2(\ell+1-k)}}\,+\,
\frac{y_{k,2(\ell+1-k)-1}}{y_{k,2(\ell+1-k)}}
\frac{y_{k+1,2(\ell-k)}}{y_{k+1,2(\ell-k)-1}}\Big[e_k^{(\ell-1)}-
 \frac{\partial}{\partial y_{k,2(\ell+1-k)-1}}\Big]\,+\nonumber \\
+\frac{y_{k+1,2(\ell-k)}}{y_{k,2(\ell+1-k)}}\left(
\frac{\partial}{\partial y_{k+1,2(\ell-k)-1}}-
\frac{\partial}{\partial y_{k+1,2(\ell-k)}}\right),\qquad
1<k<\ell,\nonumber\eqa
\bqa e_\ell^{(\ell)}=\frac{\partial}{\partial y_{\ell,2}},\eqa 
\bqa
f_1^{(\ell)}\,=\,f_1^{(\ell-1)}\,-\,y_{1,\ell}h_1^{(\ell-1)}\,+\,
y_{1,\ell}y_{2,2(\ell-1)-1}\frac{\partial}{\partial
y_{2,2(\ell-1)-1}}- y_{1,\ell}^2\frac{\partial}{\partial
y_{1,\ell}}, \nonumber\eqa \bqa f_2^{(\ell)}\,=\,f_2^{(\ell-1)}\,-\,
\Big(y_{2,2(\ell-1)-1}+y_{2,2(\ell-1)}\Big)h_2^{(\ell-1)}+
2y_{2,2(\ell-1)-1}y_{1,\ell}\frac{\partial}{\partial y_{1,\ell}}+\nonumber\\
\nonumber
+\Big(y_{2,2(\ell-1)-1}+y_{2,2(\ell-1)}\Big)y_{3,2(\ell-2)-1}
\frac{\partial}{\partial
y_{3,2(\ell-2)-1}}-\left(y_{2,2(\ell-1)-1}^2
\frac{\partial}{\partial y_{2,2(\ell-1)-1}}+\right.\\
\nonumber+\left.
2y_{2,2(\ell-1)-1}y_{2,2(\ell-1)-1}\frac{\partial}{\partial
y_{2,2(\ell-1)-1}}+y_{2,2(\ell-1)}^2\frac{\partial}{\partial
y_{2,2(\ell-1)}} \right), \eqa \bqa
f_k^{(\ell)}\,=\,f_k^{(\ell-1)}\,-\,
\Big(y_{k,2(\ell+1-k)-1}+y_{k,2(\ell+1-k)}\Big)h_k^{(\ell-1)}\,+\nonumber\\
\nonumber
+\Big(y_{k,2(\ell+1-k)-1}+y_{k,2(\ell+1-k)}\Big)y_{k+1,2(\ell-k)-1}
\frac{\partial}{\partial y_{k+1,2(\ell-k)-1}}+\\ \nonumber
+y_{k,2(\ell+1-k)}\left(y_{k-1,2(\ell+2-k)-1}\frac{\partial}{\partial
y_{k-1,2(\ell+2-k)-1}}+ y_{k-1,2(\ell+2-k)}\frac{\partial}{\partial
y_{k-1,2(\ell+2-k)}}\right)-\\ \nonumber
-\left(y_{k,2(\ell+1-k)-1}^2 \frac{\partial}{\partial
y_{k,2(\ell+1-k)-1}}+ 2y_{k,2(\ell+1-k)-1}y_{k,2(\ell+1-k)}
\frac{\partial}{\partial y_{k,2(\ell+1-k)-1}}+\right.\eqa \bqa
\nonumber+\left. y_{k,2(\ell+1-k)}^2 \frac{\partial}{\partial
y_{k,2(\ell+1-k)}} \right),\qquad 2<k<\ell, \eqa \bqa
f_\ell^{(\ell)}\,=\,(y_{\ell,1}+y_{\ell,2})\left(
-\<\mu,\alpha_\ell^\vee\>+ y_{\ell-1,1}\frac{\partial}{\partial
y_{\ell-1,1}}+
y_{\ell-1,2}\frac{\partial}{\partial y_{\ell-1,2}}\right)+\nonumber\\
\nonumber+ y_{\ell,2}\left( y_{\ell-1,3}\frac{\partial}{\partial
y_{\ell-1,3}}+y_{\ell-1,4}\frac{\partial}{\partial
y_{\ell-1,4}}\right)-\left(y_{\ell,1}^2\frac{\partial}{\partial
y_{\ell,1}}+ 2y_{\ell,1}y_{\ell,2}\frac{\partial}{\partial
y_{\ell,1}}+ y_{\ell,2}^2\frac{\partial}{\partial y_{\ell,2}}
\right).\eqa

For the Cartan generators $h_i$ we have: \bqa
h_k^{(\ell)}=\<\mu\,,\alpha_k^\vee\>\,+\,\sum_{i=1}^\ell
a_{k,i}\sum_{j=1}^{n_i}y_{i,j}\frac{\partial}{\partial y_{i,j}},\eqa
where   $n_1=\ell$ and $n_k=2(\ell+1-k)$ for $1<k\leq\ell$.

Resolving the recursion  we obtain: \bqa
e_1=\sum_{n=1}^\ell\left(\frac{\partial}{\partial y_{1,n}}-
\frac{\partial}{\partial y_{1,n+1}}\right) \prod_{j={n}}^{\ell-1}
\Big(\frac{y_{2,2j}}{y_{2,2j-1}}\Big)^2\,+\nonumber\\ \nonumber
+\sum_{n=1}^{\ell-1}\left(\frac{\partial}{\partial y_{2,2n-1}}-
\frac{\partial}{\partial y_{2,2n}}\right) \frac{y_{2,2n}}{y_{1,n}}
\Big(1+\frac{y_{2,2n}}{y_{2,2n-1}}\Big)
\prod_{j=n+1}^{\ell-1}\Big(\frac{y_{2,2j}}{y_{2,2j-1}}\Big)^2, \eqa
\bqa e_k=\sum_{n=1}^{\ell+1-k}\left(\frac{\partial}{\partial
y_{k,2n}}- \frac{\partial}{\partial y_{k,2n+1}}\right)
\prod_{i=n}^{\ell-k}\frac{y_{k+1,2j}}{y_{k+1,2j-1}}
\frac{y_{k,2(j+1)-1}}{y_{k,2(j+1)}}\,+\nonumber\eqa \bqa \nonumber
+\sum_{n=1}^{\ell-k}\left( \frac{\partial}{\partial y_{k+1,2n-1}}-
\frac{\partial}{\partial y_{k+1,2n}}\right)
\frac{y_{k+1,2n}}{y_{k,2(n-1)}}
\prod_{i=n}^{\ell-k}\frac{y_{k+1,2j}}{y_{k+1,2j-1}}
\frac{y_{k,2(j+1)-1}}{y_{k,2(j+1)}},\qquad 1<k<\ell,\eqa \bqa
\nonumber f_1=\sum_{n=1}^\ell\,y_{1,n}\Big(-\<\mu,\alpha_1^\vee\>+
\sum_{j=1}^{2(n-1)-1}y_{2,j}\frac{\partial}{\partial y_{2,j}}
-2\sum_{j=1}^{n-1}y_{1,j}\frac{\partial}{\partial y_{1,j}}-
y_{1,n}\frac{\partial}{\partial y_{1,n}}\Big), \eqa \bqa
f_2=\sum_{n=1}^{2(\ell-1)}\,y_{2,n}\Big(\<\mu,\alpha_2^\vee\>+
2\sum_{j=1}^{[n/2]+1}y_{1,j}\frac{\partial}{\partial y_{1,j}}
-2\sum_{j=1}^{n-1}y_{2,j}\frac{\partial}{\partial y_{2,j}}+\\
\nonumber +\sum_{j=1}^{2[(n+1)/2]-3} y_{3,j}\frac{\partial}{\partial
y_{3,j}}- y_{2,n}\frac{\partial}{\partial y_{2,n}}\Big),\eqa \bqa
f_k=\sum_{n=1}^{2(\ell+1-k)}\,y_{k,n}\Big(\<\mu,\alpha_k^\vee\>+
2\sum_{j=1}^{2[n/2]+1}y_{k-1,j}\frac{\partial}{\partial y_{k-1,j}}
-2\sum_{j=1}^{n-1}y_{k,j}\frac{\partial}{\partial
y_{k,j}}+\nonumber\eqa \bqa \nonumber +\sum_{j=1}^{2[(n+1)/2]-3}
y_{k+1,j}\frac{\partial}{\partial y_{k+1,j}}-
y_{k,n}\frac{\partial}{\partial y_{k,n}}\Big), \qquad 2<k<\ell.\eqa
This completes the proof of Proposition \ref{realCn}.


\subsubsection{Generators of $\mathfrak{so}_{2\ell}$:  Proof of 
Proposition \ref{GGrepD}}
\label{genso2l}

Let $e_i^{(\ell)}\,,h_i^{(\ell)}\,,f_i^{(\ell)}$ be Chevalley
generators of $\mathfrak{so}_{2\ell}$.  Below we present recursive
relations and resolved expressions for these generators.

Recursive relations are given by: \bqa
e_1^{(\ell)}=\frac{\partial}{\partial y_{1,\ell-1}}+
\frac{y_{2,\ell-1}}{y_{1,\ell-1}}
\frac{y_{3,2(\ell-2)}}{y_{3,2(\ell-2)-1}}\Big(\,e_2^{(\ell-1)}-
\frac{\partial}{\partial y_{2,\ell-1}}\Big)
+\frac{y_{3,2(\ell-2)}}{y_{1,\ell-1}}\Big(
\frac{\partial}{\partial\,y_{y_{3,2(\ell-2)-1}}}-
\frac{\partial}{\partial\,y_{y_{3,2(\ell-2)}}}\Big), \nonumber\eqa
\bqa e_2^{(\ell)}=\frac{\partial}{\partial y_{2,\ell-1}}+
\frac{y_{1,\ell-1}}{y_{2,\ell-1}}
\frac{y_{3,2(\ell-2)}}{y_{3,2(\ell-2)-1}}\Big(\,e_1^{(\ell-1)}-
\frac{\partial}{\partial y_{1,\ell-1}}\Big)+
\frac{y_{3,2(\ell-2)}}{y_{2,\ell-1}}\Big(
\frac{\partial}{\partial\,y_{y_{3,2(\ell-2)-1}}}-
\frac{\partial}{\partial\,y_{y_{3,2(\ell-2)}}}\Big),\nonumber\eqa
\bqa e_k^{(\ell)}= \frac{\partial}{\partial y_{k,2(\ell+1-k)}}+
\frac{y_{k,2(\ell+1-k)-1}}{y_{k,2(\ell+1-k)}}
\frac{y_{k+1,2(\ell-k)}}{y_{k+1,2(\ell-k)-1}}\Big(\, e_k^{(\ell-1)}-
\frac{\partial}{\partial y_{k,2(\ell+1-k)-1}}\Big)+\eqa \bqa
+\nonumber \frac{y_{k+1,2(\ell-k)}}{y_{k,2(\ell+1-k)}}\Big(
\frac{\partial}{\partial\,y_{y_{k,2(\ell+1-k)-1}}}-
\frac{\partial}{\partial\,y_{y_{k,2(\ell+1-k)}}}\Big),\qquad
2<k<\ell,\eqa \bqa e_{\ell}=\frac{\partial}{\partial y_{\ell,2}},\eqa
 \bqa f_i^{(\ell)}=f_i^{(\ell-1)}-
y_{i,\ell-1}^2\frac{\partial}{\partial y_{i,\ell-1}}+
y_{i,\ell-1}\Big(h_i^{(\ell-1)}+
y_{3,2(\ell-2)-1}\frac{\partial}{\partial
  y_{3,2(\ell-2)-1}}\Big),\qquad i=1,2,\nonumber \eqa
\bqa f_3^{(\ell)}=f_3^{(\ell-1)}-
\Big(y_{3,2(\ell-2)-1}^2\frac{\partial}{\partial y_{3,2(\ell-2)-1}}+
y_{3,2(\ell-2)}^2\frac{\partial}{\partial
  y_{3,2(\ell-2)}}\Big)+\nonumber \\
\nonumber +(y_{3,2(\ell-2)-1}+ y_{3,2(\ell-2)})
\Big(h_3^{(\ell-1)}+y_{4,2(\ell-3)-1}\frac{\partial}{\partial
y_{4,2(\ell-3)-1}}\Big)+\\ \nonumber+ y_{3,2(\ell-2)}
\Big(y_{1,\ell-1}\frac{\partial}{\partial y_{1,\ell-1}}+
y_{2,\ell-1}\frac{\partial}{\partial
y_{2,\ell-1}}\Big),\nonumber\eqa \bqa\nonumber
f_k^{(\ell)}=f_k^{(\ell-1)}- \Big(
y_{k,2(\ell+1-k)-1}^2\frac{\partial}{\partial y_{k,2(\ell+1-k)-1}}+
y_{k,2(\ell+1-k)}^2\frac{\partial}{\partial y_{k,2(\ell+1-k)}}\Big)+\\
\nonumber +(y_{k,2(\ell+1-k)-1}+
y_{k,2(\ell+1-k)})\Big(h_k^{(\ell-1)}+y_{k+1,2(\ell-k)-1}
\frac{\partial}{\partial y_{k+1,2(\ell-k)-1}}\Big)+\eqa \bqa
\nonumber
+y_{k,2(\ell+1-k)}\Big(y_{k-1,2(\ell+2-k)-1}\frac{\partial}{\partial
y_{k-1,2(\ell+2-k)-1}}+ y_{k-1,2(\ell+2-k)}\frac{\partial}{\partial
y_{k-1,2(\ell+2-k)}}\Big), \qquad 3<k\leq\ell.\eqa

For Cartan generators we have: \bqa
h_i^{(\ell)}=\<\mu,\,\alpha_i^\vee\>+\sum_{k=1}^\ell
a_{i,k}\,\sum_{j=1}^{n_k}y_{k,j}\frac{\partial}{\partial y_{k,j}},
\eqa where $n_1=n_2=\ell-1$, $n_k=2(\ell+1-k)$ for $2<k\leq\ell$.

Resolving the recursion  we obtain:

\bqa \nonumber e_1^{(\ell)}=\frac{\partial}{\partial y_{1,\ell-1}}+
\sum_{n=1}^{[\ell/2]}\Big(\frac{\partial}{\partial y_{2,\ell-n-1}}-
\frac{\partial}{\partial y_{2,\ell-n}}\Big)\prod_{k=1}^{2n-1}
\Big(\frac{y_{1,\ell-k}}{y_{2,\ell-k}}\Big)^{(-1)^k}
\frac{y_{3,2(\ell-1-k)}}{y_{3,2(\ell-1-k)-1}}+\\ \nonumber
+\sum_{n=2}^{[\ell/2]}\Big(\frac{\partial}{\partial y_{1,\ell-n-1}}-
\frac{\partial}{\partial y_{1,\ell-n}}\Big)\prod_{k=1}^{2(n-1)}
\Big(\frac{y_{1,\ell-k}}{y_{2,\ell-k}}\Big)^{(-1)^k}
\frac{y_{3,2(\ell-1-k)}}{y_{3,2(\ell-1-k)-1}}+\\ \nonumber
+\sum_{n=1}^{[\frac{\ell-1}{2}]}\Big( \frac{\partial}{\partial
y_{3,2(2n-1)-1}}- \frac{\partial}{\partial y_{3,2(2n-1)}}\Big)
\frac{y_{3,2(2n-1)}}{y_{1,2n}}\prod_{k=1}^{\ell-2n-1}
\Big(\frac{y_{1,\ell-k}}{y_{2,\ell-k}}\Big)^{(-1)^k}
\frac{y_{3,2(\ell-1-k)}}{y_{3,2(\ell-1-k)-1}}+\\ \nonumber
+\sum_{n=1}^{[\frac{\ell-2}{2}]}\Big( \frac{\partial}{\partial
y_{3,4n-1}}- \frac{\partial}{\partial y_{3,4n}}\Big)
\frac{y_{3,4n}}{y_{2,2n+1}}\prod_{k=1}^{\ell-2(n+1)}
\Big(\frac{y_{1,\ell-k}}{y_{2,\ell-k}}\Big)^{(-1)^k}
\frac{y_{3,2(\ell-1-k)}}{y_{3,2(\ell-1-k)-1}},\eqa

\bqa e_2^{(\ell)}=\frac{\partial}{\partial y_{2,\ell-1}}+
\sum_{n=2}^{[\ell/2]}\Big(\frac{\partial}{\partial y_{2,\ell-n-1}}-
\frac{\partial}{\partial y_{2,\ell-n}}\Big)\prod_{k=1}^{2(n-1)}
\Big(\frac{y_{2,\ell-k}}{y_{1,\ell-k}}\Big)^{(-1)^k}
\frac{y_{3,2(\ell-1-k)}}{y_{3,2(\ell-1-k)-1}}+\\ \nonumber
+\sum_{n=1}^{[\ell/2]}\Big(\frac{\partial}{\partial y_{1,\ell-n-1}}-
\frac{\partial}{\partial y_{1,\ell-n}}\Big)\prod_{k=1}^{2n-1}
\Big(\frac{y_{2,\ell-k}}{y_{1,\ell-k}}\Big)^{(-1)^k}
\frac{y_{3,2(\ell-1-k)}}{y_{3,2(\ell-1-k)-1}}+\\ \nonumber
+\sum_{n=1}^{[\frac{\ell-2}{2}]}\Big( \frac{\partial}{\partial
y_{3,4n-1}}- \frac{\partial}{\partial y_{3,4n}}\Big)
\frac{y_{3,4n}}{y_{1,2n+1}}\prod_{k=1}^{\ell-2(n+1)}
\Big(\frac{y_{2,\ell-k}}{y_{1,\ell-k}}\Big)^{(-1)^k}
\frac{y_{3,2(\ell-1-k)}}{y_{3,2(\ell-1-k)-1}}+\\ \nonumber
+\sum_{n=1}^{[\frac{\ell-1}{2}]}\Big( \frac{\partial}{\partial
y_{3,2(2n-1)-1}}- \frac{\partial}{\partial y_{3,2(2n-1)}}\Big)
\frac{y_{3,2(2n-1)}}{y_{1,2n}}\prod_{k=1}^{\ell-2n-1}
\Big(\frac{y_{2,\ell-k}}{y_{1,\ell-k}}\Big)^{(-1)^k}
\frac{y_{3,2(\ell-1-k)}}{y_{3,2(\ell-1-k)-1}},\eqa \bqa\nonumber
e^{(\ell)}_k=\frac{\partial}{\partial y_{k,2(\ell+1-k)}}+
\sum_{n=1}^{\ell-k}\Big(\frac{\partial}{\partial y_{k,2n}}-
\frac{\partial}{\partial y_{k,2n+1}}\Big)\prod_{i=1}^{\ell+1-n-k)}
\frac{y_{k,2(i+1)-1}}{y_{k,2(k+1)}}
+\frac{y_{k+1,2i}}{y_{k+1,2i-1}}+\eqa \bqa \nonumber
+\sum_{n=1}^{\ell-k}\Big(\frac{\partial}{\partial y_{k+1,2n-1}}-
\frac{\partial}{\partial y_{k+1,2n}}\Big)
\frac{y_{k+1,2n}}{y_{k,2(n+1)}}\prod_{i=2}^{\ell+1-n-k)}
\frac{y_{k,2(i+1)-1}}{y_{k,2(k+1)}}
\frac{y_{k+1,2i}}{y_{k+1,2i-1}},\qquad 2<k<\ell,\eqa

\bqa \nonumber
f_i^{(\ell)}=-\<\mu,\alpha_i^\vee\>\sum_{n=1}^{\ell-1}y_{i,n}-
\sum_{n=1}^{\ell-1}\Big(y_{i,n}^2\frac{\partial}{\partial y_{i,n}}+
2\sum_{k=n+1}^{\ell-1}y_{i,k}y_{i,n}
\frac{\partial}{\partial\,y_{i,n}}\Big)+\eqa\bqa \nonumber
+\sum_{n=1}^{2(\ell-2)-1}\sum_{k=[n/2]+2}y_{i,k}y_{3,n}
\frac{\partial}{\partial\,y_{3,n}}, \qquad i=1,2,\eqa

\bqa \nonumber
f_3^{(\ell)}=-\<\mu,\alpha_3^\vee\>\sum_{n=1}^{2(\ell-2)}y_{3,n}-
\sum_{n=1}^{2(\ell-2)}\Big(y_{3,n}^2
\frac{\partial}{\partial\,y_{3,n}}+
2\sum_{k=n+1}^{2(\ell-2)}y_{3,k}y_{3,n}
\frac{\partial}{\partial\,y_{3,n}}\Big)+\\ \nonumber
+\sum_{n=1}^{\ell-1}\sum_{k=2(n-1)}^{2(\ell-2)}y_{3,k}
\Big(y_{i,n}\frac{\partial}{\partial\,y_{1,n}}+ y_{2,n}
\frac{\partial}{\partial\,y_{i,n}}\Big)+
\sum_{n=1}^{2(\ell-3)-1}\sum_{k=2[n/2]+3}y_{3,k}
y_{4,n}\frac{\partial}{\partial\,y_{4,n}}, \eqa

\bqa \nonumber
f_k^{(\ell)}=-\<\mu,\alpha_k^\vee\>\sum_{n=1}^{2(\ell+1-k)}y_{k,n}-
\sum_{n=1}^{2(\ell+1-k)}\Big(y_{k,n}^2
\frac{\partial}{\partial\,y_{k,n}}+
2\sum_{i=n+1}^{2(\ell+1-k)}y_{k,i}y_{k,n}
\frac{\partial}{\partial\,y_{k,n}}\Big)+\\ \nonumber
+\sum_{n=1}^{\ell+2-k}\,\,\sum_{i=2(n-1)}^{2(\ell+1-k)}y_{k,i}
\Big(y_{k-1,2n-1}\frac{\partial}{\partial\,y_{k-1,2n-1}}+
y_{k-1,2n}\frac{\partial}{\partial\,y_{k-1,2n}}\Big)+\eqa\bqa
\nonumber +\sum_{n=1}^{2(\ell-k)-1}\sum_{i=2[n/2]+3}y_{k,i}
y_{k+1,n}\frac{\partial}{\partial\,y_{k+1,n}}, \qquad
3<k\leq\ell.\eqa
This completes the proof of Proposition \ref{GGrepD}.

\end{document}